# Yetter-Drinfel'd Hopf algebras over groups of prime order

Yorck Sommerhäuser

**Abstract**

We prove a structure theorem for Yetter-Drinfel'd Hopf algebras over groups of prime order that are nontrivial, cocommutative, and cosemisimple: Under certain assumptions on the base field, these algebras can be decomposed into a tensor product of the dual group ring of the group of prime order and an ordinary group ring of some other group. This tensor product is a crossed product as an algebra and an ordinary tensor product as a coalgebra. In particular, the dimension of such a Yetter-Drinfel'd Hopf algebra is divisible by the prime under consideration. We also find explicit examples of such Yetter-Drinfel'd Hopf algebras and apply these results to the classification program for semisimple Hopf algebras.

# Contents

# Introduction

The topic of the present investigation is a special class of Yetter-Drinfel'd Hopf algebras. Yetter-Drinfel'd Hopf algebras are Hopf algebras in a certain quasisymmetric monoidal category, namely the category of Yetter-Drinfel'd modules, which was introduced by D. N. Yetter (cf. [88]). Yetter-Drinfel'd modules are defined with respect to some underlying Hopf algebra, which in our case will be the group ring of a cyclic group of prime order. In this situation, a Yetter-Drinfel'd module is the same as a dimodule, which in turn is the same as a graded vector space together with an action that preserves the gradation.

Throughout the whole discussion, we focus on Yetter-Drinfel'd Hopf algebras that are cosemisimple and cocommutative. The consideration of such a restricted situation is motivated by the structure theory of semisimple Hopf algebras, a field that recently has experienced considerable progress, mainly through the work of A. Masuoka (cf. [47] – [54]). In trying to understand the possible structures for a semisimple Hopf algebra of a given dimension, the analysis often proceeds in the following steps: First, the possible dimensions of the simple modules are determined. This information can sometimes be used to prove the existence of nontrivial grouplike elements. If this can be shown, the existence of analogous nontrivial grouplike elements in the dual may lead to a projection onto a group ring of a cyclic group of prime order. In this case, the Radford projection theorem (cf. [65]) leads to a decomposition of the given Hopf algebra into a tensor product of a Yetter-Drinfel'd Hopf algebra and the group ring of the cyclic group. The Yetter-Drinfel'd Hopf algebras arising in this way are semisimple and tend to be commutative, because the dimensions of their simple modules are often only a fraction the dimensions of the simple modules of the initial Hopf algebra, the denominator of the fraction being the order of the cyclic group. By applying the same reasoning to the dual Hopf algebra, it is sometimes even possible to conclude that the Yetter-Drinfel'd Hopf algebra is also cocommutative.

In many cases, it is then possible to prove that this commutative and cocommutative Yetter-Drinfel'd Hopf algebra must be trivial, which leads to the conclusion that the initial Hopf algebra was a group ring or a dual group ring. It is therefore reasonable to conjecture that such a Yetter-Drinfel'd Hopf algebra over a cyclic group of prime order is always trivial, i. e., has trivial action or coaction. This conjecture was the starting point of the present considerations. However, the conjecture turns out to be false: Although rare objects, Yetter-Drinfel'd Hopf algebras of the form described above exist, although not in all dimensions. The main result of the present investigation is the precise description, under certain assumptions on the base field, of Yetter-Drinfel'd Hopf alge-bras of this type, i. e., nontrivial, cosemisimple, cocommutative Yetter-Drinfel'd Hopf algebras over groups of prime order: The structure theorem obtained in Section 7 asserts that these Yetter-Drinfel'd Hopf algebras can be decomposed into a tensor product of the dual group ring of the cyclic group under consideration and the group ring of another group. The coalgebra structure of this tensor product is the ordinary tensor product coalgebra structure, whereas the algebra structure is that of a crossed product. In particular, the dimension of such a Yetter-Drinfel'd Hopf algebra is divisible by the order of the cyclic group under consideration, a fact that in many interesting cases rules out the existence of a nontrivial Yetter-Drinfel'd Hopf algebra with these properties.

The analysis that leads to this decomposition can be reversed in order to compose nontrivial Yetter-Drinfel'd Hopf algebras. In particular, we construct in this way Yetter-Drinfel'd Hopf algebras of dimension $p^2$ that are nontrivial, commutative, cocommutative, semisimple, and cosemisimple. We prove that these examples exhaust all possibilities, and therefore see that there are, up to isomorphism, $p(p-1)$ Yetter-Drinfel'd Hopf algebras of this form. On the other hand, we show that the dimension $p^2$ is, in a sense, the minimal case; more precisely, we prove that, if a nontrivial, cocommutative, and cosemisimple Yetter-Drinfel'd Hopf algebra is also commutative, its dimension must be divisible by $p^2$. We therefore reach the conclusion that in dimension $n$ the existence of a nontrivial, commutative, cocommutative, semisimple, and cosemisimple Yetter-Drinfel'd Hopf algebra is possible if and only if $p^2$ divides $n$. Finally, we apply these results to semisimple Hopf algebras along the lines indicated above.

The presentation is organized as follows: In Section 1, we recall the basic facts on Yetter-Drinfel'd Hopf algebras, representation theory, and group cohomology that will be needed in the sequel. In Section 2, which is also preliminary, we recall the basic facts from Clifford theory. In our context, Clifford theory is used for the correspondence between the modules of the Yetter-Drinfel'd Hopf algebra and the modules of its Radford biproduct. The part of Clifford theory used in the proof of the main results is a theorem of W. Chin that sets up a correspondence between orbits of the centrally primitive idempotents of the Yetter-Drinfel'd Hopf algebra and orbits of the simple modules of the Radford biproduct under the action of the one-dimensional characters. In the applications, Clifford theory is used to prove that the Yetter-Drinfel'd Hopf algebras occurring are indeed commutative or cocommutative.

In Section 3, we construct examples of nontrivial Yetter-Drinfel'd Hopf algebras that are cocommutative and cosemisimple. We work in a framework created by N. Andruskiewitsch, or a very slight generalization thereof, that describes nicely the requirements that have to be satisfied in form of a compatibility condition. We then find various ways to satisfy this compatibility condition; in particular, we attach a Yetter-Drinfel'd Hopf algebra of this form to every group homomorphism from a finite group to the additive group of a finite ring.

In Section 4, we discuss under which circumstances the Yetter-Drinfel'd Hopf

algebras constructed in Section 3 are isomorphic. In particular, we determine the number of their isomorphism classes in dimension $p^2$.

In Section 5, we investigate the Hopf algebras arising via the Radford biproduct construction and the second construction from [78], [79] from the Yetter-Drinfel'd Hopf algebras considered in Section 3, describe their integrals, and prove that these Hopf algebras are semisimple. In this way, we construct new classes of noncommutative, noncocommutative Hopf algebras with triangular decomposition that are semisimple and cosemisimple. In particular, Masuoka's examples of semisimple Hopf algebras of dimension $p^3$ appear in these Hopf algebras as a kind of Borel subalgebra. However, these two constructions provide only one way to look at these Hopf algebras: They can also be understood from the point of view of extension theory. We show that these Hopf algebras are extensions of group rings by dual group rings, determine the corresponding groups, and find explicit normal bases for these extensions.

In Section 6, we study nontrivial Yetter-Drinfel'd Hopf algebras over groups of prime order that are commutative and semisimple. In this section, most of the technical work is done; in particular, we find a linear and colinear character of order $p$ that induces action and coaction on the other characters. This is the key step in the proof of the structure theorem in the next section.

In Section 7, we dualize the situation and consider nontrivial Yetter-Drinfel'd Hopf algebras over groups of prime order that are cocommutative and cosemisimple. From the previous section, we know that there exists an invariant, coinvariant grouplike element of order $p$ that induces action and coaction on the other grouplike elements. We then pass to a Yetter-Drinfel'd Hopf algebra quotient in which this grouplike element is equal to the unit. Now action and coaction in this quotient are trivial, and therefore this quotient is an ordinary Hopf algebra, which must, since it is cocommutative, be a group ring. It is easy to see that the initial Yetter-Drinfel'd Hopf algebra is a cleft comodule algebra over this group ring, and therefore the structure theorem stated above follows. As an application of the structure theorem, we then prove that the dimension of the Yetter-Drinfel'd Hopf algebra is divisible by $p^2$ if it is also commutative, and classify Yetter-Drinfel'd Hopf algebras of this type in dimension $p^2$.

In Section 8, these results are applied to semisimple Hopf algebras for the first time. For the prime $p$ under consideration, we show that every semisimple Hopf algebra of dimension $p^3$ is a Radford biproduct of a commutative, cocommutative Yetter-Drinfel'd Hopf algebra of dimension $p^2$, which were classified in the preceding section, with the group ring of the cyclic group of order $p$. This gives a new proof of Masuoka's classification of semisimple Hopf algebras of dimension $p^3$ (cf. [50]). In contrast, Masuoka proves that every such Hopf algebra is a Hopf algebra extension of the group ring of a cyclic group of order $p$ and the dual group ring of the elementary abelian group of order $p^2$. Of course, it follows again from the structure theorem that both approaches are related; the precise relation has been worked out in Section 5.

In Section 9, we apply the results to semisimple Hopf algebras of dimension $pq$, where $p$ and $q$ are distinct prime numbers. Although the tools developed in the preceding sections were initially forged to deal with this problem, the results presented here were obtained independently and slightly earlier by P. Etingof, S. Gelaki, and S. Westreich. We prove that a semisimple Hopf algebra of dimension $pq$ is commutative or cocommutative if it and its dual contain nontrivial grouplike elements. In Section 10, we find conditions that assure the existence of such grouplike elements. These results look slightly more general than the result proved in [21], where they were proved under the assumption that the dimensions of the simple modules divide the dimension of the Hopf algebra, i. e., are of dimension 1 or $p$, if $p$ is smaller than $q$. However, Etingof and Gelaki later proved that this is in fact always the case (cf. [18]). Their proof has already been simplified by Y. Tsang and Y. Zhu (cf. [85]), and H.-J. Schneider (cf. [74]); in addition, S. Natale has, in her investigation of semisimple Hopf algebras of dimension $pq^2$, given an extension-theoretic proof of these results (cf. [59]).

In the sequel, $K$ denotes a field. All vector spaces that we consider are defined over $K$, and all tensor products without subscripts are taken over $K$. The multiplicative group $K\setminus\{0\}$ of $K$ is denoted by $K^\times$. Although most of the results mentioned above and proved in the following need the assumption that $K$ is algebraically closed and of characteristic zero, some of the results hold in greater generality. The precise assumptions on $K$ will be stated at the beginning of every section. Likewise, $p$ denotes a natural number that in most cases, but not always, will be prime. The precise assumptions on $p$ will also be stated at the beginning of every section. The transpose of a linear map $f$ is denoted by $f^*$. Following the conventions in group theory, we denote by $\mathbb{Z}_n := \mathbb{Z}/n\mathbb{Z}$ the cyclic group of order $n$, and therefore $\mathbb{Z}_p$ should not be confused with the set of $p$-adic integers. Propositions, definitions, and similar items are referenced by the paragraph in which they occur; they are only numbered separately if this reference is ambiguous.

# 1 Preliminaries

**1.1** Suppose that $H$ is a Hopf algebra defined over the base field $K$. We shall denote its multiplication by $\mu_H$, its unit by $1_H$, its comultiplication by $\Delta_H$, its counit by $\epsilon_H$, and its antipode by $S_H$. We shall also use the following variant of the Heyneman-Sweedler sigma notation for the comultiplication:

$$\Delta_H(h) = h_{(1)} \otimes h_{(2)}$$

Recall the notion of a left Yetter-Drinfel'd module (cf. [88], [57], Def. 10.6.10, p. 213): This is a left $H$-comodule $V$ which is also a left $H$-module such that the following compatibility condition is satisfied:

$$\delta_V(h \to v) = h_{(1)}v^{(1)}S_H(h_{(3)}) \otimes (h_{(2)} \to v^{(2)})$$

for all $h \in H$ and $v \in V$. Here we have used the following Sweedler notation for the coaction: $\delta_V(v) = v^{(1)} \otimes v^{(2)} \in H \otimes V$. The arrow $\to$ denotes the module action.

We also define right Yetter-Drinfel'd modules, which are the same as left Yetter-Drinfel'd modules over the opposite and coopposite Hopf algebra. They are right comodules and right modules that satisfy:

$$\delta_V(v \leftarrow h) = (v^{(1)} \leftarrow h_{(2)}) \otimes S_H(h_{(1)})v^{(2)}h_{(3)}$$

Right Yetter-Drinfel'd modules are used in Paragraph 1.2 and in Section 5. Of course one can also define left-right and right-left Yetter-Drinfel'd modules, which are left Yetter-Drinfel'd modules over the opposite resp. coopposite Hopf algebra, but they are not used in the following. We shall use the convention that a Yetter-Drinfel'd module is a left Yetter-Drinfel'd module unless stated otherwise.

The tensor product of two Yetter-Drinfel'd modules becomes again a Yetter-Drinfel'd module if it is endowed with the diagonal module and the codiagonal comodule structure. The base field $K$ becomes a Yetter-Drinfel'd module via the trivial module structure $h \to \xi := \epsilon_H(h)\xi$ and the trivial comodule structure $\delta_K(\xi) := 1_H \otimes \xi$. Yetter-Drinfel'd modules therefore constitute a monoidal category. This category is also quasisymmetric; the quasisymmetry is given by:

$$\sigma_{V,W} : V \otimes W \longrightarrow W \otimes V$$
$$v \otimes w \mapsto (v^{(1)} \to w) \otimes v^{(2)}$$

This mapping is bijective if $H$ has a bijective antipode (cf. [57], Example 10.6.14, p. 214).

Since we have the notion of a Hopf algebra inside a quasisymmetric monoidal category (cf. [67]), we can speak of Yetter-Drinfel'd Hopf algebras. In general, these are not ordinary Hopf algebras, because the bialgebra axiom asserts that they obey the equation:

$$\Delta_A \circ \mu_A = (\mu_A \otimes \mu_A) \circ (\mathrm{id}_A \otimes \sigma_{A,A} \otimes \mathrm{id}_A) \circ (\Delta_A \otimes \Delta_A)$$

However, it may happen that the quasisymmetry is trivial in the sense that it coincides on $A$ with the usual symmetry in the category of vector spaces, i. e., we have that $\sigma_{A,A}(a \otimes a') = a' \otimes a$ for all $a, a' \in A$. We call this case the trivial case, although in several instances also examples of this type lead to interesting applications:

**Definition** A Yetter-Drinfel'd Hopf algebra $A$ over a Hopf algebra $H$ is called trivial if we have $\sigma_{A,A}(a \otimes a') = a' \otimes a$ for all $a, a' \in A$.

If $A$ is trivial, $A$ obviously is an ordinary Hopf algebra with some additional structure. However, as observed by P. Schauenburg (cf. [69], Cor. 2, p. 262), the converse of this statement is also true:

**Proposition** Suppose that $A$ is a Yetter-Drinfel'd Hopf algebra over the Hopf algebra $H$. Then the following assertions are equivalent:

1. $A$ is trivial.

2. $A$ is an ordinary Hopf algebra.

**Proof.** Obviously, the first statement implies the second. For the converse, observe that we then have:

$$a_{(1)} a'_{(1)} \otimes a_{(2)} a'_{(2)} = \Delta_A(aa') = a_{(1)} (a_{(2)}{}^{(1)} \to a'_{(1)}) \otimes a_{(2)}{}^{(2)} a'_{(2)}$$

for $a, a' \in A$. Now, convolution inversion yields:

$$S_A(a_{(1)}) a_{(2)} a'_{(1)} \otimes a_{(3)} a'_{(2)} S_A(a'_{(3)})$$
$$= S_A(a_{(1)}) a_{(2)} (a_{(3)}{}^{(1)} \to a'_{(1)}) \otimes a_{(3)}{}^{(2)} a'_{(2)} S_A(a'_{(3)})$$

and therefore $a' \otimes a = (a^{(1)} \to a') \otimes a^{(2)}$. $\square$

In the later sections, we shall consider finite-dimensional, cocommutative, cosemisimple Yetter-Drinfel'd Hopf algebras over algebraically closed fields. If these are trivial, they are, by a result of D. K. Harrison and P. Cartier (cf. [35], Thm. 3.2, p. 354, [7], p. 102, [57], Thm. 2.3.1, p. 22), group rings with some additional structure, which we consider as known objects in this context. Therefore, in these sections we shall only deal with nontrivial Yetter-Drinfel'd Hopf algebras.

**1.2** Right Yetter-Drinfel'd modules arise in a natural way as the duals of left Yetter-Drinfel'd modules:

**Proposition** If $V$ is a finite-dimensional left Yetter-Drinfel'd module, then the dual space $V^*$ is in a unique way a right Yetter-Drinfel'd module such that the natural pairing
$$\langle \cdot, \cdot \rangle : V \times V^* \to K, (v, f) \mapsto f(v)$$
is a Yetter-Drinfel'd form (cf. [79], Subsec. 2.4, p. 37), i. e. that we have:

1. $\langle h \to v, f \rangle = \langle v, f \leftarrow h \rangle$
2. $\langle v, f^1 \rangle f^2 = v^1 \langle v^2, f \rangle$

**Proof.** See [81], Lem. 2.3, p. 4. □

The transpose of an $H$-linear and colinear map between finite-dimensional Yetter-Drinfel'd modules is linear and colinear. If $V_1$ and $V_2$ are finite-dimensional left Yetter-Drinfel'd modules, then $V_1^* \otimes V_2^*$ is isomorphic to $(V_1 \otimes V_2)^*$ as a right Yetter-Drinfel'd module via the isomorphism

$$V_1^* \otimes V_2^* \to (V_1 \otimes V_2)^*, f_1 \otimes f_2 \mapsto (v_1 \otimes v_2 \mapsto f_1(v_1)f_2(v_2))$$

Up to this isomorphism, the quasisymmetry on $V_1^* \otimes V_2^*$ is the transpose of the quasisymmetry on $V_1 \otimes V_2$. One can express these facts by saying that taking the dual space is a (non-strict) quasisymmetric monoidal functor from the category of finite-dimensional left Yetter-Drinfel'd modules to the category of finite-dimensional right Yetter-Drinfel'd modules (cf. [30], Def. 2.3, p. 38). As a consequence, if $A$ is a finite-dimensional left Yetter-Drinfel'd Hopf algebra, then the dual space $A^*$ is in a unique way a right Yetter-Drinfel'd Hopf algebra such that the natural pairing described above is a bialgebra form (cf. [79], Subsec. 2.6, p. 37), i. e. we have:

1. $\langle a \otimes a', \Delta_{A^*}(b) \rangle = \langle aa', b \rangle$
2. $\langle a, bb' \rangle = \langle \Delta_A(a), b \otimes b' \rangle$
3. $\langle 1_A, b \rangle = \epsilon_B(b)$, $\langle a, 1_B \rangle = \epsilon_A(a)$

Since left Yetter-Drinfel'd modules are the same as right Yetter-Drinfel'd modules over the opposite and coopposite Hopf algebra, we have:

**Lemma**
1. If $A$ is a left Yetter-Drinfel'd Hopf algebra over $H$, then the opposite and coopposite Hopf algebra $A^{\text{op cop}}$ is a right Yetter-Drinfel'd Hopf algebra over $H^{\text{op cop}}$.

2. If $A$ is a right Yetter-Drinfel'd Hopf algebra over $H$, then $A^{\text{op cop}}$ is a left Yetter-Drinfel'd Hopf algebra over $H^{\text{op cop}}$.

This lemma and the preceding proposition imply that the dual of a finite-dimensional left Yetter-Drinfel'd Hopf algebra may be considered as a left Yetter-Drinfel'd Hopf algebra over the opposite and coopposite Hopf algebra $H^{\text{op cop}}$. By using a suitable isomorphism between $H$ and $H^{\text{op cop}}$, for example a bijective antipode, the dual may even be considered as a left Yetter-Drinfel'd Hopf algebra over $H$.

**1.3** If $H$ is finite-dimensional, the process of dualization can also be applied to $H$:

**Lemma** Suppose that $H$ is finite-dimensional. If $V$ is a right Yetter-Drinfel'd module over $H$, then $V$ becomes a left Yetter-Drinfel'd module over $H^*$ using the module structure
$$f \to v := v^{(1)} f(v^{(2)})$$
for $v \in V$ and $f \in H^*$, and the comodule structure
$$\delta_*(v) := \sum_{i=1}^n h^{i*} \otimes (h_i \to v)$$
where $h_1, \ldots, h_n$ is a basis of $H$ with dual basis $h^{1*}, \ldots, h^{n*}$.

**Proof.** This rests on direct computation (cf. [57], Lem. 1.6.4, p. 11). □

If $H$ is finite-dimensional, a mapping between two Yetter-Drinfel'd modules is linear and colinear with respect to $H$ if and only if it is linear and colinear with respect to $H^*$. The dualization process described in this lemma commutes with taking the tensor product of two Yetter-Drinfel'd modules, and the quasisymmetry on the tensor product is the same before and after the dualization. This can be expressed by saying that dualization with respect to $H$ gives rise to a strict quasisymmetric monoidal functor from the category of right Yetter-Drinfel'd modules over $H$ to left Yetter-Drinfel'd modules over $H^*$. Therefore, a Hopf algebra in the former category remains a Hopf algebra in the latter category.

**1.4** The tensor product of two ordinary Hopf algebras is again a Hopf algebra. The fact that this is not the case for Hopf algebras in quasisymmetric monoidal categories gives rise to a lot of complications, of which the following is only a small part:

**Proposition** Suppose that $A$ and $B$ are Yetter-Drinfel'd Hopf algebras over the Hopf algebra $H$. Then the following assertions are equivalent:

1. $A \otimes B$ is a Yetter-Drinfel'd Hopf algebra if endowed with the tensor product algebra and coalgebra structure.

2. For all $a \in A$ and $b \in B$, we have:
$$\sigma_{A,B}(a \otimes b) = b \otimes a \qquad \sigma_{B,A}(b \otimes a) = a \otimes b$$

**Proof.** We begin with proving that the second statement implies the first. Recall that the Yetter-Drinfel'd structure on $A \otimes B$ was defined in Paragraph 1.1. Since, under these assumptions, the tensor product multiplication $\mu_{A \otimes B}$ can be written in the form

$$\mu_{A \otimes B} = (\mu_A \otimes \mu_B) \circ (\text{id}_A \otimes \sigma_{B,A} \otimes \text{id}_B)$$

we see that $\mu_{A \otimes B}$ is a composition of $H$-linear and colinear maps, and is therefore itself a linear and colinear map. Similarly, the tensor product comultiplication $\Delta_{A \otimes B}$ can be written in the form

$$\Delta_{A \otimes B} = (\text{id}_A \otimes \sigma_{A,B} \otimes \text{id}_B) \circ (\Delta_A \otimes \Delta_B)$$

and therefore is a linear and colinear map. The unit and the counit maps are also tensor products of linear and colinear maps, and therefore are linear and colinear themselves.

Of course, $A \otimes B$ is associative and coassociative. It remains to verify the Yetter-Drinfel'd bialgebra condition. We have on the one hand:

$$\Delta_{A \otimes B}(a \otimes b) \Delta_{A \otimes B}(a' \otimes b')$$
$$= (a_{(1)} \otimes b_{(1)})(a_{(2)}{}^{(1)} b_{(2)}{}^{(1)} \to (a'_{(1)} \otimes b'_{(1)})) \otimes (a_{(2)}{}^{(2)} \otimes b_{(2)}{}^{(2)})(a'_{(2)} \otimes b'_{(2)})$$
$$= a_{(1)}(a_{(2)}{}^{(1)} b_{(2)}{}^{(1)} \to a'_{(1)}) \otimes b_{(1)}(a_{(2)}{}^{(2)} b_{(2)}{}^{(2)} \to b'_{(1)}) \otimes a_{(2)}{}^{(3)} a'_{(2)} \otimes b_{(2)}{}^{(3)} b'_{(2)}$$

and on the other hand:

$$\Delta_{A \otimes B}((a \otimes b)(a' \otimes b')) = (aa')_{(1)} \otimes (bb')_{(1)} \otimes (aa')_{(2)} \otimes (bb')_{(2)}$$
$$= a_{(1)}(a_{(2)}{}^{(1)} \to a'_{(1)}) \otimes b_{(1)}(b_{(2)}{}^{(1)} \to b'_{(1)}) \otimes a_{(2)}{}^{(2)} a'_{(2)} \otimes b_{(2)}{}^{(2)} b'_{(2)}$$

Under our assumptions, both expressions are equal.

The antipode $S_{A \otimes B} := S_A \otimes S_B$ is a tensor product of two linear and colinear maps, and therefore is itself linear and colinear. It is easy to see that it really is an antipode.

To prove that the first assertion implies the second, note that, if both expressions calculated above are equal, we get by applying $\epsilon_A \otimes \text{id}_B \otimes \text{id}_A \otimes \epsilon_B$ that:

$$b(a^{(1)} \to b') \otimes a^{(2)} a' = bb' \otimes aa'$$

For $b = 1_B$ and $a' = 1_A$, this yields $\sigma_{A,B}(a \otimes b') = b' \otimes a$. On the other hand, by applying $\text{id}_A \otimes \epsilon_B \otimes \epsilon_A \otimes \text{id}_B$ to these expressions, we get:

$$a(b^{(1)} \to a') \otimes b^{(2)}b' = aa' \otimes bb'$$

For $a = 1_A$ and $b' = 1_B$, this yields $\sigma_{B,A}(b \otimes a') = a' \otimes b$. □

The topic of tensor products of Hopf algebras in quasisymmetric categories is treated in greater detail in [60], Sec. 4.

**1.5** As for any coalgebra, we call an element $g$ of a Yetter-Drinfel'd Hopf algebra $A$ a grouplike element if $\Delta_A(g) = g \otimes g$ and $\epsilon_A(g) = 1_K$. In contrast to ordinary Hopf algebras, the grouplike elements usually do not form a group, because the product of two grouplike elements in general is not a grouplike element. However, one has the following substitute:

**Proposition 1** Suppose that $A$ is a Yetter-Drinfel'd Hopf algebra over a Hopf algebra $H$. Denote by $G(A)$ the set of grouplike elements. Denote by

$$G_I(A) := \{g \in G(A) \mid \forall h \in H : h \to g = \epsilon_H(h)g\}$$

the set of invariant grouplike elements and by

$$G_C(A) := \{g \in G(A) \mid \delta_A(g) = 1_H \otimes g\}$$

the set of coinvariant grouplike elements. Then we have:

1. Every grouplike element is invertible.

2. $G_I(A)$ and $G_C(A)$ are subgroups of the group of units.

3. $G_I(A)$ acts on $G(A)$ via right multiplication.

4. $G_C(A)$ acts on $G(A)$ via left multiplication.

**Proof.** As for ordinary Hopf algebras, $S_A(g)$ is an inverse of the grouplike element $g$. If $g$ and $g'$ are grouplike elements, we have:

$$\Delta_A(gg') = g(g^{(1)} \to g') \otimes (g^{(2)} \to g')$$

Therefore, $gg'$ is grouplike if $g$ is coinvariant or $g'$ is invariant. This proves the third and the fourth statement. The product of two invariant elements is again invariant, therefore $G_I(A)$ is a multiplicatively closed set. Similarly, the product of two coinvariant elements is again coinvariant, therefore $G_C(A)$ is a multiplicatively closed set. To prove that they are subgroups, we must prove that they contain inverses. If $g$ is a grouplike element, we have:

$$\Delta_A(S_A(g)) = S_A(g^{(1)} \to g) \otimes S_A(g^{(2)})$$

Therefore we see that, if $g$ is invariant or coinvariant, then the inverse $S_A(g)$ is again grouplike. It is also invariant or coinvariant. □

As we have seen in Paragraph 1.2, the dual of a finite-dimensional left Yetter-Drinfel'd Hopf algebra over $H$ is again a left Yetter-Drinfel'd Hopf algebra over $H^{\text{op cop}}$. The grouplike elements of the dual are precisely the one-dimensional characters, i. e., the algebra homomorphisms to the base field. We state the following proposition, which is precisely the dual of the above proposition in the finite-dimensional case, without proof. Observe that the dual of a coalgebra is always an algebra, even in the infinite-dimensional case.

**Proposition 2** Suppose that $A$ is a Yetter-Drinfel'd Hopf algebra over $H$. Denote by $G(A^*)$ the set of one-dimensional characters of $A$. Denote by $G_I(A^*)$ the subset of $H$-linear characters and by $G_C(A^*)$ the set of colinear characters. Then we have:

1. Every one-dimensional character is invertible.

2. $G_I(A^*)$ and $G_C(A^*)$ are subgroups of the group of units.

3. $G_I(A^*)$ acts on $G(A^*)$ via left multiplication.

4. $G_C(A^*)$ acts on $G(A^*)$ via right multiplication.

**1.6** Suppose that $A$ is a Yetter-Drinfel'd Hopf algebra over the Hopf algebra $H$. Since $A$ is in particular a module algebra over $H$, we can form the smash product (cf. [57], Def. 4.1.3, p. 41). This is an algebra with underlying vector space $A \otimes H$, multiplication
$$(a \otimes h)(a' \otimes h') = a(h_{(1)} \to a') \otimes h_{(2)} h'$$
and unit $1_A \otimes 1_H$.

Since $A$ is by assumption also a comodule coalgebra, we can dually form the cosmash product. This is a coalgebra with underlying vector space $A \otimes H$, comultiplication:
$$\Delta_{A \otimes H}(a \otimes h) = (a_{(1)} \otimes a_{(2)}{}^{(1)} h_{(1)}) \otimes (a_{(2)}{}^{(2)} \otimes h_{(2)})$$
and counit $\epsilon_A \otimes \epsilon_H$.

As observed by D. E. Radford (cf. [65], Thm. 1, p. 328, [57], § 10.6), the Yetter-Drinfel'd condition assures in this situation that $A \otimes H$ even becomes a Hopf algebra with these structures. This Hopf algebra is called the Radford biproduct of $A$ and $H$. Its antipode is given by:
$$S_{A \otimes H}(a \otimes h) = (1_A \otimes S_H(a^{(1)} h))(S_A(a^{(2)}) \otimes 1_H)$$

Now suppose that $A$ and $H$ are both finite-dimensional. By Paragraph 1.2, $A^*$ then is a right Yetter-Drinfel'd Hopf algebra over $H$, and thus $A^*$ is a left Yetter-Drinfel'd Hopf algebra over $H^*$ by Paragraph 1.3. Therefore, we can also form the Radford biproduct $A^* \otimes H^*$. This Hopf algebra is precisely the dual of the original one:

**Proposition** The canonical isomorphism
$$A^* \otimes H^* \to (A \otimes H)^*$$
is an isomorphism of Hopf algebras.

We leave the verification of this statement to the reader.

**1.7** Suppose that $H$ is a Hopf algebra. A module algebra over $H$ is, by definition, the same as an algebra in the monoidal category of left $H$-modules. If $A$ is a module algebra over $H$, then a module over $A$ in the category of $H$-modules is, by definition, an $A$-module $V$ for which the structure map $A \otimes V \to V$ is a morphism in the category, i. e., is an $H$-linear map. It is easy to see that $A$-modules in the category of $H$-modules are the same as modules over the smash product $A \otimes H$:

**Proposition 1**
1. Suppose that $A$ is an algebra in the category of left $H$-modules. Suppose that $V$ is an $A$-module in the category of $H$-modules. Then $V$ is a module over the smash product $A \otimes H$ via:
$$(a \otimes h) \to v := a(h \to v)$$
2. Conversely, if $V$ is an $A \otimes H$-module, $V$ is an $H$-module and an $A$-module by restriction. Then, the $A$-module structure map is $H$-linear, and therefore $V$ is an $A$-module in the category of $H$-modules.

Stating this simple fact in a more categorical language, we have constructed a functor $F : A-(H-\text{Mod}) \to (A \otimes H)-\text{Mod}$ between the category $A-(H-\text{Mod})$ of $A$-modules in the category of $H$-modules and the category of modules over the smash product $A \otimes H$, which is a category equivalence. Now suppose that $A$ is a Yetter-Drinfel'd Hopf algebra over $H$. We want to understand how the above category equivalence is compatible with tensor products.

First of all, as in all quasisymmetric monoidal categories, the tensor product of two algebras $A$ and $A'$ in the Yetter-Drinfel'd category is again an algebra with the multiplication
$$(a \otimes a')(b \otimes b') = a(a'^{(1)} \to b) \otimes a'^{(2)} b'$$

This multiplication was already used in Paragraph 1.4; it is just the left smash product of the module algebra $A$ and the comodule algebra $A'$ (cf. [4], Def. 1.2, p. 24, see also [16], Rem. (1.3), p. 374). We denote the tensor product by $A \hat{\otimes} A'$ if we want to emphasize that the tensor product carries this algebra structure. The Yetter-Drinfel'd bialgebra axiom then says precisely that $\Delta_A : A \to A \hat{\otimes} A$ is an algebra homomorphism. We now want to introduce a similar tensor product on $A$-modules.

**Proposition 2** Suppose that $A$ is a Yetter-Drinfel'd Hopf algebra over $H$ and that $V$ and $W$ are $A$-modules.

1. Suppose that $V$ is an $A$-module in the category of $H$-modules, i. e., $V$ is also an $H$-module and the module structure map $A \otimes V \to V$ is $H$-linear. Then $V \otimes W$ is an $A \hat{\otimes} A$-module via:
$$(a \otimes a')(v \otimes w) := a(a'^{(1)} \to v) \otimes a'^{(2)} w$$
We denote $V \otimes W$ by $V \hat{\otimes} W$ if we want to emphasize that it carries this module structure.

2. If also $W$ is an $A$-module in the category of $H$-modules, then $V \hat{\otimes} W$ is an $H$-module via the diagonal $H$-action, and the $A \hat{\otimes} A$-module structure map is $H$-linear.

This can be verified by direct computation. We note that in the situation of the first assertion of Proposition 2, $V \hat{\otimes} W$ also becomes an $A$-module by pulling back the $A \hat{\otimes} A$-module structure along $\Delta_A$. In the language of categories, this means that the category $A - \text{Mod}$ of left $A$-modules is a $(A \otimes H)$-category in the sense of [63], p. 351.

In the situation of the second assertion of Proposition 2, $V$ and $W$ are then both modules over the Radford biproduct. Therefore $V \otimes W$ is a module over the Radford biproduct, and we can restrict the module structure to $A$ as well as to $H$. These restricted module structures then obviously coincide with the module structures constructed in Proposition 2. Expressed in the language of categories, this means that the functor
$$F : A - (H - \text{Mod}) \to (A \otimes H) - \text{Mod}$$
defined above is a strict monoidal functor (cf. [44], Chap. VII, § 1, p. 160), i. e., we have $F(V \hat{\otimes} W) = F(V) \otimes F(W)$.

**1.8** We shall need the following variant of the Nichols-Zoeller theorem for Yetter-Drinfel'd Hopf algebras:

**Proposition** Suppose that $A$ is a finite-dimensional Yetter-Drinfel'd Hopf algebra over the finite-dimensional Hopf algebra $H$. Suppose that $B$ is a Yetter-Drinfel'd Hopf subalgebra of $A$. Then $A$ is free as a left $B$-module.

**Proof.** Let $n := \dim H$. Since the Radford biproduct $B \otimes H$ is a Hopf subalgebra of the Radford biproduct $A \otimes H$, the Nichols-Zoeller theorem (cf. [61], Thm. 7, p. 384, [57], Thm. 3.1.5, p. 30) implies that, for some $m \in \mathbb{N}$, we have $A \otimes H \cong (B \otimes H)^m$ as left $B \otimes H$-modules, and therefore also as left $B$-modules. Since we have $A \otimes H \cong A^n$ and $B \otimes H \cong B^n$ as left $B$-modules, we get:
$$A^n \cong B^{mn}$$
as left $B$-modules.

Now suppose that $B \cong \bigoplus_{i=1}^{k} P_i^{n_i}$ is a decomposition of the left regular representation of $B$ into indecomposable $B$-modules such that the indecomposable $B$-modules $P_1, \ldots, P_k$ are pairwise not isomorphic. Pick a similar decomposition

$$A \cong \bigoplus_{j=1}^{l} M_j^{m_j}$$

of the left $B$-module $A$ into indecomposable $B$-modules such that the indecomposable $B$-modules $M_1, \ldots, M_l$ are pairwise not isomorphic. Then we have:

$$\bigoplus_{j=1}^{l} M_j^{nm_j} \cong A^n \cong B^{mn} \cong \bigoplus_{i=1}^{k} P_i^{mnn_i}$$

Now the Krull-Remak-Schmidt theorem implies that, by possibly changing the enumeration of $M_1, \ldots, M_l$, we have $k = l$, $M_i \cong P_i$ and $nm_i = mnn_i$. This implies $m_i = mn_i$ and therefore we have:

$$A \cong \bigoplus_{i=1}^{l} M_i^{m_i} \cong \bigoplus_{i=1}^{k} P_i^{mn_i} \cong B^m$$

as left $B$-modules. $\square$

As observed, among others, by B. Scharfschwerdt (cf. [68]), the proof of the Nichols-Zoeller theorem carries over directly to Yetter-Drinfel'd Hopf algebras. Therefore, in the above proposition the assumption that $H$ be finite-dimensional is superfluous.

**1.9**  Now suppose that $A$ is a finite-dimensional semisimple Yetter-Drinfel'd Hopf algebra over the Hopf algebra $H$. Denote by $B := A \otimes H$ the Radford biproduct of $A$ and $H$. We choose a nonzero right integral $\rho_A \in A^*$. Observe that $A$ becomes a $B$-module via the action:

$$(a \otimes h) \to a' := a(h \to a')$$

(cf. [57], (4.5.1), p. 53). We want to study the dual module $A^*$ of $A$. To this end, we introduce the following bilinear form:

$$\langle \cdot, \cdot \rangle : A \times A \to K, (a, a') \mapsto \rho_A(S_A(a)a')$$

This bilinear form has the following properties:

**Proposition**  For $a, a', a'' \in A$ and $h \in H$, we have:

1. $\langle h \to a, a' \rangle = \langle a, S_H(h) \to a' \rangle$
2. $a^{(1)} \langle a^{(2)}, a' \rangle = S_H(a'^{(1)}) \langle a, a'^{(2)} \rangle$
3. $\langle (a \otimes h) \to a', a'' \rangle = \langle a', S_B(a \otimes h) \to a'' \rangle$

**Proof.** Since $A$ is semisimple, we have:
$$\rho_A(h \to a) = \epsilon_H(h)\rho_A(a) \qquad a^{(1)}\rho_A(a^{(2)}) = 1_H\rho_A(a)$$

(cf. [81], Prop. 2.14, p. 22). This implies:
$$\langle h_{(1)} \to a, h_{(2)} \to a'\rangle = \epsilon_H(h)\langle a, a'\rangle \qquad a^{(1)}a'^{(1)}\langle a^{(2)}, a'^{(2)}\rangle = 1_H\langle a, a'\rangle$$

which implies the first two assertions. For the third assertion, recall the formula:
$$S_A(aa') = S_A(a^{(1)} \to a')S_A(a^{(2)})$$

From the first assertion, we then get:
$$\begin{aligned}
\langle a', S_B(a \otimes h) \to a''\rangle &= \langle a', (1_A \otimes S_H(a^{(1)}h))(S_A(a^{(2)}) \otimes 1_H) \to a''\rangle \\
&= \langle a^{(1)}h \to a', S_A(a^{(2)})a''\rangle \\
&= \rho_A(S_A(a^{(1)}h \to a')S_A(a^{(2)})a'') \\
&= \rho_A(S_A(a(h \to a'))a'') \\
&= \langle (a \otimes h) \to a', a''\rangle
\end{aligned}$$

Note that the first assertion is a special case of the third assertion, namely the case $a = 1_A$. $\square$

For the dual $A^*$ of the $B$-module $A$, this proposition says that it is isomorphic to the module $A$ itself. This is because it follows directly from the third assertion of the proposition that the map
$$f : A \to A^*, a \mapsto (a' \mapsto \langle a, a'\rangle)$$

is a $B$-module homomorphism. It is an isomorphism since $A$ is a Frobenius algebra with Frobenius homomorphism $\rho_A$ (cf. [20], Cor. 5.8, p. 4885, [81], Prop. 2.10, p. 15).

**1.10** If $H$ is commutative and cocommutative, the Yetter-Drinfel'd condition reads:
$$\delta_V(h \to v) = v^{(1)} \otimes (h \to v^{(2)})$$

A module that is simultaneously a comodule such that this compatibility condition is satisfied is called a dimodule (cf. [39]). In particular, if $G$ is a finite abelian group and $H = K[G]$ is its group ring, then an $H$-comodule is the same as a $G$-graded vector space, where the homogeneous component of degree $g$ is given by
$$V_g = \{v \in V \mid \delta_V(v) = c_g \otimes v\}$$

where $c_g \in K[G]$ denotes the canonical basis element of the group ring corresponding to $g \in G$. The dimodule condition above then says precisely that the homogeneous components are submodules.

Now suppose that $G = \mathbb{Z}_p$, where $p$ is not necessarily prime. We denote by $C$ the set of grouplike elements of $H$, $C := G(H)$. These are precisely the canonical basis elements of the group ring, i. e., we have:

$$C = \{c_i \mid i \in \mathbb{Z}_p\}$$

Suppose that $K$ contains a primitive $p$-th root of unity $\zeta$, which is only possible if the characteristic of $K$ is relatively prime to $p$. Then we have a group homomorphism $\gamma : \mathbb{Z}_p \to K^\times$ that maps 1 to $\zeta$. If $\hat{\mathbb{Z}}_p$ denotes the character group of $\mathbb{Z}_p$, i. e., the group of all group homomorphisms from $\mathbb{Z}_p$ to $K^\times$ under pointwise multiplication, $\gamma$ is of order $p$ in $\hat{\mathbb{Z}}_p$, and therefore there is a group homomorphism $\mathbb{Z}_p \to \hat{\mathbb{Z}}_p$ mapping 1 to $\gamma$. It is easy to see that this is an isomorphism, which is a very special case of the Pontryagin duality theorem for finite abelian groups (cf. [28], Sec. III.2, Cor. 2, p. 118), and by linear extension we get a Hopf algebra isomorphism $H \to H^*$. Since comodules over $H$ are, in the finite-dimensional case, the same as modules over $H^*$, we see modules as well as comodules over $H = K[\mathbb{Z}_p]$ are determined by a single linear map of order $p$. If a vector space $V$ is simultaneously a module and a comodule, the Yetter-Drinfel'd (resp. dimodule) condition says precisely that these mappings have to commute. We have therefore proved the following proposition:

**Proposition** Suppose $p$ is a natural number. Suppose that $K$ contains a primitive $p$-th root of unity $\zeta$ and that $\gamma : \mathbb{Z}_p \to K^\times$ is the group homomorphism that maps 1 to $\zeta$. For a Yetter-Drinfel'd module $V$ over $H$, the mappings

$$\phi : V \to V, v \mapsto (c_1 \to v) \qquad \psi : V \to V, v \mapsto \gamma(v^{(1)})v^{(2)}$$

are commuting $K$-linear maps of order $p$. Moreover, this sets up a one-to-one correspondence between Yetter-Drinfel'd module structures over $H$ on $V$ and pairs of commuting endomorphisms of $V$ of order $p$.

**1.11** This way to express Yetter-Drinfel'd structures via pairs of commuting endomorphisms is also useful to describe algebras, coalgebras, and Hopf algebras in this category. A Yetter-Drinfel'd algebra over a cyclic group of order $p$ is an ordinary algebra such that the multiplication and the unit morphism are linear and colinear. In terms of the mappings $\phi$ and $\psi$ considered in the previous paragraph, this is equivalent to the requirement that $\phi$ and $\psi$ are algebra homomorphisms. Similarly, a Yetter-Drinfel'd coalgebra is an ordinary coalgebra together with two commuting coalgebra homomorphisms of order $p$. However, the bialgebra condition is more involved. First observe that, for a Yetter-Drinfel'd algebra $A$, the algebra structure of the tensor square $A \otimes A$ formed in the category of Yetter-Drinfel'd modules is given by the formula:

$$(a \otimes b)(a' \otimes b') := \frac{1}{p} \sum_{i,j=0}^{p-1} \zeta^{-ij} a \phi^i(a') \otimes \psi^j(b) b'$$

where $\zeta$ is the primitive $p$-th root of unity used in the definition of $\psi$. Note that this expression makes sense since the characteristic of the base field does not divide $p$. The bialgebra axiom now says that the comultiplication is an algebra homomorphism, where $A \otimes A$ is endowed with this algebra structure. An antipode for a Yetter-Drinfel'd bialgebra is a $K$-linear map that satisfies the ordinary antipode axioms and commutes with $\phi$ and $\psi$.

Note that the above algebra structure on the tensor square coincides with the usual tensor product algebra if $\phi$ or $\psi$ are equal to the identity. If $p$ is prime, this is the only possibility for $A$ to be trivial:

**Proposition** A Yetter-Drinfel'd Hopf algebra over a group of prime order is trivial if and only if the action or the coaction is trivial.

**Proof.** It is obvious that $A$ is trivial if the action or the coaction is trivial. Now suppose that $A$ is trivial. If the coaction is nontrivial, there is a nonzero element $a \in A$ which is homogeneous of a degree different from zero, i. e., we have $\delta_A(a) = c_i \otimes a$ for some nonzero $i \in \mathbb{Z}_p$. We then have for all $b \in A$:

$$b \otimes a = \sigma_{A,A}(a \otimes b) = \phi^i(b) \otimes a$$

This implies $\phi^i(b) = b$, and therefore also $\phi(b) = b$. Therefore, the action is trivial. $\Box$

**1.12** We shall use character theory for Hopf algebras. The theory of characters was invented by G. Frobenius (cf. [11]); it was first applied to Hopf algebras by R. G. Larson (cf. [36]). More detailed expositions of the results we need can be found in [71], [41], or [80].

If $H$ is any semisimple Hopf algebra and $\rho : H \to \text{End}(V)$ is a representation, we define its character to be the composition of $\rho$ and the trace map on $\text{End}(V)$:

$$\chi_V(h) := \text{Tr}(\rho(h))$$

The degree of $\chi_V$ is, by definition, the dimension of $V$.

Characters are compatible with the operations on modules by the following rules:

$$\chi_{V \oplus W} = \chi_V + \chi_W \quad \chi_{V \otimes W} = \chi_V \chi_W \quad \chi_{V^*} = S_H^*(\chi_V) \quad \chi_K = \epsilon_H$$

Isomorphic modules have the same characters. Over fields of characteristic zero, two modules are isomorphic if and only if their characters coincide. According to the above rules, the vector space spanned by the characters is actually a subalgebra of $H^*$, which is called the character ring and is denoted by $\text{Ch}(H)$. There is a unique linear map $\bar{\phantom{x}} : \text{Ch}(H) \to \text{Ch}(H)$ satisfying $\bar{\chi}_V = \chi_{V^*}$. If $V_1, \ldots, V_k$ is a system of representatives for the isomorphism types of simple $H$-modules, then, using the notation $\chi_i := \chi_{V_i}$, $\chi_1, \ldots, \chi_k$ is a basis of $\mathrm{Ch}(H)$. If $\Lambda_H$ is an integral of $H$ satisfying $\epsilon_H(\Lambda_H) = 1_K$, we can introduce the bilinear form:

$$\langle \cdot, \cdot \rangle_* : \mathrm{Ch}(H) \times \mathrm{Ch}(H) \to K, (\chi, \chi') \mapsto (\chi \overline{\chi'})(\Lambda_H)$$

We then have the important orthogonality relations:

$$\langle \chi_i, \chi_j \rangle_* = \delta_{ij}$$

It is easy to see that the adjoint of left multiplication by $\chi$ with respect to this bilinear form is left multiplication by $\overline{\chi}$:

$$\langle \chi \chi', \chi'' \rangle_* = \langle \chi', \overline{\chi} \chi'' \rangle_*$$

An important consequence of this fact, which will be used frequently in the sequel, is the observation that a one-dimensional $H$-module $W$ is isomorphic to a submodule of $V^* \otimes V$ if and only if $V \cong V \otimes W$ (cf. [62], Thm. 9, p. 303, [41], Sec. 3.1, p. 489, [80], Sec. 3.6).

**1.13** We shall need some basic facts from the cohomology theory of groups. It will be convenient to use its nonabelian variant (cf. [75], Appendix to Chap. VII, see also [5], [76], Chap. I, § 5).

Suppose that $G$ is a group. A (nonabelian) $G$-module is another group $M$ on which $G$ acts via group automorphisms. A 0-cocycle is, by definition, the same as a fixed point of $G$; the set of fixed points is denoted by $H^0(G, M)$ in this context and is called the zeroth cohomology group of $G$ with values in $M$. A 1-cocycle is a function $s : G \to M$ that satisfies:

$$s(g_1 g_2) = s(g_1)(g_1.s(g_2))$$

Note that, if the action of $G$ on $M$ is trivial, a 1-cocycle is the same as a group homomorphism. Two 1-cocycles $s$ and $s'$ are called cohomologous if there exists an element $m \in M$, called a 0-cochain in this context, such that $s'(g) = m^{-1} s(g) g.m$. This defines an equivalence relation on the set $Z^1(G, M)$ of 1-cocycles; the corresponding quotient set is called the first cohomology set of $G$ with values in $M$, denoted by $H^1(G, M)$. It is a pointed set, the distinguished element being the cohomology class of the 1-cocycle which is constantly equal to one. 1-cocycles that represent this cohomology class are called 1-coboundaries. If $M$ is abelian, then the set of 1-cocycles is a group with respect to the pointwise operations, and $H^1(G, M)$ is a factor group of this group.

We define 2-cocycles only if $M$ is abelian. A 2-cocycle is, by definition, a function $q : G \times G \to M$ that satisfies:

$$(g_1.q(g_2, g_3))q(g_1, g_2 g_3) = q(g_1 g_2, g_3)q(g_1, g_2)$$

2-cocycles form a group with respect to the pointwise operations. Two 2-cocycles $q$ and $q'$ are called cohomologous if there is a function $s : G \to M$, called a 1-cochain in this context, such that

$$q'(g_1, g_2) = q(g_1, g_2)(g_1.s(g_2))s(g_1 g_2)^{-1} s(g_1)$$

This defines an equivalence relation on the set $Z^2(G; M)$ of 2-cocycles. The corresponding quotient set is in fact a factor group; it is called the second cohomology group of $G$ with values in $M$ and is denoted by $H^2(G, M)$. A 2-cocycle representing the unit of this group is called a 2-coboundary.

In the sequel, we shall need some very basic consequences of these definitions that we summarize in the following lemma:

**Lemma**
1. Suppose that $s : G \to M$ is a 1-cocycle. Then we have:

$$s(1) = 1 \quad \text{and} \quad s(g^{-1}) = (g^{-1}.s(g))^{-1}$$

2. Suppose that $M$ is abelian and that $q : G \times G \to M$ is a 2-cocycle that satisfies $q(1, 1) = 1$. Then we have:

$$q(g, 1) = 1 = q(1, g) \quad \text{and} \quad q(g^{-1}, g) = g^{-1}.q(g, g^{-1})$$

**Proof.** We have $s(1) = s(1 \cdot 1) = s(1)(1.s(1))$, and therefore $s(1) = 1$. This implies $1 = s(g^{-1}g) = s(g^{-1})(g^{-1}.s(g))$. To prove the second statement, put $g_2 = 1$ in the definition of a 2-cocycle. We then have $(g_1.q(1, g_3))q(g_1, g_3) = q(g_1, g_3)q(g_1, 1)$. Since $M$ is commutative, this implies $g_1.q(1, g_3) = q(g_1, 1)$. For $g_1 = 1$, we get $q(1, g_3) = q(1, 1) = 1$, and therefore we also have $q(g_1, 1) = 1$.

In the definition of a 2-cocycle, we now put $g_2 = g$, $g_1 = g_3 = g^{-1}$. We then have:
$$(g^{-1}.q(g, g^{-1}))q(g^{-1}, 1) = q(1, g^{-1})q(g^{-1}, g)$$

This implies the assertion. □

A 2-cocycle that satisfies the condition $q(1, 1) = 1$ is called normalized.

If $M$ and $N$ are $G$-modules, we can introduce their tensor product $M \otimes N$ (cf. [25], Def. 25.8, p. 648). It is always abelian and is a $G$-module with respect to the diagonal module structure. If $s \in Z^1(G, M)$ and $t \in Z^1(G, N)$ are two 1-cocycles, we can define their cup product $s \cup t$ (cf. [6], Chap. V, § 3, p.110, [43], Chap. VIII, § 9, p. 244). This is a 2-cocycle of $G$ with values in the tensor product $M \otimes N$ defined by:

$$(s \cup t)(g, g') := s(g) \otimes g.t(g')$$

It is easy to see that an equivariant group homomorphism $f : M \to N$ between two $G$-modules induces mappings

$$f_0 : H^0(G, M) \to H^0(G, N) \quad f_1 : H^1(G, M) \to H^1(G, N)$$
$$f_2 : H^2(G, M) \to H^2(G, N)$$

between the various cohomology sets. Now suppose that $f : G \to G'$ is a group homomorphism. Then any $G'$-module becomes a $G$-module by pullback via $f$. $f$ induces mappings

$$f^1 : H^1(G', M) \to H^1(G, M) \qquad f^2 : H^1(G', M) \to H^2(G, M)$$

given by:

$$f^1(s)(g) := s(f(g)) \qquad f^2(q)(g_1, g_2) := q(f(g_1), f(g_2))$$

Since every fixed point of $G'$ is also a fixed point of $G$, we also have a homomorphism

$$f^0 : H^0(G', M) \to H^0(G, M)$$

which is just the inclusion map.

From the (nonabelian) cohomology theory of groups sketched above, we shall need two results that we state without proof. The first one is concerned with an exact sequence of $G$-modules:

**Proposition 1** Suppose that

$$M \xrightarrow{\iota} N \xrightarrow{\pi} P$$

is an exact sequence of $G$-modules, where $\iota$ and $\pi$ are equivariant group homomorphisms. Assume that $\iota(M)$ is contained in the center of $N$. Then there exist pointed mappings

$$\delta_0 : H^0(G, P) \to H^1(G, M) \qquad \delta_1 : H^1(G, P) \to H^2(G, M)$$

called connecting homomorphisms, such that the sequence

$$\{1\} \to H^0(G, M) \xrightarrow{\iota_0} H^0(G, N) \xrightarrow{\pi_0} H^0(G, P) \xrightarrow{\delta_0} H^1(G, M) \xrightarrow{\iota_1} H^1(G, N)$$
$$\xrightarrow{\pi_1} H^1(G, P) \xrightarrow{\delta_1} H^2(G, M)$$

is an exact sequence of pointed sets.

**Proof.** Cf. [75], Appendix to Chap. VII, Prop. 2, p. 125. □

The second result that we need is concerned with trivial modules over cyclic groups.

**Proposition 2** Suppose that $p$ is a natural number and that $M$ is a trivial abelian $\mathbb{Z}_p$-module. Then we have:

1. A 2-cocycle $q \in Z^2(\mathbb{Z}_p, M)$ satisfies $q(i,j) = q(j,i)$ for all $i, j \in \mathbb{Z}_p$.

2. $H^2(\mathbb{Z}_p, M) \cong M/pM$

3. If $i \in \mathbb{Z}_p$ and $f : \mathbb{Z}_p \to \mathbb{Z}_p, j \mapsto ij$ is multiplication by $i$, the induced endomorphism of the second cohomology group is also multiplication by $i$:

$$f^2 : H^2(\mathbb{Z}_p, M) \to H^2(\mathbb{Z}_p, M), q \mapsto iq$$

**Proof.** The last two assertions are proved in much greater generality in [43], Chap. IV, Thm. 7.1, p. 122, resp. Exerc. 7.6, p. 124 and in [22], Chap. VI, Prop. 7.1, p. 201, resp. Exerc. 7.4, p. 201. To prove the first assertion, observe that, by replacing $q$ with $q - q(0,0)$, we can assume that $q$ is normalized. By the standard description of extensions via cocycles (cf. [6], Chap. IV, § 3, p. 91, [43], Chap. IV, § 4, p. 111), we can define a group structure on the set $M \times \mathbb{Z}_p$ by:

$$(x,i)(y,j) := (x + y + q(i,j), i + j)$$

where, according to the above lemma, the unit element is $(0,0)$ and the inverse of $(x,i)$ is $(-x - q(i,-i), -i)$. The group homomorphism

$$M \to M \times \mathbb{Z}_p, x \mapsto (x,0)$$

identifies $M$ with a central subgroup of $M \times \mathbb{Z}_p$ whose corresponding quotient is isomorphic to $\mathbb{Z}_p$. Since a group that contains a central subgroup with cyclic quotient is abelian (cf. [25], Kap. III, Hilfssatz 7.1, p. 300, [34], Satz 1.19, p. 20, [82], Chap. 1, (2.26), p. 17), we see that $M \times \mathbb{Z}_p$ is abelian. Obviously, this implies that $q(i,j) = q(j,i)$. □

In a similar situation, an explicit description of cocycles that constitute a complete system of representatives of the cohomology classes can be found in [28], Sec. I.15, p. 80.

# 2 Clifford theory

**2.1** In this section, we assume that $p$ is a prime and that $K$ is an algebraically closed field whose characteristic is different from $p$. We denote by $H := K[\mathbb{Z}_p]$ the group ring of the cyclic group $\mathbb{Z}_p$ of order $p$. As in Paragraph 1.10, we denote by $C$ the group of grouplike elements of $H$, $C := G(H)$. These are precisely the canonical basis elements of the group ring, i. e., using the notation of Paragraph 1.10, we have $C = \{c_i \mid i \in \mathbb{Z}_p\}$. $A$ denotes a semisimple Yetter-Drinfel'd Hopf algebra over the group ring $H = K[\mathbb{Z}_p]$. The set of centrally primitive idempotents of $A$ is denoted by $E$. $\hat{C} \cong \hat{\mathbb{Z}}_p$ denotes the the character group of $C \cong \mathbb{Z}_p$, i. e., the set of all group homomorphisms of $C$ to the multiplicative group $K^\times$. By linear extension, we can think of elements of $\hat{C}$ as elements of the dual Hopf algebra $H^*$. As explained in Paragraph 1.10, $\hat{C}$ is a basis of $H^*$. We denote the Radford biproduct $A \otimes H$ by $B$.

The goal of this section is to understand the modules of the Radford biproduct. To achieve this, we apply Clifford theory for group-graded rings (cf. [13], [14], [15], [12], §11C) in our very special situation of Radford biproducts over groups of prime order. However, in our situation, the Hopf algebra structure of the Radford biproduct gives rise to additional properties: The character group $\hat{C}$ can be regarded as a subset of the set of irreducible characters of $B$, and therefore we get a left and a right action of $\hat{C}$ on the set of irreducible characters. We will establish a linkage principle that will yield a one-to-one correspondence between the orbits of the right action of $\hat{C}$ in the set of irreducible characters of $B$ and the orbits of centrally primitive idempotents of $A$ under the action of $C$. The results of this section can also be understood from the point of view of the theory of prime ideals in crossed products developed by M. Lorenz and D. S. Passman (cf. [40], [64], Chap. 4). We shall indicate the connections with this theory where they occur.

Throughout this section, we will constantly use the convention that indices take values between 0 and $p-1$ and are reduced modulo $p$ if they do not lie within this range. In notation, we shall not distinguish between an integer $i \in \mathbb{Z}$ and its equivalence class in $\mathbb{Z}_p := \mathbb{Z}/p\mathbb{Z}$.

**2.2** We begin by associating to every simple $B$-module a subset of the set $E$ of centrally primitive idempotents of $A$. This subset essentially characterizes which simple $A$-modules occur in the restriction of the $B$-module to $A$, since the centrally primitive idempotents in $E$ are in one-to-one correspondence with simple $A$-modules. We shall see below that a simple $A$-module occurs in a simple $B$-module at most with multiplicity one; therefore the knowledge of the simple modules occurring in the restriction already characterizes the restriction up to isomorphism. It will turn out that this set is related to the isotropy group of the character of the irreducible $B$-module under the canonical right action of $\hat{C}$:

**Definition** Suppose that $V$ is a simple $B$-module with character $\chi_V$.

1. Define:
$$\kappa(V) := \{e \in E \mid \exists v \in V \setminus \{0\} : (e \otimes 1_H)v = v\}$$

2. Define: $\kappa^*(V)$ to be the isotropy group of the canonical right action of $\hat{C}$ on the set of irreducible characters:

$$\kappa^*(V) := \{\gamma \in \hat{C} \mid \chi_V(\epsilon_A \otimes \gamma) = \chi_V\}$$

The first observation is that $\kappa(V)$ is stable under the action of $C$ on $E$. Even stronger, we have:

**Proposition** Suppose that $V$ is a simple $B$-module. Then $\kappa(V)$ is an orbit of $E$ with respect to the action of $C$.

**Proof.** If $e \in \kappa(V)$, the definition of $\kappa(V)$ implies that we can choose a nonzero vector $v \in V$ such that $(e \otimes 1_H)v = v$. We then have:

$$((c \to e) \otimes 1_H)(1_A \otimes c)v = (1_A \otimes c)(e \otimes 1_H)v = (1_A \otimes c)v$$

This implies that $c \to e \in \kappa(V)$, and therefore $\kappa(V)$ is a disjoint union of orbits of $C$. Suppose that $O$ is one of these orbits and consider for every element $e \in O$ the space
$$V_e := \{v \in V \mid (e \otimes 1_H)v = v\}$$

which is nonzero by definition. $V_e$ is obviously an $A$-submodule of $V$, and a calculation similar to the one above shows that

$$V' := \sum_{e \in O} V_e$$

is invariant with respect to $C$. This means that $V'$ is a $B$-submodule of $V$, which by simplicity implies that $V' = V$. Therefore, we have $\kappa(V) = O$. $\square$

The above proposition is a special case of [13], Exerc. 18.10, p. 288, where mostly the corresponding simple modules are emphasized over centrally primitive idempotents (cf. also [12], Prop. (11.16), p. 273). It should also be compared with [64], Lem. 14.1, p. 132. The simple module $V$ used here corresponds to a prime ideal denoted there by $P$, and $\kappa(V)$ corresponds to the ideal that is denoted there by $P \cap R$.

**2.3** The orbit $\kappa(V)$ defined in the previous paragraph obviously can consist of one or of $p$ elements. This easy observation is essential for the understanding of the simple $B$-modules; it divides these modules into two classes that behave in a fundamentally different way:

**Definition** Suppose that $e \in E$ is a centrally primitive idempotent of $A$. Denote by $C_e$ the isotropy group of $e$, i. e., $C_e := \{c_i \in C \mid c_i \to e = e\}$

1. $e$ is called purely unstable if $C_e = \{1\}$.
2. $e$ is called stable if $C_e = C$.

A simple $A$-module, or its character, will be called purely unstable, resp. stable, if the corresponding centrally primitive idempotent has this property. If $O \subset E$ is an orbit of $C$, it will be called purely unstable, resp. stable, if its elements have this property. A simple $B$-module will be called purely unstable, resp. stable if $\kappa(V)$ is purely unstable, resp. stable.

The above definition should be compared with [12], p. 269, [10], Def. 2.7, p. 356, [73], Def. 3.7, p. 278, and p. 286.

Now suppose that $O_1, \dots, O_m$ are the orbits of $E$ with respect to the action of $C$. Then it is easy to see that

$$B = \bigoplus_{i=1}^{m} AO_i \otimes H$$

is a decomposition of $B$ into two-sided ideals. However, these ideals are not always simple. But they are simple if the orbit is purely unstable:

**Proposition 1** Suppose that $O$ is a purely unstable orbit of $C$ in $E$. Then $AO \otimes H$ is a simple ideal of $B$.

**Proof.** (1) Denote by $K^{\mathbb{Z}_p}$ the algebra of functions from $\mathbb{Z}_p$ to the base field $K$, with pointwise addition and multiplication. We denote the canonical basis vector corresponding to $i \in \mathbb{Z}_p$ by $e_i$, i. e., $e_i$ is the function taking the value one on $i$ and the value zero on all other elements of $\mathbb{Z}_p$. We get a left action of $H$ on $K^{\mathbb{Z}_p}$ if we define:

$$c_i \to e_j := e_{i+j}$$

This turns $K^{\mathbb{Z}_p}$ into an $H$-module algebra. The map

$$K^{\mathbb{Z}_p} \otimes H \to \operatorname{End}_K(K^{\mathbb{Z}_p}), q \otimes c_i \mapsto (q' \mapsto q(c_i \to q'))$$

is an algebra isomorphism (cf. [57], Cor. 9.4.3, p. 162, note that the setup there is slightly different), since it maps the elements $e_i \otimes c_{i-j}$ to the matrix units. Therefore, the smash product $K^{\mathbb{Z}_p} \otimes H$ is simple.

(2) For all $e \in O$, the tensor product algebra $Ae \otimes K^{\mathbb{Z}_p}$ is therefore an $H$-module algebra if it is endowed with the action on the second tensor factor. It can be verified directly that the bijection

$$f : Ae \otimes K^{\mathbb{Z}_p} \to AO, \sum_{i=0}^{p-1} a_i \otimes e_i \mapsto \sum_{i=0}^{p-1}(c_i \to a_i)$$

is an isomorphism of $H$-module algebras. This implies that

$$f \otimes \mathrm{id}_H : Ae \otimes (K^{\mathbb{Z}_p} \otimes H) \to AO \otimes H$$

is an algebra isomorphism between the corresponding smash products. Since the left hand side is simple as the tensor product of two simple algebras, the right hand side is simple, too. □

Although, for stable orbits $O$, the ideals $AO \otimes H$ are not simple themselves, we can precisely determine their decomposition into simple ideals:

**Proposition 2** Suppose that $O = \{e\}$ is a stable orbit of $C$. Denote the corresponding two-sided ideal of $A$ by $I := Ae$ and suppose that $W$ is a simple $I$-module with corresponding representation $\rho_W : I \to \mathrm{End}_K(W)$. Then there exists a group homomorphism

$$u_e : C \to U(I)$$

from $C$ to the group of units $U(I)$ of $I$ such that, for all $\gamma \in \hat{C}$,

$$I \otimes H \to \mathrm{End}_K(W), a \otimes c_i \mapsto \rho_W(au_e(c_i))\gamma(c_i)$$

is an irreducible representation of $I \otimes H$. Every irreducible representation of $I \otimes H$ is isomorphic to precisely one of these. In particular, $I \otimes H$ is the direct sum of the direct sum of $p$ simple ideals of $B$.

**Proof.** The assumption that $O$ is stable means that $I$ is invariant with respect to the action of $C$ on $A$. By the Skolem-Noether theorem, the sequence

$$K^\times \xrightarrow{\iota} U(I) \xrightarrow{\mathrm{ad}} \mathrm{Aut}(I)$$

is exact, where $\iota(\lambda) = \lambda e$ and $\mathrm{ad}(a)(a') = aa'a^{-1}$. We regard these groups as trivial (nonabelian) $C$-modules. From Proposition 1.13.1, we then get an exact sequence of pointed sets

$$Hom(C, K^\times) \xrightarrow{\iota_*} Hom(C, U(I)) \xrightarrow{\mathrm{ad}_*} Hom(C, \mathrm{Aut}(I))$$

since $H^2(C, K^\times) = \{0\}$ by Proposition 1.13.2. This now means that the action of $C$ on $I$ can be lifted to a group homomorphism $u_e : C \to U(I)$ that satisfies $c_i \to a = u_e(c_i)au_e(c_i^{-1})$. (This can of course also be proved without the use of cohomology.) It then can be verified directly that the map

$$f : I \otimes H \to I \otimes H, a \otimes c_i \mapsto au_e(c_i) \otimes c_i$$

is an algebra isomorphism from the smash product to the ordinary tensor product. Since all irreducible representations of the ordinary tensor product are of the form $\rho_W \otimes \gamma : I \otimes H \to \operatorname{End}_K(W)$ for some $\gamma \in \hat{C}$, the irreducible representations of the smash product are of the form stated above. □

The following statement is an obvious consequence of the above proofs:

**Corollary** $\kappa$ is surjective, i. e., for all $C$-orbits $O$ in $E$ there is a simple $B$-module $V$ such that $\kappa(V) = O$.

The assertions in this paragraph should be compared with [64], Thm. 14.7, p. 138. Further references will be given below.

**2.4** Our next aim is to describe in greater detail those simple $B$-modules which are purely unstable. Recall that the function $\lambda_H : H \to K$ defined by:

$$\lambda_H(c_i) = \begin{cases} 1 & : i = 0 \\ 0 & : i \neq 0 \end{cases}$$

is a normalized integral on $H$ (cf. [57], Example 2.1.2.2, p. 17). If $W$ is an $A$-module, the set $W^p$ becomes a $B$-module by defining:

$$(a \otimes c_j).(w_i)_{i=0,\ldots,p-1} := ((c_i^{-1} \to a)w_{i-j})_{i=0,\ldots,p-1}$$

This module is isomorphic to the induced module $V := B \otimes_A W$ by the map:

$$W^p \to B \otimes_A W, (w_i)_{i=0,\ldots,p-1} \mapsto \sum_{i=0}^{p-1} (1_A \otimes c_i) \otimes_A w_i$$

(cf. [64], Lem. 3.3, p. 21).

**Proposition** Suppose that $V$ is a purely unstable simple $B$-module. Denote the character of $V$ by $\chi_V$. For any $e \in O := \kappa(V)$, choose a minimal left ideal $W_e$ of the two-sided ideal $Ae \subset A$, and denote the character of the $A$-module $W_e$ by $\eta_e$. Then we have:

1. The induced $B$-modules $B \otimes_A W_e$ are isomorphic to $V$.

2. The character of $V$ is given by the formula

$$\chi_V = \eta_O \otimes \lambda_H$$

   where $\eta_O := \sum_{e \in O} \eta_e$.

3. The restriction of $V$ to $A$ is isomorphic to $\sum_{e \in O} W_e$.

4. $\dim V = p \dim W_e$

**Proof.** Suppose that $e \in O$. It is clear from the description of the induced module $B \otimes_A W_e$ above that every centrally primitive idempotent $e' \in E \setminus O$ vanishes on $B \otimes_A W_e$. Therefore, $B \otimes_A W_e$ is a power of the unique simple module corresponding to the simple ideal $AO \otimes H$. Since $\dim AO \otimes H = p^2 \dim Ae = (\dim W_e^p)^2$, the induced module itself must be simple. Since $V$ is a simple $B$-module corresponding to the same ideal, we have $V \cong B \otimes_A W_e$. This proves the first statement. The third and the fourth statement follows from the description of the induced module given above.

For the second statement, suppose that $c_i \in C$ is not the unit element. Then the elements of $B$ of the form $a \otimes c_i$ permute the subspaces in the decomposition $V \cong \sum_{e \in O} W_e$ cyclicly. Therefore, we have:

$$\chi_V(a \otimes c_i) = 0 = \eta_O(a) \lambda_H(c_i)$$

This decomposition also implies:

$$\chi_V(a \otimes 1_H) = \eta_O(a) = \eta_O(a) \lambda_H(1_H)$$

Since $A \otimes H$ is generated by elements of the form $a \otimes c_i$, this implies the second assertion. $\square$

The above proposition is essentially a special case of [12], Prop. (11.16), p. 273 (cf. also [13], Exerc. 18.11, p. 288). A generalization of these considerations to arbitrary Hopf Galois extensions can be found in [73], Thm. 5.5, p. 288.

**2.5** The description of the characters of the purely unstable $B$-modules given in the previous paragraph is rather explicit. The analogous description for characters of stable $B$-modules is slightly less explicit, since it involves the group homomorphism from Proposition 2.3.2:

**Proposition** Suppose that $V$ is a stable simple $B$-module, i. e., $\kappa(V) = \{e\}$ has length one. Denote the character of $V$ by $\chi_V$ and the character of the restricted $A$-module $V$ by $\eta_V$.

1. The restriction of $V$ to $A$ is simple.

2. Denote the group of units of $I := Ae$ by $U(I)$. There exists a group homomorphism $u_e : C \to U(I)$, depending on $e$ but not on $V$, and a character $\gamma \in \hat{C}$ such that the character of $V$ is given by the equation:

$$\chi_V(a \otimes c_i) = \eta_V(au_e(c_i))\gamma(c_i)$$

   for all $a \in A$ and all $i = 0, \ldots, p - 1$.

3. The isotropy group $\kappa^*(V)$ is trivial.

**Proof.** If $O'$ is an orbit distinct from $\kappa(V)$, the two-sided ideal $AO' \otimes H$ annihilates $V$. Therefore, $V$ corresponds to an irreducible representation of the ideal $I \otimes H$. These representations were described in Proposition 2.3.2, and the first two assertions follow directly from this description. To prove the third assertion, assume on the contrary that the isotropy group $\kappa^*(V)$ contains nontrivial element $\gamma' \neq \epsilon_H$. Choose an element $c \in C$ that is distinct from the unit. Since $c$ generates $C$, we have $\gamma'(c) \neq 1$. We then have, for all $a \in A$, that

$$\chi_V(a \otimes c) = (\chi_V(\epsilon_A \otimes \gamma'))(a \otimes c) = \chi_V(a \otimes c)\gamma'(c)$$

and therefore

$$0 = \chi_V(a \otimes c) = \eta_V(au_e(c))\gamma(c)$$

for some $\gamma \in \hat{C}$. We therefore have to prove that, in contrast to the case of purely unstable modules, there exists an $a \in A$ such that $\eta_V(au_e(c)) \neq 0$. Denote the representation arising from the restriction of $V$ to $A$ by $\rho: A \to \mathrm{End}_K(V)$. Since $\rho(u_e(c))^p = \mathrm{id}_V$, the minimum polynomial of $\rho(u_e(c))$ divides the separable polynomial $t^p - 1$, and therefore has distinct roots. Therefore, $\rho(u_e(c))$ is diagonalizable, and we can choose a basis $v_1, \ldots, v_n$ of $V$ such that

$$\rho(u_e(c))(v_j) = \zeta_j v_j$$

for some $p$-th roots of unity $\zeta_1, \ldots, \zeta_n$. By the first part of the proposition, $\rho$ is surjective, and therefore there are elements $a_1, \ldots, a_n \in A$ such that:

$$\rho(a_i)(v_j) = \delta_{ij} v_j$$

We then have: $\eta_V(a_i u_e(c)) = \zeta_i \neq 0$ $\square$

Note that, since the isotropy group $\kappa^*(V)$ is trivial, the orbit of $\chi_V$ under the right action of $\hat{C}$ consists of $p$ characters of nonisomorphic simple $B$-modules whose restrictions to $A$ are isomorphic.

**2.6** It now turns out that there is a connection between the action of $C$ on the set $E$ of centrally primitive idempotents of $A$ and the action of $\hat{C}$ on the irreducible characters of $B$: Two simple $B$-modules have the same restriction to $A$ if and only if their characters are linked via the canonical right action of $\hat{C}$ on $\mathrm{Ch}(B)$:

**Theorem** Suppose that $V$ and $V'$ are simple $B$-modules with characters $\chi_V$ resp. $\chi_{V'}$. Then the following assertions are equivalent:

1. $\kappa(V) = \kappa(V')$

2. The restrictions of $V$ and $V'$ to $A$ are isomorphic.

3. There exists $\gamma \in \hat{C}$ such that $\chi_{V'} = \chi_V(\epsilon_A \otimes \gamma)$.

**Proof.** We distinguish two cases: Suppose first that $V$ is purely unstable. Then the restriction of $V$ to $A$ is a direct sum of $p$ nonisomorphic simple $A$-modules, each occurring with multiplicity one, that correspond to the elements of $\kappa(V)$. Therefore, the first two assertions are equivalent. But since according to Proposition 2.4, the $B$-module $V$ can be recovered from every simple $A$-submodule via induction, $V$ and $V'$ are isomorphic as $A$-modules if and only if they are isomorphic as $B$-modules, i. e., if $\chi_V = \chi_{V'}$. But since, by Proposition 2.4, $\chi_V$ has the form $\chi_V = \eta_V \otimes \lambda_H$, $\chi_V$ is invariant under the action of $\hat{C}$. Therefore, all three assertions are equivalent in this case.

Suppose now that $V$ is stable. Then the restriction of $V$ to $A$ is still simple, and $\kappa(V)$ consists precisely out of the corresponding centrally primitive idempotent of $A$. Therefore, the first two assertions are equivalent. According to Proposition 2.3.2, there are, up to isomorphism, $p$ simple $B$-modules that have the same restriction to $A$ as $V$. For $\gamma \in \hat{C}$, $\chi_V(\epsilon_A \otimes \gamma)$ is an irreducible character of $B$ that has the same restriction to $A$ as $\chi_V$. According to Proposition 2.5, these characters are all distinct, and therefore exhaust all irreducible characters of $B$ that have the same restriction to $A$ as $\chi_V$. Therefore, the third assertion is equivalent to the other two. □

The following statement is a byproduct of the above proof:

**Corollary** Suppose that $V$ is a simple $B$-module. Then we have: $\text{card}(\kappa(V)) = \text{card}(\kappa^*(V))$

The above theorem is a special case of a result of W. Chin (cf. [8], Thm. 2.1, p. 90, see also [73], Sec. 4 and [29], Bem. 1.14, p. 37) that is concerned with arbitrary smash products. Arbitrary smash products are, in contrast to Radford biproducts, not endowed with a Hopf algebra structure; nevertheless, an analogue of the right action of $\hat{C}$ can be defined also in this more general case.

**2.7** Besides the canonical right action of $\hat{C}$ on $\text{Ch}(B)$, there is of course also a canonical left action of $\hat{C}$ on $\text{Ch}(B)$, given by $\chi \mapsto (\epsilon_A \otimes \gamma)\chi$ for $\gamma \in \hat{C}$. This action takes characters of simple modules to characters of simple modules and therefore induces an action on the isomorphism classes of simple $B$-modules. We shall study this left action in conjunction with the following left action of $\hat{C}$ on orbits in $E$:

**Definition** Suppose that $\gamma \in \hat{C}$. Consider the algebra automorphism

$$\psi_\gamma : A \to A, a \mapsto \gamma(a^{(1)})a^{(2)}$$

(cf. [81], Prop. 2.6, p. 6). Denote the set of all orbits of the action of $C$ in $E$ by $\mathcal{O}$. We introduce the following left action of $\hat{C}$ on $\mathcal{O}$:

$$\hat{C} \times \mathcal{O} \to \mathcal{O}, (\gamma, O) \mapsto \gamma.O := \psi_{\gamma^{-1}}(O)$$

This is well defined because the algebra homomorphism $\psi_\gamma$ preserves $E$ and also preserves orbits since it commutes with the action of $C$ by Proposition 1.10. The fundamental property of $\kappa$ with respect to these left actions now is:

**Proposition** $\kappa$ is equivariant.

**Proof.** Suppose that $V$ is a simple $B$-module and that $\gamma$ is an element of $\hat{C}$. We can endow the base field $K$ with a $B$-module structure via the character $\epsilon_A \otimes \gamma$; we denote $K$ by $K_{\epsilon_A \otimes \gamma}$ if we want to emphasize that it is endowed with this module structure. $B$ acts on the module $\gamma.V = K_{\epsilon_A \otimes \gamma} \otimes V$ via:

$$(a \otimes h)(1_K \otimes v) = \epsilon_A(a_{(1)})\gamma(a_{(2)}{}^{(1)}h_{(1)}) \otimes (a_{(2)}{}^{(2)} \otimes h_{(2)})v$$
$$= 1_K \otimes (\psi_\gamma(a) \otimes \gamma(h_{(1)})h_{(2)})v$$

If $e \in \kappa(V)$, choose a nonzero element $v \in V$ such that $(e \otimes 1_H)v = v$. We then have:
$$(\psi_{\gamma^{-1}}(e) \otimes 1_H)(1_K \otimes v) = (1_K \otimes v)$$

and therefore $\psi_{\gamma^{-1}}(e) \in \kappa(K_{\epsilon_A \otimes \gamma} \otimes V)$ and therefore we see that $\gamma.\kappa(V) \subset \kappa(\gamma.V)$. By replacing $V$ by $\gamma^{-1}.V$, we see that the other inclusion holds, too. $\square$

# 3 Examples

**3.1** In this section, we construct examples of nontrivial cocommutative Yetter-Drinfel'd Hopf algebras. We shall see in the Section 7 that Yetter-Drinfel'd Hopf algebras of a certain type are necessarily of the form described here. The examples that we are going to construct have as underlying vector space the tensor product of a group ring and a dual group ring. Their coalgebra structure is the tensor product coalgebra structure, whereas the algebra structure is that of a crossed product. The construction of the examples proceeds in three stages: First, we present a general abstract framework that describes circumstances under which the above tensor product becomes a Yetter-Drinfel'd Hopf algebra. In this framework, the canonical basis elements of the tensor product are eigenvectors with respect to the action and homogeneous with respect to the coaction. This framework is a slight generalization of a framework created by N. Andruskiewitsch (cf. [1]) to understand a previous version of the second stage. In this second stage, we attach a Yetter-Drinfel'd Hopf algebra to a finite ring and a group homomorphism from an arbitrary finite group to the group of units of this ring. The Yetter-Drinfel'd Hopf algebra that is constructed is the tensor product of the dual group ring of the underlying additive group of the ring and the group ring of the above group; it is defined over the group ring of the additive group of the ring. The construction depends on additional data, which play the role of parameters that can be chosen freely.

In the third stage, the finite field of $p$ elements is chosen for the finite ring above, where $p$ is an odd prime. It will be shown in Section 7 that, under some restrictions on the base field, all cocommutative cosemisimple Yetter-Hopf algebras over the group ring of the group with $p$ elements are of this form. As will be seen at the end of the section, the case $p = 2$ can also be treated within this framework.

In this section, $K$ denotes a field that is not required to satisfy any restrictions except for Paragraph 3.5 and Paragraph 3.6, where we assume that $K$ contains a primitive fourth root of unity.

**3.2** We now describe the above mentioned framework. Suppose that $C$, $G$, and $P$ are finite groups. Suppose that $P$ is a (nonabelian) left $G$-module, i. e., that $G$ acts on $P$ via group automorphisms. We define $H := K[C]$, the group ring of $C$, and $A := K^P \otimes K[G]$, the tensor product of the dual group ring of $P$ and the group ring of $G$. The canonical basis elements of $K[C]$ resp. $K[G]$ are denoted by $c_b$ resp. $x_s$, where $b \in C$ and $s \in G$. The primitive idempotents of $K^P$ are denoted by $e_u$ for $u \in P$, i. e., $e_u : P \to K$ is the mapping defined by $e_u(v) := \delta_{uv}$.

By linearity, we can extend the action of $G$ on $P$ to an action of $K[G]$ on $K^P$. $K^P$ is therefore a left $K[G]$-module via:

$$x_s.e_u := e_{s.u}$$

for $s \in G$ and $u \in P$. The set of homomorphisms of any group to an abelian group is an abelian group with the pointwise operations. Therefore, the sets $\operatorname{Hom}(P, Z(C))$ and $\operatorname{Hom}(P, \hat{C})$ are abelian groups, where $Z(C)$ resp. $\hat{C} := \operatorname{Hom}(C, K^\times)$ denote the center resp. the character group of $C$. We turn these spaces into left $G$-modules by defining:

$$(s.f)(u) := f(s^{-1}.u)$$

for group homomorphisms $f : P \to Z(C)$ resp. $f : P \to \hat{C}$.

We now suppose that the following structure elements are given:

1. A 1-cocycle $z : G \to \operatorname{Hom}(P, Z(C)), s \mapsto z_s$ of $G$ with values in the $G$-module $\operatorname{Hom}(P, Z(C))$.

2. A 1-cocycle $\gamma : G \to \operatorname{Hom}(P, \hat{C}), s \mapsto \gamma_s$ of $G$ with values in the $G$-module $\operatorname{Hom}(P, \hat{C})$.

3. A normalized 2-cocycle $\sigma : G \times G \to U(K^P)$ of $G$ with values in the group of units $U(K^P)$ of $K^P$, regarded as a $G$-submodule of $K^P$. We write $\sigma$ in the form
$$\sigma(s,t) = \sum_{u \in P} \sigma_u(s,t) e_u$$
for functions $\sigma_u : G \times G \to K^\times$.

We require that these structure elements satisfy the following compatibility condition:
$$\sigma_{uv}(s,t) = ((s.\gamma_t)(u))(z_s(v))\sigma_u(s,t)\sigma_v(s,t)$$

for all $u, v \in P$ and all $s, t \in G$.

In this situation, we can construct a Yetter-Drinfel'd Hopf algebra over $H = K[C]$ as follows:

**Proposition** The vector space $A = K^P \otimes K[G]$ becomes a Yetter-Drinfel'd Hopf algebra over $H$ if it is endowed with the following structures:

1. Crossed product algebra structure:
$$(e_u \otimes x_s)(e_v \otimes x_t) := e_u(x_s.e_v)\sigma(s,t) \otimes x_s x_t = \delta_{u,s.v}\sigma_u(s,t) e_u \otimes x_{st}$$

   Unit: $1_A := \sum_{u \in P} e_u \otimes x_1$

2. Tensor product coalgebra structure:
$$\Delta_A(e_u \otimes x_s) := \sum_{v \in P} (e_v \otimes x_s) \otimes (e_{v^{-1}u} \otimes x_s)$$

Counit: $\epsilon_A(e_u \otimes x_s) := \delta_{u1}$

3. Action: $c_b \to (e_u \otimes x_s) := (\gamma_s(u))(b) e_u \otimes x_s$
4. Coaction: $\delta_A(e_u \otimes x_s) := c_{z_s(u)} \otimes (e_u \otimes x_s)$
5. Antipode: $S_A(e_u \otimes x_s) := \sigma_{u^{-1}}^{-1}(s, s^{-1}) e_{s^{-1}.u^{-1}} \otimes x_{s^{-1}}$

**Proof.** (1) The Yetter-Drinfel'd condition follows from the fact that $z_s(u)$ is contained in the center of $C$. We now prove that $\Delta_A$ is $H$-linear. We have:

$$\begin{aligned}\Delta_A(c_b \to (e_u \otimes x_s)) &= \gamma_s(u)(b) \sum_{v \in P} (e_v \otimes x_s) \otimes (e_{v^{-1}u} \otimes x_s) \\ &= \sum_{v \in P} \gamma_s(v)(b) \gamma_s(v^{-1}u)(b)(e_v \otimes x_s) \otimes (e_{v^{-1}u} \otimes x_s) \\ &= \sum_{v \in P} (c_b \to (e_v \otimes x_s)) \otimes (c_b \to (e_{v^{-1}u} \otimes x_s)) \\ &= c_b \to \Delta_A(e_u \otimes x_s)\end{aligned}$$

By a similar calculation, the equation $z_s(u) = z_s(v) z_s(v^{-1}u)$ implies that $\Delta_A$ is colinear.

(2) From Lemma 1.13, we have that $z_1 = 1 \in \text{Hom}(P, Z(C))$, i. e., we have $z_1(u) = 1_C$ for all $u \in P$. Since $z_s$ is a group homomorphism, we also have $z_s(1_P) = 1_C$ for all $s \in G$. For similar reasons, we have $\gamma_1 = 1 \in \text{Hom}(P, \hat{C})$, i. e., we have $\gamma_1(u) = \epsilon_H$ for all $u \in P$, and also $\gamma_s(1_P) = \epsilon_H$ for all $s \in G$. This implies that $1_A$ is invariant and coinvariant, since we have:

$$c_b \to 1_A = \sum_{u \in P} (\gamma_1(u))(b) e_u \otimes x_1 \qquad \delta_A(1_A) = \sum_{u \in P} c_{z_1(u)} \otimes (e_u \otimes x_1)$$

Similarly, we see that $\epsilon_A$ is $H$-linear and colinear, since we have:

$$\epsilon_A(c_b \to (e_u \otimes x_s)) = (\gamma_s(u))(b) \delta_{u1} \qquad (\text{id}_H \otimes \epsilon_A) \delta_A(e_u \otimes x_s) = \delta_{u1} c_{z_s(u)}$$

(3) We now prove that the multiplication map $\mu_A : A \otimes A \to A$ is $H$-linear and colinear. Since $\gamma$ is a 1-cocycle, we have $\gamma_{st} = \gamma_s(s.\gamma_t)$, i. e., we have $\gamma_{st}(u) = \gamma_s(u) \gamma_t(s^{-1}.u)$ for all $u \in P$. This implies:

$$\begin{aligned}(c_b \to (e_u \otimes x_s))(c_b \to (e_v \otimes x_t)) &= (\gamma_s(u))(b)(\gamma_t(v))(b)(e_u \otimes x_s)(e_v \otimes x_t) \\ &= (\gamma_s(u)\gamma_t(v))(b)\delta_{u, s.v} \sigma_u(s, t)(e_u \otimes x_{st}) \\ &= (\gamma_s(u)\gamma_t(s^{-1}.u))(b)\delta_{u, s.v} \sigma_u(s, t)(e_u \otimes x_{st}) \\ &= (\gamma_{st}(u))(b)\delta_{u, s.v} \sigma_u(s, t)(e_u \otimes x_{st}) \\ &= c_b \to ((e_u \otimes x_s)(e_v \otimes x_t))\end{aligned}$$

Similarly, since $z$ is a 1-cocycle, we have $z_{st}(u) = z_s(u)z_t(s^{-1}.u)$. Therefore we have:

$$\begin{aligned}
\delta_A(e_u \otimes x_s)\delta_A(e_v \otimes x_t) &= c_{z_s(u)}c_{z_t(v)} \otimes (e_u \otimes x_s)(e_v \otimes x_t) \\
&= c_{z_s(u)z_t(v)} \otimes \delta_{u,s.v}\sigma_u(s,t)(e_u \otimes x_{st}) \\
&= c_{z_s(u)z_t(s^{-1}.u)} \otimes \delta_{u,s.v}\sigma_u(s,t)(e_u \otimes x_{st}) \\
&= c_{z_{st}(u)} \otimes \delta_{u,s.v}\sigma_u(s,t)(e_u \otimes x_{st}) \\
&= \delta_A((e_u \otimes x_s)(e_v \otimes x_t))
\end{aligned}$$

(4) We now verify the Yetter-Drinfel'd bialgebra axiom, i. e., the fact that $\Delta_A : A \to A \hat{\otimes} A$ is an algebra homomorphism. This is the step that depends on the compatibility condition:

$$\begin{aligned}
&\Delta_A(e_u \otimes x_s)\Delta_A(e_v \otimes x_t) \\
&= \sum_{u',v' \in P} [(e_{u'} \otimes x_s) \otimes (e_{u'^{-1}u} \otimes x_s)][(e_{v'} \otimes x_t) \otimes (e_{v'^{-1}v} \otimes x_t)] \\
&= \sum_{u',v' \in P} (\gamma_t(v'))(z_s(u'^{-1}u))(e_{u'} \otimes x_s)(e_{v'} \otimes x_t) \otimes (e_{u'^{-1}u} \otimes x_s)(e_{v'^{-1}v} \otimes x_t) \\
&= \sum_{u',v' \in P} (\gamma_t(v'))(z_s(u'^{-1}u))\delta_{u',s.v'}\delta_{u'^{-1}u,s.(v'^{-1}v)}\sigma_{u'}(s,t)\sigma_{u'^{-1}u}(s,t) \\
&\qquad\qquad\qquad\qquad\qquad\qquad\qquad\qquad (e_{u'} \otimes x_{st}) \otimes (e_{u'^{-1}u} \otimes x_{st}) \\
&= \delta_{u,s.v} \sum_{u' \in P} (\gamma_t(s^{-1}.u'))(z_s(u'^{-1}u))\sigma_{u'}(s,t)\sigma_{u'^{-1}u}(s,t) \\
&\qquad\qquad\qquad\qquad\qquad\qquad\qquad\qquad (e_{u'} \otimes x_{st}) \otimes (e_{u'^{-1}u} \otimes x_{st}) \\
&= \delta_{u,s.v} \sum_{u' \in P} \sigma_u(s,t)(e_{u'} \otimes x_{st}) \otimes (e_{u'^{-1}u} \otimes x_{st}) \\
&= \Delta_A((e_u \otimes x_s)(e_v \otimes x_t))
\end{aligned}$$

$\Delta_A$ also preserves the unit.

(5) We now prove that $\epsilon_A$ is an algebra homomorphism. Inserting $u = v = 1$ into the compatibility condition, we get:

$$\sigma_1(s,t) = ((s.\gamma_t)(1))(z_s(1))\sigma_1(s,t)\sigma_1(s,t) = \sigma_1(s,t)\sigma_1(s,t)$$

and therefore $\sigma_1(s,t) = 1$. This implies:

$$\begin{aligned}
\epsilon_A((e_u \otimes x_s)(e_v \otimes x_t)) &= \epsilon_A(\delta_{u,s.v}\sigma_u(s,t)e_u \otimes x_{st}) = \delta_{u,s.v}\sigma_u(s,t)\delta_{u1} \\
&= \delta_{u1}\delta_{v1} = \epsilon_A(e_u \otimes x_s)\epsilon_A(e_v \otimes x_t)
\end{aligned}$$

It is easy to see that $\epsilon_A$ preserves the unit.

(6) We now proceed to prove that the antipode is $H$-linear and colinear. Since $\gamma$ is a 1-cocycle, we know from Lemma 1.13 that $\gamma_{s^{-1}} = (s^{-1}.\gamma_s)^{-1}$ for all $s \in G$, i. e., we have $\gamma_{s^{-1}}(u) = \gamma_s(s.u^{-1})$. This implies:

$$c_b \to S_A(e_u \otimes x_s) = \gamma_{s^{-1}}(s^{-1}.u^{-1})(b)\, S_A(e_u \otimes x_s) = \gamma_s(u)(b)\, S_A(e_u \otimes x_s)$$
$$= S_A(c_b \to e_u \otimes x_s)$$

Since $z$ is also a 1-cocycle, we know from Lemma 1.13 that $z_{s^{-1}} = (s^{-1}.z_s)^{-1}$, i. e., that $z_{s^{-1}}(u) = z_s(s.u^{-1})$. This implies:

$$\delta_A(S_A(e_u \otimes x_s)) = c_{z_{s^{-1}}(s^{-1}.u^{-1})} \otimes S_A(e_u \otimes x_s) = c_{z_s(u)} \otimes S_A(e_u \otimes x_s)$$
$$= (\mathrm{id}_H \otimes S_A)\delta_A(e_u \otimes x_s)$$

(7) At last, we show that $S_A$ really is an antipode for $A$. We begin by proving that it is a right antipode for $A$:

$$\mu_A \circ (\mathrm{id}_A \otimes S_A) \circ \Delta_A(e_u \otimes x_s) = \sum_{v \in P}(e_v \otimes x_s)S_A(e_{v^{-1}u} \otimes x_s)$$
$$= \sum_{v \in P}\sigma^{-1}_{u^{-1}v}(s, s^{-1})(e_v \otimes x_s)(e_{s^{-1}.(u^{-1}v)} \otimes x_{s^{-1}})$$
$$= \sum_{v \in P}\sigma^{-1}_{u^{-1}v}(s, s^{-1})\sigma_v(s, s^{-1})\delta_{v, u^{-1}v}(e_v \otimes x_1)$$
$$= \delta_{u1}\sum_{v \in P}e_v \otimes x_1 = \epsilon_A(e_u \otimes x_s)1_A$$

Next, we prove that $S_A$ is a left antipode for $A$. For this, we observe first that, by Lemma 1.13, we have:

$$\sigma(s^{-1}, s) = s^{-1}.\sigma(s, s^{-1})$$

This implies that $\sigma_u(s, s^{-1}) = \sigma_{s^{-1}.u}(s^{-1}, s)$. Therefore, the antipode is also given by the formula:

$$S_A(e_u \otimes x_s) = \sigma^{-1}_{s^{-1}.u^{-1}}(s^{-1}, s)e_{s^{-1}.u^{-1}} \otimes x_{s^{-1}}$$

This implies:

$$\mu_A \circ (S_A \otimes \mathrm{id}_A) \circ \Delta_A(e_u \otimes x_s) = \sum_{v \in P}S_A(e_v \otimes x_s)(e_{v^{-1}u} \otimes x_s)$$
$$= \sum_{v \in P}\sigma^{-1}_{s^{-1}.v^{-1}}(s^{-1}, s)(e_{s^{-1}.v^{-1}} \otimes x_{s^{-1}})(e_{v^{-1}u} \otimes x_s)$$
$$= \sum_{v \in P}\sigma^{-1}_{s^{-1}.v^{-1}}(s^{-1}, s)\sigma_{s^{-1}.v^{-1}}(s^{-1}, s)\delta_{s^{-1}.v^{-1}, s^{-1}.(v^{-1}u)}(e_{s^{-1}.v^{-1}} \otimes x_1)$$
$$= \delta_{u1}\sum_{v \in P}e_{s^{-1}.v^{-1}} \otimes x_1 = \epsilon_A(e_u \otimes x_s)1_A$$

This completes the proof of the proposition. $\square$

Note that $A$ is semisimple if the characteristic of the base field does not divide the order of $G$, since a crossed product of a semisimple group ring with a semisimple algebra is again semisimple (cf. [57], Thm. 7.4.2, p. 116, [64], Thm. 4.4, p. 31). $A$ is cosemisimple if the characteristic of the base field does not divide the order of $P$, since $K^P$ is cosemisimple in this case, whereas $K[G]$ is always cosemisimple.

**3.3** In the second stage, we shall construct a Yetter-Drinfel'd Hopf algebra for any group homomorphism from a finite group to the group of units of a finite ring. Suppose that $G$ is a finite group and that $R$ is a finite ring. We assume that we are given a group homomorphism

$$\nu : G \to U(R)$$

from $G$ to the multiplicative group $U(R)$ of units of $R$. We use $\nu$ to turn the additive group of $R$ into a left $G$-module in two ways. First, $R$ becomes a $G$-module via:

$$G \times R \to R, (s, u) \mapsto s \cdot u := \nu(s)u$$

We denote $R$ by ${}_G R$ if it is regarded as a left $G$-module in this way. Second, $R$ becomes a left $G$-module via:

$$G \times R \to R, (s, u) \mapsto s \, . \, u := u\nu(s^{-1})$$

We denote $R$ by $R_G$ if it is regarded as a left $G$-module in this way.

We now assume that the following additional structure elements are given:

1. Two 1-cocycles $\alpha, \beta \in Z^1(G, {}_G R)$, i. e., mappings from $G$ to $R$ satisfying

$$\alpha(st) = \alpha(s) + \nu(s)\alpha(t) \qquad \beta(st) = \beta(s) + \nu(s)\beta(t)$$

    for all $s, t \in G$.

2. A normalized 2-cocycle $q \in Z^2(G, {}_G R)$.

3. Two characters $\chi, \eta \in \hat{R}$ of the additive group of $R$, i. e., mappings from $R$ to $K^\times$ satisfying

$$\chi(u + v) = \chi(u)\chi(v) \qquad \eta(u + v) = \eta(u)\eta(v)$$

    for all $u, v \in R$.

These structure elements can be chosen freely, except that we suppose that $\chi$ satisfies the condition:
$$\chi(uvw) = \chi(vuw)$$

for all $u, v, w \in R$. For example, this condition is satisfied if $R$ is commutative. It should be observed that this condition implies

$$\chi(u_1 \cdot \ldots \cdot u_i u_{i+1} \cdot \ldots \cdot u_k) = \chi(u_i u_{i+1} \cdot \ldots \cdot u_k u_1 \cdot \ldots \cdot u_{i-1})$$
$$= \chi(u_{i+1} u_i \cdot \ldots \cdot u_k u_1 \cdot \ldots \cdot u_{i-1}) = \chi(u_1 \cdot \ldots \cdot u_{i+1} u_i \cdot \ldots \cdot u_k)$$

and therefore $\chi(u_1 \cdot \ldots \cdot u_k) = \chi(u_{\tau(1)} \cdot \ldots \cdot u_{\tau(k)})$ for all $\tau \in S_k$, since the above transpositions generate the symmetric group $S_k$. Another way of seeing this is to observe that the assumption implies that $\chi$ vanishes on the two-sided ideal generated by the additive commutators $uv - vu$, and therefore is in fact induced from a character of a commutative ring.

We define $H := K[R]$, the group ring of the additive group of $R$, and $A := K^R \otimes K[G]$, the tensor product of the dual group ring of the additive group of $R$ and the group ring of $G$. As in the previous paragraph, the canonical basis elements of $K[R]$ resp. $K[G]$ are denoted by $c_u$ resp. $x_s$, where $u \in R$ and $s \in G$. The primitive idempotents of $K^R$ are denoted by $e_u$, where $u \in R$.

**Proposition** The vector space $A = K^R \otimes K[G]$ becomes a Yetter-Drinfel'd Hopf algebra over $H$ if it is endowed with the following structures:

1. Crossed product algebra structure:
$$(e_u \otimes x_s)(e_v \otimes x_t) := \delta_{u\nu(s),v} \eta(uq(s,t)) \chi(u^2 \nu(s) \beta(s) \alpha(t)) e_u \otimes x_{st}$$

    Unit: $1_A := \sum_{u \in R} e_u \otimes x_1$

2. Tensor product coalgebra structure:
$$\Delta_A(e_u \otimes x_s) := \sum_{v \in R} (e_v \otimes x_s) \otimes (e_{u-v} \otimes x_s)$$

    Counit: $\epsilon_A(e_u \otimes x_s) := \delta_{u0}$

3. Action: $c_u \to (e_v \otimes x_s) := \chi(uv\alpha(s))^2 e_v \otimes x_s$

4. Coaction: $\delta_A(e_u \otimes x_s) := c_{u\beta(s)} \otimes (e_u \otimes x_s)$

5. Antipode: $S_A(e_u \otimes x_s) := \eta(uq(s, s^{-1})) \chi(u^2 \beta(s) \alpha(s)) e_{-u\nu(s)} \otimes x_{s^{-1}}$

**Proof.** (1) For the proof, we use the framework of Paragraph 3.2 with $C = R$ and $P = R_G$. To do this, we have to specify 1-cocycles

$$z : G \to \text{Hom}(R_G, R), s \mapsto z_s \qquad \gamma : G \to \text{Hom}(R_G, \hat{R}), s \mapsto \gamma_s$$

of $G$ with values in the $G$-modules $\text{Hom}(R_G, R)$ resp. $\text{Hom}(R_G, \hat{R})$ and a 2-cocycle $\sigma : G \times G \to U(K^R)$ of $G$ with values in the group of units $U(K^R)$ of $R$, which we write in the form

$$\sigma(s, t) = \sum_{u \in R} \sigma_u(s, t) e_u$$

for functions $\sigma_u : G \times G \to K^\times$. We define:

1. $z_s(u) := u\beta(s)$
2. $(\gamma_s(u))(v) := \chi(uv\alpha(s))^2$
3. $\sigma_u(s,t) := \eta(uq(s,t))\chi(u^2\nu(s)\beta(s)\alpha(t))$

We now have to prove that $z$, $\gamma$ and $\sigma$ are cocycles that satisfy the compatibility condition:
$$\sigma_{u+v}(s,t) = ((s.\gamma_t)(u))(z_s(v))\sigma_u(s,t)\sigma_v(s,t)$$

(2) It is obvious that $z_s : R \to R$ is a group homomorphism. $z$ is also a 1-cocycle, since we have:
$$z_{st}(u) = u\beta(st) = u\beta(s) + u\nu(s)\beta(t) = z_s(u) + z_t(s^{-1}.u)$$

It is also obvious that $\gamma_s(u)$ is a character and that $\gamma_s : R \to \hat{R}$ is a group homomorphism. $\gamma$ is a 1-cocycle, since we have:
$$(\gamma_{st}(u))(v) = \chi(vu\alpha(st))^2 = \chi(vu\alpha(s) + vu\nu(s)\alpha(t))^2$$
$$= \chi(vu\alpha(s))^2\chi(v(s^{-1}.u)\alpha(t))^2 = (\gamma_s(u))(v)(\gamma_t(s^{-1}.u))(v)$$

(3) We now prove that $\sigma$ is a 2-cocycle. We have:

$\sigma_u(r,st)\sigma_{r^{-1}.u}(s,t)$
$= \eta(uq(r,st))\chi(u^2\nu(r)\beta(r)\alpha(st))\eta((r^{-1}.u)q(s,t))\chi((r^{-1}.u)^2\nu(s)\beta(s)\alpha(t))$
$= \eta(u(q(r,st) + \nu(r)q(s,t)))$
$\qquad \chi(u^2\nu(r)\beta(r)(\alpha(s) + \nu(s)\alpha(t)) + u^2\nu(r)^2\nu(s)\beta(s)\alpha(t))$
$= \eta(u(q(rs,t) + q(r,s)))\chi(u^2\nu(r)\beta(r)\alpha(s))$
$\qquad \chi(u^2\nu(rs)\beta(r)\alpha(t) + u^2\nu(rs)\nu(r)\beta(s)\alpha(t))$
$= \sigma_u(r,s)\eta(uq(rs,t))\chi(u^2\nu(rs)(\beta(r) + \nu(r)\beta(s))\alpha(t))$
$= \sigma_u(r,s)\sigma_u(rs,t)$

It now follows from an easy calculation that $\sigma$ is a 2-cocycle.

From Lemma 1.13, we have $\alpha(1) = 0$ and $\beta(1) = 0$. Therefore, $\sigma$ is normalized since $q$ is normalized.

(4) It remains to verify the compatibility condition. We have:

$\sigma_{u+v}(s,t) = \eta((u+v)q(s,t))\chi((u+v)^2\nu(s)\beta(s)\alpha(t))$
$= \eta(uq(s,t))\eta(vq(s,t))\chi(u^2\nu(s)\beta(s)\alpha(t))\chi(v^2\nu(s)\beta(s)\alpha(t))\chi(2uv\nu(s)\beta(s)\alpha(t))$
$= \sigma_u(s,t)\sigma_v(s,t)(\gamma_t(s^{-1}.u))(v\beta(s))$

(5) It now follows from Proposition 3.2 that $A$ becomes a Yetter-Drinfel'd Hopf algebra. It is easy to see that the structure elements given there are those stated above, except for the antipode. The antipode given in Proposition 3.2 has the form

$$S_A(e_u \otimes x_s) = \sigma_{-u}^{-1}(s, s^{-1}) e_{-u\nu(s)} \otimes x_{s^{-1}}$$
$$= \eta(uq(s, s^{-1})) \chi(-u^2 \nu(s) \beta(s) \alpha(s^{-1})) e_{-u\nu(s)} \otimes x_{s^{-1}}$$

But since we have $\alpha(s^{-1}) = -\nu(s^{-1})\alpha(s)$ by Lemma 1.13, this is the form of the antipode given above. □

The Yetter-Drinfel'd Hopf algebra considered in the previous proposition will be denoted by $A_G(\alpha, \beta, q)$. This notation is ambiguous since the construction also depends on $R$, $\nu$, $\chi$, and $\eta$. We shall therefore only use this notation if the remaining structure elements are clear from the context. If $G = \mathbb{Z}_p$ is a cyclic group of prime order $p$, we shall abbreviate the notation $A_{\mathbb{Z}_p}(\alpha, \beta, q)$ to $A_p(\alpha, \beta, q)$.

In the case where $R$ is commutative, the definition of $\sigma$ can be rewritten using the cup product. The multiplication map

$$R \otimes_\mathbb{Z} R \to R, u \otimes v \mapsto uv$$

is then $G$-equivariant, provided that we consider $R$ as a $G$-module in a third way via the module structure

$$G \times R \to R, (s, u) \mapsto \nu(s)^2 u$$

The cup product $\beta \cup \alpha \in Z^2(G, R \otimes_\mathbb{Z} R)$ then is mapped to a 2-cocycle in $Z^2(G, R)$ which we, following common usage, also denote by $\beta \cup \alpha$. The definition of $\sigma$ then takes the form:

$$\sigma_u(s, t) = \eta(uq(s, t)) \chi(u^2 (\beta \cup \alpha)(s, t))$$

**3.4** We now look at the special case of the preceding construction which will be the most important one in the following, because the structure theorem mentioned above and stated in Paragraph 7.7 will say that all Yetter-Drinfel'd Hopf algebras of a certain type are necessarily of the form that we describe now. Suppose that $p$ is an odd prime and that $R = \mathbb{Z}_p$, the finite field with $p$ elements. Suppose that $\zeta \in K$ is a $p$-th root of unity, which need not be primitive here. Since $p$ is odd, 2 is an invertible element of $\mathbb{Z}_p$, and therefore the expression $i/2$ for $i \in \mathbb{Z}_p$ makes sense. We then define the characters

$$\chi : \mathbb{Z}_p \to K, i \mapsto \zeta^{i/2} \qquad \eta : \mathbb{Z}_p \to K, i \mapsto \zeta^i$$

Now, if $G$ is a finite group, $\nu : G \to \mathbb{Z}_p^\times$ is a group homomorphism, $\alpha, \beta \in Z^1(G, {}_G\mathbb{Z}_p)$ are two 1-cocycles, and $q \in Z^2(G, {}_G\mathbb{Z}_p)$ is a normalized 2-cocycle, we can construct the Yetter-Drinfel'd Hopf algebra $A_G(\alpha, \beta, q)$. Its structure elements then take the following form:

1. Crossed product algebra structure:
$$(e_i \otimes x_s)(e_j \otimes x_t) := \delta_{i\nu(s),j} \zeta^{iq(s,t)} \zeta^{i^2\nu(s)\beta(s)\alpha(t)/2} e_i \otimes x_{st}$$

Unit: $1_A := \sum_{i=0}^{p-1} e_i \otimes x_1$

2. Tensor product coalgebra structure:
$$\Delta_A(e_i \otimes x_s) := \sum_{j=0}^{p-1} (e_j \otimes x_s) \otimes (e_{i-j} \otimes x_s)$$

Counit: $\epsilon_A(e_i \otimes x_s) := \delta_{i0}$

3. Action: $c_i \to (e_j \otimes x_s) := \zeta^{ij\alpha(s)} e_j \otimes x_s$

4. Coaction: $\delta_A(e_i \otimes x_s) := c_{i\beta(s)} \otimes (e_i \otimes x_s)$

5. Antipode: $S_A(e_i \otimes x_s) := \zeta^{iq(s,s^{-1})} \zeta^{i^2\beta(s)\alpha(s)/2} e_{-i\nu(s)} \otimes x_{s^{-1}}$

**3.5** We shall now give another application of the framework considered in Paragraph 3.2. Suppose that $G$ is a finite group and that the base field $K$ contains a primitive fourth root of unity $\iota$; this can, of course, only happen if $K$ does not have characteristic 2. In the situation of Paragraph 3.2, we put $P = C := \mathbb{Z}_2$. For $i \in \mathbb{Z}_2$, we denote the corresponding primitive idempotent in $K^{\mathbb{Z}_2}$ by $e_i$ and the corresponding canonical basis vector of $K[\mathbb{Z}_2]$ by $c_i$; the canonical basis vectors of $K[G]$ are denoted by $x_s$, for $s \in G$. We regard $\mathbb{Z}_2$ as a trivial $G$-module. Suppose that $\alpha: G \to \mathbb{Z}_2$ and $\beta: G \to \mathbb{Z}_2$ are 1-cocycles; since the $G$-module structure is trivial, these are just group homomorphisms. Suppose that $q \in Z^2(G, \mathbb{Z}_4)$ is a normalized 2-cocycle of the trivial $G$-module $\mathbb{Z}_4$ such that
$$\hat{\pi} \circ q = \beta \cup \alpha$$
where $\hat{\pi}: \mathbb{Z}_4 \to \mathbb{Z}_2$ is the unique surjective group homomorphism. Here we have, as in Paragraph 3.3, used the isomorphism $\mathbb{Z}_2 \otimes_{\mathbb{Z}} \mathbb{Z}_2 \cong \mathbb{Z}_2$ to regard the cup product $\beta \cup \alpha \in Z^2(G, \mathbb{Z}_2 \otimes_{\mathbb{Z}} \mathbb{Z}_2)$ as an element of $Z^2(G, \mathbb{Z}_2)$, i. e., we have:
$$\beta \cup \alpha: G \times G \to \mathbb{Z}_2, (s,t) \mapsto \beta(s)\alpha(t)$$

Using these data, we can construct group homomorphisms $z: G \to \operatorname{Hom}(\mathbb{Z}_2, \mathbb{Z}_2)$, $s \mapsto z_s$ and $\gamma: G \to \operatorname{Hom}(\mathbb{Z}_2, \hat{\mathbb{Z}}_2), s \mapsto \gamma_s$ by:
$$z_s(i) := i\beta(s) \qquad (\gamma_s(i))(j) := (-1)^{ij\alpha(s)}$$

For $i \in \mathbb{Z}_2$, we define $\sigma_i: G \times G \to K^{\mathbb{Z}_2}$ by
$$\sigma_0(s,t) := 1 \qquad \sigma_1(s,t) := \iota^{q(s,t)}$$

and set:
$$\sigma(s,t) := \sigma_0(s,t) e_0 + \sigma_1(s,t) e_1$$

These data satisfy the requirements of Paragraph 3.2: It is obvious that $z$ and $\gamma$ are homomorphisms; $\sigma$ is a 2-cocycle, since $q$ is a 2-cocycle; we only have to verify the compatibility condition:

$$\sigma_{i+j}(s,t) = (\gamma_t(i))(z_s(j))\sigma_i(s,t)\sigma_j(s,t)$$

This equation is obvious if $i = 0$ or if $j = 0$; in the case $i = j = 1$ it says that $1 = (-1)^{\beta(s)\alpha(t)} \iota^{q(s,t)} \iota^{q(s,t)} = (-1)^{\beta(s)\alpha(t)}(-1)^{q(s,t)}$. But in this case it follows from the assumption that $\hat{\pi}(q(s,t)) = (\beta \cup \alpha)(s,t) = \beta(s)\alpha(t)$.

It now follows from Paragraph 3.2 that $A := K^{\mathbb{Z}_2} \otimes K[G]$ becomes a Yetter-Drinfel'd Hopf algebra over $H := K[\mathbb{Z}_2]$ if it is endowed with the following structures:

1. Crossed product algebra structure:

$$(e_i \otimes x_s)(e_j \otimes x_t) := \delta_{ij}\sigma_i(s,t)e_i \otimes x_{st}$$

    Unit: $1_A := \sum_{i=0}^{1} e_i \otimes x_1$

2. Tensor product coalgebra structure:

$$\Delta_A(e_i \otimes x_s) := \sum_{j=0}^{1}(e_j \otimes x_s) \otimes (e_{i-j} \otimes x_s)$$

    Counit: $\epsilon_A(e_i \otimes x_s) := \delta_{i0}$

3. Action: $c_i \to (e_j \otimes x_s) := (-1)^{ij\alpha(s)}e_j \otimes x_s$

4. Coaction: $\delta_A(e_i \otimes x_s) := c_{i\beta(s)} \otimes (e_i \otimes x_s)$

5. Antipode: $S_A(e_i \otimes x_s) := \sigma_i^{-1}(s, s^{-1})e_{-i} \otimes x_{s^{-1}}$

**3.6** The simplest case of the situation discussed in the previous paragraph is the case where also $G = \mathbb{Z}_2$. In this case, $A$ has dimension 4. If we want to exclude that $A$ is trivial, we must have that $\alpha \neq 0$ and $\beta \neq 0$; therefore, we are only left with the possibility $\alpha = \beta = \text{id}$. The 2-cocycle $q$ is then also almost determined: Since we require it to be normalized, it has to satisfy $q(i,j) = 0$ if $i = 0$ or $j = 0$; since it has to satisfy $\hat{\pi}(q(i,j)) = (\beta \cup \alpha)(i,j) = ij$, we have $\hat{\pi}(q(1,1)) = 1$ and therefore $q(1,1) = 1$ or $q(1,1) = 3$. We denote these two possibilities by $q_+$ resp. $q_-$, i. e., we define:

$$q_+(i,j) := \begin{cases} 0 & \text{if } i = 0 \text{ or } j = 0 \\ 1 & \text{if } i = 1 \text{ and } j = 1 \end{cases}$$

$$q_-(i,j) := \begin{cases} 0 & \text{if } i = 0 \text{ or } j = 0 \\ 3 & \text{if } i = 1 \text{ and } j = 1 \end{cases}$$

Since the equality

$$q_\pm(i+j, k) + q_\pm(i, j) = q_\pm(j, k) + q_\pm(i, j+k)$$

holds if $i = 0$, $j = 0$, or $k = 0$, but also in the case $i = j = k = 1$, $q_+$ resp. $q_-$ really are 2-cocycles. They are not coboundaries, because if we would have

$$q_+(i, j) = w(j) - w(i + j) + w(i)$$

for some map $w : \mathbb{Z}_2 \to \mathbb{Z}_4$, we would have $w(0) = q_+(0, 0) = 0$ and therefore $q_+(1, 1) = 2w(1)$, which is impossible. We have $q_- = -q_+$, and therefore $q_-$ is not a coboundary, too. Because $H^2(\mathbb{Z}_2, \mathbb{Z}_4) \cong \mathbb{Z}_2$ by Proposition 1.13.2, $q_+$ and $q_-$ are cohomologous. We denote the Yetter-Drinfel'd Hopf algebras arising from $q_+$ resp. $q_-$ by $A_+$ resp. $A_-$:

**Definition** $A_\pm$ denotes the Yetter-Drinfel'd Hopf algebra over $H = K[\mathbb{Z}_2]$ with underlying vector space $K^{\mathbb{Z}_2} \otimes K[\mathbb{Z}_2]$ and the following structures:

1. Crossed product algebra structure:

$$(e_i \otimes c_k)(e_j \otimes c_l) = \delta_{ij} \sigma_i^\pm(k, l) e_i \otimes c_{k+l}$$

   Unit: $1_{A_\pm} = \sum_{i=0}^1 e_i \otimes c_0$

2. Tensor product coalgebra structure:

$$\Delta_{A_\pm}(e_i \otimes c_j) = \sum_{k=0}^1 (e_k \otimes c_j) \otimes (e_{i-k} \otimes c_j)$$

   Counit: $\epsilon_{A_\pm}(e_i \otimes c_j) = \delta_{i0}$

3. Action: $c_k \to (e_i \otimes c_j) = (-1)^{ijk} e_i \otimes c_j$

4. Coaction: $\delta_{A_\pm}(e_i \otimes c_j) = c_{ij} \otimes (e_i \otimes c_j)$

5. Antipode: $S_{A_\pm}(e_i \otimes c_j) = \sigma_i^\pm(j, -j)^{-1} e_{-i} \otimes c_{-j}$

where, as in Paragraph 3.5,

$$\sigma_0^\pm(i, j) := 1 \qquad \sigma_1^\pm(i, j) := \iota^{q_\pm(i,j)}$$

Since, in the terminology of Paragraph 3.2, we are in the situation $C = G$, we have used here the notation $c_i$, and not $x_i$, for the canonical basis vectors of $K[\mathbb{Z}_2]$. It should be noted that the construction of $A$ depends on the choice of the primitive fourth root of unity $\iota$. Choosing the other possibility $-\iota$ in the above construction interchanges $A_+$ and $A_-$. Observe that both algebras are also commutative.

# 4 Isomorphisms

**4.1** In the previous section, we have constructed some classes of examples of Yetter-Drinfel'd Hopf algebras. These constructions depended on certain structure elements, and not all different choices of such structure elements lead to nonisomorphic Yetter-Drinfel'd Hopf algebras. We shall now find sufficient conditions on pairs of these structure elements to give rise to isomorphic Yetter-Drinfel'd Hopf algebras. In special cases, we shall see that these conditions are also necessary.

To construct isomorphisms, we shall more generally construct homomorphisms; it will be easy to see when these homomorphisms are actually isomorphisms. Expressed in a more categorical language, we are constructing functors; while we took care of the objects in the previous section, we will now take care of the morphisms.

In this section, $K$ denotes a field. In Paragraph 4.5 and Paragraph 4.6, we require that $K$ contains a primitive $p$-th root of unity, where $p$ is an odd prime. In Paragraph 4.7, Paragraph 4.8, and Paragraph 4.9, we require that $K$ contains a primitive fourth root of unity. Throughout this section, we will use the same notation for an integer $i$ and its equivalence class in a factor group $\mathbb{Z}_p = \mathbb{Z}/p\mathbb{Z}$.

**4.2** We now consider the question in which cases the Yetter-Drinfel'd Hopf algebras constructed in Paragraph 3.2 are isomorphic. Suppose that $C$ is a finite group. Suppose that $G$ and $G'$ are also finite groups and that $P$, resp. $P'$, is a (nonabelian) left $G$-module, resp. a $G'$-module. We denote by $H := K[C]$ the group ring of $C$. $A := K^P \otimes K[G]$, resp. $A' := K^{P'} \otimes K[G']$, denotes the tensor product of the dual group ring of $P$, resp. $P'$, and the group ring of $G$, resp. $G'$. The canonical basis elements of $K[C]$, $K[G]$ and $K[G']$ are denoted by $c_b$, $x_s$ resp. $x'_{s'}$, where $b \in C$, $s \in G$ and $s' \in G'$. The primitive idempotents of $K^P$ and $K^{P'}$ are denoted by $e_u$ resp. $e'_{u'}$, for $u \in P$ and $u' \in P'$. As in Paragraph 3.2, we turn the spaces $\mathrm{Hom}(P, Z(C))$ and $\mathrm{Hom}(P, \hat{C})$ into left $G$-modules, and similarly the spaces $\mathrm{Hom}(P', Z(C))$ and $\mathrm{Hom}(P', \hat{C})$ into left $G'$-modules. As in Paragraph 3.2, we extend the $G$-module structure of $P$ and the $G'$-module structure of $P'$ linearly to a $K[G]$-module structure of $K^P$, resp. to a $K[G']$-module structure of $K^{P'}$.

We now suppose that the following structure elements are given:

1. A 1-cocycle $z' : G' \to \mathrm{Hom}(P', Z(C)), s' \mapsto z'_{s'}$ of $G'$ with values in the $G'$-module $\mathrm{Hom}(P', Z(C))$.

2. A 1-cocycle $\gamma' : G' \to \operatorname{Hom}(P', \hat{C}), s' \mapsto \gamma'_{s'}$ of $G'$ with values in the $G'$-module $\operatorname{Hom}(P', \hat{C})$.

3. A normalized 2-cocycle $\sigma' : G' \times G' \to U(K^{P'})$ of $G'$ with values in the group of units $U(K^{P'})$ of $K^{P'}$, which we write in the form

$$\sigma'(s', t') = \sum_{u' \in P'} \sigma'_{u'}(s', t') e_{u'}$$

for functions $\sigma'_{u'} : G' \times G' \to K^\times$.

We require that these structure elements satisfy the following compatibility condition:

$$\sigma'_{u'v'}(s', t') = ((s'.\gamma'_{t'})(u'))(z'_{s'}(v'))\sigma'_{u'}(s', t')\sigma'_{v'}(s', t')$$

for all $u', v' \in P'$ and all $s', t' \in G'$. As explained in Paragraph 3.2, these structure elements can be used to define a Yetter-Drinfel'd Hopf algebra structure on the vector space $A'$.

Now suppose that

$$f_P : P \to P' \qquad f_G : G \to G'$$

are group homomorphisms. As explained in Paragraph 1.13, every $G'$-module then becomes a $G$-module by pullback via $f_G$. We assume that $f_P$ is bijective and $G$-equivariant, i. e., we have $f_G(s).f_P(u) = f_P(s.u)$, for $s \in G$ and $u \in P$. In this situation, we define $z : G \to \operatorname{Hom}(P, Z(C))$ by

$$z_s(u) := z'_{f_G(s)}(f_P(u))$$

Similarly, we define $\gamma : G \to \operatorname{Hom}(P, \hat{C})$ by

$$\gamma_s(u) := \gamma'_{f_G(s)}(f_P(u))$$

Then $z$ and $\gamma$ are 1-cocycles. We could also define $\sigma : G \times G \to K^P$ by the condition

$$f_P(\sigma(s, t)) = \sigma'(f_G(s), f_G(t))$$

where we have extended $f_P$ linearly to a map $f_P : K^P \to K^{P'}, e_u \mapsto e'_{f_P(u)}$, to obtain a 2-cocycle $\sigma$; however, we shall consider the more general case where $\sigma$ is only required to be cohomologous to the cocycle above. We therefore assume that we are given a 1-cochain $\tau : G \to U(K^{P'})$ that satisfies $\tau(1) = 1$ and define $\sigma : G \times G \to K^P$ by the condition

$$f_P(\sigma(s, t)) = \sigma'(f_G(s), f_G(t))(s.\tau(t))\tau(st)^{-1}\tau(s)$$

Obviously, $\sigma$ is then normalized. As for $\sigma'$, we write $\sigma$ and $\tau$ in the form:

$$\sigma(s, t) = \sum_{u \in P} \sigma_u(s, t) e_u \qquad \tau(s) = \sum_{u' \in P'} \tau_{u'}(s) e_{u'}$$

The definition of $\sigma$ then means for these components that we have:

$$\sigma_u(s, t) = \sigma'_{f_P(u)}(f_G(s), f_G(t))\tau_{f_P(s^{-1}.u)}(t)\tau_{f_P(u)}(st)^{-1}\tau_{f_P(u)}(s)$$

**Proposition** Suppose that we have $\tau_{u'v'}(s) = \tau_{u'}(s)\tau_{v'}(s)$ for all $s \in G$ and $u', v' \in P'$. Then we have:

1. For all $s, t \in G$ and all $u, v \in P$, the compatibility condition
$$\sigma_{uv}(s,t) = ((s.\gamma_t)(u))(z_s(v))\sigma_u(s,t)\sigma_v(s,t)$$
is satisfied. Therefore $A$, endowed with the structure elements given in Paragraph 3.2, is a Yetter-Drinfel'd Hopf algebra.

2. The map
$$f_A : A \to A', e_u \otimes x_s \mapsto \tau_{f_P(u)}(s)e'_{f_P(u)} \otimes x'_{f_G(s)}$$
is a morphism of Yetter-Drinfel'd Hopf algebras.

**Proof.** First, we verify the compatibility condition:
$$\sigma_{uv}(s,t) = \sigma'_{f_P(u)f_P(v)}(f_G(s), f_G(t))$$
$$\tau_{f_P(s^{-1}.u)f_P(s^{-1}.v)}(t)\tau_{f_P(u)f_P(v)}(st)^{-1}\tau_{f_P(u)f_P(v)}(s)$$
$$= ((f_G(s).\gamma'_{f_G(t)})(f_P(u)))(z'_{f_G(s)}(f_P(v)))$$
$$\sigma'_{f_P(u)}(f_G(s), f_G(t))\sigma'_{f_P(v)}(f_G(s), f_G(t))$$
$$\tau_{f_P(s^{-1}.u)}(t)\tau_{f_P(u)}(st)^{-1}\tau_{f_P(u)}(s)\tau_{f_P(s^{-1}.v)}(t)\tau_{f_P(v)}(st)^{-1}\tau_{f_P(v)}(s)$$
$$= ((s.\gamma_t)(u))(z_s(v))\sigma_u(s,t)\sigma_v(s,t)$$

We now prove the second assertion. Linearity and colinearity of $f_A$ follow directly from our definition of $z$ and $\gamma$. We now show that $f_A$ is an algebra homomorphism:
$$f_A(e_u \otimes x_s)f_A(e_v \otimes x_t) = \tau_{f_P(u)}(s)\tau_{f_P(v)}(t)(e'_{f_P(u)} \otimes x'_{f_G(s)})(e'_{f_P(v)} \otimes x'_{f_G(t)})$$
$$= \tau_{f_P(u)}(s)\tau_{f_P(v)}(t)\delta_{f_P(u), f_G(s).f_P(v)}\sigma'_{f_P(u)}(f_G(s), f_G(t))e'_{f_P(u)} \otimes x'_{f_G(st)}$$
$$= \delta_{u, s.v}\sigma_u(s,t)\tau_{f_P(u)}(st)e'_{f_P(u)} \otimes x'_{f_G(st)} = f_A((e_u \otimes x_s)(e_v \otimes x_t))$$

Our assumption on $\tau$ also assures that $f_A$ is a coalgebra homomorphism:
$$(f_A \otimes f_A)\Delta_A(e_u \otimes x_s)$$
$$= \sum_{v \in P} \tau_{f_P(v)}(s)\tau_{f_P(v^{-1}u)}(s)(e'_{f_P(v)} \otimes x'_{f_G(s)}) \otimes (e'_{f_P(v^{-1}u)} \otimes x'_{f_G(s)})$$
$$= \tau_{f_P(u)}(s)\Delta_{A'}(e'_{f_P(u)} \otimes x'_{f_G(s)}) = \Delta_{A'}(f_A(e_u \otimes x_s))$$

Since we have $\tau(1) = 1$, i. e., $\tau_{u'}(1) = 1$ for all $u' \in P'$, $f_A$ preserves the unit. $f_A$ preserves the counit since $u' \mapsto \tau_{u'}(s)$ is a group homomorphism, which implies that $\tau_1(s) = 1$ for all $s \in G$. $f_A$ is therefore a morphism of Yetter-Drinfel'd bialgebras. Exactly as for ordinary Hopf algebras, it can be shown that a morphism of Yetter-Drinfel'd bialgebras commutes with the antipode, and therefore is a morphism of Yetter-Drinfel'd Hopf algebras (cf. [84], Lem. 4.0.4, p. 81). $\square$

**4.3**  In Paragraph 3.3, we have explained how we can construct Yetter-Drinfel'd Hopf algebras from group homomorphisms to the additive group of a finite ring. In this paragraph, we shall find sufficient conditions under which Yetter-Drinfel'd Hopf algebras arising from this construction are isomorphic.

Suppose that $G'$ is a finite group and that $R$ is a finite ring. We assume that we are given a group homomorphism
$$\nu' : G' \to U(R)$$
from $G'$ to the multiplicative group $U(R)$ of units of $R$. As explained in Paragraph 3.3, we can use $\nu'$ to turn the additive group of $R$ into a left $G'$-module in two ways, denoted by $_{G'}R$ and $R_{G'}$ respectively.

We assume that the following structure elements are given:

1. Two 1-cocycles $\alpha', \beta' \in Z^1(G', {}_{G'}R)$.

2. A normalized 2-cocycle $q' \in Z^2(G', {}_{G'}R)$.

3. Two characters $\chi, \eta \in \hat{R}$ of the additive group of $R$.

As in Paragraph 3.3, we suppose that $\chi$ satisfies the condition $\chi(uvw) = \chi(vuw)$ for all $u, v, w \in R$. The canonical basis elements of $K[R]$ resp. $K[G']$ are denoted by $c_u$ resp. $x'_{s'}$, where $u \in R$ and $s' \in G'$. The primitive idempotents of $K^R$ are denoted by $e_u$, where $u \in R$. In Paragraph 3.3, we constructed from these data a Yetter-Drinfel'd Hopf algebra $A_{G'}(\alpha', \beta', q')$ over $H := K[R]$, the group ring of the additive group of $R$.

Now suppose that $G$ is another finite group and that $f : G \to G'$ is a group homomorphism. Suppose further that $x \in Z_{U(R)}(\nu'(f(G)))$ is an element of the centralizer $Z_{U(R)}(\nu'(f(G)))$ of $\nu'(f(G))$ in the group of units $U(R)$ of $R$, i. e., an element $x \in U(R)$ satisfying $x\nu'(f(s)) = \nu'(f(s))x$ for all $s \in G$. In this situation, we define:

1. $\nu : G \to R, s \mapsto \nu(s) := \nu'(f(s))$

2. $\alpha : G \to R, s \mapsto \alpha(s) := x\alpha'(f(s))$

3. $\beta : G \to R, s \mapsto \beta(s) := x\beta'(f(s))$

If $_GR$ denotes $R$ regarded as a $G$-module via $\nu$, then $\alpha$ and $\beta$ are 1-cocycles in $_GR$, since they are the images of the 1-cocycles $f^1(\alpha')$ resp. $f^1(\beta')$ with respect to the $G$-equivariant homomorphism $u \mapsto xu$. We could also define for all $s, t \in G$
$$q(s,t) := xq'(f(s), f(t))$$

to get a 2-cocycle; however, we shall consider the more general case where $q$ is only cohomologous to this cocycle. We therefore assume that we are given a 1-cochain $w : G \to {}_G R$ satisfying $w(1) = 0$ and define:

$$q : G \times G \to R, \, (s,t) \mapsto q(s,t) := x(q'(f(s), f(t)) + \nu(s)w(t) - w(st) + w(s))$$

**Proposition** The map

$$f_A : A_G(\alpha, \beta, q) \to A_{G'}(\alpha', \beta', q'), e_u \otimes x_s \mapsto \eta(uxw(s)) e_{ux} \otimes x'_{f(s)}$$

is a morphism of Yetter-Drinfel'd Hopf algebras.

**Proof.** By definition, $A_{G'}(\alpha', \beta', q')$ is the Yetter-Drinfel'd Hopf algebra arising from the construction described in Paragraph 3.2, using the structure elements $z' \in Z^1(G', \mathrm{Hom}(R_{G'}, R))$, $\gamma' \in Z^1(G', \mathrm{Hom}(R_{G'}, \hat{R}))$ and $\sigma' \in Z^2(G', U(R))$ defined by:

1. $z'_{s'}(u) := u\beta'(s')$

2. $(\gamma'_{s'}(u))(v) := \chi(uv\alpha'(s'))^2$

3. $\sigma'_u(s', t') := \eta(uq'(s', t'))\chi(u^2 \nu'(s')\beta'(s')\alpha'(t'))$

with $\sigma'(s', t') = \sum_{u \in R} \sigma'_u(s', t') e_u$. Using the maps $f_R : R \to R, u \mapsto ux$, $f_G := f$, and $\tau : G \to K^R$ defined by $\tau(s) := \sum_{u \in R} \tau_u(s) e_u$, where

$$\tau_u(s) := \eta(uw(s))$$

we get new structure elements $z \in Z^1(G, \mathrm{Hom}(R_G, R))$, $\gamma \in Z^1(G, \mathrm{Hom}(R_G, \hat{R}))$ and $\sigma \in Z^2(G, U(R))$ via the construction considered in Paragraph 4.2. The cocycle $z$ is determined by:

$$z_s(u) = z'_{f(s)}(f_R(u)) = ux\beta'(f(s)) = u\beta(s)$$

Similarly, $\gamma$ is determined by:

$$(\gamma_s(u))(v) = (\gamma'_{f(s)}(f_R(u)))(v) = \chi(vux\alpha'(f(s)))^2 = \chi(vu\alpha(s))^2$$

Finally $\sigma$, if written in the form $\sigma(s,t) = \sum_{u \in P} \sigma_u(s,t) e_u$, is determined by:

$$\sigma_u(s,t) = \sigma'_{f_P(u)}(f_G(s), f_G(t))\tau_{f_P(s^{-1}.u)}(t)\tau_{f_P(u)}(st)^{-1}\tau_{f_P(u)}(s)$$
$$= \eta(uxq'(f(s), f(t)))\chi((ux)^2 \nu'(f(s))\beta'(f(s))\alpha'(f(t)))$$
$$\quad \eta(u\nu(s)xw(t))\eta(-uxw(st))\eta(uxw(s))$$
$$= \eta(uxq'(f(s), f(t)) + u\nu(s)xw(t) - uxw(st) + uxw(s))\chi(u^2 \nu(s)\beta(s)\alpha(t))$$
$$= \eta(uq(s,t))\chi(u^2 \nu(s)\beta(s)\alpha(t))$$

But these are precisely the structure elements defining $A_G(\alpha, \beta, q)$. Therefore, the assertion follows from Proposition 4.2. $\square$

In the special case $G = G'$, $f = \mathrm{id}_G$, and $x = 1$, this proposition yields:

**Corollary** If $q$ and $q'$ are cohomologous, $A_G(\alpha, \beta, q)$ and $A_G(\alpha, \beta, q')$ are isomorphic.

**4.4** As already noted in Paragraph 4.1, to attach a Yetter-Drinfel'd Hopf algebra to a set of structure elements means to construct a functor, at least if one also cares about morphisms, what we do in this section. Although we are not going to pursue this topic further, we shall, as an example, make this viewpoint explicit in the situation of Paragraph 4.3.

Suppose that $R$ is a finite ring and that $\chi, \eta \in \hat{R}$ are two characters of the additive group of $R$. We suppose that $\chi$ satisfies the condition $\chi(uvw) = \chi(vuw)$ for all $u, v, w \in R$. We define the category $\mathcal{C}_{\chi,\eta}(R)$ as follows: The objects of $\mathcal{C}_{\chi,\eta}(R)$ are quintuples $(G, \nu, \alpha, \beta, q)$ consisting of a finite group $G$, a group homomorphism $\nu : G \to U(R)$ from $G$ to the multiplicative group $U(R)$ of units of $R$, two 1-cocycles $\alpha, \beta \in Z^1(G, {}_G R)$, and a normalized 2-cocycle $q \in Z^2(G, {}_G R)$. Here, ${}_G R$ denotes the additive group of $R$ considered as a left $G$-module as explained in Paragraph 3.3. If $(G', \nu', \alpha', \beta', q')$ is another object of $\mathcal{C}_{\chi,\eta}(R)$, a morphism from $(G, \nu, \alpha, \beta, q)$ to $(G', \nu', \alpha', \beta', q')$ is a triple $(f, x, w)$ consisting of a group homomorphism $f : G \to G'$, an element $x \in Z_{U(R)}(\nu'(f(G)))$ of the centralizer $Z_{U(R)}(\nu'(f(G)))$ of $\nu'(f(G))$ in $U(R)$, and a 1-cochain $w : G \to {}_G R$ satisfying $w(1) = 0$ such that

1. $\forall s \in G : \nu(s) = \nu'(f(s))$
2. $\forall s \in G : \alpha(s) = x\alpha'(f(s))$
3. $\forall s \in G : \beta(s) = x\beta'(f(s))$
4. $\forall s, t \in G : q(s,t) := x(q'(f(s), f(t)) + \nu(s)w(t) - w(st) + w(s))$

If $(f', x', w') : (G', \nu', \alpha', \beta', q') \to (G'', \nu'', \alpha'', \beta'', q'')$ is another morphism, we define the composition as:
$$(f', x', w') \circ (f, x, w) := (f' \circ f, xx', x'^{-1}w + f^1(w'))$$

This is in fact a morphism, since we obviously have $\nu(s) = \nu''(f'(f(s)))$, $\alpha(s) = xx'\alpha''(f'(f(s)))$, and $\beta(s) = xx'\beta''(f'(f(s)))$, but also
$$q(s,t) := xx'(q''(f'(f(s)), f'(f(t))) + \nu(s)w'(f(t)) - w'(f(st)) + w'(f(s))$$
$$+ \nu(s)x'^{-1}w(t) - x'^{-1}w(st) + x'^{-1}w(s))$$

Note that, since $x' \in Z_{U(R)}(\nu'(G')) \subset Z_{U(R)}(\nu(G))$, we have $xx' \in Z_{U(R)}(\nu(G))$. Since
$$(f'', x'', w'') \circ ((f', x', w') \circ (f, x, w))$$
$$= (f'' \circ f' \circ f, xx'x'', x''^{-1}x'^{-1}w + x''^{-1}f^1(w') + f^1(f'^1(w'')))$$
$$= ((f'', x'', w'') \circ (f', x', w')) \circ (f, x, w)$$

this composition law is associative. The identity homomorphism of the object $(G, \nu, \alpha, \beta, q)$ is $(\mathrm{id}_G, 1_R, 0)$.

Now consider the Hopf algebra $H := K[R]$, i. e., the group ring of the additive group of $R$. For $u \in R$, we denote the corresponding canonical basis element of $K[R]$ by $c_u$ and the corresponding primitive idempotent of $K^R$ by $e_u$. If $\mathcal{YDH}(H)$ denotes the category of Yetter-Drinfel'd Hopf algebras over $H$, the construction carried out in Paragraph 3.3 yields a functor

$$\mathcal{C}_{\chi,\eta}(R) \to \mathcal{YDH}(H), (G, \nu, \alpha, \beta, q) \mapsto A_G(\alpha, \beta, q)$$

if we assign to the morphism

$$(f, x, w) : (G, \nu, \alpha, \beta, q) \to (G', \nu', \alpha', \beta', q')$$

in $\mathcal{C}_{\chi,\eta}(R)$ the morphism

$$f_A : A_G(\alpha, \beta, q) \to A_{G'}(\alpha', \beta', q'), e_u \otimes x_s \mapsto \eta(uxw(s))e_{ux} \otimes x'_{f(s)}$$

in $\mathcal{YDH}(H)$ constructed in Paragraph 4.3. Here $x_s$, resp. $x'_{s'}$, denotes the canonical basis elements of $K[G]$, resp. $K[G']$, where $s \in G$ and $s' \in G'$. This assignment is compatible with composition, because, if $f_{A'}$ denotes the morphism assigned to $(f', x', w') : (G', \nu', \alpha', \beta', q') \to (G'', \nu'', \alpha'', \beta'', q'')$, we have:

$$f_{A'} \circ f_A(e_u \otimes x_s) = \eta(uxx'(x'^{-1}w(s) + w'(f(s))))e_{uxx'} \otimes x''_{f'(f(s))}$$

and therefore $f_{A'} \circ f_A$ is the morphism assigned to $(f', x', w') \circ (f, x, w)$. Obviously, identity morphisms are assigned to identity morphisms.

**4.5** Now suppose that $p$ is an odd prime and that $R = \mathbb{Z}_p$. We shall prove now that isomorphisms between Yetter-Drinfel'd Hopf algebras of the type $A_G(\alpha, \beta, q)$ over $H := K[\mathbb{Z}_p]$ described in Paragraph 3.4 are necessarily of the form presented in the previous two paragraphs.

To consider two such algebras, we suppose that we are given two finite groups $G$ and $G'$, two group homomorphisms $\nu : G \to \mathbb{Z}_p^\times$ and $\nu' : G' \to \mathbb{Z}_p^\times$, two cocycles $\alpha, \beta \in Z^1(G, {}_G\mathbb{Z}_p)$ as well as two cocycles $\alpha', \beta' \in Z^1(G', {}_{G'}\mathbb{Z}_p)$, and two normalized cocycles $q \in Z^2(G, {}_G\mathbb{Z}_p), q' \in Z^2(G', {}_{G'}\mathbb{Z}_p)$. We define the characters

$$\chi : \mathbb{Z}_p \to K, i \mapsto \zeta^{i/2} \qquad \eta : \mathbb{Z}_p \to K, i \mapsto \zeta^i$$

where $\zeta$ is a fixed primitive $p$-th root of unity. We note that, since $p$ is odd, 2 is an invertible element of $\mathbb{Z}_p$, and therefore the expression $i/2$ for $i \in \mathbb{Z}_p$ makes sense. We denote by $e_i$, for $i \in \mathbb{Z}_p$, the primitive idempotents of $K^{\mathbb{Z}_p}$, whereas $c_i$, $x_s$, resp. $x'_{s'}$ denote the canonical basis elements of $K[\mathbb{Z}_p]$, $K[G]$, resp. $K[G']$.

As in Paragraph 1.10, we introduce the notation

$$\phi(e_i \otimes x_s) := c_1 \to (e_i \otimes x_s) \qquad \psi(e_i \otimes x_s) := \gamma((e_i \otimes x_s)^{(1)})(e_i \otimes x_s)^{(2)}$$

where $\gamma \in G(H^*)$ is the character which is uniquely determined by the condition $\gamma(c_1) = \zeta$. In a similar way, we introduce endomorphisms $\phi'$ and $\psi'$ of $A_{G'}(\alpha', \beta', q')$.

Now suppose that $f_A : A_G(\alpha, \beta, q) \to A_{G'}(\alpha', \beta', q')$ is an isomorphism.

**Proposition** Suppose that $\alpha \neq 0$ or $\beta \neq 0$. Then there is an element $k \in \mathbb{Z}_p^\times$, a group isomorphism $f : G \to G'$, and a 1-cochain $w : G \to {}_G\mathbb{Z}_p$ satisfying $w(1) = 0$ such that

1. $\forall s \in G : \nu(s) = \nu'(f(s))$
2. $\forall s \in G : \alpha(s) = k\alpha'(f(s))$
3. $\forall s \in G : \beta(s) = k\beta'(f(s))$
4. $\forall s, t \in G : q(s,t) = k(q'(f(s), f(t)) + \nu(s)w(t) - w(st) + w(s))$
5. $\forall i \in \mathbb{Z}_p \, \forall s \in G : f_A(e_i \otimes x_s) = \zeta^{kiw(s)} e_{ki} \otimes x'_{f(s)}$

In particular, $q$ and $kf^2(q')$ are cohomologous.

**Proof.** (1) In the algebras $A_G(\alpha, \beta, q)$ resp. $A_{G'}(\alpha', \beta', q')$, we consider the elements
$$u := \sum_{j=0}^{p-1} \zeta^j e_j \otimes x_1 \qquad u' := \sum_{j=0}^{p-1} \zeta^j e_j \otimes x'_1$$

We then have
$$u(e_i \otimes x_s) = \sum_{j=0}^{p-1} \zeta^j (e_j \otimes x_1)(e_i \otimes x_s) = \sum_{j=0}^{p-1} \zeta^j \delta_{j\nu(1),i}(e_j \otimes x_s) = \zeta^i e_i \otimes x_s$$

and similarly $u'(e_i \otimes x'_{s'}) = \zeta^i e_i \otimes x'_{s'}$. Since $\phi(e_i \otimes x_s) = \zeta^{i\alpha(s)} e_i \otimes x_s$ and $\psi(e_i \otimes x_s) = \zeta^{i\beta(s)} e_i \otimes x_s$, we see that:
$$\phi(e_i \otimes x_s) = u^{\alpha(s)} e_i \otimes x_s \qquad \psi(e_i \otimes x_s) = u^{\beta(s)} e_i \otimes x_s$$

Similarly, we have:
$$\phi'(e_i \otimes x'_{s'}) = u'^{\alpha'(s')} e_i \otimes x'_{s'} \qquad \psi'(e_i \otimes x'_{s'}) = u'^{\beta'(s')} e_i \otimes x'_{s'}$$

(2) The grouplike elements of the dual group ring $K^{\mathbb{Z}_p}$ are the Fourier transformed elements $\sum_{j=0}^{p-1} \zeta^{ij} e_j$ of the idempotents $e_j$. Since the coalgebra structure of $A_G(\alpha, \beta, q)$ is the ordinary tensor product coalgebra structure, the elements $\sum_{j=0}^{p-1} \zeta^{ij} e_j \otimes x_s$ form a basis of $A_G(\alpha, \beta, q)$ consisting of grouplike elements. Since the powers of $u$ are given by the formula $u^i = \sum_{j=0}^{p-1} \zeta^{ij} e_j \otimes x_1$, these elements can be written in the form:

$$\sum_{j=0}^{p-1} \zeta^{ij} e_j \otimes x_s = u^i(1 \otimes x_s)$$

$A_{G'}(\alpha', \beta', q')$ has a similar basis consisting of grouplike elements. Since $f_A$ takes grouplike elements to grouplike elements, there exists, for every $s \in G$, an element $w(s) \in \mathbb{Z}_p$ and an element $f(s) \in G'$ such that

$$f_A(1 \otimes x_s) = u'^{w(s)}(1 \otimes x'_{f(s)})$$

Now the linearity of $f_A$ over $H$ implies:

$$f_A(u^{\alpha(s)})u'^{w(s)}(1 \otimes x'_{f(s)}) = f_A(u^{\alpha(s)}(1 \otimes x_s)) = f_A(\phi(1 \otimes x_s))$$
$$= \phi'(f_A(1 \otimes x_s)) = u'^{\alpha'(f(s))}u'^{w(s)}(1 \otimes x'_{f(s)})$$

Since grouplike elements are invertible, this implies:

$$f_A(u^{\alpha(s)}) = u'^{\alpha'(f(s))}$$

Similarly, the colinearity of $f_A$ implies $f_A(u^{\beta(s)}) = u'^{\beta'(f(s))}$. Now suppose that $\alpha \neq 0$. Then there exists an element $s \in G$ such that $\alpha(s) \neq 0$, and, since $f_A$ is injective, we also have $\alpha'(f(s)) \neq 0$. Since $\alpha'(f(s))$ generates $\mathbb{Z}_p$, there exists $k \in \mathbb{Z}_p^\times$ such that $f_A(u^k) = u'$. By a similar reasoning, this also holds if $\beta \neq 0$. The above equations now yield:

$$\alpha(s) = k\alpha'(f(s)) \qquad \beta(s) = k\beta'(f(s))$$

This proves the second and the third assertion.

(3) We have:

$$(e_i \otimes x_s)u = \sum_{j=0}^{p-1} \zeta^j(e_i \otimes x_s)(e_j \otimes x_1) = \sum_{j=0}^{p-1} \zeta^j \delta_{i\nu(s),j} e_i \otimes x_s$$
$$= \zeta^{i\nu(s)} e_i \otimes x_s = u^{\nu(s)}(e_i \otimes x_s)$$

By a similar calculation, we have $(e_i \otimes x'_{s'})u' = u'^{\nu'(s')}(e_i \otimes x'_{s'})$. The first equation implies $(e_i \otimes x_s)u^k = u^{k\nu(s)}(e_i \otimes x_s)$. Summing over $i$ and applying $f_A$, we get $u'^{w(s)}(1 \otimes x'_{f(s)})u' = u'^{\nu(s)}u'^{w(s)}(1 \otimes x'_{f(s)})$. Since on the other hand $u'^{w(s)}(1 \otimes x'_{f(s)})u' = u'^{\nu'(f(s))}u'^{w(s)}(1 \otimes x'_{f(s)})$, we see that $u'^{\nu(s)} = u'^{\nu'(f(s))}$ and therefore $\nu(s) = \nu'(f(s))$. This proves the first assertion.

(4) By inverting the discrete Fourier transform above, we get:
$$e_i \otimes x_s = \frac{1}{p}\sum_{j=0}^{p-1} \zeta^{-ij} u^j (1 \otimes x_s)$$

Therefore, we have:
$$f_A(e_i \otimes x_s) = \frac{1}{p}\sum_{j=0}^{p-1} \zeta^{-ij} f_A(u^j(1 \otimes x_s)) = \frac{1}{p}\sum_{j=0}^{p-1} \zeta^{-ij} u'^{j/k} u'^{w(s)} (1 \otimes x'_{f(s)})$$
$$= u'^{w(s)} \frac{1}{p}\sum_{j=0}^{p-1} \zeta^{-kij} u'^j (1 \otimes x'_{f(s)}) = u'^{w(s)} e_{ki} \otimes x'_{f(s)}$$
$$= \zeta^{kiw(s)} e_{ki} \otimes x'_{f(s)}$$

This proves that $f_A$ has the form given in the fifth assertion.

(5) We have:
$$f_A(e_i \otimes x_s) f_A(e_j \otimes x_t) = \zeta^{kiw(s)+kjw(t)}(e_{ki} \otimes x'_{f(s)})(e_{kj} \otimes x'_{f(t)})$$
$$= \zeta^{kiw(s)+kjw(t)} \delta_{ki\nu'(f(s)),kj} \zeta^{kiq'(f(s),f(t))} \zeta^{k^2 i^2 \nu'(f(s))\beta'(f(s))\alpha'(f(t))/2}$$
$$(e_{ki} \otimes x'_{f(s)f(t)})$$
$$= \zeta^{ki(w(s)+\nu(s)w(t))} \delta_{i\nu(s),j} \zeta^{kiq'(f(s),f(t))} \zeta^{i^2\nu(s)\beta(s)\alpha(t)/2}(e_{ki} \otimes x'_{f(s)f(t)})$$

On the other hand, we have:
$$f_A((e_i \otimes x_s)(e_j \otimes x_t)) = \delta_{i\nu(s),j} \zeta^{iq(s,t)} \zeta^{i^2\nu(s)\beta(s)\alpha(t)/2} \zeta^{kiw(st)}(e_{ki} \otimes x'_{f(st)})$$

Since $f_A$ is an algebra homomorphism, both expressions must be equal. Therefore, we see that $f(st) = f(s)f(t)$, i. e., $f$ is a group homomorphism, and we have $q(s,t) = k(q'(f(s),f(t)) + \nu(s)w(t) - w(st) + w(s))$. This proves the fourth assertion. Since $q$ and $q'$ are normalized, we get, by inserting $s = t = 1$ in the previous equation, the fact that $w(1) = 0$. $\square$

**4.6** In this paragraph, we continue our analysis of the case $R = \mathbb{Z}_p$ for an odd prime $p$, but specialize the situation further to the case where also $G = \mathbb{Z}_p$. We want to determine the isomorphism classes of nontrivial Yetter-Drinfel'd Hopf algebras among the algebras $A_p(\alpha, \beta, q)$. We shall use the notation of the previous paragraph.

The algebras $A_p(\alpha, \beta, q)$ are defined with respect to the characters
$$\chi : \mathbb{Z}_p \to K, i \mapsto \zeta^{i/2} \qquad \eta : \mathbb{Z}_p \to K, i \mapsto \zeta^i$$

and a group homomorphism
$$\nu : \mathbb{Z}_p \to \mathbb{Z}_p^\times$$

Obviously, the only group homomorphism of this type is the trivial homomorphism that is identically equal to one. Therefore, in this situation, $\mathbb{Z}_p$ is a trivial $\mathbb{Z}_p$-module, and the 1-cocycles $\alpha$ and $\beta$ are ordinary group homomorphisms. If we require that $A_p(\alpha, \beta, q)$ be nontrivial, $\alpha$ and $\beta$ have to be nonzero by Proposition 1.11.

For $m \in \mathbb{Z}_p^\times$, we denote by $\alpha_m$ the group homomorphism

$$\alpha_m : \mathbb{Z}_p \to \mathbb{Z}_p, i \mapsto mi$$

Then, $\alpha_1, \ldots, \alpha_{p-1}$ are all nonzero group homomorphisms from $\mathbb{Z}_p$ to itself. From Proposition 1.13.2, we know that the second cohomology group of the trivial $\mathbb{Z}_p$-module $H^2(\mathbb{Z}_p, \mathbb{Z}_p)$ is isomorphic to $\mathbb{Z}_p$. We choose a complete system of representatives for these cohomology classes, i. e., cocycles $q_0, \ldots, q_{p-1} \in Z^2(\mathbb{Z}_p, \mathbb{Z}_p)$ such that $H^2(\mathbb{Z}_p, \mathbb{Z}_p) = \{\bar{q}_0, \ldots, \bar{q}_{p-1}\}$.

**Proposition**
1. Suppose that $\alpha, \beta \in \operatorname{Hom}(\mathbb{Z}_p, \mathbb{Z}_p)$ are nonzero group homomorphisms and that $q \in Z^2(\mathbb{Z}_p, \mathbb{Z}_p)$ is a 2-cocycle. Define $m \in \mathbb{Z}_p^\times$ by $m := \alpha(1)/\beta(1)$. Suppose that $q/\beta(1)$ is cohomologous to $q_n$, for $n \in \mathbb{Z}_p$. Then $A_p(\alpha, \beta, q)$ is isomorphic to $A_p(\alpha_m, \operatorname{id}, q_n)$.

2. The Yetter-Drinfel'd Hopf algebras $A_p(\alpha_m, \operatorname{id}, q_n)$, for $m = 1, \ldots, p-1$ and $n = 0, \ldots, p-1$, are mutually nonisomorphic.

**Proof.** To prove the first statement, define $k := \beta(1)$. Since group endomorphisms of $\mathbb{Z}_p$ are determined by their value on 1, we then have $\beta = k \operatorname{id}$ and, since $k\alpha_m(1) = km = \alpha(1)$, we also have $\alpha = k\alpha_m$. Since $q$ is cohomologous to $kq_n$, there is a 1-cochain $w : \mathbb{Z}_p \to \mathbb{Z}_p$ such that

$$q(i, j) = kq_n(i, j) + kw(j) - kw(i+j) + kw(i)$$

It now follows from Proposition 4.3 that the map

$$f_A : A_p(\alpha, \beta, q) \to A_p(\alpha_m, \operatorname{id}, q_n), e_i \otimes c_j \mapsto \zeta^{kiw(j)} e_{ki} \otimes c_j$$

is an isomorphism.

To prove the second statement, suppose that

$$f_A : A_p(\alpha_m, \operatorname{id}, q_n) \to A_p(\alpha_{m'}, \operatorname{id}, q_{n'})$$

is an isomorphism. From Proposition 4.5, we know that there exists an element $k \in \mathbb{Z}_p^\times$, a group isomorphism $f : \mathbb{Z}_p \to \mathbb{Z}_p$, and a 1-cochain $w : \mathbb{Z}_p \to \mathbb{Z}_p$ satisfying $w(0) = 0$ such that $f_A$ has the form

$$f_A(e_i \otimes c_j) = \zeta^{kiw(j)} e_{ki} \otimes c_{f(j)}$$

Moreover, the 1-cocycles are related by

$$\alpha_m(i) = k\alpha_{m'}(f(i)) \qquad i = kf(i)$$

for all $i \in \mathbb{Z}_p$, wheras the 2-cocycles are related by

$$q_n(i,j) = k(q_{n'}(f(i), f(j))) + w(j) - w(ij) + w(i))$$

for all $i, j \in \mathbb{Z}_p$. From the first relation, we see that $\alpha_m = \alpha_{m'}$, and therefore we have $m = m'$. From the second relation, we see that $q_n$ and $kf^2(q_{n'})$ are cohomologous. But this implies, by Proposition 1.13.2, that $q_n$ and $q_{n'}$ are cohomologous, and therefore we have $n = n'$. $\square$

Therefore, the algebras $A_p(\alpha, \beta, q)$ fall into $p(p-1)$ isomorphism classes, which are represented by the algebras $A_p(\alpha_m, \mathrm{id}, q_n)$, for $m = 1, \ldots, p-1$ and $n = 0, \ldots, p-1$. Note that, by Proposition 1.13.2, these algebras are also commutative.

**4.7** The framework considered in Paragraph 4.3 can also be applied to the case $p = 2$. Suppose that $\iota$ is a primitive fourth root of unity. We denote by $\hat{\iota} : \mathbb{Z}_2 \to \mathbb{Z}_4$ the unique injective group homomorphism and by $\hat{\pi} : \mathbb{Z}_4 \to \mathbb{Z}_2$ the unique surjective group homomorphism.

Suppose that $G'$ is a finite group. We regard $\mathbb{Z}_2$ as a trivial $G'$-module. Suppose that $\alpha' : G' \to \mathbb{Z}_2$ and $\beta' : G' \to \mathbb{Z}_2$ are 1-cocycles; since the $G'$-module structure is trivial, these are just group homomorphisms. In addition, suppose that $q' \in Z^2(G', \mathbb{Z}_4)$ is a normalized 2-cocycle of the trivial $G'$-module $\mathbb{Z}_4$ satisfying

$$\hat{\pi} \circ q' = \beta' \cup \alpha'$$

For $i \in \mathbb{Z}_2$, we denote the corresponding primitive idempotentin $K^{\mathbb{Z}_2}$ by $e_i$ and the corresponding canonical basis vector of $K[\mathbb{Z}_2]$ by $c_i$; the canonical basis vectors of $K[G']$ are denoted by $x'_{s'}$, for $s' \in G'$. We define $\sigma'_i : G' \times G' \to K^{\mathbb{Z}_2}$ by

$$\sigma'_0(s', t') := 1 \qquad \sigma'_1(s', t') := \iota^{q'(s', t')}$$

and set:

$$\sigma'(s', t') := \sigma'_0(s', t') e_0 + \sigma'_1(s', t') e_1$$

As explained in Paragraph 3.5, these data can be used to construct a Yetter-Drinfel'd Hopf algebra over $H := K[\mathbb{Z}_2]$, which we denote by $A'$.

Now suppose that $G$ is another finite group and that $f : G \to G'$ is a group homomorphism. We define:

$$\alpha : G \to \mathbb{Z}_2, s \mapsto \alpha(s) := \alpha'(f(s)) \qquad \beta : G \to \mathbb{Z}_2, s \mapsto \beta(s) := \beta'(f(s))$$

We could also define for all $s, t \in G$

$$q(s, t) := q'(f(s), f(t))$$

to get a 2-cocycle; however, we shall consider a more general case where $q$ is only cohomologous to this cocycle. Before we really define $q$, we remark that any 2-cocycle $q \in Z^2(G, \mathbb{Z}_4)$ that satisfies $\hat{\pi}_2(q) = \beta \cup \alpha$ obviously has the property that:

$$\hat{\pi}_2(q) = f^1(\beta') \cup f^1(\alpha') = f^2(\beta' \cup \alpha') = f^2(\hat{\pi}_2(q')) = \hat{\pi}_2(f^2(q'))$$

Therefore, we see that there is a 2-cocycle $q'' \in Z^2(G, \mathbb{Z}_2)$ such that $q - f^2(q') = \hat{\imath}_2(q'')$. The cocycles $q$ that we want to consider are those for which $q''$ is even a coboundary.

We therefore assume that we are given a 1-cochain $w : G \to \mathbb{Z}_2$ satisfying $w(1) = 0$ and define

$$q : G \times G \to \mathbb{Z}_2, \ (s, t) \mapsto q(s, t) := q'(f(s), f(t)) + \hat{\imath}(w(t) - w(st) + w(s))$$

From $\alpha$, $\beta$, and $q$ we can construct another Yetter-Drinfel'd Hopf algebra over $H$, which we denote by $A$.

**Proposition** The map

$$f_A : A \to A', e_i \otimes x_s \mapsto (-1)^{iw(s)} e_i \otimes x'_{f(s)}$$

is a morphism of Yetter-Drinfel'd Hopf algebras.

**Proof.** By definition, $A'$ is the Yetter-Drinfel'd Hopf algebra arising from the construction described in Paragraph 3.2, using the structure elements $z' \in Z^1(G', \text{Hom}(\mathbb{Z}_2, \mathbb{Z}_2))$ and $\gamma' \in Z^1(G', \text{Hom}(\mathbb{Z}_2, \hat{\mathbb{Z}}_2))$ defined as:

$$z'_{s'}(i) := i\beta'(s') \qquad (\gamma'_{s'}(i))(j) := (-1)^{ij\alpha'(s')}$$

and $\sigma' \in Z^2(G', U(K^{\mathbb{Z}_2}))$ as defined above. Define, for $i \in \mathbb{Z}_2$ and $s \in G$, $\tau_i(s) := (-1)^{iw(s)}$. Using the maps $f_{\mathbb{Z}_2} := \text{id}$ and $f_G := f$, we get new structure elements $z \in Z^1(G, \text{Hom}(\mathbb{Z}_2, \mathbb{Z}_2))$, $\gamma \in Z^1(G, \text{Hom}(\mathbb{Z}_2, \hat{\mathbb{Z}}_2))$ and $\sigma \in Z^2(G, U(K^{\mathbb{Z}_2}))$ via the construction considered in Paragraph 4.2. These structure elements can be expressed as follows:

1. $z_s(i) := i\beta(s)$
2. $(\gamma_s(i))(j) := (-1)^{ij\alpha(s)}$
3. $\sigma_0(s, t) = 1 \qquad \sigma_1(s, t) = \iota^{q(s,t)}$

where $\sigma(s, t) = \sum_{i=0}^{1} \sigma_i(s, t) e_i$.

Here, the last equation follows from the definition of $\sigma_i$:

$$\sigma_i(s, t) := \sigma'_i(f(s), f(t)) \tau_i(t) \tau_i(st)^{-1} \tau_i(s)$$

This obviously yields $\sigma_0(s, t) = 1$, whereas for $i = 1$ it reads:

$$\sigma_1(s, t) = \iota^{q'(f(s), f(t))} (-1)^{w(t) - w(st) + w(s)}$$
$$= \iota^{q'(f(s), f(t)) + \hat{\imath}(w(t) - w(st) + w(s))} = \iota^{q(s,t)}$$

Now, the assertion follows from Proposition 4.2. $\square$

**4.8**  In the preceding paragraph, we have constructed homomorphisms of Yetter-Drinfel'd Hopf algebras in the case $p = 2$, which, in special cases, may be isomorphisms. We now consider the question whether all isomorphisms are necessarily of this form. We continue to use the notation of the preceding paragraph.

Suppose that $G$ and $G'$ are finite groups. We regard $\mathbb{Z}_2$ as a trivial $G$-module, resp. as a trivial $G'$-module. Suppose that $\alpha : G \to \mathbb{Z}_2$ and $\beta : G \to \mathbb{Z}_2$ are 1-cocycles, i. e., group homomorphisms. Similarly, suppose that $\alpha' : G' \to \mathbb{Z}_2$ and $\beta' : G' \to \mathbb{Z}_2$ are 1-cocycles. In addition, suppose that $q \in Z^2(G, \mathbb{Z}_4)$ and $q' \in Z^2(G', \mathbb{Z}_4)$ are normalized 2-cocycles satisfying:

$$\hat{\pi} \circ q = \beta \cup \alpha \qquad \hat{\pi} \circ q' = \beta' \cup \alpha'$$

We define $\sigma_i : G \times G \to K^{\mathbb{Z}_2}$ by

$$\sigma_0(s,t) := 1 \qquad \sigma_1(s,t) := \iota^{q(s,t)}$$

and set:

$$\sigma(s,t) := \sigma_0(s,t)e_0 + \sigma_1(s,t)e_1$$

In a similar way, we also define $\sigma'_0$, $\sigma'_1$, and $\sigma'$. As explained in Paragraph 3.5, these data can be used to construct two Yetter-Drinfel'd Hopf algebras over $H := K[\mathbb{Z}_2]$, which we denote by $A$ and $A'$.

As in Paragraph 1.10, we introduce the notation

$$\phi(e_i \otimes x_s) := c_1 \to (e_i \otimes x_s) \qquad \psi(e_i \otimes x_s) := \gamma((e_i \otimes x_s)^{(1)})(e_i \otimes x_s)^{(2)}$$

where $\gamma \in G(H^*)$ is the unique nontrivial character of $H$; it satisfies $\gamma(c_1) = -1$. Similarly, we introduce endomorphisms $\phi'$ and $\psi'$ of $A'$.

Now suppose that $f_A : A \to A'$ is an isomorphism.

**Proposition**  Suppose that $\alpha \neq 0$ or $\beta \neq 0$. Then there is a group isomorphism $f : G \to G'$ and a 1-cochain $w : G \to \mathbb{Z}_2$ satisfying $w(1) = 0$ such that

1. $\forall s \in G : \alpha(s) = \alpha'(f(s))$
2. $\forall s \in G : \beta(s) = \beta'(f(s))$
3. $\forall s, t \in G : q(s,t) - q'(f(s), f(t)) = \hat{\iota}(w(t) - w(st) + w(s))$
4. $\forall i \in \mathbb{Z}_2 \, \forall s \in G : f_A(e_i \otimes x_s) = (-1)^{iw(s)} e_i \otimes x'_{f(s)}$

**Proof.**  (1) We follow the line of reasoning in Paragraph 4.5. In the algebras $A$ resp. $A'$, we consider the elements

$$u := \sum_{j=0}^{1}(-1)^j e_j \otimes x_1 \qquad u' := \sum_{j=0}^{1}(-1)^j e_j \otimes x'_1$$

We then have

$$u(e_i \otimes x_s) = \sum_{j=0}^{1}(-1)^j(e_j \otimes x_1)(e_i \otimes x_s)$$

$$= \sum_{j=0}^{1}(-1)^j \delta_{ji}\sigma_j(1,s)(e_j \otimes x_s) = (-1)^i e_i \otimes x_s$$

and similarly $u'(e_i \otimes x'_{s'}) = (-1)^i e_i \otimes x'_{s'}$. Since $\phi(e_i \otimes x_s) = (-1)^{i\alpha(s)} e_i \otimes x_s$ and $\psi(e_i \otimes x_s) = (-1)^{i\beta(s)} e_i \otimes x_s$, we see that:

$$\phi(e_i \otimes x_s) = u^{\alpha(s)} e_i \otimes x_s \qquad \psi(e_i \otimes x_s) = u^{\beta(s)} e_i \otimes x_s$$

Similarly, we have:

$$\phi'(e_i \otimes x'_{s'}) = u'^{\alpha'(s')} e_i \otimes x'_{s'} \qquad \psi'(e_i \otimes x'_{s'}) = u'^{\beta'(s')} e_i \otimes x'_{s'}$$

(2) The grouplike elements of the dual group ring $K^{\mathbb{Z}_2}$ are the Fourier transformed elements $\sum_{j=0}^{1}(-1)^{ij}e_j$ of the idempotents $e_j$. Since the coalgebra structure of $A$ is the ordinary tensor product coalgebra structure, the elements of the form $\sum_{j=0}^{1}(-1)^{ij}e_j \otimes x_s$ constitute a basis of $A$ consisting of grouplike elements. These elements can be written in the form:

$$\sum_{j=0}^{1}(-1)^{ij}e_j \otimes x_s = u^i(1 \otimes x_s)$$

$A'$ has a similar basis consisting of grouplike elements. Since $f_A$ takes grouplike elements to grouplike elements, there exists, for every $s \in G$, an element $w(s) \in \mathbb{Z}_2$ and an element $f(s) \in G'$ such that

$$f_A(1 \otimes x_s) = u'^{w(s)}(1 \otimes x'_{f(s)})$$

Now the linearity of $f_A$ over $H$ implies:

$$f_A(u^{\alpha(s)})u'^{w(s)}(1 \otimes x'_{f(s)}) = f_A(u^{\alpha(s)}(1 \otimes x_s)) = f_A(\phi(1 \otimes x_s))$$
$$= \phi'(f_A(1 \otimes x_s)) = u'^{\alpha'(f(s))} u'^{w(s)}(1 \otimes x'_{f(s)})$$

Since grouplike elements are invertible, this implies:

$$f_A(u^{\alpha(s)}) = u'^{\alpha'(f(s))}$$

Similarly, the colinearity of $f_A$ implies $f_A(u^{\beta(s)}) = u'^{\beta'(f(s))}$. Now suppose that $\alpha \neq 0$. Then there exists an element $s \in G$ such that $\alpha(s) = 1$, and, since $f_A$ is injective, we also have $\alpha'(f(s)) = 1$. Therefore, we have $f_A(u) = u'$. By a similar reasoning, this also holds if $\beta \neq 0$. The above equations now yield:

$$\alpha(s) = \alpha'(f(s)) \qquad \beta(s) = \beta'(f(s))$$

This proves the first and the second assertion.

(3) By inverting the discrete Fourier transform above, we get:

$$e_i \otimes x_s = \frac{1}{2}\sum_{j=0}^{1}(-1)^{ij}u^j(1\otimes x_s)$$

Therefore, we have:

$$f_A(e_i \otimes x_s) = \frac{1}{2}\sum_{j=0}^{1}(-1)^{ij}f_A(u^j(1\otimes x_s)) = \frac{1}{2}\sum_{j=0}^{1}(-1)^{ij}u'^j u'^{w(s)}(1\otimes x'_{f(s)})$$

$$= u'^{w(s)}\frac{1}{2}\sum_{j=0}^{1}(-1)^{ij}u'^j(1\otimes x'_{f(s)}) = u'^{w(s)}e_i \otimes x'_{f(s)}$$

$$= (-1)^{iw(s)}e_i \otimes x'_{f(s)}$$

This proves that $f_A$ has the form given in the fourth assertion.

(4) We have:

$$f_A(e_1 \otimes x_s)f_A(e_1 \otimes x_t) = (-1)^{w(s)+w(t)}(e_1 \otimes x'_{f(s)})(e_1 \otimes x'_{f(t)})$$
$$= (-1)^{w(s)+w(t)}\iota^{q'(f(s),f(t))}(e_1 \otimes x'_{f(s)f(t)})$$

On the other hand, we have:

$$f_A((e_1 \otimes x_s)(e_1 \otimes x_t)) = (-1)^{w(st)}\iota^{q(s,t)}(e_1 \otimes x'_{f(st)})$$

Since $f_A$ is an algebra homomorphism, both expressions must be equal. Therefore, we see that $f(st) = f(s)f(t)$, i. e., $f$ is a group homomorphism, and we have $q(s,t) - q'(f(s),f(t)) = \hat{\iota}(w(t) - w(st) + w(s))$. This proves the third assertion. Since $q$ and $q'$ are normalized, we get, by inserting $s = t = 1$ in the previous equation, the fact that $w(1) = 0$. □

**4.9** As an application of the preceding considerations, we consider the Yetter-Drinfel'd Hopf algebras $A_+$ and $A_-$ constructed in Paragraph 3.6. We keep the notation of Paragraph 4.7 and Paragraph 4.8.

**Proposition** $A_+$ and $A_-$ are not isomorphic.

**Proof.** As in Paragraph 3.6, we denote the cocycles used in the definition of $A_+$ resp. $A_-$ by $q_+$ resp. $q_-$. We have $q_- - q_+ = \hat{\iota}_2(q)$, where $q \in Z^2(\mathbb{Z}_2,\mathbb{Z}_2)$ is the 2-cocycle of the trivial $\mathbb{Z}_2$-module $\mathbb{Z}_2$ defined as:

$$q(i,j) := \begin{cases} 0 & \text{if } i = 0 \text{ or } j = 0 \\ 1 & \text{if } i = 1 \text{ and } j = 1 \end{cases}$$

If $f_A : A_- \to A_+$ were an isomorphism, we know from Proposition 4.8 that it would have the form

$$f_A(e_i \otimes c_j) = (-1)^{iw(j)} e_i \otimes c_j$$

for some 1-cochain $w : \mathbb{Z}_2 \to \mathbb{Z}_2$ satisfying $w(0) = 0$ and

$$q(i,j) = w(j) - w(i+j) + w(i)$$

for all $i, j \in \mathbb{Z}_2$. Therefore, $q$ were the coboundary arising from $w$. However, this is not the case, since from $w(0) = 0$ we get, as in Paragraph 3.6, that $1 = q(1,1) = w(1) + w(1) = 0$, which is a contradiction. $\square$

# 5 Constructions

**5.1** Yetter-Drinfel'd Hopf algebras can be used to construct ordinary Hopf algebras. We have already seen one of these constructions in Paragraph 1.6, namely the Radford biproduct construction. The second construction from [78] is another example of such a construction; it yields a Hopf algebra in which the Radford biproduct appears as a kind of Borel subalgebra.

In this section, we apply these constructions to the Yetter-Drinfel'd Hopf algebra $A_G(\alpha, \beta, q)$ considered in Paragraph 3.3. We describe the resulting Hopf algebras by exhibiting a basis for which the structure elements can be written down explicitly. Afterwards, we determine when the resulting Hopf algebras are semisimple.

In the whole section, we assume that $G$ is a finite group and that $R$ is a finite ring. As in Paragraph 3.3, we assume that we are given a group homomorphism

$$\nu : G \to U(R)$$

from $G$ to the multiplicative group $U(R)$ of units of $R$, and use it to introduce the left $G$-module structures $_GR$ and $R_G$ on $R$ described there. Also, we assume that the following additional structure elements are given:

1. Two 1-cocycles $\alpha, \beta \in Z^1(G, {}_GR)$.

2. A normalized 2-cocycle $q \in Z^2(G, {}_GR)$.

3. Two characters $\chi, \eta \in \hat{R}$ of the additive group of $R$, where $\chi$ is required to satisfy $\chi(uvw) = \chi(vuw)$ for all $u, v, w \in R$.

From Paragraph 5.5 on, where we begin to consider the second construction, we require that $\alpha$ and $\beta$ are compatible in the sense that we have

$$\chi(u\alpha(s)\beta(t)) = \chi(u\beta(s)\alpha(t))$$

for all $s, t \in G$ and all $u \in R$.

We use the notation $H := K[R]$ for the group ring of the additive group of $R$ and $A := A_G(\alpha, \beta, q) = K^R \otimes K[G]$ for the tensor product of the dual group ring of the additive group of $R$ and the group ring of $G$, considered as a Yetter-Drinfel'd Hopf algebra with the structure elements described in Proposition 3.3. The canonical basis elements of $K[R]$ resp. $K[G]$ are denoted by $c_u$ resp. $x_s$, where $u \in R$ and $s \in G$, and the primitive idempotents of $K^R$ resp. $K^G$ are denoted by $e_u$ resp. $d_s$.

**5.2** In this paragraph, we consider the Radford biproduct $B := A \otimes H$. We introduce the basis
$$b_{uv}(s) := e_u \otimes x_s \otimes c_v$$
of the Radford biproduct, where $u, v \in R$ are elements of the finite ring and $s \in G$ is an element of the finite group under consideration. With respect to this basis, the structure elements of $B$ take the following form:

**Proposition**
1. Multiplication: $b_{uv}(s) b_{u'v'}(s') =$
$$\delta_{u\nu(s), u'} \eta(uq(s, s')) \chi(2vu'\alpha(s') + u^2 \nu(s) \beta(s) \alpha(s')) b_{u, v+v'}(ss')$$
2. Unit: $1_B = \sum_{u \in R} b_{u0}(1)$
3. Comultiplication: $\Delta_B(b_{uv}(s)) = \sum_{w \in R} b_{u-w, w\beta(s)+v}(s) \otimes b_{wv}(s)$
4. Counit: $\epsilon_B(b_{uv}(s)) = \delta_{u0}$
5. Antipode: $S_B(b_{uv}(s)) =$
$$\eta(uq(s, s^{-1})) \chi(-u^2 \beta(s) \alpha(s) - 2uv\beta(s)\alpha(s)) b_{-u\nu(s), -u\beta(s)-v}(s^{-1})$$

**Proof.** This follows by direct computation:

(1) Multiplication:
$$b_{uv}(s) b_{u'v'}(s') = (e_u \otimes x_s)(c_v \to (e_{u'} \otimes x_{s'})) \otimes c_v c_{v'}$$
$$= \chi(vu'\alpha(s'))^2 \delta_{u\nu(s), u'} \eta(uq(s, s')) \chi(u^2 \nu(s) \beta(s) \alpha(s')) e_u \otimes x_{ss'} \otimes c_{v+v'}$$
$$= \delta_{u\nu(s), u'} \eta(uq(s, s')) \chi(2vu'\alpha(s') + u^2 \nu(s) \beta(s) \alpha(s')) b_{u, v+v'}(ss')$$

(2) Comultiplication:
$$\Delta_B(e_u \otimes x_s \otimes c_v) = \sum_{w \in R} (e_{u-w} \otimes x_s \otimes (e_w \otimes x_s)^{(1)} c_v) \otimes ((e_w \otimes x_s)^{(2)} \otimes c_v)$$
$$= \sum_{w \in R} (e_{u-w} \otimes x_s \otimes c_{w\beta(s)+v}) \otimes (e_w \otimes x_s \otimes c_v)$$

(3) Antipode:
$$S_B(e_u \otimes x_s \otimes c_v) = (1_A \otimes S_H((e_u \otimes x_s)^{(1)} c_v))(S_A((e_u \otimes x_s)^{(2)}) \otimes 1_H)$$
$$= \eta(uq(s, s^{-1})) \chi(u^2 \beta(s) \alpha(s))(1_A \otimes c_{-u\beta(s)-v})(e_{-u\nu(s)} \otimes x_{s^{-1}} \otimes 1_H)$$
$$= \eta(uq(s, s^{-1})) \chi(u^2 \beta(s) \alpha(s)) \chi(u\nu(s)(u\beta(s) + v)\alpha(s^{-1}))^2$$
$$(e_{-u\nu(s)} \otimes x_{s^{-1}} \otimes c_{-u\beta(s)-v})$$
$$= \eta(uq(s, s^{-1})) \chi(u^2 \beta(s) \alpha(s) - 2u(u\beta(s) + v)\alpha(s))$$
$$(e_{-u\nu(s)} \otimes x_{s^{-1}} \otimes c_{-u\beta(s)-v})$$
$$= \eta(uq(s, s^{-1})) \chi(-u^2 \beta(s) \alpha(s) - 2uv\beta(s)\alpha(s))(e_{-u\nu(s)} \otimes x_{s^{-1}} \otimes c_{-u\beta(s)-v})$$

where we have used the equation $\alpha(s) = -\nu(s)\alpha(s^{-1})$ from Lemma 1.13.

We leave the verification of the formulas for unit and counit to the reader. $\square$

The reader should compare the above description with the description in Paragraph 8.5, where a slightly different basis for the Radford biproduct is used.

**5.3** It is known that the crossed product of a semisimple algebra and a semisimple group ring is semisimple (cf. [64], Thm. 4.4, p. 31, [57], Thm. 7.4.2. p. 116). Therefore, the algebra $A$ constructed in Paragraph 3.3 is semisimple if the characteristic of $K$ does not divide the cardinality of $G$. Here we approach this issue in a different way that provides slightly more information:

**Proposition 1**
1. $\Lambda_A := \sum_{s \in G} e_0 \otimes x_s$ is a two-sided integral of $A$. It is invariant and coinvariant.

2. $A$ is semisimple if and only if the characteristic of $K$ does not divide the cardinality of $G$.

3. $B$ is semisimple if and only if the characteristic of $K$ neither divides the cardinality of $G$ nor the cardinality of $R$.

**Proof.** $\Lambda_A$ is a left integral since we have:

$$(e_u \otimes x_s)\Lambda_A = \sum_{t \in G}(e_u \otimes x_s)(e_0 \otimes x_t) = \sum_{t \in G}\delta_{u0}e_0 \otimes x_{st} = \epsilon_A(e_u \otimes x_s)\Lambda_A$$

The fact that $\Lambda_A$ is also a right integral follows from a similar calculation. It is obvious that $\Lambda_A$ is invariant and coinvariant. Since $\epsilon_A(\Lambda_A) = \text{card}(G)$, the second assertion follows from Maschke's theorem for Yetter-Drinfel'd Hopf algebras (cf. [20], Cor. 5.8, p. 4885, [81], Prop. 2.14, p. 22). The third assertion on the Radford biproduct $B$ follows from [65], Prop. 3, p. 333. $\square$

Since the coalgebra structure of $A$ is the ordinary tensor product coalgebra structure, it is easy to describe when $A$ is cosemisimple:

**Proposition 2**
1. The linear form $\lambda_A : A \to K$ determined by

$$\lambda_A(e_u \otimes x_s) = \delta_{s1}$$

is a two-sided integral of $A^*$. It is $H$-linear and colinear.

2. $A$ is cosemisimple if and only if the characteristic of $K$ does not divide the cardinality of $R$.

3. $B$ is cosemisimple if and only if the characteristic of $K$ does not divide the cardinality of $R$.

**Proof.** The linear form $\lambda_A$ is the tensor product of the integral $\lambda_R$ on $K^R$ and the integral $\lambda_G$ on $K[G]$ that are determined by the conditions:

$$\lambda_R(e_u) = 1 \qquad \lambda_G(x_s) = \delta_{s1}$$

Therefore, it is a two-sided integral itself; the fact that it is linear and colinear follows from the fact that $\alpha(1) = \beta(1) = 0$, which we have established in Lemma 1.13. Since $\lambda_A(1_A) = \text{card}(R)$, the second assertion follows from the dual of Maschke's theorem for Yetter-Drinfel'd Hopf algebras (cf. [20], Cor. 5.8, p. 4885, [81], Cor. 2.14, p. 23). The third assertion on the Radford biproduct $B$ follows from [65], Prop. 4, p. 335. $\square$

**5.4** The vector space underlying the Radford biproduct is $K^R \otimes K[G] \otimes K[R]$, where the first two tensor factors constitute the Yetter-Drinfel'd Hopf algebra and the last tensor factor represents its base Hopf algebra. A slight shift of the viewpoint yields a rather different picture: Dividing the triple tensor product into two parts consisting of the first tensor factor on the one hand and the last two tensor factors on the other hand, we get a Hopf algebra extension: The first tensor factor $K^R$ is a Hopf subalgebra of the Radford biproduct, whereas the last two tensor factors appear as a Hopf algebra quotient. We make these assertions precise in the following:

**Proposition** Define the linear mappings

$$\iota : K^R \to B, e_u \mapsto b_{u0}(1) = e_u \otimes x_1 \otimes c_0$$

and

$$\pi : B \to K[G] \otimes K[R], b_{uv}(s) = e_u \otimes x_s \otimes c_v \mapsto \delta_{u0} x_s \otimes c_v$$

Then $\iota$ and $\pi$ are Hopf algebra homomorphisms and

$$K^R \xrightarrow{\iota} B \xrightarrow{\pi} K[G] \otimes K[R]$$

is a short exact sequence of Hopf algebras.

**Proof.** Using the fact that $\alpha(1) = \beta(1) = 0$ from Lemma 1.13, it is easy to see that $\iota$ is an injective Hopf algebra homomorphism. We have:

$$\pi(b_{uv}(s)b_{u'v'}(s'))$$
$$= \delta_{u\nu(s),u'} \eta(uq(s,s'))\chi(2vu'\alpha(s') + u^2\nu(s)\beta(s)\alpha(s'))\pi(b_{u,v+v'}(ss'))$$
$$= \delta_{u0}\delta_{u'0}(x_s \otimes c_v)(x_{s'} \otimes c_{v'}) = \pi(b_{uv}(s))\pi(b_{u'v'}(s'))$$

and therefore $\pi$ is an algebra homomorphism, since it also preserves the unit. The fact that $\pi$ is a coalgebra homomorphism follows from a similar direct computation.

To prove the exactness of the above sequence, it suffices to show that the space of coinvariants
$$\{b \in B \mid (\mathrm{id}_B \otimes \pi)\Delta_B(b) = b \otimes 1\}$$
coincides with $\iota(K^R)$; in fact, this is the definition of a short exact sequence of finite-dimensional Hopf algebras (cf. [46], Def. 1.3, p. 821, [55], Def. 5.6, p. 129, see also [70], p. 3338). Now, if $b = \sum_{u,v \in R, s \in G} \mu_{uvs} b_{uv}(s)$, we have:
$$(\mathrm{id}_B \otimes \pi)\Delta_B(b) = \sum_{u,v \in R, s \in G} \mu_{uvs} b_{uv}(s) \otimes (x_s \otimes c_v)$$

Therefore, $b$ is coinvariant if and only if the coefficients $\mu_{uvs}$ vanish for all indices with the property $(v, s) \neq (0, 1)$, which means precisely that $b \in \iota(K^R)$. □

In the general theory of Hopf algebra extensions, one variant of the normal basis theorem asserts in our situation that $B$ is a crossed product of $K[G] \otimes K[R]$ and $K^R$, or of $G \times R$ and $K^R$ in the more traditional terminology (cf. [72], Thm. 2.2, p. 299, [57], Thm. 8.4.6, p. 141). It is not difficult to determine the cocycle and the corresponding action for $B$ explicitly.

**Corollary** Consider $K^R$ as a $G \times R$–module via
$$(s, u).e_v := e_{v\nu(s^{-1})}$$

For $s, t \in G$ and $u, v, w \in R$, define:
$$\tau_u(s, v; t, w) := \eta(uq(s, t))\chi(2uv\nu(s)\alpha(t) + u^2\nu(s)\beta(s)\alpha(t))$$

and set:
$$\tau(s, v; t, w) := \sum_{u \in R} \tau_u(s, v; t, w) e_u$$

Then $B$ is isomorphic to the crossed product of $G \times R$ and $K^R$ with respect to the 2-cocycle $\tau$ and the specified action.

**Proof.** In this crossed product, two basis elements are multiplied as follows:
$$(e_u \otimes x_s \otimes c_v)(e_{u'} \otimes x_{s'} \otimes c_{v'}) = e_u(x_s \otimes c_v).e_{u'}\tau(s, v; s', v') \otimes x_s x_{s'} \otimes c_v c_{v'}$$
$$= \eta(uq(s, s'))\chi(2uv\nu(s)\alpha(s') + u^2\nu(s)\beta(s)\alpha(s'))\delta_{u,u'\nu(s^{-1})} e_u \otimes x_{ss'} \otimes c_{vv'}$$

But this is precisely the formula for the product of two of the basis elements $b_{uv}(s) = e_u \otimes x_s \otimes c_v$ in $B$. This implies the assertion, since the cocycle condition is equivalent to the associativity of the crossed product, which follows from the associativity of $B$. □

We note that the above module structure can also be written in the form $(s, u).e_v := e_{s.v}$, using the module structure $R_G$ defined in Paragraph 3.3.

**5.5** We now turn to another construction, namely the second construction from [78], Sec. 3, resp. [79], Sec. 4. This construction, which generalizes the construction of deformed enveloping algebras, yields a Hopf algebra structure on the vector space $B := A \otimes H \otimes A^*$. This construction does not apply to arbitrary Yetter-Drinfel'd Hopf algebras, but only to those that satisfy the following main assumption on $A$ (cf. [79], Subsec. 4.1), namely the requirement that we have

$$(a^{(1)} \to a') \otimes a^{(2)} = a'^{(2)} \otimes (a'^{(1)} \to a)$$

for all $a, a' \in A$. To satisfy this main assumption, we impose from now on the condition that we have

$$\chi(u\alpha(s)\beta(t)) = \chi(u\beta(s)\alpha(t))$$

for all $s, t \in G$ and all $u \in R$, as already noted in Paragraph 5.1. This condition then assures that the main assumption is satisfied, because we have for $a = e_u \otimes x_s$ and $a' = e_v \otimes x_t$ that:

$$\begin{aligned}(a^{(1)} \to a') \otimes a^{(2)} &= (c_{u\beta(s)} \to (e_v \otimes x_t)) \otimes (e_u \otimes x_s) \\ &= \chi(u\beta(s)v\alpha(t))^2 (e_v \otimes x_t) \otimes (e_u \otimes x_s)\end{aligned}$$

whereas:

$$\begin{aligned}a'^{(2)} \otimes (a'^{(1)} \to a) &= (e_v \otimes x_t) \otimes (c_{v\beta(t)} \to (e_u \otimes x_s)) \\ &= \chi(v\beta(t)u\alpha(s))^2 (e_v \otimes x_t) \otimes (e_u \otimes x_s)\end{aligned}$$

Since we have allowed arbitrary permutations among the arguments of $\chi$ (cf. Paragraph 3.3), both expressions are equal.

**5.6** We have explained in Paragraph 1.2 that the dual of a finite-dimensional left Yetter-Drinfel'd Hopf algebra is a right Yetter-Drinfel'd Hopf algebra. Here we use a modified dual instead. First, we observe that the vector space $A^* := K[R] \otimes K^G$ may be considered as the dual of $A$ with respect to the nondegenerate bilinear form:

$$\langle e_u \otimes x_s, c_v \otimes d_t \rangle_A = \delta_{uv} \delta_{st}$$

Therefore, $A^*$ is a right Yetter-Drinfel'd Hopf algebra over $H$ by Paragraph 1.2. However, we want to consider a different right Yetter-Drinfel'd Hopf algebra structure on $A^*$, namely the coopposite structure in the categorical sense: Instead of $\Delta_{A^*}$, we use $\sigma^{-1}_{A^*,A^*} \circ \Delta_{A^*}$ as a comultiplication; instead of $\delta_{A^*}$, we use $(\operatorname{id}_{A^*} \otimes S_H) \circ \delta_{A^*}$ as a coaction; and instead of $S_{A^*}$, we use $S^{-1}_{A^*}$ as an antipode. With these structures, $A^*$ becomes again a right Yetter-Drinfel'd Hopf algebra over $H$ (cf. [79], Subsec. 4.4). Since we will not need the earlier structures anymore, we will denote the new structures again by $\Delta_{A^*}$, $\delta_{A^*}$, and $S_{A^*}$.

With respect to the basis $c_u \otimes d_s$ of $A^*$, these structure elements take the following form:

**Proposition**

1. Multiplication: $(c_u \otimes d_s)(c_v \otimes d_t) = \delta_{st} c_{u+v} \otimes d_t$

2. Unit: $1_{A^*} = \sum_{s \in G} c_0 \otimes d_s$

3. Comultiplication: $\Delta_{A^*}(c_u \otimes d_s) =$
$$\sum_{t \in G} \eta(uq(t, t^{-1}s))\chi(-u^2\nu(t)\beta(t)\alpha(t^{-1}s))(c_{u\nu(t)} \otimes d_{t^{-1}s}) \otimes (c_u \otimes d_t)$$

4. Counit: $\epsilon_{A^*}(c_u \otimes d_s) = \delta_{s1}$

5. Antipode: $S_{A^*}(c_u \otimes d_s) = \eta(-uq(s, s^{-1}))\chi(-u^2\beta(s)\alpha(s))c_{-u\nu(s)} \otimes d_{s^{-1}}$

6. Action: $(c_u \otimes d_s) \leftarrow c_v = \chi(vu\alpha(s))^2 c_u \otimes d_s$

7. Coaction: $\delta_{A^*}(c_u \otimes d_s) = (c_u \otimes d_s) \otimes c_{-u\beta(s)}$

**Proof.** This rests on straightforward verification. We therefore only prove the more complicated parts concerning the comultiplication and the antipode. From the bilinear form $\langle \cdot, \cdot \rangle_A$, we get as in [79], Subsec. 2.5 a nondegenerate bilinear pairing between $A \otimes A$ and $A^* \otimes A^*$, which we also denote by $\langle \cdot, \cdot \rangle_A$. If we denote the multiplication mapping of $A$ by $\mu_A$, we have by [79], Subsec. 4.4:

$$\langle (e_v \otimes x_p) \otimes (e_w \otimes x_r), \Delta_{A^*}(c_u \otimes d_s) \rangle_A$$
$$= \langle \mu_A \circ \sigma_{A,A}^{-1}((e_v \otimes x_p) \otimes (e_w \otimes x_r)), c_u \otimes d_s \rangle_A$$
$$= \chi(-wv\beta(r)\alpha(p))^2 \langle (e_w \otimes x_r)(e_v \otimes x_p), c_u \otimes d_s \rangle_A$$
$$= \chi(-wv\beta(r)\alpha(p))^2 \delta_{w\nu(r),v} \eta(wq(r,p))\chi(w^2\nu(r)\beta(r)\alpha(p))\langle e_w \otimes x_{rp}, c_u \otimes d_s \rangle_A$$
$$= \sum_{t \in G} \chi(-u^2\nu(t)\beta(t)\alpha(t^{-1}s))^2 \eta(uq(t, t^{-1}s))\chi(u^2\nu(t)\beta(t)\alpha(t^{-1}s))$$
$$\langle (e_v \otimes x_p) \otimes (e_w \otimes x_r), (c_{u\nu(t)} \otimes d_{t^{-1}s}) \otimes (c_u \otimes d_t) \rangle_A$$
$$= \langle (e_v \otimes x_p) \otimes (e_w \otimes x_r),$$
$$\sum_{t \in G} \eta(uq(t, t^{-1}s))\chi(-u^2\nu(t)\beta(t)\alpha(t^{-1}s))(c_{u\nu(t)} \otimes d_{t^{-1}s}) \otimes (c_u \otimes d_t) \rangle_A$$

Since the pairing between $A \otimes A$ and $A^* \otimes A^*$ is nondegenerate, this establishes the form of the comultiplication.

To establish the form of the antipode, we take the formula given above as the definition of a linear endomorphism $S_{A^*}$ of $A^*$. It is then easy to verify that this endomorphism satisfies:

$$\langle S_A(e_u \otimes x_s), S_{A^*}(c_v \otimes d_t) \rangle_A = \langle e_u \otimes x_s, c_v \otimes d_t \rangle_A$$

This implies that we have $\langle S_A^{-1}(a), b \rangle_A = \langle a, S_{A^*}(b) \rangle_A$ for all $a \in A$ and $b \in A^*$. Therefore, it follows from [79], Subsec. 4.4 that $S_{A^*}$ really is the antipode of $A^*$. $\square$

**5.7**  The adjoint action of a Hopf algebra on itself has various generalizations for Yetter-Drinfel'd Hopf algebras, depending on how the necessary interchanging of the tensor factors is accomplished and whether the antipode or its inverse is used in the definition. As in [79], Subsec. 4.6, we here consider the left adjoint action constructed with the inverse quasisymmetry and the inverse antipode, which we denote by $\rightharpoonup$:

$$a \rightharpoonup a' := \mu_A(\mu_A \otimes S_A^{-1})\sigma_{A \otimes A, A}^{-1}(\Delta_A \otimes \mathrm{id}_A)(a \otimes a')$$

Via the left adjoint action, $A$ becomes a left $A$-module. Similarly, $A^*$ becomes a right $A^*$-module via the right adjoint action:

$$b' \leftharpoonup b := \mu_{A^*}(S_{A^*}^{-1} \otimes \mu_{A^*})\sigma_{A^*, A^* \otimes A^*}^{-1}(\mathrm{id}_{A^*} \otimes \Delta_{A^*})(b' \otimes b)$$

where we use the modified structures from the previous paragraph. With respect to the bases $e_u \otimes x_s$ resp. $c_u \otimes d_s$ of $A$ resp. $A^*$, the adjoint actions take the following form:

**Proposition**  For $u, v \in R$ and $s, t \in G$, we have:

1. $(e_u \otimes x_s) \rightharpoonup (e_v \otimes x_t) =$
   $\delta_{u\nu(s), v-v\nu(t)} \eta(v\nu(s^{-1})(q(s,t) + q(st, s^{-1}) - \nu(sts^{-1})q(s, s^{-1})))$
   $\chi(v^2 \nu(s^{-2})(\nu(s)\beta(t) + \nu(t)\beta(st) - \nu(t^2)\beta(s))\alpha(s)) e_{v\nu(s^{-1})} \otimes x_{sts^{-1}}$

2. $(c_v \otimes d_t) \leftharpoonup (c_u \otimes d_s) = \delta_{s1}\chi(-vu\nu(t)\beta(t)\alpha(t))^2 c_v \otimes d_t$

**Proof.**  Using Heyneman-Sweedler sigma notation, the left adjoint action can also be written in the form:

$$a \rightharpoonup a' = a_{(2)}{}^{(2)} a'^{(2)} S_A^{-1}(S_H(a_{(2)}{}^{(1)} a'^{(1)}) \rightharpoonup a_{(1)})$$

With Lemma 1.13, we see that the inverse of the antipode is given by the equation:

$$S_A^{-1}(e_u \otimes x_s) = \eta(uq(s, s^{-1}))\chi(-u^2 \beta(s)\alpha(s)) e_{-u\nu(s)} \otimes x_{s^{-1}}$$

Using this lemma again, together with our assumption that $\chi(u\alpha(s)\beta(t)) = \chi(u\beta(s)\alpha(t))$, we get:

$(e_u \otimes x_s) \rightharpoonup (e_v \otimes x_t)$
$= \sum_{w \in R} (e_w \otimes x_s)^{(2)} (e_v \otimes x_t)^{(2)}$
$\qquad\qquad S_A^{-1}(S_H((e_w \otimes x_s)^{(1)} (e_v \otimes x_t)^{(1)}) \rightharpoonup (e_{u-w} \otimes x_s))$
$= \sum_{w \in R} (e_w \otimes x_s)(e_v \otimes x_t) S_A^{-1}(c_{-w\beta(s)-v\beta(t)} \rightharpoonup (e_{u-w} \otimes x_s))$

$$\begin{aligned}
&= \sum_{w \in R} \delta_{w\nu(s),v} \eta(wq(s,t))\chi(w^2\nu(s)\beta(s)\alpha(t))\\
&\qquad \chi(-(w\beta(s)+v\beta(t))(u-w)\alpha(s))^2(e_w \otimes x_{st})S_A^{-1}(e_{u-w} \otimes x_s)\\
&= \eta(v\nu(s^{-1})q(s,t))\chi(v^2\nu(s^{-1})\beta(s)\alpha(t))\\
&\qquad \chi(v(\nu(s^{-1})\beta(s)+\beta(t))(v\nu(s^{-1})-u)\alpha(s))^2\\
&\qquad (e_{v\nu(s^{-1})} \otimes x_{st})S_A^{-1}(e_{u-v\nu(s^{-1})} \otimes x_s)\\
&= \eta(v\nu(s^{-1})q(s,t))\chi(v^2\nu(s^{-1})\beta(s)\alpha(t))\\
&\qquad \chi(v\nu(s^{-1})(\beta(s)+\nu(s)\beta(t))(v\nu(s^{-1})-u)\alpha(s))^2\\
&\qquad \eta((u-v\nu(s^{-1}))q(s,s^{-1}))\chi(-(u-v\nu(s^{-1}))^2\beta(s)\alpha(s))\\
&\qquad (e_{v\nu(s^{-1})} \otimes x_{st})(e_{v-u\nu(s)} \otimes x_{s^{-1}})\\
&= \eta(v\nu(s^{-1})q(s,t))\chi(v^2\nu(s^{-1})\beta(s)\alpha(t))\\
&\qquad \chi(v\nu(s^{-1})\beta(st)(v\nu(s^{-1})-u)\alpha(s))^2 \eta((u-v\nu(s^{-1}))q(s,s^{-1}))\\
&\qquad \chi(-(u-v\nu(s^{-1}))^2\beta(s)\alpha(s))\delta_{v\nu(t),v-u\nu(s)}\eta(v\nu(s^{-1})q(st,s^{-1}))\\
&\qquad \chi(v^2\nu(s^{-2})\nu(st)\beta(st)\alpha(s^{-1}))e_{v\nu(s^{-1})} \otimes x_{sts^{-1}}\\
&= \delta_{u\nu(s),v-v\nu(t)}\eta(v\nu(s^{-1})q(s,t))\eta(-v\nu(ts^{-1})q(s,s^{-1}))\eta(v\nu(s^{-1})q(st,s^{-1}))\\
&\qquad \chi(v^2\nu(s^{-1})\beta(t)\alpha(s))\chi(v^2\nu(s^{-1})\beta(st)\nu(ts^{-1})\alpha(s))^2\\
&\qquad \chi(-v^2\nu(ts^{-1})^2\beta(s)\alpha(s))\chi(-v^2\nu(s^{-2}t)\beta(st)\alpha(s))e_{v\nu(s^{-1})} \otimes x_{sts^{-1}}\\
&= \delta_{u\nu(s),v-v\nu(t)}\eta(v\nu(s^{-1})(q(s,t)+q(st,s^{-1})-\nu(sts^{-1})q(s,s^{-1})))\\
&\qquad \chi(v^2\nu(s^{-2})(\nu(s)\beta(t)+\nu(t)\beta(st)-\nu(t^2)\beta(s))\alpha(s))e_{v\nu(s^{-1})} \otimes x_{sts^{-1}}
\end{aligned}$$

This proves the first statement; we now turn to the second. Using Heyneman-Sweedler sigma notation, the right adjoint action can also be written in the form:
$$b' \leftharpoonup b = S_{A^*}^{-1}(b_{(2)} \leftharpoonup S_H(b'^{(2)}b_{(1)}{}^{(2)}))b'^{(1)}b_{(1)}{}^{(1)}$$

Again with Lemma 1.13, we see that the inverse of the antipode of $A^*$ is given by the equation:
$$S_{A^*}^{-1}(c_u \otimes d_s) = \eta(-uq(s,s^{-1}))\chi(u^2\beta(s)\alpha(s))c_{-u\nu(s)} \otimes d_{s^{-1}}$$

We therefore have:

$$(c_v \otimes d_t) \leftharpoonup (c_u \otimes d_s)$$
$$= \sum_{r \in G} \eta(uq(r,r^{-1}s))\chi(-u^2\nu(r)\beta(r)\alpha(r^{-1}s))$$
$$\qquad S_{A^*}^{-1}((c_u \otimes d_r) \leftharpoonup S_H((c_v \otimes d_t)^{(2)}(c_{u\nu(r)} \otimes d_{r^{-1}s})^{(2)}))$$
$$\qquad\qquad\qquad (c_v \otimes d_t)^{(1)}(c_{u\nu(r)} \otimes d_{r^{-1}s})^{(1)}$$
$$= \sum_{r \in G} \eta(uq(r,r^{-1}s))\chi(-u^2\nu(r)\beta(r)\alpha(r^{-1}s))$$

$$S_{A^*}^{-1}((c_u \otimes d_r) \leftharpoonup c_{v\beta(t)+u\nu(r)\beta(r^{-1}s)})(c_v \otimes d_t)(c_{u\nu(r)} \otimes d_{r^{-1}s})$$

$$= \sum_{r \in G} \eta(uq(r,r^{-1}s))\chi(-u^2\nu(r)\beta(r)\alpha(r^{-1}s))$$

$$\chi((v\beta(t)+u\nu(r)\beta(r^{-1}s))u\alpha(r))^2 \delta_{t,r^{-1}s} S_{A^*}^{-1}(c_u \otimes d_r)(c_{v+u\nu(r)} \otimes d_t)$$

$$= \eta(uq(st^{-1},t))\chi(-u^2\nu(st^{-1})\beta(st^{-1})\alpha(t))$$

$$\chi((v\beta(t)+u\nu(st^{-1})\beta(t))u\alpha(st^{-1}))^2 \eta(-uq(st^{-1},ts^{-1}))$$

$$\chi(u^2\beta(st^{-1})\alpha(st^{-1}))(c_{-u\nu(st^{-1})} \otimes d_{ts^{-1}})(c_{v+u\nu(st^{-1})} \otimes d_t)$$

$$= \delta_{s1} \eta(uq(t^{-1},t))\chi(-u^2\nu(t^{-1})\beta(t^{-1})\alpha(t))$$

$$\chi((v\beta(t)+u\nu(t^{-1})\beta(t))u\alpha(t^{-1}))^2 \eta(-uq(t^{-1},t))\chi(u^2\beta(t^{-1})\alpha(t^{-1}))c_v \otimes d_t$$

$$= \delta_{s1}\chi(u^2\beta(t^{-1})\alpha(t^{-1}))^2 \chi((v\beta(t)-u\beta(t^{-1}))u\alpha(t^{-1}))^2 c_v \otimes d_t$$

$$= \delta_{s1}\chi(vu\beta(t)\alpha(t^{-1}))^2 c_v \otimes d_t = \delta_{s1}\chi(-vu\nu(t)\beta(t)\alpha(t))^2 c_v \otimes d_t \qquad \square$$

**5.8** We now dualize the adjoint actions to get the coadjoint actions (cf. [79], Subsec. 4.7). Of course, this dualization depends on the bilinear form that is used. To define the right coadjoint action $\leftharpoonup$ of $A$ on $A^*$, we use the bilinear form $\langle \cdot, \cdot \rangle_A$ from Paragraph 5.6; i. e., we define the right coadjoint action by the condition:

$$\langle a', b \leftharpoonup a \rangle_A = \langle a \rightharpoonup a', b \rangle_A$$

For the dualization of the right adjoint action of $A^*$ on itself, we use a different bilinear form $\langle \cdot, \cdot \rangle_{A^*}$ that is defined as:

$$\langle \cdot, \cdot \rangle_{A^*} : A \otimes A^* \to K, a \otimes b \mapsto \langle S_A^{-1}(a), b \rangle_A$$

This bilinear form is the convolution inverse of the bilinear form $\langle \cdot, \cdot \rangle_A$; i. e., we have:

$$\langle a_{(1)}, b_{(1)} \rangle_A \langle a_{(2)}, b_{(2)} \rangle_{A^*} = \epsilon_A(a)\epsilon_{A^*}(b) = \langle a_{(1)}, b_{(1)} \rangle_{A^*} \langle a_{(2)}, b_{(2)} \rangle_A$$

for all $a \in A$ and $b \in A^*$.

We then define the left coadjoint action $\rightharpoonup$ of $A^*$ on $A$ by the condition:

$$\langle b \rightharpoonup a, b' \rangle_{A^*} = \langle a, b' \leftharpoonup b \rangle_{A^*}$$

Finally, we introduce the mapping

$$\sharp : A^* \otimes A \to H, b \otimes a \mapsto b \sharp a := \langle a_{(1)}, b_{(1)}{}^{(1)} \rangle_{A^*} b_{(1)}{}^{(2)} a_{(2)}{}^{(1)} \langle a_{(2)}{}^{(2)}, b_{(2)} \rangle_A$$

With respect to the bases $e_u \otimes x_s$ of $A$, resp. $c_u \otimes d_s$ of $A^*$, these structure elements take the following form:

**Proposition** For $u, v \in R$ and $s, t \in G$, we have:

1. $(c_u \otimes d_s) \rightharpoonup (e_v \otimes x_t) = \delta_{s1} \chi(vu\nu(t^{-2})\beta(t)\alpha(t))^2 e_v \otimes x_t$

2. $(c_v \otimes d_t) \leftharpoonup (e_u \otimes x_s) = \delta_{u,v-v\nu(t)} \eta(v(q(s, s^{-1}ts) + q(ts, s^{-1}) - \nu(t)q(s, s^{-1})))$
   $\chi(v^2(\nu(s)\beta(s^{-1}ts) + \nu(t)\beta(ts) - \nu(t^2)\beta(s))\alpha(s))^2 c_{v\nu(s)} \otimes d_{s^{-1}ts}$

3. $(c_u \otimes d_s) \sharp (e_v \otimes x_t) = \delta_{v0} \delta_{s1} c_{2u\beta(t)}$

**Proof.** Using Lemma 1.13, we have for the left coadjoint action:

$$\langle (c_u \otimes d_s) \rightharpoonup (e_v \otimes x_t), c_w \otimes d_r \rangle_{A^*} = \langle e_v \otimes x_t, (c_w \otimes d_r) \leftharpoonup (c_u \otimes d_s) \rangle_{A^*}$$
$$= \delta_{s1} \chi(-wu\nu(r)\beta(r)\alpha(r))^2 \langle e_v \otimes x_t, c_w \otimes d_r \rangle_{A^*}$$
$$= \delta_{s1} \chi(-wu\nu(r)\beta(r)\alpha(r))^2 \eta(vq(t, t^{-1})) \chi(-v^2 \beta(t)\alpha(t))$$
$$\langle e_{-v\nu(t)} \otimes x_{t^{-1}}, c_w \otimes d_r \rangle_A$$
$$= \delta_{s1} \chi(v\nu(t)u\nu(t^{-1})\beta(t^{-1})\alpha(t^{-1}))^2 \eta(vq(t, t^{-1})) \chi(-v^2 \beta(t)\alpha(t))$$
$$\langle e_{-v\nu(t)} \otimes x_{t^{-1}}, c_w \otimes d_r \rangle_A$$
$$= \delta_{s1} \chi(vu\beta(t^{-1})\alpha(t^{-1}))^2 \langle e_v \otimes x_t, c_w \otimes d_r \rangle_{A^*}$$
$$= \delta_{s1} \chi(vu\nu(t^{-2})\beta(t)\alpha(t))^2 \langle e_v \otimes x_t, c_w \otimes d_r \rangle_{A^*}$$

We leave the verification of the formula for the right coadjoint action to the reader. The third formula follows, using again Lemma 1.13, from the following calculation:

$(c_u \otimes d_s) \sharp (e_v \otimes x_t)$
$$= \sum_{w \in R, r \in G} \eta(uq(r, r^{-1}s)) \chi(-u^2 \nu(r)\beta(r)\alpha(r^{-1}s))$$
$$\langle e_{v-w} \otimes x_t, (c_{u\nu(r)} \otimes d_{r^{-1}s})^{(1)} \rangle_{A^*} (c_{u\nu(r)} \otimes d_{r^{-1}s})^{(2)} (e_w \otimes x_t)^{(1)}$$
$$\langle (e_w \otimes x_t)^{(2)}, c_u \otimes d_r \rangle_A$$
$$= \eta(uq(t, t^{-1}s)) \chi(-u^2 \nu(t)\beta(t)\alpha(t^{-1}s))$$
$$\langle e_{v-u} \otimes x_t, c_{u\nu(t)} \otimes d_{t^{-1}s} \rangle_{A^*} c_{-u\nu(t)\beta(t^{-1}s)} c_{u\beta(t)}$$
$$= \eta(uq(t, t^{-1}s)) \chi(-u^2 \nu(t)\beta(t)\alpha(t^{-1}s)) \eta((v-u)q(t, t^{-1})) \chi(-(v-u)^2 \beta(t)\alpha(t))$$
$$\langle e_{(u-v)\nu(t)} \otimes x_{t^{-1}}, c_{u\nu(t)} \otimes d_{t^{-1}s} \rangle_A c_{-u\nu(t)\beta(t^{-1}s)} c_{u\beta(t)}$$
$$= \delta_{v0} \delta_{s1} c_{2u\beta(t)} \qquad \square$$

**5.9** The second construction described in [78], Sec. 3, resp. [79], Sec. 4, now enables us to build a second ordinary Hopf algebra from the Yetter-Drinfel'd Hopf algebra $A$ considered in Paragraph 3.3. The underlying vector space of this Hopf algebra is $A \otimes H \otimes A^*$; it contains the Radford biproduct $A \otimes H$ as a kind of Borel subalgebra.

We describe the Hopf algebra arising from the second construction with respect to the basis
$$z_{uvw}(s,t) := e_u \otimes x_s \otimes c_v \otimes c_w \otimes d_t$$
of $A \otimes H \otimes A^*$, where $u, v, w \in R$ and $s, t \in G$. With respect to this basis, the structure elements of $A \otimes H \otimes A^*$ take the following form:

**Proposition**
1. Multiplication:
$$z_{uvw}(s,t) z_{u'v'w'}(s',t') =$$
$$\delta_{u'-u\nu(s), w-w\nu(t)} \delta_{ts', s't'}$$
$$\eta(uq(s,s') + w(q(s',t') + q(s't', s'^{-1}) - \nu(t)q(s', s'^{-1})))$$
$$\chi(2uw\nu(sts'^{-2})\beta(s')\alpha(s') + 2vu\nu(s)\alpha(s') + 2v'w\nu(s')\alpha(t')$$
$$+ 2w^2(\nu(s')\beta(t') + \nu(t)\beta(t))\alpha(s') + u^2\nu(s)\beta(s)\alpha(s'))$$
$$z_{u, v+v'+w\beta(s')+w\nu(t)\beta(s'), w\nu(s')+w'}(ss', t')$$

2. Unit: $1 = \sum_{u \in R, s \in G} z_{u00}(1, s)$

3. Comultiplication:
$$\Delta(z_{uvw}(s,t)) =$$
$$\sum_{r \in G, k \in R} \eta(wq(r, r^{-1}t))\chi(-w^2\nu(r)\beta(r)\alpha(r^{-1}t))$$
$$z_{u-k, v+k\beta(s), w\nu(r)}(s, r^{-1}t) \otimes z_{k, v-w\nu(r)\beta(r^{-1}t), w}(s, r)$$

4. Counit: $\epsilon(z_{uvw}(s,t)) = \delta_{u0}\delta_{t1}$

5. Antipode:
$$S(z_{uvw}(s,t)) =$$
$$\eta((u - w\nu(ts^{-1}))q(s, s^{-1}) - wq(t, t^{-1})$$
$$+ w\nu(ts^{-1})(q(st^{-1}, s^{-1}) + q(s, t^{-1})))$$
$$\chi((-u^2 + 2uw\nu(s) - 2w^2\nu(t) + 2w^2 + 2uw\nu(s^{-1}t) - 2uw\nu(s^{-1})$$
$$+ 2w^2\nu(s^{-2}t) - 2w^2\nu(s^{-2}t^2))\beta(s)\alpha(s) + 4w^2\nu(s^{-1})\beta(t)\alpha(s)$$
$$+ w^2\beta(t)\alpha(t) - 2(uv - vw\nu(s^{-1}t) + vw\nu(s^{-1}))\alpha(s) - 2vw\alpha(t))$$
$$z_{-u\nu(s)+w\nu(t)-w, w\beta(t)-v-u\beta(s)-w\nu(t)\beta(s^{-1})-w\beta(s^{-1}), -w\nu(ts^{-1})}(s^{-1}, st^{-1}s^{-1})$$

**Proof.** We first establish the form of the multiplication. We have:
$$(\Delta_A \otimes \mathrm{id}_A) \circ \Delta_A(e_u \otimes x_s) = \sum_{k,l \in R} (e_k \otimes x_s) \otimes (e_{l-k} \otimes x_s) \otimes (e_{u-l} \otimes x_s)$$

Similarly, we have:

$$(\mathrm{id}_{A^*} \otimes \Delta_{A^*}) \circ \Delta_{A^*}(c_w \otimes d_t) =$$
$$\sum_{p,r \in G} \eta(w(q(r, r^{-1}t) + q(p, p^{-1}r)))$$
$$\chi(-w^2(\nu(r)\beta(r)\alpha(r^{-1}t) + \nu(p)\beta(p)\alpha(p^{-1}r)))$$
$$(c_{w\nu(r)} \otimes d_{r^{-1}t}) \otimes (c_{w\nu(p)} \otimes d_{p^{-1}r}) \otimes (c_w \otimes d_p)$$

This implies:

$$(e_u \otimes x_s \otimes c_v \otimes c_w \otimes d_t)(e_{u'} \otimes x_{s'} \otimes c_{v'} \otimes c_{w'} \otimes d_{t'})$$
$$= \sum_{k,l \in R, p,r \in G} \eta(w(q(r, r^{-1}t) + q(p, p^{-1}r)))$$
$$\chi(-w^2(\nu(r)\beta(r)\alpha(r^{-1}t) + \nu(p)\beta(p)\alpha(p^{-1}r)))$$
$$(e_u \otimes x_s)(c_v \to [(c_{w\nu(r)} \otimes d_{r^{-1}t})^{(1)} \rightharpoonup (e_k \otimes x_{s'})])$$
$$\otimes c_v(c_{w\nu(r)} \otimes d_{r^{-1}t})^{(2)}[c_{w\nu(p)} \otimes d_{p^{-1}r} \sharp e_{l-k} \otimes x_{s'}](e_{u'-l} \otimes x_{s'})^{(1)} c_{v'}$$
$$\otimes ([(c_w \otimes d_p) \leftharpoonup (e_{u'-l} \otimes x_{s'})^{(2)}] \leftharpoonup c_{v'})(c_{w'} \otimes d_{t'})$$
$$= \sum_{k,l \in R, p,r \in G} \eta(w(q(r, r^{-1}t) + q(p, p^{-1}r)))$$
$$\chi(-w^2(\nu(r)\beta(r)\alpha(r^{-1}t) + \nu(p)\beta(p)\alpha(p^{-1}r)))$$
$$\delta_{rt}\delta_{u'-l,w-w\nu(p)}\delta_{lk}\delta_{rp}\chi(kw\nu(r)\nu(s'^{-2})\beta(s')\alpha(s'))^2$$
$$\eta(w(q(s', s'^{-1}ps') + q(ps', s'^{-1}) - \nu(p)q(s', s'^{-1})))$$
$$\chi(w^2(\nu(s')\beta(s'^{-1}ps') + \nu(p)\beta(ps') - \nu(p^2)\beta(s'))\alpha(s'))^2$$
$$(e_u \otimes x_s)(c_v \to (e_k \otimes x_{s'})) \otimes c_v c_{-w\nu(r)\beta(r^{-1}t)} c_{2w\nu(p)\beta(s')} c_{(u'-l)\beta(s')} c_{v'}$$
$$\otimes ((c_{w\nu(s')} \otimes d_{s'^{-1}ps'}) \leftharpoonup c_{v'})(c_{w'} \otimes d_{t'})$$
$$= \chi((u' - w + w\nu(t))w\nu(t)\nu(s'^{-2})\beta(s')\alpha(s'))^2$$
$$\eta(w(q(s', s'^{-1}ts') + q(ts', s'^{-1}) - \nu(t)q(s', s'^{-1})))$$
$$\chi(w^2(\nu(s')\beta(s'^{-1}ts') + \nu(t)\beta(ts') - \nu(t^2)\beta(s'))\alpha(s'))^2$$
$$\chi((u' - w + w\nu(t))v\alpha(s'))^2 \chi(v'w\nu(s')\alpha(s'^{-1}ts'))^2$$
$$(e_u \otimes x_s)(e_{u'-w+w\nu(t)} \otimes x_{s'}) \otimes c_v c_{2w\nu(t)\beta(s')} c_{(w-w\nu(t))\beta(s')} c_{v'}$$
$$\otimes (c_{w\nu(s')} \otimes d_{s'^{-1}ts'})(c_{w'} \otimes d_{t'})$$
$$= \delta_{u\nu(s),u'-w+w\nu(t)} \delta_{s'^{-1}ts',t'} \chi(u\nu(s)w\nu(t)\nu(s'^{-2})\beta(s')\alpha(s'))^2$$
$$\eta(w(q(s', t') + q(s't', s'^{-1}) - \nu(t)q(s', s'^{-1})))$$
$$\chi(w^2(\nu(s')\beta(t') + \nu(t)\beta(s't') - \nu(t^2)\beta(s'))\alpha(s'))^2$$
$$\chi(vu\nu(s)\alpha(s'))^2 \chi(v'w\nu(s')\alpha(t'))^2 \eta(uq(s, s')) \chi(u^2\nu(s)\beta(s)\alpha(s'))$$
$$e_u \otimes x_{ss'} \otimes c_{v+v'+w\beta(s')+w\nu(t)\beta(s')} \otimes c_{w\nu(s')+w'} \otimes d_{t'}$$

$$= \delta_{u'-u\nu(s),w-w\nu(t)}\delta_{ts',s't'}$$
$$\eta(uq(s,s') + w(q(s',t') + q(s't',s'^{-1}) - \nu(t)q(s',s'^{-1})))$$
$$\chi(uw\nu(sts'^{-2})\beta(s')\alpha(s') + vu\nu(s)\alpha(s') + v'w\nu(s')\alpha(t'))^2$$
$$\chi(w^2(\nu(s')\beta(t') + \nu(t)\beta(s't') - \nu(t^2)\beta(s'))\alpha(s'))^2 \chi(u^2\nu(s)\beta(s)\alpha(s'))$$
$$e_u \otimes x_{ss'} \otimes c_{v+v'+w\beta(s')+w\nu(t)\beta(s')} \otimes c_{w\nu(s')+w'} \otimes d_{t'}$$

If $s't' = ts'$, we get from the cocycle identity that:
$$\nu(t)\beta(s't') - \nu(t^2)\beta(s') = \nu(t)\beta(ts') - \nu(t^2)\beta(s') = \nu(t)\beta(t)$$

Therefore, the second argument of $\chi$ in the above expression for the product is equal to $w^2(\nu(s')\beta(t') + \nu(t)\beta(t))\alpha(s')$. This implies the form of the multiplication stated above.

The form of the comultiplication follows easily from the definition (cf. [79], Subsec. 3.2, p. 39):

$$\Delta(e_u \otimes x_s \otimes c_v \otimes c_w \otimes d_t) =$$
$$(e_u \otimes x_s)_{(1)} \otimes (e_u \otimes x_s)_{(2)}{}^{(1)}c_v \otimes (c_w \otimes d_t)_{(1)}{}^{(1)}$$
$$\otimes (e_u \otimes x_s)_{(2)}{}^{(2)} \otimes c_v(c_w \otimes d_t)_{(1)}{}^{(2)} \otimes (c_w \otimes d_t)_{(2)}$$
$$\sum_{r \in G, k \in R} \eta(wq(r,r^{-1}t))\chi(-w^2\nu(r)\beta(r)\alpha(r^{-1}t))$$
$$(e_{u-k} \otimes x_s) \otimes (e_k \otimes x_s)^{(1)}c_v \otimes (c_{w\nu(r)} \otimes d_{r^{-1}t})^{(1)}$$
$$\otimes (e_k \otimes x_s)^{(2)} \otimes c_v(c_{w\nu(r)} \otimes d_{r^{-1}t})^{(2)} \otimes (c_w \otimes d_r)$$
$$\sum_{r \in G, k \in R} \eta(wq(r,r^{-1}t))\chi(-w^2\nu(r)\beta(r)\alpha(r^{-1}t))$$
$$(e_{u-k} \otimes x_s \otimes c_{v+k\beta(s)} \otimes c_{w\nu(r)} \otimes d_{r^{-1}t})$$
$$\otimes (e_k \otimes x_s \otimes c_{v-w\nu(r)\beta(r^{-1}t)} \otimes c_w \otimes d_r)$$

For the antipode, we have:

$$S(e_u \otimes x_s \otimes c_v \otimes c_w \otimes d_t)$$
$$= (1_A \otimes 1_H \otimes S_{A^*}((c_w \otimes d_t)^{(1)}))(1_A \otimes S_H((e_u \otimes x_s)^{(1)}c_v(c_w \otimes d_t)^{(2)}) \otimes 1_{A^*})$$
$$(S_A((e_u \otimes x_s)^{(2)}) \otimes 1_H \otimes 1_{A^*})$$
$$= \eta(-wq(t,t^{-1}))\chi(-w^2\beta(t)\alpha(t))\eta(uq(s,s^{-1}))\chi(u^2\beta(s)\alpha(s))$$
$$(1_A \otimes 1_H \otimes c_{-w\nu(t)} \otimes d_{t^{-1}})(1_A \otimes c_{w\beta(t)-v-u\beta(s)} \otimes 1_{A^*})$$
$$(e_{-u\nu(s)} \otimes x_{s^{-1}} \otimes 1_H \otimes 1_{A^*})$$
$$= \eta(uq(s,s^{-1}) - wq(t,t^{-1}))\chi(u^2\beta(s)\alpha(s) - w^2\beta(t)\alpha(t))$$
$$\chi(-u\nu(s)(w\beta(t) - v - u\beta(s))\alpha(s^{-1}))^2$$
$$(1_A \otimes 1_H \otimes c_{-w\nu(t)} \otimes d_{t^{-1}})(e_{-u\nu(s)} \otimes x_{s^{-1}} \otimes c_{w\beta(t)-v-u\beta(s)} \otimes 1_{A^*})$$

$$\begin{aligned}
&= \sum_{k \in R, r \in G} \eta(uq(s,s^{-1}) - wq(t,t^{-1}))\chi(u^2\beta(s)\alpha(s) - w^2\beta(t)\alpha(t)) \\
&\quad \chi(-u\nu(s)(w\beta(t) - v - u\beta(s))\alpha(s^{-1}))^2 \\
&\quad (e_k \otimes x_1 \otimes c_0 \otimes c_{-w\nu(t)} \otimes d_{t^{-1}})(e_{-u\nu(s)} \otimes x_{s^{-1}} \otimes c_{w\beta(t)-v-u\beta(s)} \otimes c_0 \otimes d_r) \\
&= \eta(uq(s,s^{-1}) - wq(t,t^{-1}))\chi(u^2\beta(s)\alpha(s) - w^2\beta(t)\alpha(t)) \\
&\quad \chi(-u\nu(s)(w\beta(t) - v - u\beta(s))\alpha(s^{-1}))^2 \\
&\quad \eta(-w\nu(t)(q(s^{-1}, st^{-1}s^{-1}) + q(t^{-1}s^{-1}, s) - \nu(t^{-1})q(s^{-1}, s))) \\
&\quad \chi(w\nu(t)(u\nu(s) - w\nu(t) + w)\nu(t^{-1}s^2)\beta(s^{-1})\alpha(s^{-1}) \\
&\qquad\qquad - (w\beta(t) - v - u\beta(s))w\nu(t)\nu(s^{-1})\alpha(st^{-1}s^{-1}))^2 \\
&\quad \chi(w^2\nu(t)^2(\nu(s^{-1})\beta(st^{-1}s^{-1}) + \nu(t^{-1})\beta(t^{-1}s^{-1}) - \nu(t^{-2})\beta(s^{-1}))\alpha(s^{-1}))^2 \\
&\quad e_{-u\nu(s)+w\nu(t)-w} \otimes x_{s^{-1}} \otimes c_{w\beta(t)-v-u\beta(s)-w\nu(t)\beta(s^{-1})-w\beta(s^{-1})} \otimes \\
&\qquad\qquad\qquad\qquad\qquad\qquad\qquad\qquad c_{-w\nu(t)\nu(s^{-1})} \otimes d_{st^{-1}s^{-1}}
\end{aligned}$$

To simplify this expression, we treat the arguments of $\chi$ and $\eta$ separately. As already pointed out in Paragraph 3.3, $\chi$ vanishes by assumption on the two-sided ideal $I$ generated by the additive commutators $uv - vu$ and the elements of the form $\beta(s)\alpha(t) - \beta(t)\alpha(s)$, for $u, v \in R$ and $s, t \in G$. Therefore, denoting congruence modulo $I$ by $\equiv$, we can rewrite the argument of $\chi$ as follows:

$$\begin{aligned}
&u^2\beta(s)\alpha(s) - w^2\beta(t)\alpha(t) - 2u\nu(s)(w\beta(t) - v - u\beta(s))\alpha(s^{-1}) \\
&\quad + 2w\nu(t)(u\nu(s) - w\nu(t) + w)\nu(t^{-1}s^2)\beta(s^{-1})\alpha(s^{-1}) \\
&\quad - 2(w\beta(t) - v - u\beta(s))w\nu(t)\nu(s^{-1})\alpha(st^{-1}s^{-1}) \\
&\quad + 2w^2\nu(t)^2(\nu(s^{-1})\beta(st^{-1}s^{-1}) + \nu(t^{-1})\beta(t^{-1}s^{-1}) - \nu(t^{-2})\beta(s^{-1}))\alpha(s^{-1}) \\
&\equiv u^2\beta(s)\alpha(s) - w^2\beta(t)\alpha(t) + 2u(w\beta(t) - v - u\beta(s))\alpha(s) \\
&\quad + 2w(u\nu(s) - w\nu(t) + w)\beta(s)\alpha(s) - 2(w\beta(t) - v - u\beta(s))w\nu(s^{-1}t)\alpha(s) \\
&\quad - 2(w\beta(t) - v - u\beta(s))w\nu(t)\alpha(t^{-1}) + 2(w\beta(t) - v - u\beta(s))w\nu(s^{-1})\alpha(s) \\
&\quad + 2(w^2\nu(s^{-2}t) - w^2\nu(s^{-2}t^2))\beta(s)\alpha(s) + 2(w^2\nu(s^{-1}) + w^2\nu(s^{-1}t))\beta(t)\alpha(s) \\
&\equiv (-u^2 + 2uw\nu(s) - 2w^2\nu(t) + 2w^2 + 2uw\nu(s^{-1}t) - 2uw\nu(s^{-1}) \\
&\qquad\qquad\qquad\qquad\qquad\qquad + 2w^2\nu(s^{-2}t) - 2w^2\nu(s^{-2}t^2))\beta(s)\alpha(s) \\
&\quad + 4w^2\nu(s^{-1})\beta(t)\alpha(s) + w^2\beta(t)\alpha(t) - 2(uv - vw\nu(s^{-1}t) + vw\nu(s^{-1}))\alpha(s) \\
&\quad - 2vw\alpha(t)
\end{aligned}$$

Here we have used the equality
$$\alpha(sts^{-1}) \equiv \alpha(s) + \nu(s)\alpha(t) - \nu(t)\alpha(s)$$

which follows easily from the definition of a 1-cocycle and Lemma 1.13, and,

based on this equality, the expansion

$$2w^2\nu(t)^2(\nu(s^{-1})\beta(st^{-1}s^{-1}) + \nu(t^{-1})\beta(t^{-1}s^{-1}) - \nu(t^{-2})\beta(s^{-1}))\alpha(s^{-1})$$
$$\equiv 2(w^2\nu(s^{-2}t) - w^2\nu(s^{-2}t^2))\beta(s)\alpha(s) + 2(w^2\nu(s^{-1}) + w^2\nu(s^{-1}t))\beta(t)\alpha(s)$$

for the last summand in the first expression above.

To simplify the argument of $\eta$, we first note that we have from the definition of a normalized 2-cocycle that:

$$q(s^{-1}, st^{-1}s^{-1}) + q(t^{-1}s^{-1}, s) = \nu(s^{-1})q(st^{-1}s^{-1}, s) + q(s^{-1}, st^{-1})$$
$$= \nu(s^{-1})q(st^{-1}s^{-1}, s) - \nu(s^{-1})q(st^{-1}, s^{-1}s) + q(s^{-1}, st^{-1}) - q(s^{-1}s, t^{-1})$$
$$= \nu(s^{-1})\nu(st^{-1})q(s^{-1}, s) - \nu(s^{-1})q(st^{-1}, s^{-1}) + q(s^{-1}, s) - \nu(s^{-1})q(s, t^{-1})$$
$$= \nu(t^{-1})q(s^{-1}, s) + q(s^{-1}, s) - \nu(s^{-1})q(st^{-1}, s^{-1}) - \nu(s^{-1})q(s, t^{-1})$$

Therefore, the argument of $\eta$ in the above expression is given by the formula:

$$uq(s, s^{-1}) - wq(t, t^{-1})$$
$$\qquad - w\nu(t)(q(s^{-1}, st^{-1}s^{-1}) + q(t^{-1}s^{-1}, s) - \nu(t^{-1})q(s^{-1}, s))$$
$$= uq(s, s^{-1}) - wq(t, t^{-1})$$
$$\qquad - w\nu(t)(q(s^{-1}, s) - \nu(s^{-1})q(st^{-1}, s^{-1}) - \nu(s^{-1})q(s, t^{-1}))$$
$$= (u - w\nu(ts^{-1}))q(s, s^{-1}) - wq(t, t^{-1}) + w\nu(ts^{-1})(q(st^{-1}, s^{-1}) + q(s, t^{-1}))$$

This implies the asserted form of the antipode.

The formulas for the unit and the counit follow immediately from their respective definition. □

The reader is invited to check directly that the above structures make $A \otimes H \otimes A^*$ into a Hopf algebra, which is not entirely obvious.

**5.10** We now want to determine when the above Hopf algebra is semisimple. As in Paragraph 5.3, we approach this problem by explicitly exhibiting an integral.

**Proposition 1**
1. $\Lambda := \sum_{s \in G, u, v \in R} z_{0uv}(s, 1)$ is a two-sided integral in $A \otimes H \otimes A^*$.

2. $A \otimes H \otimes A^*$ is semisimple if and only if the characteristic of $K$ neither divides the cardinality of $R$ nor the cardinality of $G$.

**Proof.** We have seen in Proposition 5.3.1 that $\Lambda_A := \sum_{s \in G} e_0 \otimes x_s$ is a two-sided integral of $A$ that is invariant and coinvariant. The integral character therefore coincides with the counit, and the integral group element coincides with the unit (cf. [81], Prop. 2.10, p. 15). Since, as an algebra, $A^*$ is the ordinary tensor product of a group ring and a dual group ring, it is easy to see that $\Lambda_{A^*} := \sum_{v \in R} c_v \otimes d_1$ is a two-sided integral of $A^*$, which is, as a consequence of Lemma 1.13, also invariant and coinvariant. If $\Lambda_H := \sum_{u \in R} c_u$ denotes the integral of the group ring $H$, we have by [81], Thm. 5.4, p. 62 that $\Lambda = \Lambda_A \otimes \Lambda_H \otimes \Lambda_{A^*}$ is a two-sided integral of $A \otimes H \otimes A^*$. Since $\epsilon(\Lambda) = \text{card}(G) \text{card}(R)^2$, the second assertion follows from Maschke's theorem for ordinary Hopf algebras (cf. [57], Thm. 2.2.1, p. 20). □

We now turn to the question of cosemisimplicity. In general, we obtain a left integral $\lambda$ on the Hopf algebra $A \otimes H \otimes A^*$ emerging from the second construction by the formula:
$$\lambda(a \otimes h \otimes b) = \lambda_A(a)\lambda_H(g_A h)\lambda_B(b)$$
where $\lambda_A \in A^*$, $\lambda_H \in H^*$, and $\lambda_B \in B^*$ are left integrals and $g_A \in H$ is the integral group element of $A$. A right integral $\rho$ on $A \otimes H \otimes A^*$ is given by the formula:
$$\rho(a \otimes h \otimes b) = \rho_A(a)\rho_H(hg_A^{-1})\rho_B(b)$$
where $\rho_A \in A^*$, $\rho_H \in H^*$, and $\rho_B \in B^*$ are right integrals. In our situation, this yields the following:

**Proposition 2**
1. The linear form $\lambda : A \otimes H \otimes A^* \to K$ determined by
$$\lambda(z_{uvw}(s,t)) = \delta_{s1}\delta_{v0}\delta_{w0}$$
is a two-sided integral on $A \otimes H \otimes A^*$.

2. $A \otimes H \otimes A^*$ is cosemisimple if and only if the characteristic of $K$ neither divides the cardinality of $R$ nor the cardinality of $G$.

**Proof.** We have seen in Proposition 5.3.2 that the linear form $\lambda_A : A \to K$ determined by
$$\lambda_A(e_u \otimes x_s) = \delta_{s1}$$
is a two-sided integral of $A^*$. If $\Lambda_A := \sum_{s \in G} e_0 \otimes x_s$ denotes the two-sided integral of $A$ obtained in Proposition 5.3.1, we know from [81], Lem. 5.4, p. 61 that the linear form
$$\lambda_{A^*} : A^* \to K, b \mapsto \langle \Lambda_A, b \rangle_B$$
is a two-sided integral on $A^*$. Since the antipode fixes $\Lambda_A$, it is easy to see that we have:
$$\lambda_{A^*}(c_u \otimes d_s) = \delta_{u0}$$
As already noted in the proof of Proposition 1, the integral group element $g_A$ is equal to one. Using the two-sided integral
$$\lambda_H : H \to K, c_u \mapsto \delta_{u0}$$

on the group ring $H$, we therefore get from the formulas above that the linear form $\lambda := \lambda_A \otimes \lambda_H \otimes \lambda_{A^*}$ is a two-sided integral on $A \otimes H \otimes A^*$. It is obvious that this integral satisfies:
$$\lambda(z_{uvw}(s,t)) = \delta_{s1}\delta_{v0}\delta_{w0}$$
Since $\lambda(1) = \operatorname{card}(G)\operatorname{card}(R)$, the second assertion follows again from Maschke's theorem. □

**5.11** As we have seen in Paragraph 5.4, the Radford biproduct can also be understood from the point of view of Hopf algebra extensions. In this paragraph, we shall obtain a similar result for the Hopf algebra arising from the second construction. For this, we need some preparation.

Consider the right action of $G$ on the additive group $R \times R$ defined by:
$$(v,w)s := (v + 2w\beta(s), w\nu(s))$$
The fact that this really is a right action follows from the cocycle condition for $\beta$:
$$(v,w)st = (v + 2w\beta(s) + 2w\nu(s)\beta(t), w\nu(s)\nu(t)) = ((v,w)s)t$$
It is also easy to see that $G$ acts on $R \times R$ via group automorphisms. We can therefore form the corresponding semidirect product $T := G \ltimes (R \times R)$. In this semidirect product, the multiplication is determined by the formula:
$$(s,v,w)(s',v',w') = (ss', v + v' + 2w\beta(s'), w\nu(s') + w')$$
We denote the canonical basis elements of the group ring $K[T]$ by $y_{vw}(s)$, where $v, w \in R$ and $s \in G$.

Besides this semidirect product, we will use the dual group ring $K^{G_{\mathrm{op}}}$ of the opposite group $G_{\mathrm{op}}$, in which the elements are multiplied in the reverse order. The coproduct in $K^{G_{\mathrm{op}}}$ is determined by the formula:
$$\Delta_{K^{G_{\mathrm{op}}}}(d_s) = \sum_{t \in G} d_{t^{-1}s} \otimes d_t$$
Of course, $K^{G_{\mathrm{op}}}$ is just the coopposite Hopf algebra $(K^G)^{\mathrm{cop}}$ of $K^G$.

**Proposition** Define the linear mappings
$$\iota : K^R \otimes K^{G_{\mathrm{op}}} \to B, e_u \otimes d_t \mapsto z_{u00}(1,t)$$
and
$$\pi : B \to K[T], z_{uvw}(s,t) \mapsto \delta_{u0}\delta_{t1}y_{vw}(s)$$
Then $\iota$ and $\pi$ are Hopf algebra homomorphisms and
$$K^R \otimes K^{G_{\mathrm{op}}} \overset{\iota}{\rightarrowtail} B \overset{\pi}{\twoheadrightarrow} K[T]$$
is a short exact sequence of Hopf algebras.

**Proof.** $\iota$ is an algebra homomorphism since we have:

$$z_{u00}(1,t)z_{u'00}(1,t') = \delta_{uu'}\delta_{tt'}z_{u00}(1,t')$$

and since it preserves the unit. $\iota$ is a coalgebra homomorphism since we have:

$$\Delta(z_{u00}(1,t)) = \sum_{s\in G, v\in R} z_{u-v,0,0}(1,s^{-1}t) \otimes z_{v00}(1,s)$$

and since it preserves the counit. The assertion that $\pi$ is an algebra homomorphism follows from:

$$\pi(z_{uvw}(s,t)z_{u'v'w'}(s',t')) = \delta_{u0}\delta_{u'0}\delta_{t1}\delta_{t'1}y_{v+v'+2w\beta(s'),w\nu(s')+w'}(ss')$$
$$= \pi(z_{uvw}(s,t))\pi(z_{u'v'w'}(s',t'))$$

To show that $\pi$ is a coalgebra homomorphism, we first observe that:

$$(\mathrm{id}\otimes\pi)\Delta(z_{uvw}(s,t)) = z_{uvw}(s,t) \otimes y_{v-w\beta(t),w}(s)$$

Since $\beta(1) = 0$, this implies immediately that:

$$(\pi\otimes\pi)\Delta(z_{uvw}(s,t)) = \delta_{u0}\delta_{t1}y_{vw}(s) \otimes y_{vw}(s)$$

and therefore $\pi$ is a coalgebra homomorphism. But the above equation also shows that an element

$$b = \sum_{\substack{u,v,w\in R \\ s,t\in G}} \mu_{uvwst}z_{uvw}(s,t)$$

is coinvariant if and only if the coefficients $\mu_{uvwst}$ vanish for all indices with the property $(s, v - w\beta(t), w) \neq (1,0,0)$. Since this property is equivalent to $(s,v,w) \neq (1,0,0)$, this means precisely that $b \in \iota(K^R \otimes K^{G_{\mathrm{op}}})$. As explained in Paragraph 5.4, this means that the above sequence is exact. $\square$

**5.12** As we have also pointed out in Paragraph 5.4, the middle term of an exact sequence of finite-dimensional Hopf algebras is a crossed product of the outer terms. However, in the present case the determination of the action and the cocycle is more difficult than in the case of the Radford biproduct in Paragraph 5.4. The reason for this is that the basis $z_{uvw}(s,t)$ is not a normal basis of $B$ in the sense that the map

$$K^R \otimes K^{G_{\mathrm{op}}} \otimes K[T] \to B, e_u \otimes d_t \otimes y_{vw}(s) \mapsto z_{uvw}(s,t)$$

is in general neither $K^R \otimes K^{G_{\mathrm{op}}}$-linear nor $K[T]$-colinear, as we have, at least in part, already seen in the proof of Proposition 5.11. Therefore, we introduce a new basis of $B$, which is defined as:

$$z'_{uvw}(s,t) := z_{u,v+w\beta(s^{-1}ts),w}(s,s^{-1}ts)$$

The product of two of these basis elements is given by the formula:

$$z'_{uvw}(s,t)z'_{u'v'w'}(s',t') =$$
$$\delta_{u'-u\nu(s),w-w\nu(s^{-1}ts)}\delta_{s^{-1}ts,t'}$$
$$\eta(uq(s,s') + w(q(s',s'^{-1}t's') + q(t's',s'^{-1}) - \nu(s^{-1}ts)q(s',s'^{-1})))$$
$$\chi(2uw\nu(tss'^{-2})\beta(s')\alpha(s') + 2vu\nu(s)\alpha(s') + 2w\beta(s^{-1}ts)u\nu(s)\alpha(s')$$
$$+ 2v'w\nu(s')\alpha(s'^{-1}t's') + 2ww'\nu(s')\beta(s'^{-1}t's')\alpha(s'^{-1}t's')$$
$$+ 2w^2(\nu(s')\beta(s'^{-1}t's') + \nu(s^{-1}ts)\beta(s^{-1}ts))\alpha(s') + u^2\nu(s)\beta(s)\alpha(s'))$$
$$z_{u,v+v'+w\beta(s^{-1}ts)+w'\beta(s'^{-1}t's')+w\beta(s')+w\nu(s^{-1}ts)\beta(s'),w\nu(s')+w'}(ss', s'^{-1}t's')$$

We want to rewrite this formula slightly. First, since

$$\beta(t) + \beta(s) + \nu(t)\beta(s) - \nu(s)\beta(s^{-1}ts)$$
$$= \beta(s) + \beta(ts) - \nu(s)\beta(s^{-1}) - \beta(ts) = 2\beta(s)$$

by Lemma 1.13, we have for $t' = s^{-1}ts$:

$$z_{u,v+v'+w\beta(s^{-1}ts)+w'\beta(s'^{-1}t's')+w\beta(s')+w\nu(s^{-1}ts)\beta(s'),w\nu(s')+w'}(ss', s'^{-1}t's')$$
$$= z_{u,v+v'+w\beta(t')+w'\beta(s'^{-1}t's')+w\beta(s')+w\nu(t')\beta(s'),w\nu(s')+w'}(ss', s'^{-1}t's')$$
$$= z'_{u,v+v'+w\beta(t')+w\beta(s')+w\nu(t')\beta(s')-w\nu(s')\beta(s'^{-1}t's'),w\nu(s')+w'}(ss', st's^{-1})$$
$$= z'_{u,v+v'+2w\beta(s'),w\nu(s')+w'}(ss', t)$$

We also have:

$$q(s,s^{-1}ts) + q(ts,s^{-1}) - \nu(t)q(s,s^{-1})$$
$$= q(s,s^{-1}ts) - q(ss^{-1},ts) + q(ts,s^{-1}) - \nu(t)q(s,s^{-1})$$
$$= q(s,s^{-1}) - \nu(s)q(s^{-1},ts) + q(t,ss^{-1}) - q(t,s)$$

and therefore, if $t' = s^{-1}ts$, the argument of $\eta$ in the above equation can be written in the form:

$$uq(s,s') + w(q(s',s'^{-1}t's') + q(t's',s'^{-1}) - \nu(s^{-1}ts)q(s',s'^{-1}))$$
$$= uq(s,s') + w(q(s',s'^{-1}) - \nu(s')q(s'^{-1},t's') - q(t',s'))$$

We now turn to the argument of $\chi$. As pointed out in Paragraph 5.9, $\chi$ vanishes on the two-sided ideal $I$ generated by the additive commutators $uv - vu$ and the elements of the form $\beta(s)\alpha(t) - \beta(t)\alpha(s)$, for $u, v \in R$ and $s, t \in G$. Since we have, using Lemma 1.13, that:

$$\beta(s^{-1}ts) = \nu(s^{-1})(-\beta(s) + \beta(t) + \nu(t)\beta(s)) = \nu(s^{-1})(\nu(t)-1)\beta(s) + \nu(s^{-1})\beta(t)$$

we get, denoting congruence modulo $I$ by $\equiv$, that:

$$\beta(s^{-1}ts)\alpha(s^{-1}ts) \equiv \nu(s^{-2})(1 - 2\nu(t) + \nu(t^2))\beta(s)\alpha(s)$$
$$+ 2\nu(s^{-2})(\nu(t) - 1)\beta(s)\alpha(t) + \nu(s^{-2})\beta(t)\alpha(t)$$

Applying this to the primed arguments instead and collecting terms, the above formula can be rewritten in the form:

$$z'_{uvw}(s,t)z'_{u'v'w'}(s',t') =$$
$$\delta_{u'-u\nu(s),w-w\nu(t')}\delta_{s^{-1}ts,t'}$$
$$\eta(uq(s,s') + w(q(s',s'^{-1}) - \nu(s')q(s'^{-1},t's') - q(t',s')))$$
$$\chi([2u\nu(tss'^{-2}) + 2w'\nu(s'^{-1}) - 4w'\nu(s'^{-1}t') + 2w'\nu(s'^{-1}t'^2)$$
$$+ 2w\nu(t') - 2w]w\beta(s')\alpha(s') + 2ww'\nu(s'^{-1})\beta(t')\alpha(t')$$
$$+ [2u\nu(s) + 4w'\nu(s'^{-1}t') - 4w'\nu(s'^{-1}) + 2w + 2w\nu(t')]w\beta(s')\alpha(t')$$
$$+ 2vu\nu(s)\alpha(s') + 2v'w\nu(s')\alpha(s'^{-1}t's') + u^2\nu(s)\beta(s)\alpha(s'))$$
$$z'_{u,v+v'+2w\beta(s'),w\nu(s')+w'}(ss',t)$$

We now can determine the cocycle and the corresponding action for $B$ explicitly. Consider $K^R \otimes K^{G_{\mathrm{op}}}$ as a $T$–module via:

$$(s,v,w).(e_u \otimes d_t) := e_{u\nu(s^{-1})-w\nu(s^{-1})+w\nu(ts^{-1})} \otimes d_{sts^{-1}}$$

We extend this action to a $K[T]$–module structure by linearity. Furthermore, we define for $s, s', t \in G$ and $u, v, w, v', w' \in R$:

$$\rho_{ut}(s,v,w;s',v',w') :=$$
$$\eta(uq(s,s') + w(q(s',s'^{-1}) - \nu(s')q(s'^{-1},s^{-1}tss') - q(s^{-1}ts,s')))$$
$$\chi([2u\nu(tss'^{-2}) + 2w'\nu(s'^{-1}) - 4w'\nu(s'^{-1}s^{-1}ts) + 2w'\nu(s'^{-1}s^{-1}t^2s)$$
$$+ 2w\nu(s^{-1}ts) - 2w]w\beta(s')\alpha(s') + 2ww'\nu(s'^{-1})\beta(s^{-1}ts)\alpha(s^{-1}ts)$$
$$+ [2u\nu(s) + 4w'\nu(s'^{-1}s^{-1}ts) - 4w'\nu(s'^{-1}) + 2w + 2w\nu(s^{-1}ts)]w\beta(s')\alpha(s^{-1}ts)$$
$$+ 2vu\nu(s)\alpha(s') + 2v'w\nu(s')\alpha(s'^{-1}s^{-1}tss') + u^2\nu(s)\beta(s)\alpha(s'))$$

and set:

$$\rho(s,v,w;s',v',w') := \sum_{u \in R, t \in G} \rho_{ut}(s,v,w;s',v',w')e_u \otimes d_t$$

**Proposition** The map

$$K^R \otimes K^{G_{\mathrm{op}}} \otimes K[T] \to B, e_u \otimes d_t \otimes y_{vw}(s) \mapsto z'_{uvw}(s,t)$$

is an isomorphism between $B$ and the crossed product of $T$ and $K^R \otimes K^{G_{\mathrm{op}}}$ with respect to the 2-cocycle $\rho$ and the specified action. $z'_{uvw}(s,t)$ is a normal basis of $B$ in the sense that this mapping is $K^R \otimes K^{G_{\mathrm{op}}}$-linear and $K[T]$-colinear, where $B$ is a $K^R \otimes K^{G_{\mathrm{op}}}$-module via $\iota$ and a $K[T]$-comodule via $\pi$.

**Proof.** In this crossed product, two basis elements are multiplied as follows:

$$(e_u \otimes d_t \otimes y_{vw}(s))(e_{u'} \otimes d_{t'} \otimes y_{v'w'}(s'))$$
$$= (e_u \otimes d_t)(y_{vw}(s).(e_{u'} \otimes d_{t'}))\rho(s,v,w;s',v',w') \otimes y_{vw}(s)y_{v'w'}(s')$$
$$= (e_u \otimes d_t)(e_{u'\nu(s^{-1})-w\nu(s^{-1})+w\nu(t's^{-1})} \otimes d_{st's^{-1}})\rho_{ut}(s,v,w;s',v',w')$$
$$\otimes y_{v+v'+2w\beta(s'),w\nu(s')+w'}(ss')$$
$$= \delta_{u'-u\nu(s),w-w\nu(t')}\delta_{s^{-1}ts,t'}$$
$$\eta(uq(s,s') + w(q(s',s'^{-1}) - \nu(s')q(s'^{-1},t's') - q(t',s')))$$
$$\chi([2u\nu(tss'^{-2}) + 2w'\nu(s'^{-1}) - 4w'\nu(s'^{-1}t') + 2w'\nu(s'^{-1}t'^2)$$
$$+ 2w\nu(t') - 2w]w\beta(s')\alpha(s') + 2ww'\nu(s'^{-1})\beta(t')\alpha(t')$$
$$+ [2u\nu(s) + 4w'\nu(s'^{-1}t') - 4w'\nu(s'^{-1}) + 2w + 2w\nu(t')]w\beta(s')\alpha(t')$$
$$+ 2vu\nu(s)\alpha(s') + 2v'w\nu(s')\alpha(s'^{-1}t's') + u^2\nu(s)\beta(s)\alpha(s'))$$
$$(e_u \otimes d_t) \otimes y_{v+v'+2w\beta(s'),w\nu(s')+w'}(ss')$$

Now the formula for the product of two of the basis elements $z'_{uvw}(s,t)$ implies that the above map is an algebra homomorphism, which is obviously bijective. Since the cocycle condition is equivalent to the associativity of the crossed product, this also implies that $\rho$ is a 2-cocycle. The fact that the above map is $K^R \otimes K^{G_{\text{op}}}$-linear follows from the formula

$$z'_{u00}(1,t)z'_{u'v'w'}(s',t') = \delta_{uu'}\delta_{tt'}z'_{u'v'w'}(s',t')$$

whereas the formula

$$(\text{id} \otimes \pi)\Delta(z'_{uvw}(s,t)) = z'_{uvw}(s,t) \otimes y_{vw}(s)$$

which is a consequence of the corresponding formula for $z_{uvw}(s,t)$ obtained in the proof of Proposition 5.11, implies that this map is $K[T]$-colinear. □

# 6 Commutative Yetter-Drinfel'd Hopf algebras

**6.1** In this section, we assume that $p$ is a prime and that $K$ is an algebraically closed field whose characteristic is different from $p$. We can therefore choose a primitive $p$-th root of unity that we denote by $\zeta$. The group ring of the cyclic group $\mathbb{Z}_p$ of order $p$ will be denoted by $H := K[\mathbb{Z}_p]$; the canonical basis element of the group ring corresponding to $i \in \mathbb{Z}_p$ will be denoted by $c_i$. The group of grouplike elements of $H$, which consists precisely of these canonical basis elements, will be denoted by $C$. $A$ denotes a Yetter-Drinfel'd Hopf algebra over $H$ that is commutative and semisimple. $A$ therefore has a unique basis that consists of primitive idempotents, which is denoted by $E$. If $\gamma : C \to K^\times$ is the group homomorphism that maps $c_1$ to $\zeta$, we introduce as in Paragraph 1.10 the mappings

$$\phi : V \to V, v \mapsto (c_1 \to v) \qquad \psi : V \to V, v \mapsto \gamma(v^{(1)})v^{(2)}$$

where we extend $\gamma$ to $H$ by linearity. (Note that the symbol $\gamma$, in contrast to $\psi$, will also be used for other characters below.) Since $\phi$ and $\psi$ are algebra automorphisms, they induce permutations of $E$. The Radford biproduct $A \otimes H$ will be denoted by $B$.

Throughout the section, we will constantly use the convention that indices take values between 0 and $p-1$ and are reduced modulo $p$ if they do not lie within this range. In notation, we shall not distinguish between an integer $i \in \mathbb{Z}$ and its equivalence class in $\mathbb{Z}_p := \mathbb{Z}/p\mathbb{Z}$.

**6.2** In the commutative case considered here, the results of Section 2 yield the following:

**Proposition** Suppose that $V$ is a simple $B$-module.

1. $V$ is purely unstable if and only if $\dim V = p$. In this case, $\text{Span}(\kappa(V))$ is a $B$-submodule of $A$ that is isomorphic to $V$.

2. $V$ is stable if and only if $\dim V = 1$.

In particular, we have $\dim V = p$ or $\dim V = 1$.

**Proof.** If $V$ is purely unstable, we have $\dim V = p$ by Proposition 2.4. Then $\text{Span}(\kappa(V))$ is an $H$-submodule that is simultaneously and ideal of $A$, and therefore a $B$-submodule of $A$. To see that $V \cong \text{Span}(\kappa(V))$, we compare characters.

Denote the character of $\mathrm{Span}(\kappa(V))$ by $\chi$. If $e \in E$ and $i \in \mathbb{Z}_p$, $i \neq 0$, $e \otimes c_i$ induces a fixed-point free map of $\kappa(V)$ to itself, which implies $\chi(e \otimes c_i) = 0$. If $i = 0$, we have $\chi(e \otimes c_i) = 1$. Therefore, $\chi$ coincides with the character of $V$, which was given in Proposition 2.4.

If $V$ is stable, the restriction of $V$ to $A$ is still simple by Proposition 2.5, and is therefore one-dimensional. $\square$

**6.3** The correspondence between simple modules and orbits afforded by Clifford theory has an additional aspect that we have not considered so far: To every simple $B$-module, we can associate its dual. An interesting question is which orbit corresponds to the dual module. We here consider the more important case where the module is purely unstable.

**Lemma** Suppose that $V$ is a purely unstable simple $B$-module. Then we have: $\mathrm{Span}(\kappa(V^*)) = S_A^{-1}(\mathrm{Span}(\kappa(V)))$.

**Proof.** By Proposition 6.2, we can then assume that $V = \mathrm{Span}(\kappa(V))$. We define $V' := S_A^{-1}(V)$. $V'$ is an $H$-submodule of $A$. For $v \in V$ and $a \in A$, we have:
$$S_A^{-1}(a) S_A^{-1}(v) = S_A^{-1}((a^{(1)} \to v) a^{(2)}) \in S_A^{-1}(V)$$
and therefore $V'$ is even an ideal of $A$. Therefore, $V'$ is a $B$-submodule of the $B$-module $A$.

Now choose a nonzero right integral $\rho_A \in A^*$ and define, as in Paragraph 1.9, the bilinear form:
$$\langle \cdot, \cdot \rangle : A \times A \to K, (a, a') \mapsto \rho_A(S_A(a) a')$$

Since $A$ is a Frobenius algebra with Frobenius homomorphism $\rho_A$, this bilinear form is nondegenerate and restricts to a nondegenerate bilinear form
$$\langle \cdot, \cdot \rangle : V' \times V \to K$$

Now Proposition 1.9 implies that $V' \cong V^*$ as $B$-modules. Therefore, we have $\mathrm{Span}(\kappa(V^*)) = V' = S_A^{-1}(V)$. $\square$

The above lemma has the following consequence that will be important in the sequel:

**Proposition** Suppose that $V$ is a simple $B$-module.

1. If $\kappa(V)$ is $\psi$-invariant, then $\kappa(V^*)$ is $\psi$-invariant.

2. If $\kappa(V)$ consists of fixed points of $\psi$, then $\kappa(V^*)$ consists of fixed points of $\psi$.

**Proof.** We consider the case where $V$ is purely unstable first. Since $S_A$ is colinear, it commutes with $\psi$, and therefore $\text{Span}(\kappa(V)) = S_A(\text{Span}(\kappa(V^*)))$ is $\psi$-invariant if and only if $\text{Span}(\kappa(V^*))$ is $\psi$-invariant, which is obviously equivalent to the $\psi$-invariance of $\kappa(V^*)$. Also, $\text{Span}(\kappa(V)) = S_A(\text{Span}(\kappa(V^*)))$ consists of $\psi$-fixed points if and only if $\text{Span}(\kappa(V^*))$ consists of $\psi$-fixed points.

We now consider the case where $V$ is stable. In this case, the character $\chi_V$ of $V$ is an algebra homomorphism to the base field, and therefore has the form $\chi_V = \eta \otimes \gamma$ for some $H$-linear algebra homomorphism $\eta : A \to K$ and some algebra homomorphism $\gamma : H \to K$. If $e \in E$ is the primitive idempotent that satisfies $\eta(e) = 1$, we have $\kappa(V) = \{e\}$. If $\kappa(V)$ is $\psi$-invariant, $e$ is a fixed point of $\psi$, which means that $\eta$ is colinear. This implies that $\eta \otimes \epsilon_H$ commutes with $\epsilon_A \otimes \gamma$. Then $\eta^{-1} := \eta \circ S_A$ is also an $H$-linear and colinear algebra homomorphism to the base field, and we have

$$\chi_{V^*} = S_{B^*}(\chi_V) = (\epsilon_A \otimes \gamma^{-1})(\eta^{-1} \otimes \epsilon_H) = \eta^{-1} \otimes \gamma^{-1}$$

The primitive idempotent corresponding to $\eta^{-1}$ is $S_A^{-1}(e)$: Since $e$ is invariant and coinvariant, $S_A^{-1}(e)$ is an idempotent, and we can prove as in the preceding lemma that $S_A^{-1}(e)$ generates a one-dimensional ideal of $A$, which means that $S_A^{-1}(e)$ is primitive. Therefore, we have $\kappa(V^*) = \{S_A^{-1}(e)\}$, which is also $\psi$-invariant. □

**6.4** The possibility to consider dual modules leads to a connection between the left and right actions of $\hat{C}$ on the set of irreducible characters of $B$ considered in Section 2. This fact has the surprising consequence that purely unstable orbits are also invariant with respect to $\psi$:

**Proposition** Suppose that $O$ is a purely unstable orbit of $C$ in the set $E$ of primitive idempotents of $A$. Then we have: $\psi(O) \subset O$.

**Proof.** Consider the simple $B$-submodule $V := \text{Span}(O)$ of $A$ and denote its character by $\chi_V$. By Corollary 2.6, the character $\chi_{V^*}$ of the dual module $V^*$ satisfies $\chi_{V^*} = \chi_{V^*}(\epsilon_A \otimes \gamma)$ for all $\gamma \in \hat{C}$. Dualizing this, we get for the character of $V$ that $(\epsilon_A \otimes \gamma^{-1})\chi_V = \chi_V$. Applying $\kappa$, the assertion now follows from Proposition 2.7. □

**6.5** In Paragraph 1.7, we have already investigated modules over the algebra $A \hat{\otimes} A$. Now we want to determine, in our situation, all modules over this algebra. We shall achieve this by describing explicitly the decomposition of $A \hat{\otimes} A$ into simple two-sided ideals.

Let $O_1, \ldots, O_l$ denote the $C$-orbits in the set $E$ of primitive idempotents. Similarly, we denote the orbits of $\psi$ by $O'_1, \ldots, O'_m$. We define $V_i := \text{Span}(O_i)$ and $W_j := \operatorname{Span}(O'_j)$. Then $V_i$ is an $H$-submodule of $A$ and $W_j$ is an $H$-subcomodule of $A$. It is easy to see that

$$A \hat{\otimes} A = \bigoplus_{i,j=1}^{l,m} V_i \hat{\otimes} W_j$$

is a decomposition of $A \hat{\otimes} A$ into two-sided ideals. However, not all of these are simple. There are four possibilities:

**Proposition** Suppose that $O$ is an orbit of $\phi$ and that $O'$ is an orbit of $\psi$. Define $V := \operatorname{Span}(O)$ and $W := \operatorname{Span}(O')$.

1. If $O = \{e\}$ and $O' = \{e'\}$ are both orbits of length one, $V \hat{\otimes} W$ is a one-dimensional two-sided ideal and $e \otimes e'$ is a centrally primitive idempotent of $A \hat{\otimes} A$.

2. If $O = \{e\}$ is an orbit of length one, i. e., a stable orbit, and $O' = \{e'_0, \ldots, e'_{p-1}\}$ is an orbit of length $p$, i. e., a purely unstable orbit, $V \hat{\otimes} W = \bigoplus_{j=0}^{p-1} K e \otimes e'_j$ is a decomposition of $V \hat{\otimes} W$ into simple two-sided ideals.

3. If $O = \{e_0, \ldots, e_{p-1}\}$ is an orbit of length $p$, i. e., a purely unstable orbit, and $O' = \{e'\}$ is an orbit of length one, $V \hat{\otimes} W = \bigoplus_{i=0}^{p-1} K e_i \otimes e'$ is a decomposition of $V \hat{\otimes} W$ into simple two-sided ideals.

4. If $O = \{e_0, \ldots, e_{p-1}\}$ and $O' = \{e'_0, \ldots, e'_{p-1}\}$ are both orbits of length $p$, $V \hat{\otimes} W$ is a simple two-sided ideal.

**Proof.** The first three statements are obvious. We now prove the fourth. $V$ is a module algebra over $H$, whereas $W$ is a comodule algebra over $H$. For each such pair, we can form the left smash product $V \otimes W$ already considered in Paragraph 1.7, where the multiplication is defined as:

$$(v \otimes w)(v' \otimes w') := v(w^{(1)} \to v') \otimes w^{(2)} w'$$

This is precisely the algebra structure of the ideal $V \hat{\otimes} W$ in $A \hat{\otimes} A$. We now perform a discrete Fourier transform to pass to a new basis $c'_0, \ldots, c'_{p-1}$ of $W$ defined as:

$$c'_j := \sum_{i=0}^{p-1} \zeta^{-ij} e'_i$$

It is easy to see that these basis elements satisfy $c'_i c'_j = c'_{i+j}$ and $\psi(c'_j) = \zeta^{-j} c'_j$, i. e., we have:

$$\delta_A(c'_j) = c_j \otimes c'_j$$

This proves that the mapping $H \to W, c_j \mapsto c'_j$ is an isomorphism of comodule algebras, where $H$ is considered as a comodule algebra via the left regular coaction, i. e., the comodule structure whose coaction is equal to the comultiplication. It is obvious that $V$ is isomorphic, as a module algebra, to $K^{\mathbb{Z}_p}$, where $K^{\mathbb{Z}_p}$ is endowed with the module algebra structure considered in the first step of the proof of Proposition 2.3.1. Therefore, $V \hat{\otimes} W$ is isomorphic to the smash product $K^{\mathbb{Z}_p} \otimes H$, and thus is simple as seen in the proof of Proposition 2.3.1. □

This proposition has the following corollary that will be useful later on:

**Corollary** Suppose that $V$ is a simple $B$-module of dimension $p$ and that $e, e' \in E$ are primitive idempotents. If $e' = \psi(e)$, the $A$-modules $V \hat{\otimes} Ae$ and $V \hat{\otimes} Ae'$ are isomorphic.

**Proof.** According to Paragraph 1.7, the modules $V \hat{\otimes} Ae$ and $V \hat{\otimes} Ae'$ are even $A \hat{\otimes} A$-modules, and we shall prove the stronger statement that they are isomorphic as such. By Proposition 6.2, we can assume that $V = \mathrm{Span}(O)$ for some purely unstable $\phi$-orbit $O$. Denote the $\psi$-orbit of $e$ by $O'$; it is also the $\psi$-orbit of $e'$ and, in the relevant case where $e \neq e'$, it has length $p$. Let $W := \mathrm{Span}(O')$. By the preceding proposition, $V \hat{\otimes} W$ is a simple two-sided ideal of $A \hat{\otimes} A$ of dimension $p^2$; $V \hat{\otimes} Ae$ and $V \hat{\otimes} Ae'$ are $p$-dimensional left ideals of $A \hat{\otimes} A$, which are contained in $V \hat{\otimes} W$. Since all simple modules of $V \hat{\otimes} W$ are isomorphic and of dimension $p$, the assertion follows. □

**6.6** We have seen in the preceding paragraph that the product of two one-dimensional characters $\eta$ and $\eta'$ is again a one-dimensional character if $\eta$ is linear or $\eta'$ is colinear over $H$, but is in general not a character in other cases. We now investigate what can be said instead.

**Proposition** Suppose that $\eta : A \to K$ is a character that is not $H$-linear and that $\eta' : A \to K$ is a character that is not colinear. Then there are distinct characters $\omega_0, \ldots, \omega_{p-1}$ with the following properties:

1. If $\eta_0, \ldots, \eta_{p-1}$ are the conjugates of $\eta$ with respect to $\phi^*$ and $\eta'_0, \ldots, \eta'_{p-1}$ are the conjugates of $\eta'$ with respect to $\psi^*$, we have

$$\eta_i \eta'_j \in \mathrm{Span}(\{\omega_0, \ldots, \omega_{p-1}\})$$

   for all $i, j \in \mathbb{Z}_p$. In particular, we have $\eta \eta' \in \mathrm{Span}(\{\omega_0, \ldots, \omega_{p-1}\})$.

2. $\{\omega_0, \ldots, \omega_{p-1}\}$ is invariant under $\phi^*$ and $\psi^*$.

3. If $\eta$ is colinear, $\{\omega_0, \ldots, \omega_{p-1}\}$ is an orbit with respect to $\psi^*$.

4. If $\eta'$ is $H$-linear, $\{\omega_0, \ldots, \omega_{p-1}\}$ is an orbit with respect to $\phi^*$.

**Proof.** (1) We can enumerate $\eta_0, \ldots, \eta_{p-1}$ in such a way that $\eta_i = \phi^{*i}(\eta)$. Similarly, we can enumerate $\eta'_0, \ldots, \eta'_{p-1}$ in such a way that $\eta'_j = \psi^{*j}(\eta')$. Denote by $e_i$ resp. $e'_j$ the primitive idempotents in $E$ satisfying $\eta_i(e_i) = 1$ resp. $\eta'_j(e'_j) = 1$. Then $O := \{e_0, \ldots, e_{p-1}\}$ is an orbit with respect to $\phi$ and $O' = \{e'_0, \ldots, e'_{p-1}\}$ is an orbit with respect to $\psi$. Define $V := \text{Span}(O)$ and $W := \text{Span}(O')$. For the primitive idempotents belonging to $\eta$ resp. $\eta'$, we use the abbreviation $e := e_0$ resp. $e' := e'_0$.

Since $A$ is commutative, the $A$-module $V \hat{\otimes} Ae'$ can be decomposed into a direct sum of one-dimensional $A$-modules. Therefore, there is a basis $v_0, \ldots, v_{p-1}$ of $V$ such that
$$a(v_k \otimes e') = \omega_k(a) v_k \otimes e'$$
for characters $\omega_0, \ldots, \omega_{p-1}$. Since $a(v_k \otimes e') = a_{(1)}(a_{(2)}{}^{(1)} \to v_k) \otimes a_{(2)}{}^{(2)} e'$, we get by applying $\text{id}_V \otimes \eta'$:
$$\omega_k(a) v_k = a_{(1)}(a_{(2)}{}^{(1)} \to v_k) \eta'(a_{(2)}{}^{(2)})$$

(2) Since $v_0, \ldots, v_{p-1}$ is a basis of $V$, the invariant idempotent $e_I := \sum_{k=0}^{p-1} e_k$ can be written in the form $e_I = \sum_{k=0}^{p-1} \lambda_k v_k$ for some $\lambda_0, \ldots, \lambda_{p-1} \in K$. We then have:
$$a_{(1)} e_I \eta'(a_{(2)}) = a_{(1)}(a_{(2)}{}^{(1)} \to e_I) \eta'(a_{(2)}{}^{(2)})$$
$$= \sum_{k=0}^{p-1} \lambda_k a_{(1)}(a_{(2)}{}^{(1)} \to v_k) \eta'(a_{(2)}{}^{(2)}) = \sum_{k=0}^{p-1} \lambda_k \omega_k(a) v_k$$

Applying $\eta_i$ to this equation, we get:
$$(\eta_i \eta')(a) = \sum_{k=0}^{p-1} \lambda_k \eta_i(v_k) \omega_k(a)$$

Therefore, we have $\eta_i \eta' \in \text{Span}(\{\omega_0, \ldots, \omega_{p-1}\})$. In particular, $\omega_0, \ldots, \omega_{p-1}$ must be distinct, since $\eta_0 \eta', \ldots, \eta_{p-1} \eta'$ are linearly independent.

(3) It can be verified directly that we have
$$(\psi \otimes \psi)(a(v \otimes w)) = \psi(a)(\psi(v) \otimes \psi(w))$$
for $a \in A$, $v \in V$ and $w \in W$. Since $\psi(e'_0) = e'_{p-1}$, we get by applying $\psi \otimes \psi$ to the equation $a(v_i \otimes e'_0) = \omega_i(a) v_i \otimes e'_0$ that:
$$a(\psi(v_i) \otimes e'_{p-1}) = \omega_i(\psi^{-1}(a)) \psi(v_i) \otimes e'_{p-1}$$

According to Proposition 6.4, $\psi(v_0), \ldots, \psi(v_{p-1})$ is again a basis of $V$. Since the $A$-modules $V \otimes Ae'_0$ and $V \otimes Ae'_{p-1}$ are isomorphic by Corollary 6.5, the sets of one-dimensional characters occurring in both modules must coincide, i. e., we have:
$$\{\psi^{-1*}(\omega_0), \ldots, \psi^{-1*}(\omega_{p-1})\} = \{\omega_0, \ldots, \omega_{p-1}\}$$

Therefore, $\{\omega_0, \ldots, \omega_{p-1}\}$ is invariant under $\psi^*$.

(4) As for $\psi$, we have for $\phi$ that

$$(\phi \otimes \phi)(a(v \otimes w)) = \phi(a)(\phi(v) \otimes \phi(w))$$

for $a \in A$, $v \in V$, and $w \in W$. By Proposition 6.4, we have that $\phi(e_0') = e_j'$ for some $j \in \mathbb{Z}_p$. Applying $\phi \otimes \phi$ to the equation $a(v_i \otimes e_0') = \omega_i(a)v_i \otimes e_0'$, we get:

$$a(\phi(v_i) \otimes e_j') = \omega_i(\phi^{-1}(a))\phi(v_i) \otimes e_j'$$

Since $V \otimes Ae_0' \cong V \otimes Ae_j'$ by Corollary 6.5, we get as above that $\{\omega_0, \ldots, \omega_{p-1}\}$ is invariant under $\phi^*$.

(5) Since $\{\omega_0, \ldots, \omega_{p-1}\}$ is invariant under $\psi^*$, we have that:

$$\psi^{*j}(\eta_i)\psi^{*j}(\eta') = \psi^{*j}(\eta_i\eta') \in \text{Span}(\{\omega_0, \ldots, \omega_{p-1}\})$$

By Proposition 6.4, we have $\{\psi^{*j}(\eta_0), \ldots, \psi^{*j}(\eta_{p-1})\} = \{\eta_0, \ldots, \eta_{p-1}\}$. This implies the first assertion of the proposition.

(6) Now suppose that $\eta$ is colinear. If $\{\omega_0, \ldots, \omega_{p-1}\}$ is not an orbit with respect to $\psi^*$, it consists of colinear characters. Since $\eta\eta' \in \text{Span}(\{\omega_0, \ldots, \omega_{p-1}\})$, we have $\eta' \in \text{Span}(\{\eta^{-1}\omega_0, \ldots, \eta^{-1}\omega_{p-1}\})$. We have seen in Paragraph 1.5 that colinear characters form a group, therefore the linear functions $\eta^{-1}\omega_k$ are colinear characters. Since distinct characters are linearly independent, we have $\eta' = \eta^{-1}\omega_k$ for some $k \in \{0, \ldots, p-1\}$. Therefore, $\eta'$ is colinear, which contradicts our assumptions.

(7) Now suppose that $\eta'$ is $H$-linear. If $\{\omega_0, \ldots, \omega_{p-1}\}$ is not an orbit with respect to $\phi^*$, it consists of $H$-linear characters. As in the previous step, we get that $\eta \in \text{Span}(\{\omega_0\eta'^{-1}, \ldots, \omega_{p-1}\eta'^{-1}\})$ is a linear combination of $H$-linear characters, and since we then have $\eta = \omega_k\eta'^{-1}$ for some $k \in \{0, \ldots, p-1\}$, $\eta$ is $H$-linear itself, which contradicts our assumptions. $\square$

**Corollary** If all primitive idempotents are invariant or coinvariant, then all primitive idempotents are invariant or all primitive idempotents are coinvariant.

**Proof.** If this is not the case, there is a primitive idempotent $e$ that is not invariant and a primitive idempotent $e'$ that is not coinvariant. By assumption, we then have that $e$ is coinvariant and that $e'$ is invariant. We denote the characters corresponding to $e$ resp. $e'$ by $\eta$ resp. $\eta'$. By the preceding proposition, we then have characters $\{\omega_0, \ldots, \omega_{p-1}\}$ satisfying $\psi^*(\omega_j) = \omega_{j+1}$ such that $\eta\eta' \in \text{Span}(\omega_0, \ldots, \omega_{p-1})$. Since $\omega_j$ is not colinear, it must be linear, i. e., we have $\phi^*(\omega_j) = \omega_j$ and therefore $\phi^*(\eta\eta') = \eta\eta'$. This implies:

$$\phi^*(\eta)\eta' = \phi^*(\eta)\phi^*(\eta') = \phi^*(\eta\eta') = \eta\eta'$$

Since $\eta'$ is an invertible element by Proposition 1.5.2, we have $\phi^*(\eta) = \eta$, which is a contradiction. $\square$

**6.7** The proposition proved in the previous paragraph has particularly striking consequences if applied to the situation where the two characters $\eta$ and $\eta'$ correspond to dual $B$-modules:

**Proposition** Suppose that $\eta : A \to K$ is a character which is neither $H$-linear nor colinear. Then there exists a character $\omega : A \to K$ with the following properties:

1. $\omega : A \to K$ is an $H$-linear and colinear character of order $p$.

2. We have: $\phi^*(\eta) = \omega\eta = \eta\omega^j$ for some $j \in \mathbb{Z}_p$.

3. If $e$ is the primitive idempotent corresponding to $\eta$, $O$ the orbit of $e$ under the action of $C$, and $V := \mathrm{Span}(O)$ the corresponding simple $B$-module, the $B$-module $V \otimes V^*$ decomposes into $p^2$ one-dimensional modules, whose characters are $\omega^i \otimes \gamma^j$, for $i, j = 0, \ldots, p-1$ and a generator $\gamma$ of $\hat{C}$.

In particular, $p$ divides $\dim A$.

**Proof.** (1) Since $V^*$ is also a simple $B$-module, there is a $C$-orbit $O' = \{e'_0, \ldots, e'_{p-1}\}$ such that $V^* \cong W := \mathrm{Span}(O')$. By Proposition 6.4, we have that $O$ and $O'$ are invariant under $\psi$. The fact that $\eta$ is neither linear nor colinear means that $O$ contains no fixed points of $\phi$ and $\psi$, and therefore, according to Proposition 6.3, $O'$ contains no fixed points of $\phi$ and $\psi$, and thus is an orbit with respect to $\phi$ and $\psi$.

As in the proof of the preceding proposition, there exists a basis $v_0, \ldots, v_{p-1}$ of $V$ and distinct characters $\omega_0, \ldots, \omega_{p-1}$ such that:

$$a(v_i \otimes e'_0) = \omega_i(a) v_i \otimes e'_0$$

(2) From the decomposition

$$V \otimes V^* \cong V \otimes W = \bigoplus_{i=0}^{p-1} V \otimes Ae'_i$$

we see, since the $A$-modules $V \otimes Ae'_i$ are all isomorphic by Corollary 6.5, that the $A$-module $V \otimes V^*$ decomposes into one-dimensional modules corresponding to the characters $\omega_0, \ldots, \omega_{p-1}$, each occurring with multiplicity $p$. By Schur's lemma, one of these characters must be $\epsilon_A$, which is a fixed point of $\phi^*$ and $\psi^*$. Therefore, the set $\{\omega_0, \ldots, \omega_{p-1}\}$, which is invariant under $\phi^*$ and $\psi^*$ by Proposition 6.6, cannot be an orbit of length $p$ with respect to these maps, and therefore must contain only fixed points. This implies:

$$\phi^*(\omega_i) = \omega_i \qquad \psi^*(\omega_i) = \omega_i$$

for all $i = 0, \ldots, p-1$. Therefore, the $B$-module $V \otimes V^*$ cannot contain a purely unstable submodule, since its restriction to $A$ would not correspond to $H$-linear characters. Therefore, $V \otimes V^*$ decomposes into stable, i. e., one-dimensional, submodules.

(3) We denote by $\chi_V$, resp. $\chi_{V^*}$, the characters of $V$, resp. $V^*$. We then know from Theorem 2.6 that $\chi_{V^*}(\epsilon_A \otimes \gamma') = \chi_{V^*}$ for all characters $\gamma'$ of $H$. Therefore, we also have $(\epsilon_A \otimes \gamma')\chi_V = \chi_V$. Now, for $i = 0, \ldots, p-1$, we can find a one-dimensional submodule of the $B$-module $V \otimes V^*$ whose restriction to $A$ has the character $\omega_i$. The character of the submodule considered as a $B$-module then has the form $\omega_i \otimes \gamma'$ for some character $\gamma'$ of $H$. From Paragraph 1.12, we then know that $(\omega_i \otimes \gamma')\chi_V = \chi_V$, and therefore we have $(\omega_i \otimes \epsilon_H)\chi_V = \chi_V$. By reversing the argument, we see that

$$\{\omega_0, \ldots, \omega_{p-1}\} = \{\omega' \in G_I(A^*) \mid (\omega' \otimes \epsilon_H)\chi_V = \chi_V\}$$

where, as in Paragraph 1.5, $G_I(A^*)$ denotes the set of $H$-linear characters of $A$. In particular, $\{\omega_0, \ldots, \omega_{p-1}\}$ is a subgroup of $G_I(A^*)$ of order $p$. By changing the enumeration of the characters, we can therefore assume that $\omega_i = \omega_1^i$ for $i = 0, \ldots, p-1$. In particular, all these characters, except $\epsilon_A$, have order $p$. In addition, we that see that $(\omega_i \otimes \gamma^j)\chi_V = \chi_V$ for all $i, j = 0, \ldots, p-1$, and therefore all the distinct characters $\omega_i \otimes \gamma^j$ are contained in $V \otimes V^*$. Since a sum of nonisomorphic simple modules is direct, we see by counting dimensions that $V \otimes V^*$ is the direct sum of one-dimensional submodules corresponding to the characters $\omega_i \otimes \gamma^j$.

(4) We denote the characters that correspond to the idempotents in $O$ by $\eta_0, \ldots, \eta_{p-1}$, with the convention that $\eta_0 = \eta$. From Proposition 2.4, we have $\chi_V(a \otimes 1_H) = \sum_{j=0}^{p-1} \eta_j(a)$. Therefore, the fact that $(\omega_i \otimes \epsilon_H)\chi_V = \chi_V$ implies that $\sum_{j=0}^{p-1} \omega_i \eta_j = \sum_{j=0}^{p-1} \eta_j$. Since $\omega_i \eta_j$ is again a character by Proposition 1.5.2, this implies $\{\omega_i \eta_0, \ldots, \omega_i \eta_{p-1}\} = \{\eta_0, \ldots, \eta_{p-1}\}$. But this means that we have $\{\omega_0 \eta_j, \ldots, \omega_{p-1} \eta_j\} = \{\eta_0, \ldots, \eta_{p-1}\}$, and therefore we must have $\phi^*(\eta) = \omega_i \eta$ for some $i \in \mathbb{Z}_p$. Since $\eta$ is not $H$-linear, we cannot have $\omega_i = \epsilon_A$. If we define $\omega := \omega_i$, $\omega$ therefore is an $H$-linear and colinear character of order $p$ satisfying $\phi^*(\eta) = \omega\eta$.

(5) Now denote the character corresponding to $e_i'$ by $\eta_i'$ and set $\eta' := \eta_0'$. By applying the results proved so far to $\eta'$ instead of $\eta$, we see that there is an $H$-linear and colinear character $\omega' : A \to K$ of order $p$ such that $\phi^*(\eta') = \omega'\eta'$. From the previous paragraph, we know that $\eta_0 \eta', \ldots, \eta_{p-1}\eta'$ as well as $\eta\eta_0', \ldots, \eta\eta_{p-1}'$ are linearly independent subsets of $\mathrm{Span}(\omega^0, \omega^1, \ldots, \omega^{p-1})$. This implies that

$$\mathrm{Span}(\eta_0\eta', \ldots, \eta_{p-1}\eta') = \mathrm{Span}(\eta\omega'^0\eta', \ldots, \eta\omega'^{p-1}\eta')$$

and therefore $\mathrm{Span}(\eta_0, \ldots, \eta_{p-1}) = \mathrm{Span}(\eta\omega'^0, \ldots, \eta\omega'^{p-1})$. Therefore, we get that $\phi^*(\eta) = \eta\omega'^j$ for some $j \in \{1, \ldots, p-1\}$.

By a similar argument, we see that

$$\mathrm{Span}(\eta_0\eta, \ldots, \eta_{p-1}\eta) = \mathrm{Span}(\eta\eta_0, \ldots, \eta\eta_{p-1})$$

and therefore $\mathrm{Span}(\eta\omega'^0\eta, \ldots, \eta\omega'^{p-1}\eta) = \mathrm{Span}(\eta\omega^0\eta, \ldots, \eta\omega^{p-1}\eta)$. Multiplying on the left and on the right by $\eta^{-1}$, we deduce that:

$$\mathrm{Span}(\omega'^0, \ldots, \omega'^{p-1}) = \mathrm{Span}(\omega^0, \ldots, \omega^{p-1})$$

This implies that $\omega'$ is a power of $\omega$. Replacing $j$ by a different one, we arrive at $\phi^*(\eta) = \eta\omega^j$.

(6) From Paragraph 1.2, we know that $A^*$ is a (right) Yetter-Drinfel'd Hopf algebra over $H$. The subspace spanned by the powers of $\omega$ is obviously a Yetter-Drinfel'd Hopf subalgebra of $A^*$. Since, by Proposition 1.8 and Lemma 1.2, $A^*$ is free as a right module over this subalgebra, its dimension divides the dimension of $A^*$. Therefore, $p$ divides $\dim A$. $\square$

In the case where $A$ is nontrivial, we know from Paragraph 1.11 that the action and the coaction must be nontrivial. Therefore, it follows from Corollary 6.6 that there is an idempotent which is neither invariant nor coinvariant, i. e., a character which is neither linear nor colinear. By the preceding proposition, this implies that $p \mid \dim A$. Therefore, we have the following corollary:

**Corollary** If $A$ is a nontrivial, $p$ divides $\dim A$.

This corollary will be substantially sharpened by the structure theorems that will be obtained in Paragraph 7.7 and Paragraph 7.8.

**6.8** We have seen in the preceding paragraph that the action and the coaction on a given nonlinear and noncolinear character is induced by some other character of order $p$. We shall prove now that the action and the coaction are induced universally by a single character of order $p$.

**Proposition** Suppose that $A$ is nontrivial. Then there exists an $H$-linear and colinear character $\omega : A \to K$ of order $p$ such that for all characters $\eta : A \to K$ we have
$$\phi^*(\eta) = \omega^i\eta = \eta\omega^k \qquad \psi^*(\eta) = \omega^j\eta = \eta\omega^l$$
for some $i, j, k, l \in \mathbb{Z}_p$.

**Proof.** (1) By Corollary 6.6, there exists a character $\eta : A \to K$ which is neither linear nor colinear over $H$. By Proposition 6.7, there exists an $H$-linear and colinear character $\omega : A \to K$ of order $p$ such that $\phi^*(\eta) = \omega\eta = \eta\omega^m$ for some $m \in \mathbb{Z}_p$. Now suppose that $\eta'$ is any one-dimensional character of $A$. We have to prove that a similar equation holds with $\eta'$ instead of $\eta$. For this, we distinguish various cases.

(2) First, suppose that $\eta'$ is neither linear nor colinear over $H$. Since, by Proposition 6.4, $\psi^*(\eta') = \phi^{*r}(\eta')$ for some $r \in \mathbb{Z}_p$, Proposition 6.7 yields that there exists an $H$-linear and colinear character $\omega' : A \to K$ of order $p$ such that $\psi^*(\eta') = \omega'\eta' = \eta'\omega'^n$ for some $n \in \mathbb{Z}_p$. We then know from Proposition 6.6 that there are characters $\eta''_0, \ldots, \eta''_{p-1}$ such that the elements $\eta\omega^0\eta', \ldots, \eta\omega^{p-1}\eta'$ as well as the elements $\eta\omega'^0\eta', \ldots, \eta\omega'^{p-1}\eta'$ are contained in $\mathrm{Span}(\eta_0'', \ldots, \eta_{p-1}'')$. Since both families consist of linear independent vectors, we have

$$\mathrm{Span}(\eta\omega^0\eta', \ldots, \eta\omega^{p-1}\eta') = \mathrm{Span}(\eta\omega'^0\eta', \ldots, \eta\omega'^{p-1}\eta')$$

which, by multiplication by $\eta^{-1}$ and $\eta'^{-1}$, implies that $\mathrm{Span}(\omega^0, \ldots, \omega^{p-1}) = \mathrm{Span}(\omega'^0, \ldots, \omega'^{p-1})$. Therefore, $\omega'$ is a power of $\omega$, and the assertion follows in this case.

(3) Now suppose that $\eta'$ is colinear, but not linear over $H$. Consider the character $\eta'' := \eta\eta'$. Then $\eta''$ is not colinear, because in this case, by Proposition 1.5.2, $\eta = \eta''\eta'^{-1}$ would be colinear, too. We now treat separately the two cases where $\eta''$ is $H$-linear and where it is not $H$-linear. Suppose first that $\eta''$ is $H$-linear. We then have:

$$\eta\omega^m \phi^*(\eta') = \phi^*(\eta)\phi^*(\eta') = \phi^*(\eta'') = \eta'' = \eta\eta'$$

and therefore $\phi^*(\eta') = \omega^{-m}\eta'$.

Now suppose that $\eta''$ is not $H$-linear. Then we can apply the above result to conclude that $\phi^*(\eta'') = \omega^n \eta''$ for some other $n \in \mathbb{Z}_p$ that is not necessarily equal to the one above. This implies:

$$\omega^n \eta\eta' = \phi^*(\eta'') = \phi^*(\eta)\phi^*(\eta') = \eta\omega^m \phi^*(\eta')$$

Since $\omega\eta = \eta\omega^m$, this yields:

$$\eta\omega^{nm}\eta' = \eta\omega^m \phi^*(\eta')$$

and therefore $\phi^*(\eta') = \omega^{(n-1)m}\eta'$.

Since $\eta'^{-1}$ is also a character which is colinear, but not linear over $H$, we have $\phi^*(\eta'^{-1}) = \omega^k \eta'^{-1}$ for some $k \in \mathbb{Z}_p$. This implies $\phi^*(\eta') = \eta'\omega^{-k}$.

(4) Now suppose that $\eta'$ is linear, but not colinear over $H$. Consider the character $\eta'' := \eta'\eta$. Then $\eta''$ is not $H$-linear, because in this case, by Proposition 1.5.2, $\eta = \eta'^{-1}\eta''$ would be $H$-linear, too. We now treat separately the two cases where $\eta''$ is colinear and where it is not colinear. Suppose first that $\eta''$ is colinear. If $\psi^*(\eta) = \omega^s \eta = \eta\omega^{sm}$, we have

$$\psi^*(\eta')\omega^s \eta = \psi^*(\eta')\psi^*(\eta) = \psi^*(\eta'') = \eta'' = \eta'\eta$$

and therefore $\psi^*(\eta') = \eta'\omega^{-s}$.

Now suppose that $\eta''$ is not colinear. Then we can apply the above result to conclude that $\psi^*(\eta'') = \eta''\omega^k$ for some other $k \in \mathbb{Z}_p$. This implies:

$$\eta'\eta\omega^k = \psi^*(\eta'') = \psi^*(\eta')\psi^*(\eta) = \psi^*(\eta')\omega^s \eta$$

Since $\eta\omega = \omega^t \eta$, where $t \in \mathbb{Z}_p$ is the multiplicative inverse of $m$, this yields:

$$\eta'\omega^{kt}\eta = \psi^*(\eta')\omega^s \eta$$

and therefore $\psi^*(\eta') = \eta'\omega^{kt-s}$.

Since $\eta'^{-1}$ is also a character which is linear, but not colinear over $H$, we have $\psi^*(\eta'^{-1}) = \eta'^{-1}\omega^r$ for some other $r \in \mathbb{Z}_p$. This implies $\psi^*(\eta') = \omega^{-r}\eta'$.

(5) The remaining case is that $\eta'$ is linear as well as colinear over $H$. This is the trivial case, since the assertion is satisfied with $i, j, k, l = 0$. □

# 7 Cocommutative Yetter-Drinfel'd Hopf algebras

**7.1** In this section, we assume that $p$ is a prime and that $K$ is an algebraically closed field whose characteristic is different from $p$. We can therefore choose a primitive $p$-th root of unity that we denote by $\zeta$. The group ring of the cyclic group $\mathbb{Z}_p$ of order $p$ will be denoted by $H := K[\mathbb{Z}_p]$; the canonical basis element of the group ring corresponding to $i \in \mathbb{Z}_p$ will be denoted by $c_i$. The group of grouplike elements of $H$, which consists precisely of these canonical basis elements, will be denoted by $C$. $A$ denotes a nontrivial Yetter-Drinfel'd Hopf algebra over $H$ that is cocommutative and cosemisimple. $A$ therefore has a unique basis that consists of grouplike elements, which is denoted by $G(A)$. If $\gamma : C \to K^\times$ is the group homomorphism that maps $c_1$ to $\zeta$, we introduce as in Paragraph 1.10 the mappings

$$\phi : V \to V, v \mapsto (c_1 \to v) \qquad \psi : V \to V, v \mapsto \gamma(v^{(1)})v^{(2)}$$

where we extend $\gamma$ to $H$ by linearity. (Note that the symbol $\gamma$, in contrast to $\psi$, will also be used for other characters below.) Since $\phi$ and $\psi$ are coalgebra automorphisms, they induce permutations of $G(A)$. The Radford biproduct $A \otimes H$ will be denoted by $B$.

Throughout the section, we will constantly use the convention that indices take values between 0 and $p-1$ and are reduced modulo $p$ if they do not lie within this range. In notation, we shall not distinguish between an integer $i \in \mathbb{Z}$ and its equivalence class in $\mathbb{Z}_p := \mathbb{Z}/p\mathbb{Z}$.

The aim of the section is to prove that $A$ arises from the construction explained in Section 3, i. e., that $A$ is, as an algebra, a crossed product of the dual group ring $K^{\mathbb{Z}_p}$ and the group ring of a certain group, and as a coalgebra the ordinary tensor product coalgebra of these spaces.

**7.2** First, we dualize the results of Section 6:

**Proposition** There exists a grouplike element which is not invariant and not coinvariant. Furthermore, there exists a grouplike element $u$ of order $p$ which is invariant and coinvariant such that, for all grouplike elements $g \in G(A)$, there are numbers $i, j, k, l \in \mathbb{Z}_p$ such that:

$$\phi(g) = u^i g = g u^k \qquad \psi(g) = u^j g = g u^l$$

In particular, $p$ divides $\dim A$.

**Proof.** From Paragraph 1.2, we know that $A^*$ is a right Yetter-Drinfel'd Hopf algebra over $H$, and therefore, by Lemma 1.2, we know that $A^{*\text{op cop}}$ is a left Yetter-Drinfel'd Hopf algebra over $H$. The existence of grouplike elements that are not invariant and not coinvariant now follows by applying Corollary 6.6 to $A^{*\text{op cop}}$, and the existence of $u$ follows by applying Proposition 6.8 to $A^{*\text{op cop}}$. $\square$

For the rest of this section, we fix a grouplike element $u$ of order $p$ that has the properties stated in the preceding proposition. The subspace spanned by the powers of $u$ will be denoted by $U$; it is obviously a Yetter-Drinfel'd Hopf subalgebra.

**7.3** In this paragraph, we shall construct a Yetter-Drinfel'd Hopf algebra quotient of $A$ in which the equivalence class of $u$ is equal to the unit. Since action and coaction are induced by $u$, they become trivial in this quotient, and therefore this quotient is an ordinary Hopf algebra which is, in addition, cocommutative. Therefore, it is isomorphic to a group ring.

We shall use the usual notation $U^+ := \ker(\epsilon_U) = U \cap \ker \epsilon_A$. Observe that, since the product of $u$ and a grouplike element $g$ is again a grouplike element, $G(A)$ can be decomposed into orbits with respect to left multiplication by $u$.

**Proposition**
1. $u - 1, u^2 - 1, \ldots, u^{p-1} - 1$ is a basis of $U^+$.

2. $AU^+$ is a two-sided ideal, a two-sided coideal, an $H$-submodule and an $H$-subcomodule of $A$ that is invariant with respect to the antipode $S_A$. $A/AU^+$ is therefore a Yetter-Drinfel'd Hopf algebra.

3. If $g_1, \ldots, g_n \in G(A)$ is a system of representatives for the orbits of the action of $u$ on $G(A)$, the equivalence classes $\bar{g}_1, \ldots, \bar{g}_n$ form a basis of $A/AU^+$ consisting of grouplike elements. Therefore, the dimension of $A/AU^+$ is $\frac{1}{p} \dim A$.

4. Action and coaction of $H$ on $A/AU^+$ are trivial. $A/AU^+$ is therefore an ordinary cocommutative Hopf algebra, $G(A/AU^+) = \{\bar{g}_1, \ldots, \bar{g}_n\}$ is a group, and we have $A/AU^+ \cong K[G(A/AU^+)]$.

**Proof.** The first statement holds because $1, u - 1, u^2 - 1, \ldots, u^{p-1} - 1$ is a basis of $U$. It is obvious that $AU^+$ is a left ideal of $A$. We now prove that $AU^+ = U^+A$, which implies that $AU^+$ is also a right ideal. Since the grouplike elements form a basis of $A$, this will follow if we can prove that $Ug = gU$ for all grouplike elements $g \in G(A)$. If $g$ is not invariant or not coinvariant, we have by Proposition 7.2 that

$$\phi(g) = u^i g = g u^j \quad \text{or} \quad \psi(g) = u^i g = g u^j$$

for some nonzero elements $i, j \in \mathbb{Z}_p$, which implies the assertion. If $g$ is invariant and coinvariant, choose a grouplike element $g'$ that is neither invariant nor coinvariant, which exists by Proposition 7.2. Since $gg'$ is also neither invariant nor coinvariant, we then have
$$g'u = u^k g' \quad \text{and} \quad ugg' = gg'u^l$$
for some nonzero elements $k, l \in \mathbb{Z}_p$. This implies $ugg' = gu^{kl}g'$ and therefore $ug = gu^{kl}$.

We have $\Delta_A(U^+) \subset \ker(\epsilon_U \otimes \epsilon_U) = U^+ \otimes U + U \otimes U^+$. If $a \in A$ and $a' \in U^+$, we have
$$\Delta_A(aa') = a_{(1)}(a_{(2)}{}^{(1)} \to a'_{(1)}) \otimes a_{(2)}{}^{(2)} a'_{(2)}$$
$$= a_{(1)} a'_{(1)} \otimes a_{(2)} a'_{(2)} \in AU^+ \otimes A + A \otimes AU^+$$

Therefore, $AU^+$ is a two-sided coideal. Since $u$ is invariant and coinvariant, it is an $H$-submodule and an $H$-subcomodule. Since we have
$$S_A(aa') = S_A(a^{(1)} \to a') S_A(a^{(2)}) = S_A(a') S_A(a) \in U^+ A = AU^+$$
for $a' \in U^+$, it is invariant with respect to the antipode.

To prove the third statement, observe that, since $A = \bigoplus_{i=1}^n g_i U$, we have $AU^+ = \bigoplus_{i=1}^n g_i U^+$, and therefore:
$$A/AU^+ = \bigoplus_{i=1}^n g_i U / g_i U^+ = \bigoplus_{i=1}^n K \bar{g}_i$$

Since the $n$ orbits all consist of $p$ elements, we have $\dim A = pn$.

To prove the last statement, observe that, if $g \in G(A)$ is an arbitrary grouplike element, we have by Proposition 7.2 that
$$\phi(g) = gu^i \qquad \psi(g) = gu^j$$
for some $i, j \in \mathbb{Z}_p$. This implies that $\phi(g) - g = g(u^i - 1) \in AU^+$, which means that $\phi$, and similarly $\psi$, induces the identity on $A/AU^+$. Therefore, action and coaction on $A/AU^+$ are trivial. The remaining assertions are obvious. $\square$

**7.4** We introduce some more notation. We denote by $G := G(A/AU^+)$ the group of grouplike elements of $A/AU^+$. In the proof of Proposition 7.3, we have seen that, for all $g \in G(A)$, there exists $j \in \mathbb{Z}_p$ such that $gu = u^j g$. The number $j$ obviously only depends on the $u$-orbit of $g$. We therefore have a map $\nu : G \to \mathbb{Z}_p^\times$ such that
$$gu = u^{\nu(\bar{g})} g$$
where $\mathbb{Z}_p^\times = \mathbb{Z}_p \backslash \{0\}$ denotes the multiplicative group of the finite field $\mathbb{Z}_p$.

By Proposition 7.2, we have, for every grouplike element $g \in G(A)$, numbers $i, j \in \mathbb{Z}_p$ such that:
$$\phi(g) = u^i g \qquad \psi(g) = u^j g$$

Again the numbers $i$ and $j$ only depend on the orbit of $g$ with respect to left multiplication by $u$, since $u$ is invariant and coinvariant. We therefore have mappings $\alpha : G \to \mathbb{Z}_p$ and $\beta : G \to \mathbb{Z}_p$ such that:
$$\phi(g) = u^{\alpha(\bar{g})} g \qquad \psi(g) = u^{\beta(\bar{g})} g$$

**Proposition**
1. $\nu : G \to \mathbb{Z}_p^\times$ is a group homomorphism. Therefore, $\mathbb{Z}_p$ is a $G$-module via $g.i := \nu(g)i$.

2. $\alpha$ and $\beta$ are 1-cocycles with respect to this module structure, i. e., we have:
$$\alpha(st) = \alpha(s) + \nu(s)\alpha(t) \quad \beta(st) = \beta(s) + \nu(s)\beta(t)$$
for all $s, t \in G$.

**Proof.** For $g, g' \in G(A)$, we have:
$$u^{\nu(\bar{g}\bar{g}')} g g' = g g' u = g u^{\nu(\bar{g}')} g' = u^{\nu(\bar{g})\nu(\bar{g}')} g g'$$

This implies that $\nu(\bar{g}\bar{g}') = \nu(\bar{g})\nu(\bar{g}')$. Similarly, we have:
$$u^{\alpha(\bar{g}\bar{g}')} g g' = \phi(gg') = \phi(g)\phi(g') = u^{\alpha(\bar{g})} g u^{\alpha(\bar{g}')} g'$$
$$= u^{\alpha(\bar{g})} u^{\nu(\bar{g})\alpha(\bar{g}')} g g' = u^{\alpha(\bar{g}) + \nu(\bar{g})\alpha(\bar{g}')} g g'$$

which implies $\alpha(\bar{g}\bar{g}') = \alpha(\bar{g}) + \nu(\bar{g})\alpha(\bar{g}')$. The proof for $\beta$ is similar. $\square$

As in Paragraph 3.3, we will denote $\mathbb{Z}_p$ by ${}_G\mathbb{Z}_p$ if it is considered as a $G$-module via $\nu$ as in the preceding proposition.

**7.5** The fact that $A/AU^+$ is a quotient coalgebra of $A$ leads to a comodule structure of $A$ over $A/AU^+$. If $\pi : A \to A/AU^+$ denotes the canonical projection, we introduce the right coaction
$$\delta_G : A \to A \otimes A/AU^+, a \mapsto a_{(1)} \otimes \bar{a}_{(2)}$$

which should be distinguished from the left coaction $\delta_A$ that is part of the Yetter-Drinfel'd structure. Since $A/AU^+ \cong K[G]$ is a group ring, $A$ becomes a $G$-graded vector space, where, for a grouplike element $g \in G(A)$, the homogeneous component corresponding to $\bar{g} \in G$ is:
$$A_{\bar{g}} := \{a \in A \mid \delta_G(a) = a \otimes \bar{g}\}$$

Now we select for every $s \in G$ a representative $g_s \in G(A)$ that satisfies $\bar{g}_s = s$, where we choose $g_1 = 1_A$. We then have $G(A) = \{u^j g_s \mid j \in \mathbb{Z}_p, s \in G\}$.

**Proposition**

1. $A$ is a right comodule algebra with respect to $A/AU^+$.

2. For $s \in G$, the homogeneous component $A_s$ is:
$$A_s = \operatorname{Span}(g_s, ug_s, u^2 g_s, \ldots, u^{p-1} g_s)$$
In particular, we have $A_1 = U$.

3. The extension $U \subset A$ is cleft over $K[G]$ with respect to the map
$$K[G] \to A, s \mapsto g_s$$

4. $A$ is a crossed product of $U$ and $K[G]$, with cocycle
$$\sigma : G \times G \to U, (s,t) \mapsto g_s g_t g_{st}^{-1}$$
with respect to the $G$-module structure on $U$ determined by $s.u = u^{\nu(s)}$.

**Proof.** The first statement follows from the fact that $A/AU^+$ is a trivial $H$-comodule: For $a, a' \in A$, we have:

$$\delta_G(aa') = a_{(1)}(a_{(2)}{}^{(1)} \to a'_{(1)}) \otimes \bar{a}_{(2)}{}^{(2)} \bar{a}'_{(2)} = a_{(1)} a'_{(1)} \otimes \bar{a}_{(2)} \bar{a}'_{(2)} = \delta_G(a)\delta_G(a')$$

To prove the second statement, observe that, since $u$ is homogeneous of degree 1, we have $u^i g_s \in A_s$. Since we already have:
$$A = \bigoplus_{s \in G} \operatorname{Span}(g_s, ug_s, \ldots, u^{p-1} g_s)$$
and the decomposition into homogeneous components is also direct, we must have $A_s = \operatorname{Span}(g_s, ug_s, \ldots, u^{p-1} g_s)$.

The third assertion is obvious (cf. [17], p. 806, [57], Def. 7.2.1, p. 105 for the definition of cleft extensions). Cleft extensions are always crossed products (cf. [17], Thm. 11, p. 815, [57], Thm. 7.2.2, p. 106, see also [64], Chap. 1, p. 2); the cocycle arises from the cleaving map in the way described in the proposition, whereas the action is determined by the condition
$$s.u = g_s u g_s^{-1} = u^{\nu(s)}$$
since $u$ generates $U$. □

**7.6** We now introduce a new basis of $A$ by performing an inverse discrete Fourier transform: For $i \in \mathbb{Z}_p$ and $s \in G$, we define the elements
$$e_i(s) := \frac{1}{p} \sum_{j=0}^{p-1} \zeta^{-ij} u^j g_s$$

where $\zeta$ is the primitive $p$-th root of unity fixed at the beginning of the section. Then $e_0(s), \ldots, e_{p-1}(s)$ is a basis of $A_s$. These basis elements satisfy:

$$ue_i(s) = \zeta^i e_i(s)$$

We therefore see that these elements explicitly depend on the chosen representative, but only up to a root of unity: To replace $g_s$ by another representative $u^k g_s$ means to replace $e_i(s)$ by $\zeta^{ik} e_i(s)$.

The cocycle $\sigma : G \times G \to U$ exhibited in Proposition 7.5 now can be written in the form

$$\sigma(s,t) = \sum_{i=0}^{p-1} \sigma_i(s,t) e_i(1)$$

where $\sigma_i(s,t)$ is an element of the base field. Since the cocycle is invertible, $\sigma_i(s,t)$ is nonzero.

We now try to describe the structure elements of $A$ with respect to the basis $e_i(s)$. In doing so, we will find again the compatibility condition first stated in Paragraph 3.2. The possible solutions of this compatibility condition will heavily depend on the parity of $p$, as we will see in the next paragraphs.

**Proposition**
1. Action and coaction are, with respect to this basis, determined by:

$$\phi(e_i(s)) = \zeta^{i\alpha(s)} e_i(s) \qquad \psi(e_i(s)) = \zeta^{i\beta(s)} e_i(s)$$

We therefore have: $\delta_A(e_i(s)) = c_{i\beta(s)} \otimes e_i(s)$

2. The coalgebra structure has, with respect to this basis, the following form:

$$\Delta_A(e_i(s)) = \sum_{j=0}^{p-1} e_j(s) \otimes e_{i-j}(s)$$

and $\epsilon_A(e_i(s)) = \delta_{i0}$.

3. For all $i, j \in \mathbb{Z}_p$ and $s, t \in G$, we have:

$$e_i(s) e_j(t) = \delta_{i\nu(s),j} \sigma_i(s,t) e_i(st)$$

4. The functions $\sigma_i$ satisfy the compatibility condition:

$$\sigma_{i+j}(s,t) = \zeta^{ij\nu(s)\beta(s)\alpha(t)} \sigma_i(s,t) \sigma_j(s,t)$$

**Proof.** The form of the action and the coaction follows directly from the equations $ue_i(s) = \zeta^i e_i(s)$, $\phi(g_s) = u^{\alpha(s)} g_s$, and $\psi(g_s) = u^{\beta(s)} g_s$. We leave the verification of the form of the coalgebra structure to the reader.

The form of the multiplication is a variant of the formula for the multiplication in a crossed product: First, observe that we have:

$$s.e_i(1) = \frac{1}{p}\sum_{j=0}^{p-1} \zeta^{-ij} u^{j\nu(s)} = \frac{1}{p}\sum_{j=0}^{p-1} \zeta^{-ij\nu(s^{-1})} u^j = e_{i\nu(s^{-1})}(1)$$

Therefore, since $e_i(s) = e_i(1)g_s$, we have:

$$e_i(s)e_j(t) = e_i(1)(s.e_j(1))\sigma(s,t)g_{st} = \delta_{i\nu(s),j}\sigma_i(s,t)e_i(st)$$

It remains to establish the asserted compatibility condition for the functions $\sigma_i$. For this, we look at the implications of the condition:

$$\Delta_A(e_i(s)e_{i\nu(s)}(t)) = \Delta_A(e_i(s))\Delta_A(e_{i\nu(s)}(t))$$

We have:

$$\Delta_A(e_i(s))\Delta_A(e_{i\nu(s)}(t)) = \sum_{m,n=0}^{p-1} (e_m(s) \otimes e_{i-m}(s))(e_{n\nu(s)}(t) \otimes e_{(i-n)\nu(s)}(t))$$

$$= \sum_{m,n=0}^{p-1} e_m(s)(c_{(i-m)\beta(s)} \to e_{n\nu(s)}(t)) \otimes e_{i-m}(s)e_{(i-n)\nu(s)}(t)$$

$$= \sum_{m,n=0}^{p-1} \zeta^{n\nu(s)(i-m)\beta(s)\alpha(t)} e_m(s)e_{n\nu(s)}(t) \otimes e_{i-m}(s)e_{(i-n)\nu(s)}(t)$$

$$= \sum_{m,n=0}^{p-1} \delta_{m\nu(s),n\nu(s)} \zeta^{n(i-m)\nu(s)\beta(s)\alpha(t)} \sigma_m(s,t)\sigma_{i-m}(s,t)e_m(st) \otimes e_{i-m}(st)$$

$$= \sum_{m=0}^{p-1} \zeta^{m(i-m)\nu(s)\beta(s)\alpha(t)} \sigma_m(s,t)\sigma_{i-m}(s,t)e_m(st) \otimes e_{i-m}(st)$$

On the other hand, we have:

$$\Delta_A(e_i(s)e_{i\nu(s)}(t)) = \Delta_A(\sigma_i(s,t)e_i(st)) = \sum_{m=0}^{p-1} \sigma_i(s,t)e_m(st) \otimes e_{i-m}(st)$$

By comparing coefficients, we get:

$$\sigma_i(s,t) = \zeta^{m(i-m)\nu(s)\beta(s)\alpha(t)} \sigma_m(s,t)\sigma_{i-m}(s,t)$$

By replacing $i$ by $i+m$, we arrive at the assertion. $\square$

**7.7** We shall now, in the case where $p$ is odd, prove the main result, namely that $A$ is isomorphic to a Yetter-Drinfel'd Hopf algebra of the form considered in Paragraph 3.4. The remaining task is to determine the solutions of the above compatibility condition. Note that, in the odd case, 2 is an invertible element of $\mathbb{Z}_p$, and therefore the expression $i/2$ for $i \in \mathbb{Z}_p$ makes sense. As in Paragraph 3.4, we work with respect to the characters:

$$\chi : \mathbb{Z}_p \to K, i \mapsto \zeta^{i/2} \qquad \eta : \mathbb{Z}_p \to K, i \mapsto \zeta^i$$

**Theorem** Suppose that $p$ is odd. Then there exists a normalized 2-cocycle $q \in Z^2(G, {}_G\mathbb{Z}_p)$ of the $G$-module ${}_G\mathbb{Z}_p$ such that

$$f_A : A_G(\alpha, \beta, q) \to A, e_i \otimes x_s \mapsto e_i(s)$$

is an isomorphism of Yetter-Drinfel'd Hopf algebras.

**Proof.** We use the notation of the previous proposition. Fix $s$ and $t$ and define $\tilde{\sigma}_i := \zeta^{-i^2 \nu(s)\beta(s)\alpha(t)/2} \sigma_i(s, t)$. We then have from the preceding proposition:

$$\tilde{\sigma}_{i+j} = \zeta^{-(i+j)^2 \nu(s)\beta(s)\alpha(t)/2} \sigma_{i+j}(s, t)$$
$$= \zeta^{-i^2 \nu(s)\beta(s)\alpha(t)/2} \zeta^{-ij\nu(s)\beta(s)\alpha(t)} \zeta^{-j^2 \nu(s)\beta(s)\alpha(t)/2} \zeta^{ij\nu(s)\beta(s)\alpha(t)} \sigma_i(s,t)\sigma_j(s,t)$$
$$= \tilde{\sigma}_i \tilde{\sigma}_j$$

This means that $\tilde{\sigma}$ defines a character to the base field. Therefore, there exists an element $q(s, t) \in \mathbb{Z}_p$, depending on $s$ and $t$, such that $\tilde{\sigma}_i = \zeta^{iq(s,t)}$, i. e., we have $\sigma_i(s, t) = \zeta^{iq(s,t)} \zeta^{i^2 \nu(s)\beta(s)\alpha(t)/2}$.

We now prove that the so-defined $q$ is a normalized 2-cocycle. The cocycle condition for $\sigma$, i. e., the equality $r.\sigma(s,t)\sigma(r, st) = \sigma(rs, t)\sigma(r, s)$, implies for the components that:

$$\sigma_{i\nu(r)}(s,t)\sigma_i(r, st) = \sigma_i(rs, t)\sigma_i(r, s)$$

As explained in Paragraph 3.3, we can use the isomorphism $\mathbb{Z}_p \otimes_{\mathbb{Z}} \mathbb{Z}_p \cong \mathbb{Z}_p$ to regard the cup product $\beta \cup \alpha \in Z^2(G, \mathbb{Z}_p \otimes_{\mathbb{Z}} \mathbb{Z}_p)$ as an element of $Z^2(G, \mathbb{Z}_p)$, if the $G$-module structure on $\mathbb{Z}_p$ is chosen correctly. We then have:

$$\sigma_i(s, t) = \zeta^{iq(s,t)} \zeta^{i^2 (\beta \cup \alpha)(s,t)/2}$$

The above condition then reads:

$$\zeta^{i\nu(r)q(s,t)} \zeta^{i^2 \nu(r)^2 (\beta \cup \alpha)(s,t)/2} \zeta^{iq(r,st)} \zeta^{i^2 (\beta \cup \alpha)(r,st)/2}$$
$$= \zeta^{iq(rs,t)} \zeta^{i^2 (\beta \cup \alpha)(rs,t)/2} \zeta^{iq(r,s)} \zeta^{i^2 (\beta \cup \alpha)(r,s)/2}$$

Since we already know that $\beta \cup \alpha$ is a 2-cocycle, this implies:

$$\nu(r)q(s, t) + q(r, st) = q(r, s) + q(rs, t)$$

and therefore $q$ is a cocycle. Since $\sigma$ is normalized, we have $\sigma_i(1,1) = 1$ for all $i \in \mathbb{Z}_p$. From Lemma 1.13, we have $\beta(1) = 0$. Therefore, we must have $q(1,1) = 0$, which means that $q$ is normalized.

It is obvious from the description of the structure elements of $A_G(\alpha, \beta, q)$ in Paragraph 3.4 and the description of the structure elements of $A$ in Paragraph 7.6 that $f_A$ is an $H$-linear and colinear coalgebra homomorphism. The cocycle $q$ is constructed in such a way that $f_A$ is also algebra homomorphism, and therefore a morphism of Yetter-Drinfel'd bialgebras. As in Paragraph 4.2, it is therefore a morphism of Yetter-Drinfel'd Hopf algebras. $\square$

**7.8** We shall now prove a similar result in the case $p = 2$, namely, we shall prove in this case that $A$ is isomorphic to a Yetter-Drinfel'd Hopf algebra of the form constructed in Paragraph 3.5.

Suppose that $\iota$ is a primitive fourth root of unity. First, it should be observed that, since the group of units of $\mathbb{Z}_2$ is trivial, we have $\nu(s) = 1$ for all $s \in G$, and therefore $\mathbb{Z}_2$ is a trivial $G$-module. From Paragraph 7.6, we have:
$$\sigma_{i+j}(s,t) = (-1)^{ij(\beta \cup \alpha)(s,t)} \sigma_i(s,t) \sigma_j(s,t)$$
Here we have, as in Paragraph 3.5, used the isomorphism $\mathbb{Z}_2 \otimes_\mathbb{Z} \mathbb{Z}_2 \cong \mathbb{Z}_2$ to regard the cup product $\beta \cup \alpha \in Z^2(G, \mathbb{Z}_2 \otimes_\mathbb{Z} \mathbb{Z}_2)$ as an element of $Z^2(G, \mathbb{Z}_2)$, i. e., we have:
$$\beta \cup \alpha : G \times G \to \mathbb{Z}_2, (s,t) \mapsto \beta(s)\alpha(t)$$
For $i = j = 0$, this implies that $\sigma_0(s,t) = 1$, since $\sigma_i(s,t) \neq 0$. For $i = j = 1$, we then have:
$$\sigma_1(s,t)^2 = (-1)^{(\beta \cup \alpha)(s,t)}$$
This implies $\sigma_1(s,t) = \iota^{q(s,t)}$, where $q(s,t) \in \mathbb{Z}_4$ satisfies:
$$q(s,t) \in \begin{cases} \{0,2\} & \text{if } (\beta \cup \alpha)(s,t) = 0 \\ \{1,3\} & \text{if } (\beta \cup \alpha)(s,t) = 1 \end{cases}$$
This condition means precisely that $\hat{\pi} \circ q = \beta \cup \alpha$, where $\hat{\pi} : \mathbb{Z}_4 \to \mathbb{Z}_2$ is the unique surjective group homomorphism. Since $\sigma$ is a normalized 2-cocycle, $q$ is a normalized 2-cocycle, too.

As described in Paragraph 3.5, the structure elements $\alpha$, $\beta$, and $q$ now can be used to define a new Yetter-Drinfel'd Hopf algebra, which we denote by $A'$.

**Theorem** The mapping
$$f_A : A' \to A, e_i \otimes x_s \mapsto e_i(s)$$
is an isomorphism of Yetter-Drinfel'd Hopf algebras.

**Proof.** From the description of the structure elements in Paragraph 3.5, it is clear that $f_A$ is an isomorphism of Yetter-Drinfel'd bialgebras. As in the proof of Proposition 4.2, it is therefore an isomorphism of Yetter-Drinfel'd Hopf algebras. $\square$

**7.9** We have already seen in Corollary 6.7 that $p$ divides $\dim A$. If $A$ is also commutative, this result can be sharpened:

**Theorem** If $A$ is commutative, $p^2$ divides $\dim A$.

**Proof.** If $A$ is commutative, we have $u^{\nu(s)}g_s = g_s u = u g_s$ for all $s \in G$, and therefore $\nu(s) = 1$. This means that the $G$-module structure on $\mathbb{Z}_p$ is trivial, and therefore $\alpha$ and $\beta$ are ordinary group homomorphisms from $G$ to $\mathbb{Z}_p$. Since $A$ is nontrivial by hypothesis, $\alpha$ and $\beta$ are nonzero by Proposition 1.11, and are therefore surjective. This implies that the order of $G$ is divisible by $p$, and therefore the dimension of $A$ is divisible by $p^2$. $\square$

Since we have already constructed Yetter-Drinfel'd Hopf algebras of dimension $p^2$ with these properties in Paragraph 3.6 and Paragraph 4.6, we have the following corollary:

**Corollary** For a natural number $n$, the following assertions are equivalent:

1. The exists a nontrivial, commutative, cocommutative, cosemisimple Yetter-Drinfel'd Hopf algebra over $K[\mathbb{Z}_p]$ that has dimension $n$.

2. $p^2$ divides $n$.

**Proof.** By the above theorem, the first assertion implies the second. For the converse, choose a nontrivial, commutative, cocommutative, cosemisimple Yetter-Drinfel'd Hopf algebra $A$ over $K[\mathbb{Z}_p]$ of dimension $p^2$. Determine $m$ satisfying $n = p^2 m$. The group ring $K[\mathbb{Z}_m]$ of the cyclic group of order $m$ becomes a Yetter-Drinfel'd Hopf algebra if endowed with the trivial module and the trivial comodule structure. From Proposition 1.4, we know that $A \otimes K[\mathbb{Z}_m]$, endowed with the tensor product algebra and the tensor product coalgebra structure, is a Yetter-Drinfel'd Hopf algebra, which is nontrivial by Proposition 1.11, and obviously is commutative, cocommutative, cosemisimple, and $n$-dimensional. $\square$

**7.10** As another application of the structure theorem, we classify the nontrivial Yetter-Drinfel'd Hopf algebras of dimension $p^2$ that are cocommutative and cosemisimple. Recall that we are still assuming that $p$ is a prime and that $K$ is an algebraically closed field whose characteristic is different from $p$.

**Theorem** Suppose that $A$ is a nontrivial Yetter-Drinfel'd Hopf algebra over $H = K[\mathbb{Z}_p]$ of dimension $p^2$ that is cocommutative and cosemisimple. Then $A$ is commutative.

1. If $p$ is odd, there are nonzero group homomorphisms $\alpha, \beta \in \text{Hom}(\mathbb{Z}_p, \mathbb{Z}_p)$ and a cocycle $q \in Z^2(\mathbb{Z}_p, \mathbb{Z}_p)$ such that $A \cong A_p(\alpha, \beta, q)$.

2. If $p = 2$, $A$ is isomorphic to $A_+$ or $A_-$.

3. In any case, there are $p(p-1)$ isomorphism classes of nontrivial, cocommutative, cosemisimple Yetter-Drinfel'd Hopf algebras of dimension $p^2$ over $K[\mathbb{Z}_p]$.

**Proof.** We first consider the case where $p$ is odd. Using the notation of the previous paragraphs, we must have $\text{card}(G) = p$, and therefore we have that $G \cong \mathbb{Z}_p$. As already explained in Paragraph 4.6, $\nu$ must, since it is a group homomorphism from $G$ to $\mathbb{Z}_p^\times$, be constantly equal to one, and therefore the $G$-module structure on $\mathbb{Z}_p$ is trivial. By Theorem 7.7, $A$ is isomorphic to some algebra of the type $A_p(\alpha, \beta, q)$, for group homomorphisms $\alpha, \beta \in \text{Hom}(\mathbb{Z}_p, \mathbb{Z}_p)$ and a normalized 2-cocycle $q \in Z^2(\mathbb{Z}_p, \mathbb{Z}_p)$. From Proposition 1.13.2, we see that $A$ is commutative. As explained in Paragraph 4.6, there are $p(p-1)$ isomorphism classes of these Yetter-Drinfel'd Hopf algebras.

In the case $p = 2$, we have $\text{card}(G) = 2$, and therefore $G \cong \mathbb{Z}_2$. From Theorem 7.8 and the discussion in Paragraph 3.6, we see that $A$ is isomorphic to $A_+$ or $A_-$. Therefore, $A$ is commutative. Since, by Proposition 4.9, $A_+$ and $A_-$ are not isomorphic, the assertions follow. $\square$

# 8 Semisimple Hopf algebras of dimension $p^3$

**8.1** In this section, we assume that $K$ is an algebraically closed field of characteristic zero, and that $p$ is a prime number. $B$ denotes a semisimple Hopf algebra of dimension $p^3$ over $K$ that is neither commutative nor cocommutative. We want to prove that $B$ is a Radford biproduct of a group ring of a group of order $p$ and a Yetter-Drinfel'd Hopf algebra of the form described in Section 3. We then apply this result to give a new proof of the theorem of A. Masuoka that, if $p$ is odd, there are $p+1$ isomorphism classes of semisimple Hopf algebras of dimension $p^3$ that are neither commutative nor commutative (cf. [50]). In Masuoka's approach, these Hopf algebras are not constructed via the Radford biproduct construction, but rather via the construction of Hopf algebra extensions described by W. M. Singer, M. Feth, and I. Hofstetter (cf. [77], [19], [23], [24]). As explained towards the end of the section, the case $p=2$ can also be treated from the point of view of the Radford biproduct construction in a rather analogous fashion. This gives a new proof of the classification result in dimension 8 obtained earlier by R. Williams and A. Masuoka (cf. [87], [47]).

As in the previous sections, we shall not, in notation, distinguish between an integer $i \in \mathbb{Z}$ and its equivalence class in a quotient group $\mathbb{Z}_n$.

**8.2** First, we determine the group $G(B)$ of grouplike elements of $B$.

**Proposition** The group $G(B)$ of grouplike elements of $B$ is isomorphic to $\mathbb{Z}_p \times \mathbb{Z}_p$.

**Proof.** (1) By a result of G. I. Kac and A. Masuoka ([33], Cor. 2, p. 159, [49], Thm. 1, p. 736), $B$ contains a nontrivial central grouplike element $g$. By the Nichols-Zoeller theorem, the order of $g$ is a power of $p$; we can therefore assume by Cauchy's theorem (cf. [3], Chap. 2, Exerc. 3, p. 20, [25], Kap. I, Satz 7.4, p. 34, [34], Satz 3.9, p. 41, [82], Chap. 2, § 2, p. 97) that its order is exactly $p$. Denote by $R$ the group ring spanned by the powers of $g$; $R$ is then a normal Hopf subalgebra of dimension $p$. By the normal basis theorem (cf. [72], Thm. 2.4, p. 300, [57], Cor. 8.4.7, p. 142), the corresponding Hopf algebra quotient $B/BR^+$ has dimension $p^2$. By a theorem of A. Masuoka (cf. [49], Thm. 2, p. 736), Hopf algebras of dimension $p^2$ are group rings.

(2) Recall that groups of order $p^2$ are isomorphic to $\mathbb{Z}_p \times \mathbb{Z}_p$ or $\mathbb{Z}_{p^2}$ (cf. [3], Chap. 2, Exerc. 4, p. 20, [25], Kap. I, Satz 6.10, p. 31, [34], Satz 4.3, p. 57, [82], Chap. 1, § 3, p. 27). Assume that the quotient $B/BR^+$ is isomorphic to the group ring of $\mathbb{Z}_{p^2}$. If $\pi : B \to B/BR^+$ denotes the corresponding projection, the set of coinvariant elements with respect to this projection is precisely $R$:

$$R := \{b \in B \mid (\mathrm{id}_B \otimes \pi)\Delta_B(b) = b \otimes 1\}$$

(cf. [57], Prop. 3.4.3, p. 34). By a different version of the normal basis theorem (cf. [72], Thm. 2.2, p. 299, [57], Thm. 8.4.6, p. 141), $B$ is, as an algebra, isomorphic to a crossed product of $R$ and $K[\mathbb{Z}_{p^2}]$, and in this crossed product the multiplication has the form:

$$(r \otimes x_i)(s \otimes x_j) = r(x_i.s)\sigma(i,j) \otimes x_{i+j}$$

for $r, s \in R$ and $i, j \in \mathbb{Z}_{p^2}$, where, for $i \in \mathbb{Z}_{p^2}$, $x_i$ denotes the corresponding canonical basis vector of $K[\mathbb{Z}_{p^2}]$, and $\sigma \in Z^2(\mathbb{Z}_{p^2}, U(R))$ is a normalized 2-cocycle. With respect to this isomorphism, the inclusion $R \subset B$ corresponds to the embedding $R \to R \otimes K[\mathbb{Z}_{p^2}], r \mapsto r \otimes 1$. Therefore, the fact that $R$ is central in $B$ implies that the action above is trivial. As in Proposition 1.13.2, we can then show that $\sigma(i,j) = \sigma(j,i)$. This implies that $B$ is commutative. Since this is not the case by assumption, we have ruled out the case that $B/BR^+$ is isomorphic to the group ring of $\mathbb{Z}_{p^2}$, and therefore it must be isomorphic to the group ring of $\mathbb{Z}_p \times \mathbb{Z}_p$.

(3) We have just established the existence of a surjective Hopf algebra morphism from $B$ to the group ring of $\mathbb{Z}_p \times \mathbb{Z}_p$. Applying this to $B^*$ and using the fact that the group ring of $\mathbb{Z}_p \times \mathbb{Z}_p$ is self-dual (cf. Paragraph 1.10), we get a Hopf algebra injection from this group ring to $B$. Therefore, $G(B)$ contains a subgroup that is isomorphic to $\mathbb{Z}_p \times \mathbb{Z}_p$. By the Nichols-Zoeller theorem, the order of $G(B)$ divides $p^3$, and must be strictly smaller since $B$ is not cocommutative. Therefore, $G(B)$ is isomorphic to $\mathbb{Z}_p \times \mathbb{Z}_p$. $\square$

The argument used in the proof is taken from work of A. Masuoka, where it appears in several places (cf. the proofs in [47], Prop. 2.3, p. 368, [48], Lem. 2, p. 1933, [49], Thm. 2, p. 737). A similar argument will be used in Paragraph 9.2.

**8.3** In order to invoke the methods developed so far, we have to show that $B$ is a Radford biproduct:

**Proposition** $B$ is isomorphic to a Radford biproduct:

$$B \cong A \otimes H$$

Here $H := K[\mathbb{Z}_p]$ is the group ring of the cyclic group of order $p$, and $A$ is a Yetter-Drinfel'd Hopf algebra over $H$.

**Proof.** As we have seen in the preceding paragraph, the dualization of the last proposition says that there is a surjective Hopf algebra homomorphism

$$\pi : B \to K[\mathbb{Z}_p \times \mathbb{Z}_p]$$

Restricting $\pi$ to $G(B)$, we get a group homomorphism

$$f : G(B) \to \mathbb{Z}_p \times \mathbb{Z}_p$$

It is impossible that $f$ is trivial, i. e., identically equal to the unit element, because in this case $G(B)$ would be contained in the set of coinvariant elements $\{b \in B \mid (\mathrm{id} \otimes \pi)\Delta_B(b) = b \otimes 1\}$, which has dimension $p$ by the normal basis theorem.

The group $G(B) \cong \mathbb{Z}_p \times \mathbb{Z}_p$ is elementary abelian, and therefore may be considered as a vector space over the field with $p$ elements. Considered in this way, $f$ is a linear mapping of vector spaces. We have just seen that the kernel of $f$ has dimension at most one, and therefore there is a one-dimensional subspace $U$ that $f$ maps isomorphically to a one-dimensional subspace of $\mathbb{Z}_p \times \mathbb{Z}_p$. By projecting to $f(U)$, we can construct group homomorphisms from $\mathbb{Z}_p$ to $G(B)$ and from $\mathbb{Z}_p \times \mathbb{Z}_p$ to $\mathbb{Z}_p$ such that the composite

$$\mathbb{Z}_p \to G(B) \xrightarrow{f} \mathbb{Z}_p \times \mathbb{Z}_p \to \mathbb{Z}_p$$

is the identity. Extending these group homomorphisms linearly, we can construct Hopf algebra homomorphisms from $K[\mathbb{Z}_p]$ to $B$ and back such that the composite

$$K[\mathbb{Z}_p] \to B \to K[\mathbb{Z}_p]$$

is the identity. By the Radford projection theorem (cf. [65], Thm. 3, p. 336), this implies the assertion. □

**8.4** From the preceding proposition, we know that we can assume that $B$ is a Radford biproduct $A \otimes H$, and we shall do so in the following. It then follows from many results in the literature that $A$ is semisimple (cf. [45], Thm. 3.1, p. 1368, [65], Prop. 3, p. 333, [20], Cor. 5.8, p. 4885, [81], Prop. 2.14, p. 22). We now prove:

**Proposition** $A$ is commutative.

**Proof.** Suppose that $W$ is a simple $A$-module, and denote the corresponding centrally primitive idempotent by $e$. We have to prove that $\dim W = 1$. By Corollary 2.3, there is a simple $B$-module $V$ such that $e \in \kappa(V)$ (cf. Definition 2.2). By a result of S. Montgomery and S. J. Witherspoon (cf. [58], Cor. 3.6, p. 325), $V$ has dimension $1$, $p$, $p^2$, or $p^3$. Since $B$ cannot contain two-sided ideals of dimension $p^4$ resp. $p^6$, the last two cases are impossible, and $V$ has dimension $1$ or $p$. Since $e \in \kappa(V)$, $W$ is a submodule of the restriction of $V$ to $A$, and therefore $\dim V = 1$ implies that $\dim W = 1$. We therefore may assume that the dimension of $V$ is $p$.

Suppose first that $V$ is stable (cf. Definition 2.3). We then know from Paragraph 2.5 that there are $p$ nonisomorphic simple $B$-modules whose restriction to $A$ is isomorphic to $W$. This implies that $B$ contains $p$ simple two-sided ideals of dimension $p^2$, whose direct sum has dimension $p^3$ and is therefore equal to $B$. This would imply that all simple $B$-modules have dimension $p$, which is obviously not the case for the trivial module. Therefore, this case is impossible, and $V$ must be purely unstable. Therefore, we know from Proposition 2.4 that $\dim V = p \dim W$, which means that $\dim W = 1$. Since a semisimple algebra whose simple modules are all one-dimensional is commutative, we have established the assertion. $\square$

From Paragraph 1.6, we have that $B^* \cong A^* \otimes H^*$ is again a Radford biproduct. Since $B$ is also cosemisimple (cf. [37], Thm. 3.3, p. 276) and $H = K[\mathbb{Z}_p]$ is self-dual (cf. Paragraph 1.10), the above proposition therefore also yields:

**Corollary** $A$ is cocommutative.

For similar reasons as above, $A$ is also cosemisimple (cf. [65], Prop. 4, p. 335, [81], Cor. 2.14, p. 23). If the action of $H$ on $A$ were trivial, it would follow directly from the formula for the multiplication of the Radford biproduct that $B$ were commutative. Similarly, if the coaction were trivial, it would follow from the comultiplication of the Radford biproduct that $B$ were cocommutative. Therefore, action and coaction are nontrivial, and therefore $A$ is nontrivial by Proposition 1.11.

**8.5** Now suppose that $p$ is odd and that $\zeta$ is a primitive $p$-th root of unity. Nontrivial, cocommutative, cosemisimple Yetter-Drinfel'd Hopf algebras over $H$ of dimension $p^2$ were classified in Theorem 7.10. We therefore know that $A$ appears in the list given in Paragraph 4.6. We now introduce a basis for these algebras and describe their structure elements with respect to this basis. The basis used here is slightly different from the basis used in Paragraph 5.2.

**Definition** Suppose that $a, b \in \mathbb{Z}_p^\times$ are nonzero elements of $\mathbb{Z}_p$ and that $q \in Z^2(\mathbb{Z}_p, \mathbb{Z}_p)$ is a normalized 2-cocycle. Denote by $\alpha$, resp. $\beta$, the group automorphism of $\mathbb{Z}_p$ given by multiplication by $a$, resp. $b$. In the Radford biproduct $B_p(a, b, q) := A_p(\alpha, \beta, q) \otimes H$, we introduce the basis

$$b_{ijk} := e_i \otimes c_j \otimes d_k$$

where $i, j, k \in \mathbb{Z}_p$. Here $e_i$ denotes, as in Paragraph 2.3, the $i$-th primitive idempotent in $K^{\mathbb{Z}_p}$, $c_j$, and not $x_j$, denotes the canonical basis vector of the group ring $K[\mathbb{Z}_p]$, as in Paragraph 1.10, and $d_k$ is the primitive idempotent

$$d_k := \frac{1}{p} \sum_{l=0}^{p-1} \zeta^{-kl} c_l$$

obtained by an inverse discrete Fourier transform from the canonical basis vectors. Under a Hopf algebra isomorphism between $K[\mathbb{Z}_p]$ and $K^{\mathbb{Z}_p}$, $d_0, \ldots, d_{p-1}$ would correspond to $e_0, \ldots, e_{p-1}$.

With respect to this basis, the structure maps of $B_p(a,b,q)$ take the form:

1. Multiplication: $b_{ijk} b_{lmn} = \delta_{k-n,alm} \delta_{il} \zeta^{iq(j,m)} \zeta^{abjmi^2/2} b_{i,j+m,n}$

2. Unit: $1 = \sum_{i,k=0}^{p-1} b_{i,0,k}$

3. Comultiplication: $\Delta(b_{ijk}) = \sum_{l,m=0}^{p-1} \zeta^{b(i-l)jm} b_{ljm} \otimes b_{i-l,j,k-m}$

4. Counit: $\epsilon(b_{ijk}) = \delta_{i0} \delta_{k0}$

5. Antipode: $S(b_{ijk}) := \zeta^{bijk} \zeta^{abi^2 j^2/2} \zeta^{iq(j,-j)} b_{-i,-j,-k-aij}$

**8.6** To complete the classification of semisimple Hopf algebras of dimension $p^3$, we have to describe which of the Hopf algebras $B_p(a,b,q)$ are isomorphic. For this, we must first understand the structure of $B_p(a,b,q)$ in greater detail. We maintain the assumptions that $p$ is odd and that $\zeta$ is a primitive $p$-th root of unity.

In the Yetter-Drinfel'd Hopf algebra $A_p(\alpha, \beta, q)$ used in the previous paragraph, we can consider the element

$$u := \sum_{i=0}^{p-1} \zeta^i e_i \otimes c_0$$

Since $\sum_{i=0}^{p-1} \zeta^i e_i$ is a grouplike element of $K^{\mathbb{Z}_p}$, and since the coalgebra structure of $A_p(\alpha, \beta, q)$ is the tensor product coalgebra structure, $u$ is a grouplike element. Since $\alpha(0) = \beta(0) = 0$, $u$ is invariant and coinvariant.

Analogously, we can define a linear form $\omega$ on $A_p(\alpha, \beta, q)$ by:

$$\omega(e_i \otimes c_j) := \delta_{i0} \zeta^j$$

From the description of the structure elements of $A_p(\alpha, \beta, q)$ given in Paragraph 3.4, it is easy to verify that $\omega$ is an $H$-linear and colinear character.

Using these elements, we can define the following elements of the Radford biproduct:

**Definition**
1. Define the elements $g_Z, g_N \in B_p(a,b,q)$ as: $g_Z := u \otimes c_0 \quad g_N := 1 \otimes c_1$

2. Define the linear form $\chi \in B_p(a,b,q)^*$ as: $\chi := \omega \otimes \epsilon_H$

3. Define the linear map $\Omega \in \text{End}(B_p(a,b,q))$ as:

$$\Omega : B_p(a,b,q) \to B_p(a,b,q), b \mapsto (\text{id} \otimes \chi) \Delta(b)$$

With respect to the basis defined in Paragraph 8.5, these elements can be expressed in the following way:

$$g_Z = \sum_{i,k=0}^{p-1} \zeta^i b_{i,0,k} \qquad g_N = \sum_{i,k=0}^{p-1} \zeta^k b_{i,0,k} \qquad \chi(b_{ijk}) = \delta_{i0}\delta_{k0}\zeta^j$$

Now we can characterize the basis elements $b_{ijk}$ in terms of the elements $g_Z$, $g_N$, and $\chi$:

**Proposition**
1. The element $g_Z$ is a central grouplike element.
2. The element $g_N$ is a grouplike element which is not central.
3. The element $\chi$ is a central grouplike element of $B_p(a,b,q)^*$.
4. The basis element $b_{ijk}$ satisfies

$$g_Z b_{ijk} = \zeta^i b_{ijk} \qquad b_{ijk} g_N = \zeta^k b_{ijk} \qquad \Omega(b_{ijk}) = \zeta^j b_{ijk}$$

Every other element satisfying these equations is a scalar multiple of $b_{ijk}$.

**Proof.** $g_Z$ is grouplike since $u$ is coinvariant, and it is central since $u$ is invariant. Similarly, $\chi$ is grouplike, i. e., a character, since $\omega$ is $H$-linear, and it is central since $\omega$ is colinear. $g_N$ is not central since conjugation by $g_N$ induces the action of $c_1$ on $A_p(\alpha, \beta, q)$, which is nontrivial. This proves the first three assertions.

It is easy to verify that $b_{ijk}$ in fact satisfies the equations stated in the fourth assertion. In particular, we see that the operators 'Left multiplication by $g_Z$', 'Right multiplication by $g_N$', and $\Omega$ commute. This also follows from the fact that $\Omega$ is an algebra homomorphism that, since $\chi(g_Z) = \omega(u) = 1$ and also $\chi(g_N) = \omega(1) = 1$, satisfies $\Omega(g_Z) = g_Z \chi(g_Z) = g_Z$ and $\Omega(g_N) = g_N \chi(g_N) = g_N$. Since the basis $b_{ijk}$ is a basis of simultaneous eigenvectors for these operators, the simultaneous eigenspaces cannot have a dimension greater than one. $\square$

Since we know from Proposition 8.2 that $B_p(a,b,q)$ contains $p^2$ grouplike elements, we see that $g_Z$ and $g_N$ generate $G(B_p(a,b,q))$. In addition, since not all grouplike elements are central, the subgroup of central grouplike elements must have order $p$, and is therefore generated by $g_Z$.

**8.7** We have seen in Theorem 7.10 that, up to isomorphism, there are $p(p-1)$ isomorphism classes of nontrivial, cocommutative, cosemisimple Yetter-Drinfel'd Hopf algebras of dimension $p^2$ over $K[\mathbb{Z}_p]$. However, two nonisomorphic Yetter-Drinfel'd Hopf algebras may have isomorphic Radford biproducts, and this turns out to be the case in our situation. The following lemma gives a first description of the form of such an isomorphism:

**Lemma** Suppose that $a, b, a', b' \in \mathbb{Z}_p^\times$ and that $q, q' \in Z^2(\mathbb{Z}_p, \mathbb{Z}_p)$ are normalized 2-cocycles. Suppose that

$$f : B_p(a, b, q) \to B_p(a', b', q')$$

is a Hopf algebra isomorphism. Then there are $r, s, t \in \mathbb{Z}_p^\times$, $u \in \mathbb{Z}_p$ and a map $\tau : \mathbb{Z}_p \times \mathbb{Z}_p \times \mathbb{Z}_p \to K^\times$ such that

$$f(b_{ijk}) = \tau(i, j, k) b'_{ri, sj, tk+ui}$$

for all $i, j, k \in \mathbb{Z}_p$. Here the basis elements in $B_p(a', b', q')$ are denoted by $b'_{ijk}$.

**Proof.** Denote the corresponding objects of $B_p(a', b', q')$ by $g'_Z, g'_N, \chi'$, and $\Omega'$. $g'_Z$ must be the image of a central grouplike element in $B_p(a, b, q)$. Therefore, there exists $r \in \mathbb{Z}_p^\times$ such that:

$$f(g_Z^r) = g'_Z$$

Similarly, $g'_N$ must be the image of a grouplike element in $B_p(a, b, q)$ that is not central. Therefore, there exist $t \in \mathbb{Z}_p^\times$ and $u \in \mathbb{Z}_p$ such that:

$$f(g_Z^u g_N^t) = g'_N$$

Of course, this $u$ is different from the grouplike element used in Paragraph 8.6. Since the transpose $f^*$ of maps central grouplikes to central grouplikes, there exists $s \in \mathbb{Z}_p^\times$ such that:

$$f^*(\chi') = \chi^s$$

We therefore get that:

$$\Omega' \circ f = (\mathrm{id} \otimes \chi') \circ \Delta \circ f = (\mathrm{id} \otimes \chi') \circ (f \otimes f) \circ \Delta = f \circ (\mathrm{id} \otimes \chi^s) \circ \Delta = f \circ \Omega^s$$

We now have the following equations:

$$\begin{aligned}
g'_Z f(b_{ijk}) &= f(g_Z^r b_{ijk}) = \zeta^{ri} f(b_{ijk}) \\
f(b_{ijk}) g'_N &= f(b_{ijk} g_Z^u g_N^t) = \zeta^{tk+ui} f(b_{ijk}) \\
\Omega'(f(b_{ijk})) &= f(\Omega^s(b_{ijk})) = \zeta^{sj} f(b_{ijk})
\end{aligned}$$

By Proposition 8.6, $f(b_{ijk})$ must be proportional to $b'_{ri, sj, tk+ui}$. □

**Proposition** Suppose that $a, b, a', b' \in \mathbb{Z}_p^\times$ and that $q, q' \in Z^2(\mathbb{Z}_p, \mathbb{Z}_p)$ are normalized 2-cocycles. Suppose that $B_p(a, b, q)$ and $B_p(a', b', q')$ are isomorphic. Then there are $r, t \in \mathbb{Z}_p^\times$ such that $a' = \frac{ta}{r}$, $b' = \frac{b}{rt}$, and $q$ and $rq'$ are cohomologous.

**Proof.** (1) Suppose that $f : B_p(a,b,q) \to B_p(a',b',q')$ is a Hopf algebra isomorphism. By the preceding lemma, there are $r,s,t \in \mathbb{Z}_p^\times$, $u \in \mathbb{Z}_p$, and a map $\tau : \mathbb{Z}_p \times \mathbb{Z}_p \times \mathbb{Z}_p \to K^\times$ such that $f(b_{ijk}) = \tau(i,j,k) b'_{ri,sj,tk+ui}$. We have:

$$\Delta(f(b_{ijk})) = \tau(i,j,k) \sum_{l,m=0}^{p-1} \zeta^{b'r(i-l)sjm} b'_{rl,sj,m} \otimes b'_{r(i-l),sj,tk+ui-m}$$

$$= \tau(i,j,k) \sum_{l,m=0}^{p-1} \zeta^{b'r(i-l)sj(tm+ul)} b'_{rl,sj,tm+ul} \otimes b'_{r(i-l),sj,t(k-m)+u(i-l)}$$

On the other hand, we have:

$$(f \otimes f)\Delta(b_{ijk}) =$$
$$\sum_{l,m=0}^{p-1} \zeta^{b(i-l)jm} \tau(l,j,m) \tau(i-l, j, k-m) b'_{rl,sj,tm+ul} \otimes b'_{r(i-l),sj,t(k-m)+u(i-l)}$$

By comparing coefficients, we get:

$$\tau(i,j,k) \zeta^{b'r(i-l)sj(tm+ul)} = \tau(i-l, j, k-m) \tau(l,j,m) \zeta^{b(i-l)jm}$$

or, by replacing $i$ with $i+l$ and $k$ with $k+m$:

$$\tau(i+l, j, k+m) \zeta^{b'risj(tm+ul)} = \tau(i,j,k) \tau(l,j,m) \zeta^{bijm}$$

(2) The last equation can be rewritten in the form

$$\frac{\tau(i,j,k)\tau(l,j,m)}{\tau(i+l, j, k+m)} = \zeta^{(b'rst-b)ijm} \zeta^{b'rsuijl}$$

The left hand side of this equation remains invariant if we exchange simultaneously $i$ and $l$ as well as $k$ and $m$. Therefore, we have $\zeta^{(b'rst-b)ijm} = \zeta^{(b'rst-b)ljk}$. This yields for $i = j = m = 1$ and $k = l = 0$ that $\zeta^{b'rst-b} = 1$, and therefore $b = b'rst$.

(3) We now argue as in the proof of Theorem 7.7. Fix $j \in \mathbb{Z}_p$. For $i, k \in \mathbb{Z}_p$, we define:
$$\tilde{\tau}(i,k) := \tau(i,j,k) \zeta^{b'rsuji^2/2}$$

Note that, since $p$ is odd, 2 is an invertible element of $\mathbb{Z}_p$, and therefore the expression $i/2$ in the above formula makes sense. We then have:

$$\tilde{\tau}(i+l, k+m) = \tau(i+l, j, k+m) \zeta^{b'rsuj(i+l)^2/2}$$
$$= \tau(i+l, j, k+m) \zeta^{b'rsuji^2/2} \zeta^{b'rsujil} \zeta^{b'rsujl^2/2}$$
$$= \tau(i,j,k) \tau(l,j,m) \zeta^{b'rsuji^2/2} \zeta^{b'rsujl^2/2}$$
$$= \tilde{\tau}(i,k) \tilde{\tau}(l,m)$$

This shows that $\tilde{\tau}: \mathbb{Z}_p \times \mathbb{Z}_p \to K^\times$ is a group homomorphism. Therefore, there exist, for every $j \in \mathbb{Z}_p$, elements $v(j), w(j) \in \mathbb{Z}_p$ such that $\tilde{\tau}(i,k) = \zeta^{v(j)i+w(j)k}$, i. e.:
$$\tau(i,j,k) = \zeta^{v(j)i+w(j)k}\zeta^{-b'rsuji^2/2}$$

(4) We have:
$$f(b_{ijk})f(b_{ilm}) = \zeta^{v(j)i+w(j)k+v(l)i+w(l)m}\zeta^{-b'rsu(j+l)i^2/2}b'_{ri,sj,tk+ui}b'_{ri,sl,tm+ui}$$
$$= \delta_{t(k-m),a'risl}\zeta^{i(v(j)+v(l))+w(j)k+w(l)m}\zeta^{-b'rsu(j+l)i^2/2}$$
$$\zeta^{riq'(sj,sl)}\zeta^{a'b'sjsl(ri)^2/2}b'_{ri,s(j+m),tm+ui}$$

On the other hand, we have:
$$f(b_{ijk}b_{ilm}) = \delta_{k-m,ail}\zeta^{iq(j,l)}\zeta^{abjli^2/2}f(b_{i,j+l,m})$$
$$= \delta_{k-m,ail}\zeta^{iq(j,l)}\zeta^{abjli^2/2}\zeta^{iv(j+l)+mw(j+l)}\zeta^{-b'rsu(j+l)i^2/2}b'_{ri,s(j+m),tm+ui}$$

Now suppose that $k - m = ail$. Then we have $f(b_{ijk}b_{ilm}) \neq 0$, and therefore also $f(b_{ijk})f(b_{ilm}) \neq 0$. This implies that $t(k-m) = a'risl$. In particular, in the case $k = a, m = 0, i = l = 1$, we have that $at = a'rs$.

(5) Inserting this result into the first calculation above, we get:
$$f(b_{ijk})f(b_{ilm})$$
$$= \delta_{k-m,ail}\zeta^{i(v(j)+v(l))+w(j)k+w(l)m}\zeta^{-b'rsu(j+l)i^2/2}\zeta^{riq'(sj,sl)}\zeta^{abjli^2/2}$$
$$b'_{ri,s(j+m),tm+ui}$$

In the case $k - m = ail$, a comparison with the expression for $f(b_{ijk}b_{ilm})$ yields:
$\zeta^{iq(j,l)}\zeta^{iv(j+l)+mw(j+l)} = \zeta^{i(v(j)+v(l))+kw(j)+mw(l)}\zeta^{irq'(sj,sl)}$ Therefore, if we set $k = m + ail$, we get:
$$iq(j,l) - irq'(sj,sl) = i(v(j)+v(l)-v(j+l))+ailw(j)+m(w(j)+w(l)-w(j+l))$$

for all $i, j, l, m \in \mathbb{Z}_p$. For $i = 0$ and $m = 1$, this yields $w(j) + w(l) - w(j+l) = 0$. Therefore, we have $w(j) = wj$ for some $w \in \mathbb{Z}_p$, for which we use the same notation as for the function $w$. On the other hand, we get for $i = 1$ and $m = 0$ that:
$$q(j,l) - rq'(sj,sl) = v(j) + v(l) - v(j+l) + awjl$$

(6) Now define $\tilde{v}(j) := v(j) - \frac{aw}{2}j^2$. Then we have:
$$\tilde{v}(j) + \tilde{v}(l) - \tilde{v}(j+l) = v(j) + v(l) - v(j+l) + \frac{aw}{2}(j+l)^2 - \frac{aw}{2}j^2 - \frac{aw}{2}l^2$$
$$= v(j) + v(l) - v(j+l) + awjl$$
$$= q(j,l) - rq'(sj,sl)$$

Therefore, $q$ is cohomologous to the cocycle $(j,l) \mapsto rq'(sj,sl)$. By Proposition 1.13.2, we now have that this cocycle, and therefore $q$, is cohomologous to $rsq'$. Now the assertion follows if we change notation and denote the product $rs$ by $r$. $\square$

**8.8** In the preceding paragraph, we have obtained a necessary condition for $B_p(a, b, q)$ and $B_p(a', b', q')$ to be isomorphic. However, this condition is also sufficient. For, if there are $r, t \in \mathbb{Z}_p^\times$ such that $q$ and $rq'$ are cohomologous, $a' = \frac{ta}{r}$, and $b' = \frac{b}{rt}$, it can be verified directly that

$$f : B_p(a, b, q) \to B_p(a', b', q'), b_{ijk} \to \zeta^{iv(j)} b'_{ri, j, tk}$$

is a Hopf algebra isomorphism, where $v : \mathbb{Z}_p \to \mathbb{Z}_p$ is a 1-cochain whose coboundary is $rq' - q$, i. e., we have:

$$rq'(j, l) - q(j, l) = v(j + l) - v(j) - v(l)$$

The remaining task is to determine the number of isomorphism classes among the $B_p(a, b, q)$. For this, the following lemma is helpful:

**Lemma** Suppose that $p$ is odd. Consider the action of $\mathbb{Z}_p^\times \times \mathbb{Z}_p^\times$ on $M := \mathbb{Z}_p^\times \times \mathbb{Z}_p^\times \times \mathbb{Z}_p$ defined by

$$(r, t).(a, b, q) := (\frac{ta}{r}, \frac{b}{rt}, \frac{q}{r})$$

Then $M$ can be decomposed into $p + 1$ orbits with respect to this action.

**Proof.** We denote the same set $\mathbb{Z}_p^\times \times \mathbb{Z}_p^\times \times \mathbb{Z}_p$ by $N$ if endowed with the action:

$$(r, t).(a, b, q) := (t^2 a, rtb, \frac{q}{r})$$

of $\mathbb{Z}_p^\times \times \mathbb{Z}_p^\times$. It is then easy to see that

$$f : M \to N, (a, b, q) \mapsto (\frac{a}{b}, \frac{1}{b}, q)$$

is an equivariant bijection. It therefore suffices to count the orbits of $N$. It is obvious that the isotropy group of the element $(a, 1, 1) \in N$ is trivial, and any two such element cannot belong to the same orbit. Therefore, we have found $p - 1$ distinct orbits of length $(p - 1)^2$.

Now fix an element $a_0 \in \mathbb{Z}_p^\times$ that is not a square. Then the elements $(1, 1, 0) \in N$ and $(a_0, 1, 0) \in N$ are not conjugate under the action of $\mathbb{Z}_p^\times \times \mathbb{Z}_p^\times$. If $(r, t)$ is contained in the isotropy group of $(a_0, 1, 0)$, we have $t^2 a_0 = a_0$ and $rt = 1$, i. e., $t = \pm 1$ and $r = t$. Therefore, the isotropy group of $(a_0, 1, 0)$ consists of two, the orbit of $(a_0, 1, 0)$ of $\frac{1}{2}(p - 1)^2$ elements. Since the same applies to $(1, 1, 0)$, we have found two additional orbits of length $\frac{1}{2}(p-1)^2$. Since the combined lengths of these orbits satisfy

$$(p-1)^3 + 2\frac{(p-1)^2}{2} = (p-1)^2 p = \text{card}(N)$$

these are already all orbits, which means that we have $p + 1$ orbits in total. □

We summarize the obtained results in the following theorem:

**Theorem** Suppose that $K$ is an algebraically closed field of characteristic zero and that $p$ is an odd prime.

1. Every semisimple, noncommutative, noncocommutative Hopf algebra of dimension $p^3$ is isomorphic to a Hopf algebra $B_p(a, b, q)$ for some $a, b \in \mathbb{Z}_p^\times$ and a normalized 2-cocycle $q \in Z^2(\mathbb{Z}_p, \mathbb{Z}_p)$.

2. $B_p(a, b, q)$ and $B_p(a', b', q')$ are isomorphic if and only if there are $r, t \in \mathbb{Z}_p^\times$ such that $a' = \frac{ta}{r}$, $b' = \frac{b}{rt}$, and $q$ and $rq'$ are cohomologous.

3. There are $p+1$ isomorphism classes of semisimple, noncommutative, noncocommutative Hopf algebras of dimension $p^3$.

**Proof.** The first result was obtained in Paragraph 8.5. One part of the second assertion was proved in Proposition 8.7, the other part was proved at the beginning of this paragraph. Since we know from Proposition 1.13.2 that $H^2(\mathbb{Z}_p, \mathbb{Z}_p) \cong \mathbb{Z}_p$, and since every cocycle is cohomologous to a normalized cocycle, the third assertion follows from the preceding lemma. □

Those semisimple Hopf algebras of dimension $p^3$ that are commutative or cocommutative are described below in Paragraph 8.10.

**8.9** We now consider the case $p = 2$. Fix a primitive fourth root of unity $\iota$. In Paragraph 7.10 and Paragraph 4.9, we have described completely the four-dimensional Yetter-Drinfel'd Hopf algebras over $K[\mathbb{Z}_2]$ that are nontrivial, cocommutative, and cosemisimple: There are two isomorphism types that are represented by the two nonisomorphic Yetter-Drinfel'd Hopf algebras that were denoted $A_+$ and $A_-$. We therefore can say that any eight-dimensional, noncommutative, noncocommutative, semisimple Hopf algebra $B$ over the algebraically closed field $K$ of characteristic zero is isomorphic either to the Radford biproduct $B_+ := A_+ \otimes H$ or to the Radford biproduct $B_- := A_- \otimes H$. We now introduce a basis for these algebras and describe their structure elements with respect to this basis:

**Definition** In the Hopf algebra $B_+ := A_+ \otimes H$, we introduce the basis

$$b_{ijk}^+ := e_i \otimes c_j \otimes d_k$$

where $d_k$ is the idempotent

$$d_k := \frac{1}{2} \sum_{l=0}^{1} (-1)^{kl} c_l$$

constructed in analogy to Definition 8.5. Here the indices $i, j, k$ take the values 0 and 1. We introduce an analogous basis $b_{ijk}^-$ in $B_-$.

As in Paragraph 3.6, we define cocycles $q_+, q_- \in Z^2(\mathbb{Z}_2, \mathbb{Z}_4)$ by:

$$q_+(i,j) := \begin{cases} 0 & \text{if } i = 0 \text{ or } j = 0 \\ 1 & \text{if } i = 1 \text{ and } j = 1 \end{cases} \qquad q_-(i,j) := \begin{cases} 0 & \text{if } i = 0 \text{ or } j = 0 \\ 3 & \text{if } i = 1 \text{ and } j = 1 \end{cases}$$

Also as in Paragraph 3.6, we use the notation

$$\sigma_0^\pm(i,j) := 1 \qquad \sigma_1^\pm(i,j) := \iota^{q_\pm(i,j)}$$

With this notation, the structure maps of $B_\pm$ can be expressed with respect to this basis:

1. Multiplication: $b_{ijk}^\pm b_{lmn}^\pm = \delta_{k-n,lm}\delta_{il}\sigma_i^\pm(j,m)b_{i,j+m,n}^\pm$

2. Unit: $1 = \sum_{i,k=0}^{1} b_{i,0,k}^\pm$

3. Comultiplication: $\Delta(b_{ijk}^\pm) = \sum_{l,m=0}^{1}(-1)^{(i-l)jm} b_{ljm}^\pm \otimes b_{i-l,j,k-m}^\pm$

4. Counit: $\epsilon(b_{ijk}^\pm) = \delta_{i0}\delta_{k0}$

5. Antipode: $S(b_{ijk}) = (-1)^{ijk}\sigma_i^\pm(j,-j)^{-1} b_{-i,-j,-k-ij}^\pm$

The remaining open question in our treatment of eight-dimensional semisimple Hopf algebras is whether $B_+$ and $B_-$ are isomorphic. We have seen in Proposition 4.9 that $A_+$ and $A_-$ are not isomorphic. Nevertheless, $B_+$ and $B_-$ are isomorphic:

**Proposition** The linear map

$$f: B_+ \to B_-, \ b_{ijk}^+ \mapsto \sigma_1^+(i,j) b_{i,j,k+i}^-$$

is a Hopf algebra isomorphism.

**Proof.** It is easy to verify the equations:

$$\sigma_i^-(j,k) = (-1)^{ijk}\sigma_i^+(j,k) \qquad \sigma_1^+(i,j+k) = (-1)^{ijk}\sigma_1^+(i,j)\sigma_1^+(i,k)$$

Therefore, $f$ is an algebra homomorphism:

$$\begin{aligned} f(b_{ijk}^+ b_{lmn}^+) &= \delta_{k-n,lm}\delta_{il}\sigma_i^+(j,m) f(b_{i,j+m,n}^+) \\ &= \delta_{k-n,lm}\delta_{il}\sigma_i^+(j,m)\sigma_1^+(i,j+m) b_{i,j+m,n+i}^- \\ &= \delta_{k-n,lm}\delta_{il}(-1)^{ijm}(-1)^{ijm}\sigma_i^-(j,m)\sigma_1^+(i,j)\sigma_1^+(i,m) b_{i,j+m,n+i}^- \\ &= \delta_{k-n,lm}\delta_{il}\sigma_i^-(j,m)\sigma_1^+(i,j)\sigma_1^+(l,m) b_{i,j+m,n+i}^- \\ &= \sigma_1^+(i,j)\sigma_1^+(l,m) b_{i,j,k+i}^- b_{l,m,n+l}^- = f(b_{ijk}^+)f(b_{lmn}^+) \end{aligned}$$

We can also deduce from these equations that $f$ is a coalgebra homomorphism:

$$\begin{aligned}
(f \otimes f)\Delta(b_{ijk}^+) &= \sum_{l,m=0}^{1} (-1)^{(i-l)jm} \sigma_1^+(l,j)\sigma_1^+(i-l,j) b_{l,j,m+l}^- \otimes b_{i-l,j,k-m+i-l}^- \\
&= \sum_{l,m=0}^{1} (-1)^{(i-l)j(m-l)} \sigma_1^+(l,j)\sigma_1^+(i-l,j) b_{l,j,m}^- \otimes b_{i-l,j,k-m+i}^- \\
&= \sum_{l,m=0}^{1} (-1)^{(i-l)jm} \sigma_1^+(i,j) b_{ljm}^- \otimes b_{i-l,j,k-m+i}^- \\
&= \sigma_1^+(i,j)\Delta(b_{i,j,k+i}^-) = \Delta(f(b_{ijk}^+))
\end{aligned}$$

It is easy to see that $f$ preserves the unit and the counit, and therefore also commutes with the antipodes (cf. [84], Lem. 4.0.4, p. 81). □

We therefore see that, in dimension eight, there is only one isomorphism class of noncommutative, noncocommutative, semisimple Hopf algebras. This was first proved by R. Williams and A. Masuoka (cf. [87], [47], Thm. 2.13, p. 371), who, however, used comparatively different methods. The Hopf algebra itself was first constructed by G. I. Kac and V. G. Paljutkin (cf. [31], [32]).

**8.10** We have described above the semisimple Hopf algebras of dimension $p^3$ that are neither commutative nor cocommutative. We now describe the semisimple Hopf algebras of dimension $p^3$ that are commutative or cocommutative. We maintain our assumption that the base field $K$ is algebraically closed of characteristic zero. In this case, finite-dimensional Hopf algebras that are commutative or cocommutative are automatically semisimple; this follows, for example, from an even more general result of R. G. Larson and D. E. Radford (cf. [37], Cor. 2.6, p. 275). By a result of D. K. Harrison and P. Cartier (cf. [35], Thm. 3.2, p. 354, [7], p. 102, [57], Thm. 2.3.1, p. 22), a finite-dimensional cocommutative Hopf algebra is isomorphic to a group ring. The description of all commutative or cocommutative Hopf algebras of dimension $p^3$ therefore reduces to the description of all groups of order $p^3$. This is, of course, well known, and here we shall only describe the results (cf. [25], Kap. I, Satz 14.10, p. 93, [34], Kap. IV, § 1, p. 61/62, [83], Chap. 4, (4.13), p. 67); nevertheless, it is not a totally easy task. By the structure theorem for finite abelian groups (cf. [82], Chap. 2, Thm. 5.2, p. 145, [34], Kap. II, § 2, p. 30), the abelian groups of order $p^3$ are:

$$\mathbb{Z}_{p^3} \qquad \mathbb{Z}_{p^2} \times \mathbb{Z}_p \qquad \mathbb{Z}_p \times \mathbb{Z}_p \times \mathbb{Z}_p$$

The nonabelian groups of order 8 are the dihedral group $D_4$ and the quaternion group $Q$. Now suppose that $p \neq 2$ and that $G$ is a nonabelian group of order $p^3$. It is easy to see that the center of $G$ is of order $p$ and that $G$ contains a subgroup $N$ of order $p^2$, which is therefore normal. Now two cases may occur: Either $G$ contains an element of order $p^2$ or all nontrivial elements are of order $p$,

i. e., we have $N \cong \mathbb{Z}_{p^2}$ or $N \cong \mathbb{Z}_p \times \mathbb{Z}_p$. In any case, we have a short exact sequence:
$$N \rightarrowtail G \twoheadrightarrow \mathbb{Z}_p$$
The nontrivial task is to prove that this sequence is split, i. e., $G$ is a semidirect product in both cases. Since $\text{card}(\text{Aut}(\mathbb{Z}_{p^2})) = p(p-1)$, there is an essentially unique automorphism of order $p$ in the first case, which can be chosen to be multiplication by $p+1$. In the second case, $N$ is elementary abelian, and therefore group homomorphisms are linear maps over the field with $p$ elements. The automorphism in the semidirect product can be chosen in such a way that it is described by the matrix
$$\begin{pmatrix} 1 & 1 \\ 0 & 1 \end{pmatrix}$$
which is of order $p$. This yields:

**Proposition** Suppose that $p$ is an odd prime. Then the semidirect product
$$G_1 := \mathbb{Z}_{p^2} \rtimes \mathbb{Z}_p$$
with respect to the action of $\mathbb{Z}_p$ on $\mathbb{Z}_{p^2}$ given by multiplication by $p+1$, and the semidirect product
$$G_2 := (\mathbb{Z}_p \times \mathbb{Z}_p) \rtimes \mathbb{Z}_p$$
with respect to the action of $\mathbb{Z}_p$ on $\mathbb{Z}_p \times \mathbb{Z}_p$ given by the matrix
$$\begin{pmatrix} 1 & 1 \\ 0 & 1 \end{pmatrix}$$
are, up to isomorphism, the only nonabelian groups of order $p^3$.

We therefore have the following corollary:

**Corollary** Suppose that $K$ is an algebraically closed field of characteristic zero. Suppose that $p$ is a prime number and that $B$ is a semisimple Hopf algebra of dimension $p^3$.

1. Suppose that $B$ is commutative and cocommutative. Then $B$ is isomorphic to $K[\mathbb{Z}_{p^3}]$, $K[\mathbb{Z}_{p^2} \times \mathbb{Z}_p]$, or $K[\mathbb{Z}_p \times \mathbb{Z}_p \times \mathbb{Z}_p]$.

2. Suppose that $B$ is cocommutative, but not commutative. Then $B$ is isomorphic to $K[D_4]$ or $K[Q]$ if $p = 2$, and to $K[G_1]$ or $K[G_2]$ if $p \neq 2$.

3. Suppose that $B$ is commutative, but not cocommutative. Then $B$ is isomorphic to $K^{D_4}$ or $K^Q$ if $p = 2$, and to $K^{G_1}$ or $K^{G_2}$ if $p \neq 2$, where $K^G := \text{Map}(G, K)$ is the Hopf algebra of functions on the group $G$.

The reader should note that the notation $G_1$ and $G_2$ is interchanged in comparison to [50], Example 2.6, p. 796.

# 9 Semisimple Hopf algebras of dimension $pq$

**9.1** In this section, we assume that $K$ is an algebraically closed field of characteristic zero, and that $p$ and $q$ are two distinct prime numbers. $B$ denotes a semisimple Hopf algebra of dimension $pq$ over $K$ that is neither commutative nor cocommutative. Note that $B$ is then also cosemisimple (cf. [37], Thm. 3.3, p. 276). We assume that $B$ contains a nontrivial grouplike element $g \neq 1_B$ and that $B^*$ contains a nontrivial grouplike element $\gamma \neq \epsilon_B$. Our goal is to prove that these assumptions are contradictory. $H$ denotes the Hopf subalgebra of $B$ spanned by the powers of $g$. Observe that, by the Nichols-Zoeller theorem, the orders of $g$ and $\gamma$ must divide $\dim B$. Since $B$ is neither commutative nor cocommutative, these orders cannot be equal to $pq$, and therefore are equal to $p$ or $q$. By convention, we suppose that the order of $g$ is $p$.

**9.2** The first step in this investigation is similar to the first step of the previous section in Paragraph 8.2. It does not rely on the existence of $\gamma$.

**Proposition** $g$ is not central.

**Proof.** Assume on the contrary that $g$ is central. Then $H$ is a normal Hopf subalgebra of $B$. The corresponding Hopf algebra quotient $B/BH^+$ has dimension $q$ by the normal basis theorem (cf. [72], Thm. 2.4, p. 300, [57], Cor. 8.4.7, p. 142). By Zhu's theorem (cf. [90], Thm. 2, p. 57), this quotient is isomorphic to the group ring of the cyclic group of order $q$. By a different version of the normal basis theorem (cf. [72], Thm. 2.2, p. 299, [57], Thm. 8.4.6, p. 141), we can conclude that, as an algebra, $B$ is isomorphic to a crossed product of $H$ and the group ring $K[\mathbb{Z}_q]$, since the set of coinvariant elements with respect to $B/BH^+$ is precisely $H$ (cf. [57], Prop. 3.4.3, p. 34). The fact that $g$ is central implies that the corresponding action is trivial (cf. [57], Prop. 7.2.3, p. 106).

Now we know from Proposition 1.13.2 that $H^2(\mathbb{Z}_q, U(H)) \cong U(H)/U(H)^q$, where $U(H)$ denotes the group of units of $H$. Since $H \cong K^p$ as an algebra and since $K$ is algebraically closed, every unit is a $q$-th power, and we see that $H^2(\mathbb{Z}_q, U(H)) = \{1\}$. The cocycle involved in the crossed product is therefore trivial. This implies that $B \cong H \otimes K[\mathbb{Z}_q]$ as an algebra (cf. [57], Thm. 7.3.4, p. 113); in particular, $B$ is commutative. This is a contradiction. □

**9.3** We now rule out the case that the orders of $g$ and $\gamma$ are different.

**Proposition** $g$ and $\gamma$ have the same order $p$.

**Proof.** (1) Suppose that the order of $\gamma$ is $q$. By exchanging $B$ and $B^*$, we could exchange $g$ and $\gamma$, and therefore we can assume that $p < q$. Since then $B^*$ contains a grouplike element of order $q$, we get a Hopf algebra injection $K[\mathbb{Z}_q] \to B^*$. Since $K[\mathbb{Z}_q]$ is self-dual (cf. Paragraph 1.10), we get by dualization a Hopf algebra surjection
$$\pi : B \to K[\mathbb{Z}_q]$$
Then $\pi(g)$ is a grouplike element whose order simultaneously divides $p$ and $q$. Therefore, we have $\pi(g) = 1$. Consider the space of coinvariant elements:
$$A := \{b \in B \mid (\mathrm{id}_B \otimes \pi)\Delta_B(b) = b \otimes 1\}$$
By the normal basis theorem, we know that $\dim A = p$. Since $A$ contains $g$ and its powers, we see that $A$ is precisely the vector space spanned by the powers of $g$, i. e., $A$ is isomorphic to $K[\mathbb{Z}_p]$ as an algebra. In particular, $A$ is commutative.

(2) By the second version of the normal basis theorem mentioned in the previous paragraph, we know that $B$ is isomorphic to a crossed product of $A$ and $K[\mathbb{Z}_q]$. Since $A$ is commutative and $K[\mathbb{Z}_q]$ is cocommutative, $A$ is not only a twisted $K[\mathbb{Z}_q]$-module (cf. [57], Lem. 7.1.2, p. 101), but rather an ordinary $K[\mathbb{Z}_q]$-module algebra. As in Paragraph 6.1, $\mathbb{Z}_q$ acts by permutations on the set $E$ of primitive idempotents of $A$. The orbits of this action have length 1 or $q$. Since $\dim A = \mathrm{card}(E) = p < q$, orbits of length $q$ cannot occur. Therefore, all orbits are of length 1, which means that the action is trivial. Now the argument used in the preceding paragraph shows that $B \cong A \otimes K[\mathbb{Z}_q]$ as an algebra (cf. [57], Thm. 7.3.4, p. 113); in particular, $B$ is commutative, which is a contradiction. □

**9.4** Suppose that $\lambda_B \in B^*$ is an integral that satisfies $\lambda_B(1_B) = 1$. In consistency with our conventions in Paragraph 1.12, we denote the character ring of $B$ by $\mathrm{Ch}(B)$. If $e$ is a primitive idempotent in $\mathrm{Ch}(B)$, the class equation of G. I. Kac and Y. Zhu (cf. [33], Thm. 2, p. 158, [90], Thm. 1, p. 56, [42], p. 2842) says that $\dim B^*e$ divides $\dim B$. Since $\dim B^*e$ cannot be equal to $pq$, it must be 1, $p$, or $q$.

The left coregular representation of $B$ induces an isomorphism
$$Z(B) \to \mathrm{Ch}(B), b \mapsto (b \to \lambda_B)$$
which restricts to a bijection between the central grouplike elements of $B$ and those idempotents $e$ of $\mathrm{Ch}(B)$ that generate one-dimensional left ideals of $B^*$, i. e., satisfy $\dim B^*e = 1$ (cf. [71], Lem. 4.14, p. 50, [66], Prop. 6, p. 598). Note that, by the Wedderburn structure theorem, all one-dimensional left ideals of $B^*$ are already two-sided ideals and contain a unique idempotent, which is central.

Since, by Proposition 9.2, the unit element is the only central grouplike element of $B$, we see that $\lambda_B$ is the only primitive idempotent of $\text{Ch}(B)$ that generates a one-dimensional left ideal of $B^*$. Therefore, the following elementary lemma allows us to determine the number of primitive idempotents $e$ in $\text{Ch}(B)$ that satisfy $\dim B^*e = p$ resp. $\dim B^*e = q$:

**Lemma** There are unique nonnegative integers $n_p$ and $n_q$ such that:

$$1 + pn_p + qn_q = pq$$

They satisfy the inequalities $0 < n_p < q$ and $0 < n_q < p$.

**Proof.** To prove existence, let $0 \le n_p \le q-1$ be the unique number that satisfies $n_p p \equiv -1 \mod q$. Similarly, let $0 \le n_q \le p-1$ be the unique number that satisfies $n_q q \equiv -1 \mod p$. Since $q \mid n_p p + 1$ and $p \mid n_q q + 1$, we see that:

$$pq \mid (n_p p + 1)(n_q q + 1) = n_p n_q pq + n_p p + n_q q + 1$$

Therefore, we have $pq \mid n_p p + n_q q + 1$, which means that there is an integer $m$ such that $pqm = n_p p + n_q q + 1$. But then we have:

$$pqm \le (q-1)p + (p-1)q + 1 = 2pq - (p+q) + 1 < 2pq$$

We conclude that $0 < m < 2$, i. e., $m = 1$, and therefore $pq = n_p p + n_q q + 1$.

To prove uniqueness, observe that we obviously have $n_p p \equiv -1 \mod q$ and $n_q q \equiv -1 \mod p$; in particular, $n_p$ and $n_q$ are nonzero. We have $n_p p < pq$ and therefore $n_p < q$. Similarly, we have $n_q < p$. □

**9.5** The grouplike element $\gamma$ gives rise to idempotents in the character ring by performing an inverse discrete Fourier transform as in Paragraph 7.6. For $j = 0, \ldots, p-1$, we introduce the elements:

$$e_j := \frac{1}{p} \sum_{i=0}^{p-1} \zeta^{-ij} \gamma^i$$

where $\zeta$ is a primitive $p$-th root of unity.

**Proposition**
1. $p$ divides $q - 1$.

2. The unique nonnegative integers $n_p$ and $n_q$ satisfying $1 + pn_p + qn_q = pq$ are $n_p = \frac{q-1}{p}$ and $n_q = p - 1$.

3. For all $j = 0, \ldots, p-1$, we have $\dim B^* e_j = q$.

4. For $j = 1, \ldots, p-1$, $e_j$ is primitive in $\text{Ch}(B)$.

5. The idempotent $e_0$ can be decomposed in the form:

$$e_0 = \lambda_B + \sum_{i=1}^{n_p} d_i$$

where $\lambda_B \in B^*$ is an integral satisfying $\lambda_B(1_B) = 1$ and $d_1, \ldots, d_{n_p}$ are primitive idempotents in $\mathrm{Ch}(B)$ satisfying $\dim B^* d_i = p$.

**Proof.** (1) As we have already noted several times, $e_0, \ldots, e_{p-1}$ is a complete system of primitive, orthogonal idempotents in the group ring $K[G(B^*)]$, which is spanned by the powers of $\gamma$. By the Nichols-Zoeller theorem, $B^*$ is free as a right $K[G(B^*)]$-module, i. e., we have:

$$B^* \cong K[G(B^*)]^q$$

This implies:

$$B^* e_j \cong (K[G(B^*)]e_j)^q$$

Since $\dim K[G(B^*)]e_j = 1$, this implies $\dim B^* e_j = q$. This establishes the third assertion.

(2) Since $\epsilon_{B^*}(e_0) = 1$, we have $\lambda_B e_0 = \lambda_B$. Therefore $e_0 = \lambda_B + (e_0 - \lambda_B)$ is a nontrivial decomposition of $e_0$ into orthogonal idempotents, which means that $e_0$ cannot be primitive. Suppose that

$$e_0 - \lambda_B = \sum_{i=1}^{m} d_i$$

is a decomposition of $e_0 - \lambda_B$ into primitive idempotents of $\mathrm{Ch}(B)$. Then we have:

$$B^* e_0 = B^* \lambda_B \oplus \bigoplus_{i=1}^{m} B^* d_i$$

As explained in Paragraph 9.4, we have $\dim B^* d_i = p$ or $\dim B^* d_i = q$. Since we have $\dim B^* d_i < \dim B^* e_0 = q$, we must have $\dim B^* d_i = p$ and $p < q$. Now we have:

$$q = \dim B^* e_0 = 1 + \sum_{i=1}^{m} \dim B^* d_i = 1 + mp$$

We therefore see that $p \mid pm = q - 1$. Since we have:

$$pq = 1 + \frac{q-1}{p} p + (p-1)q$$

we must have $n_p = m = \frac{q-1}{p}$ and $n_q = p - 1$. This proves the first, second, and fifth assertion.

(3) It remains to prove the fourth assertion. Assume on the contrary that, for some $i \in \{1, \ldots, p-1\}$, $e_i$ were not primitive. Then it could be decomposed into primitive idempotents:
$$e_i = \sum_{j=1}^{n} e'_j$$

This would lead to the vector space decomposition:
$$B^* e_i = \bigoplus_{j=1}^{n} B^* e'_j$$

As above, we must have $\dim B^* e'_j = p$, and therefore we have:
$$q = \dim B^* e_i = \sum_{j=1}^{n} \dim B^* e'_j = np$$

This is obviously impossible. $\square$

**9.6** We will prove now that the dimensions of the simple modules of $B$ that are not 1 are divisible by $p$. Recall the definition of the left adjoint representation (cf. [57], Def. 3.4.1, p. 33):
$$\mathrm{ad}_B : B \to \mathrm{End}(B), b \mapsto \mathrm{ad}_B(b)$$

where $\mathrm{ad}_B(b)(b') := b_{(1)} b' S_B(b_{(2)})$. We denote its character by $\chi_A$:
$$\chi_A(b) = \mathrm{Tr}(\mathrm{ad}_B(b))$$

**Lemma**
1. $\zeta := \gamma(g)$ is a primitive $p$-th root of unity.

2. We have: $\chi_A(g) = p$

3. If $\chi$ is the character of a simple module of dimension $n > 1$, we have $\chi(g) = 0$ and $p \mid n$.

**Proof.** (1) The proof of the first statement follows the proof of [47], Prop. 1.2, p. 362: In any case, $\gamma(g)$ is a $p$-th root of unity; we have to show that it is primitive. So assume on the contrary that $\gamma(g) = 1$. As in the proof of Proposition 9.3, we have a Hopf algebra injection $K[\mathbb{Z}_p] \to B^*$. Dually, we get a Hopf algebra surjection:
$$\pi : B \to K^{\mathbb{Z}_p}$$

We can identify $K^{\mathbb{Z}_p}$ with $K^p$ by mapping the basis of primitive idempotents to the canonical basis of $K^p$. Under this identification, the mapping $\pi$ takes the form:
$$\pi : B \to K^p, b \mapsto (\epsilon_B(b), \gamma(b), \ldots, \gamma^{p-1}(b))$$

Therefore, we have $\pi(g) = 1$. Consider the space of coinvariant elements:
$$A := \{b \in B \mid (\mathrm{id}_B \otimes \pi)\Delta_B(b) = b \otimes 1\}$$

By the normal basis theorem, we know that $\dim A = q$. Since $A$ contains $g$, it is a left relative $K[\mathbb{Z}_p]$-$B$-Hopf module with respect to the action of $g$ and the regular coaction of $B$. By the Nichols-Zoeller theorem, $A$ is a free $K[\mathbb{Z}_p]$-module; in particular, $p$ divides $q$, which is a contradiction.

(2) Recall from Paragraph 9.5 the definition of the idempotents $e_0, \ldots, e_{p-1}$, where we define $\zeta := \gamma(g)$. The mappings
$$\xi_j : \mathrm{Ch}(B) \to K, \chi \mapsto \chi(g^j)$$
are one-dimensional characters of the character ring that satisfy $\xi_j(e_i) = \delta_{ji}$. Since the idempotents $e_1, \ldots, e_{p-1}$ are primitive by Proposition 9.5, these must be the central idempotents of $\mathrm{Ch}(B)$ that generate the one-dimensional two-sided ideals corresponding to $\xi_1, \ldots, \xi_{p-1}$.

On the other hand, we also have a different description of these idempotents. If $\{\chi_1, \ldots, \chi_k\}$ is the set of irreducible characters of $B$, we know from [80], Prop. 3.5, p. 211 that $\sum_{i=1}^{k} \chi_i \otimes \chi_{\bar{i}}$ is a Casimir element of the Frobenius algebra $\mathrm{Ch}(B)$, i. e., we have:
$$\sum_{i=1}^{k} \chi\chi_i \otimes \chi_{\bar{i}} = \sum_{i=1}^{k} \chi_i \otimes \chi_{\bar{i}}\chi$$
for all $\chi \in \mathrm{Ch}(B)$. This implies that we have:
$$\sum_{i=1}^{k} \chi\chi_i \xi_j(\chi_{\bar{i}}) = \xi_j(\chi) \sum_{i=1}^{k} \chi_i \xi_j(\chi_{\bar{i}})$$

In particular, since $\chi_A = \sum_{i=1}^{k} \chi_i \chi_{\bar{i}}$ (cf. [80], Par. 3.3, p. 208), we have:
$$(\sum_{i=1}^{k} \chi_i \xi_j(\chi_{\bar{i}}))^2 = \xi_j(\sum_{i=1}^{k} \chi_i \xi_j(\chi_{\bar{i}})) \sum_{i=1}^{k} \chi_i \xi_j(\chi_{\bar{i}}) = \chi_A(g^j) \sum_{i=1}^{k} \chi_i \xi_j(\chi_{\bar{i}})$$

Since $\chi_A$ is invertible (cf. [80], Thm. 3.8, p. 215), $\chi_A(g^j)$ is nonzero. Therefore, for $j = 1, \ldots, p-1$, the element $\frac{1}{\chi_A(g^j)} \sum_{i=1}^{k} \chi_i \xi_j(\chi_{\bar{i}})$ is an idempotent that must coincide with the idempotent $e_j$. This means that we have for $j = 1, \ldots, p-1$:
$$\frac{1}{\chi_A(g^j)} \sum_{i=1}^{k} \chi_i \xi_j(\chi_{\bar{i}}) = \frac{1}{p} \sum_{i=0}^{p-1} \zeta^{-ij} \gamma^i$$

Since the grouplike elements $\gamma^i$ appear among all characters, we can compare coefficients. We then get that $\chi_A(g^j) = p$ and $\xi_j(\chi_i) = \chi_i(g^j) = 0$ if $\chi_i \notin \{\epsilon_B, \gamma, \ldots, \gamma^{p-1}\}$, i. e., if the degree of $\chi_i$ is greater than 1.

(3) If $V$ is a simple module of dimension greater than 1 with corresponding character $\chi$, we look at the eigenspaces of $g$:

$$V^j := \{v \in V \mid gv = \zeta^j v\}$$

Then we have:

$$\sum_{j=0}^{p-1} (\dim V^j)\zeta^j = \chi(g) = 0$$

$\zeta$ is therefore a zero of the polynomial $\sum_{j=0}^{p-1}(\dim V^j)t^j \in \mathbb{Q}[t]$, and therefore the minimum polynomial $\sum_{j=0}^{p-1} t^j$ of $\zeta$ divides this polynomial (cf. [28], Sec. III.1, p. 112). Since both polynomials have the same degree, they must be equal up to a scalar, and we get:

$$\dim V^0 = \dim V^1 = \ldots = \dim V^{p-1}$$

Therefore, $p$ divides $\dim V = p \dim V^0$ (cf. [2], Lem. 2.6, p. 433 for a similar argument). $\square$

**9.7** We are now in a position that allows us to determine precisely the dimensions of the simple modules of $B$.

**Proposition**
1. Up to isomorphism, $B$ has $p$ simple modules of dimension 1 corresponding to the characters $\epsilon_B, \gamma, \gamma^2, \ldots, \gamma^{p-1}$, and $\frac{q-1}{p}$ simple modules of dimension $p$.

2. The character ring $\mathrm{Ch}(B)$ is commutative.

**Proof.** $B$ has $p$ one-dimensional representations. Denote the number of irreducible representations that are not one-dimensional by $m$. Since these are of dimension at least $p$, we have $pq = \dim B \geq p + p^2 m$, and equality holds if and only if all representations are of dimension 1 or $p$. This implies $q \geq 1 + pm$ resp. $m \leq \frac{q-1}{p}$. Therefore, we have:

$$\dim \mathrm{Ch}(B) = p + m \leq p + \frac{q-1}{p}$$

On the other hand, we know from Proposition 9.5 that a complete system of primitive, orthogonal idempotents of the character ring contains $1 + n_p + n_q = 1 + \frac{q-1}{p} + p - 1$ elements. Therefore, we have:

$$\dim \mathrm{Ch}(B) \geq p + \frac{q-1}{p}$$

This implies that equality holds, i. e., we have $m = \frac{q-1}{p}$ irreducible representations that are not one-dimensional, and all of these are of dimension $p$. Moreover, we see that the dimension of the character ring is equal to the cardinality of a complete system of primitive, orthogonal idempotents. For a semisimple algebra over an algebraically closed field (cf. [90], Lem. 2, p. 55, [80], Thm. 3.8, p. 215), this is only possible if it is commutative. □

**9.8** We have proved in Lemma 9.6 that $\gamma(g)$ is a primitive $p$-th root of unity. If
$$\pi : B \to K^p, b \mapsto (\epsilon_B(b), \gamma(b), \dots, \gamma^{p-1}(b))$$
is the Hopf algebra homomorphism considered in the proof of that result, $\pi(g)$ is a nontrivial grouplike element of order $p$, and therefore the Hopf subalgebra $\pi(H)$ is equal to $K^p$. Since $K^p$ is isomorphic to $H$, we get by composition with such an isomorphism a Hopf algebra map from $B$ to $H$ that restricts to the identity on $H$. By the Radford projection theorem (cf. [65], Thm. 3, p. 336), the map
$$A \otimes H \to B, a \otimes h \mapsto ah$$
is an isomorphism from the Radford biproduct $A \otimes H$ to $B$, where the subspace of coinvariant elements with respect to $\pi$ is denoted by
$$A := \{b \in B \mid (\mathrm{id} \otimes \pi)\Delta_B(b) = b \otimes 1\}$$
and is regarded as a Yetter-Drinfel'd Hopf algebra over $H$ in a suitable way (cf. [65], Eq. (3.3b), p. 337). Therefore, we shall assume from now on that $B = A \otimes H$ is a Radford biproduct, where $H = K[\mathbb{Z}_p]$ is the group ring of the cyclic group of order $p$ and $A$ is a left Yetter-Drinfel'd Hopf algebra over $H$.

As in Paragraph 8.4, we have that $A$ is semisimple. We now prove:

**Proposition** $A$ is commutative.

**Proof.** Suppose that $W$ is a simple $A$-module which is not the trivial module, and denote the corresponding centrally primitive idempotent by $e$. We have to prove that $\dim W = 1$. By Corollary 2.3, there is a simple $B$-module $V$ such that $e \in \kappa(V)$. Denote the character of $V$ by $\chi_V$. By Proposition 9.7, $V$ has dimension 1 or $p$. Since $e \in \kappa(V)$, $W$ is a submodule of the restriction of $V$ to $A$, and therefore $\dim V = 1$ implies that $\dim W = 1$. We therefore may assume that the dimension of $V$ is $p$. The $B$-module $V^* \otimes V$ can be decomposed into simple modules of dimension 1 or $p$. If $m_1$ is the number of one-dimensional modules and $m_p$ is the number of $p$-dimensional modules occurring in this decomposition, we have:
$$p^2 = m_1 + pm_p$$
Since by Schur's lemma the trivial module appears exactly once in this decomposition, we have $m_1 \geq 1$, and since $m_1$ is divisible by $p$, this implies that a nontrivial one-dimensional module appears in this decomposition. As noted in Paragraph 1.12, this means that we have $\chi_V(\epsilon_A \otimes \gamma') = \chi_V$ for some nontrivial element $\gamma' \in \hat{\mathbb{Z}}_p$, i. e., the isotropy group $\kappa^*(V)$ is nontrivial, and therefore equal to $\hat{\mathbb{Z}}_p$. This implies by Corollary 2.6 that the $\mathbb{Z}_p$-orbit $\kappa(V)$ that contains $e$ has length $p$, which means that $V$ is purely unstable. Therefore, we have by Proposition 2.4 that $p = \dim V = p \dim W$. This implies that $\dim W = 1$. □

As in Paragraph 8.4, we have the following corollary:

**Corollary**  $A$ is cocommutative.

**9.9**   We summarize the results of this section in the following theorem:

**Theorem**   Suppose that $K$ is an algebraically closed field of characteristic zero. Suppose that $p$ and $q$ are distinct prime numbers. Suppose that $B$ is a semisimple Hopf algebra over $K$ of dimension $pq$ such that both $B$ and $B^*$ contain nontrivial grouplike elements, i. e, grouplike elements different from the unit. Then $B$ is commutative or cocommutative.

**Proof.**   We have shown that, if this is not the case, $B$ is isomorphic to a Radford biproduct $A \otimes H$, where $A$ is a semisimple, commutative, cocommutative Yetter-Drinfel'd Hopf algebra over $H = K[\mathbb{Z}_p]$. Since $p^2$ does not divide $\dim A = q$, we get from Theorem 7.9, or already from Corollary 6.7, that $A$ is trivial. From Proposition 1.11, we therefore see that either the action or the coaction of $H$ on $A$ is trivial. Therefore, the Radford biproduct $A \otimes H$ is either as an algebra or as a coalgebra the ordinary tensor product of $A$ and $H$. This implies that $B$ is commutative or cocommutative, which is a contradiction. □

**9.10**   We have already explained in Paragraph 8.10 why the determination of Hopf algebras of a given dimension that are commutative or cocommutative reduces to the determination of finite groups of that order. However, the determination of groups of order $pq$ is considerably simpler than the determination of groups of order $p^3$. Suppose that $G$ is a group of order $pq$ for two distinct primes $p$ and $q$, where $p < q$. The Sylow subgroups $G_p$ and $G_q$ are cyclic; the number of $q$-Sylow subgroups is 1 or $p$, and is congruent to 1 modulo $q$. Therefore, $G_q$ is normal (cf. [34], Kap. III.3, Aufg. 4, p. 47), and $G$ is the semidirect product of $G_p$ and $G_q$. Since $\text{Aut}(G_q)$ is cyclic of order $q-1$, $\text{Aut}(G_q)$ does not contain elements of order $p$ if $p$ does not divide $q-1$, and therefore the action determining the semidirect product is trivial in this case. If $p$ does divide $q-1$, $\text{Aut}(G_q)$ contains a unique subgroup of order $p$, and the action can be nontrivial. A detailed analysis of the isomorphism classes arising in this latter case yields (cf. [3], Chap. 3, Exerc. 4, p. 33, [25], Kap. I, Satz 8.10, p. 40, [34], Kap. II.3,

Aufg. 1, p. 35, [82], Chap. 2, § 2, Example 1, p. 103):

**Proposition** Suppose that $p$ and $q$ are distinct prime numbers with $p < q$.

1. If $p \nmid q - 1$, then all groups of order $pq$ are isomorphic to $\mathbb{Z}_{pq}$.

2. If $p \mid q - 1$, then besides $\mathbb{Z}_{pq}$ there is a unique isomorphism type of nonabelian groups of order $pq$, which is a semidirect product $\mathbb{Z}_p \ltimes \mathbb{Z}_q$.

As in Paragraph 8.10, we have the following corollary:

**Corollary** Suppose that $K$ is an algebraically closed field of characteristic zero. Suppose that $p$ and $q$ are distinct prime numbers with $p < q$ and that $B$ is a semisimple Hopf algebra of dimension $pq$.

1. Suppose that $B$ is commutative and cocommutative. Then $B$ is isomorphic to $K[\mathbb{Z}_{pq}]$.

2. Suppose that $B$ is cocommutative, but not commutative. Then we have $q \equiv 1 \bmod p$, and $B$ is isomorphic to $K[\mathbb{Z}_p \ltimes \mathbb{Z}_q]$, where $\mathbb{Z}_p \ltimes \mathbb{Z}_q$ is a nonabelian semidirect product of $\mathbb{Z}_p$ and $\mathbb{Z}_q$, which is unique up to isomorphism.

3. Suppose that $B$ is commutative, but not cocommutative. Then we have $q \equiv 1 \bmod p$, and $B$ is isomorphic to $K^{\mathbb{Z}_p \ltimes \mathbb{Z}_q}$, where $\mathbb{Z}_p \ltimes \mathbb{Z}_q$ is the semidirect product above.

# 10 Applications

**10.1** In this section, we assume that $K$ is an algebraically closed field of characteristic zero, and that $p$ and $q$ are two distinct prime numbers. Since the case of the even prime 2 has been treated by A. Masuoka (cf. [48]), we shall assume that $p$ and $q$ are odd. $B$ denotes a semisimple Hopf algebra of dimension $pq$ over $K$ that is neither commutative nor cocommutative. Our goal is to find sufficient conditions for the existence of nontrivial grouplike elements, thereby arriving at a contradiction to Theorem 9.9. Among the sufficient conditions exhibited below will be the condition that the dimensions of the simple modules of $B$ divide the dimension of $B$. We note that most of the results derived in this section have been established in a similar form independently and slightly earlier by S. Gelaki and S. Westreich (cf. [21]). Although the results of their article are slightly less general, it has been shown even more recently by P. Etingof and S. Gelaki that the dimensions of the simple $B$-modules divide, in this situation, the dimension of $B$, thereby establishing that $B$ is always commutative or cocommutative (cf. [18]). The proof of this result has already been simplified by Y. Tsang and Y. Zhu (cf. [85]), and by H.-J. Schneider (cf. [74]).

Throughout this section, $n_p$ and $n_q$ will denote the unique nonnegative integers that satisfy $1 + pn_p + qn_q = pq$ considered in Lemma 9.4.

**10.2** The basic idea for the proofs of the results of this section is the following: The class equation can be used to determine the cardinality of a complete system of primitive, orthogonal idempotents in the character ring. This cardinality is a lower bound for the dimension of the character ring, and therefore for the number of irreducible representations of $B$. Together with some knowledge about the dimensions of the irreducible representations, this can be turned into a lower bound for the dimension of $B$. This idea, together with some technical improvements, leads to the following inequality:

**Proposition** Suppose that $B^*$ does not contain nontrivial grouplike elements. Then we have:
$$(q^2 - 9q + 9n_p)p \geq 24q + 9qn_p - 9$$

**Proof.** (1) Suppose that $\lambda_B \in B^*$ is an integral that satisfies $\lambda_B(1_B) = 1$. We have already used in Paragraph 9.4 that, for a primitive idempotent $e$ of the character ring $\mathrm{Ch}(B)$, $\dim B^*e$ divides $\dim B$. We have also explained there that $\dim B^*e$ must be 1, $p$, or $q$ and that the idempotents $e$ of $\mathrm{Ch}(B)$ that satisfy $\dim B^*e = 1$ are in bijection with the central grouplike elements of $B$. Since by Proposition 9.2 the unit element is the only central grouplike element of $B$, we see that $\lambda_B$ is the only primitive idempotent of $\mathrm{Ch}(B)$ that generates a one-dimensional left ideal of $B^*$.

Now suppose that
$$\epsilon_B = \lambda_B + \sum_{i=1}^{m_p} e_i + \sum_{j=1}^{m_q} e'_j$$
is a decomposition of the unit of the character ring into primitive, orthogonal idempotents such that we have $\dim B^* e_i = p$ for $i = 1, \ldots, m_p$ and $\dim B^* e'_j = q$ for $j = 1, \ldots, m_q$. We then get the decomposition
$$B^* = B^* \lambda_B \oplus \bigoplus_{i=1}^{m_p} B^* e_i \oplus \bigoplus_{j=1}^{m_q} B^* e'_j$$
and therefore, by comparing dimensions, we have $pq = 1 + pm_p + qm_q$. We therefore have $m_p = n_p$ and $m_q = n_q$ by the uniqueness of $n_p$ and $n_q$. Since orthogonal idempotents are linearly independent, we have that $\dim \mathrm{Ch}(B) \geq 1 + n_p + n_q$. Therefore, $B$ has at least $1 + n_p + n_q$ nonisomorphic irreducible representations.

(2) We have assumed that the trivial representation is the only one-dimensional representation of $B$. By a result of W. Nichols and M. B. Richmond, $B$ cannot have a two-dimensional simple module, since its dimension is odd (cf. [62], Cor. 12, p. 306). Therefore, all simple $B$-modules have at least dimension three. By a result of S. Zhu (cf. [89], Lem. 11, p. 3879), there are at least four distinct dimensions of simple $B$-modules. Therefore, there is one simple $B$-module of dimension greater than three, and another one of dimension greater than four. This implies that:
$$pq = \dim B \geq 1 + 3^2(n_p + n_q - 2) + 4^2 + 5^2$$

(3) We have $pq \geq 9n_p + 9n_q + 24$, and therefore $pq^2 \geq 9qn_p + 9qn_q + 24q$. Since $pq = 1 + pn_p + qn_q$, this implies $pq^2 \geq 9qn_p + 9(pq - pn_p - 1) + 24q$, and therefore:
$$pq^2 - 9pq + 9pn_p \geq 9qn_p + 24q - 9$$

This implies the assertion. $\square$

**10.3** Semisimple Hopf algebras of dimension $3p$ have been classified by M. Izumi and H. Kosaki (cf. [26], p. 369), and by A. Masuoka (cf. [51]). These results can also be obtained along the lines discussed below, since the present methods should be viewed as a generalization of Masuoka's approach; however, we shall not discuss this in detail, but rather consider now the case $q = 5$. In this case, the above argument can be slightly refined:

**Proposition** Suppose that $B$ is a semisimple Hopf algebra of dimension $5p$ that does not contain nontrivial grouplike elements. Then we have:
$$(61n_p - 180)p \geq 225n_p + 139$$

**Proof.** We keep the notation of the proof of the preceding proposition. For $j = 1, \ldots, n_5$, consider the left modules $B^* e'_j$ of dimension 5. By assumption, the trivial representation is the only one-dimensional representation of $B^*$; it occurs with multiplicity one. By the result of Nichols and Richmond also used in the preceding paragraph, $B^*$ cannot have a two-dimensional simple module. Therefore, the modules $B^* e'_j$ cannot contain a simple submodule of dimension 1 or 2. But then they also do not contain a simple submodule of dimension 3 or 4, because the complements of these modules would have dimension 1 or 2. This means that they are simple themselves.

If two of the left modules $B^* e'_j$ are isomorphic, they are contained in the same two-sided ideal of $B^*$. Such a two-sided ideal is isomorphic to the ring of $5 \times 5$-matrices, and therefore can contain at most five primitive, orthogonal idempotents. This implies that the isomorphism classes of the left modules $B^* e'_j$ contain at most five elements, and therefore $B^*$ has at least $\frac{n_5}{5}$ nonisomorphic simple modules of dimension 5.

Denote the number of isomorphism classes of 5-dimensional simple $B^*$-modules by $k_5$. By the result of S. Zhu mentioned above, there must be another simple module of dimension larger than three. By Proposition 9.2, $B^*$ does also not contain a central grouplike element, and therefore we conclude from Paragraph 9.4 that $B^*$ has at least $1 + n_p + n_5$ irreducible representations, too. We therefore have:

$$5p \geq 1 + 3^2(n_p + n_5 - k_5 - 1) + 5^2 k_5 + 4^2 = 8 + 9n_p + 9n_5 + 16k_5$$
$$\geq 8 + 9n_p + 9n_5 + \frac{16}{5}n_5$$

Multiplying this by 25, we get:

$$125p \geq 200 + 225n_p + 225n_5 + 80n_5 = 225n_p + 305n_5 + 200$$
$$= 225n_p + 61(5p - pn_p - 1) + 200$$

This implies that $61pn_p - 180p \geq 225n_p + 139$. $\square$

The above proposition can be used to determine certain cases in which semisimple Hopf algebras of dimension $5p$ are commutative or cocommutative:

**Corollary** Suppose that $B$ is a semisimple Hopf algebra of dimension $5p$ over an algebraically closed field of characteristic zero, where $p$ is an odd prime.

1. If $p \equiv 2 \bmod 5$ or $p \equiv 4 \bmod 5$, then $B$ is commutative and cocommutative.

2. If $p \in \{3, 13, 23, 43, 53, 73, 83, 103, 113, 163, 173, 193, 223, 233, 263\}$, then $B$ is commutative and cocommutative.

3. If $p = 11$, then $B$ is commutative or cocommutative.

**Proof.** Suppose that $B$ does not contain a nontrivial grouplike element. By Lemma 9.4, we know that $pn_p \equiv -1 \mod 5$. By the preceding proposition, we therefore have the following cases:

$$p \equiv 1 \mod 5 \Rightarrow n_p = 4 \Rightarrow 64p \geq 1039 \Rightarrow p \neq 11$$
$$p \equiv 2 \mod 5 \Rightarrow n_p = 2 \Rightarrow -58p \geq 589$$
$$p \equiv 3 \mod 5 \Rightarrow n_p = 3 \Rightarrow 3p \geq 814 \Rightarrow$$
$$p \notin \{3, 13, 23, 43, 53, 73, 83, 103, 113, 163, 173, 193, 223, 233, 263\}$$
$$p \equiv 4 \mod 5 \Rightarrow n_p = 1 \Rightarrow -119p \geq 364$$

In the cases stated above, $B$ therefore contains a nontrivial grouplike element. Dually, $B^*$ contains a nontrivial grouplike element, and therefore $B$ is commutative or cocommutative by Theorem 9.9. $B$ is therefore a group ring or a dual group ring. Since groups of these orders are commutative by Proposition 9.10, except for the case $p = 11$, the assertion follows. $\square$

**10.4** We now consider the case $q = 7$. In this case, the inequality derived in Proposition 10.2 yields the following results:

**Corollary** Suppose that $B$ is a semisimple Hopf algebra of dimension $7p$ over an algebraically closed field of characteristic zero, where $p$ is an odd prime.

1. If $p \equiv 6 \mod 7$, then $B$ is commutative and cocommutative.

2. If $p \in \{5, 11, 17, 23, 31, 59\}$, then $B$ is commutative and cocommutative.

**Proof.** If $B^*$ does not contain a nontrivial grouplike element, the inequality derived in Proposition 10.2 implies in this case:

$$(9n_p - 14)p \geq 63n_p + 159$$

We therefore have the following cases:

$$p \equiv 1 \mod 7 \Rightarrow n_p = 6 \Rightarrow 40p \geq 537$$
$$p \equiv 2 \mod 7 \Rightarrow n_p = 3 \Rightarrow 13p \geq 348 \Rightarrow p \neq 23$$
$$p \equiv 3 \mod 7 \Rightarrow n_p = 2 \Rightarrow 4p \geq 285 \Rightarrow p \notin \{3, 17, 31, 59\}$$
$$p \equiv 4 \mod 7 \Rightarrow n_p = 5 \Rightarrow 31p \geq 474 \Rightarrow p \neq 11$$
$$p \equiv 5 \mod 7 \Rightarrow n_p = 4 \Rightarrow 22p \geq 411 \Rightarrow p \neq 5$$
$$p \equiv 6 \mod 7 \Rightarrow n_p = 1 \Rightarrow -5p \geq 222$$

In the cases stated above, $B^*$ therefore contains a nontrivial grouplike element. The assertions now follow as in the preceding paragraph. $\square$

**10.5** We now look at the situation where the dimensions of the simple $B$-modules divide the dimension of $B$. These dimensions then must be 1, $p$, $q$, or $pq$. It is obviously impossible that $B$ has a simple module of dimension $pq$, because then it would contain a matrix ring of dimension $p^2q^2$. It therefore follows directly from the result of S. Zhu already used above, or alternatively from the fact that we cannot have $pq = 1 + m_p p^2 + m_q q^2$ for nonnegative integers $m_p$ and $m_q$, that $B^*$ contains a nontrivial grouplike element. Combining this with Theorem 9.9, we get:

**Corollary** Suppose that $B$ is a semisimple Hopf algebra of dimension $pq$ over an algebraically closed field of characteristic zero, where $p$ and $q$ are distinct primes. If the dimensions of the simple $B$-modules and the simple $B$-comodules divide the dimension of $B$, then $B$ is commutative or cocommutative.

As we already noted at the beginning of this section, the premise of this corollary has been recently established by P. Etingof and S. Gelaki (cf. [18]). Their proof has been simplified by Y. Tsang and Y. Zhu (cf. [85]), and by H.-J. Schneider (cf. [74]).

**Acknowledgement** The author thanks N. Andruskiewitsch, S. Natale, B. Pareigis, and P. Schauenburg for interesting discussions. He thanks M. Lorenz for helpful comments on the literature, in particular for pointing out reference [40]. He also thanks N. Andruskiewitsch for the kind permission to use his approach to the construction of Yetter-Drinfel'd Hopf algebras (cf. [1]).

Part of the results were presented at the conference 'Hopf algebras and quantum groups', June 15-18, 1998, Brussels, Belgium. The structure theorem was first presented at the 'International conference on algebra and its applications', March 25-28, 1999, Athens, United States. The author thanks the Lautrach foundation for financial support making the visit to Brussels possible, and the 'Graduiertenkolleg: Mathematik im Bereich ihrer Wechselwirkung mit der Physik' of the Deutsche Forschungsgemeinschaft for financial support making the visit to Athens possible.

Typeset using $\mathcal{AMS}$ - LaTeX

# Subject index

# Symbol index